\documentclass[11pt,english]{article}
\usepackage{amssymb, amsmath, color}
\usepackage{epsfig,epstopdf,color}
\usepackage{multicol}
\usepackage{graphicx}
\usepackage{amsfonts}
\usepackage{babel}
\newcommand{\bfi}{\bfseries\itshape}
\newcommand{\rem}[1]{}
\newcommand{\remfigure}[1]{}
\makeatletter
\@addtoreset{equation}{section}

\makeatother

\definecolor{Violet}{cmyk}{0.79,0.88,0,0}
\definecolor{Plum}{cmyk}{0.50,1,0,0}
\definecolor{Periwinkle}{cmyk}{0.57,0.55,0,0}
\definecolor{ForestGreen}{cmyk}{0.91,0,0.88,0.12}
\definecolor{OliveGreen}{cmyk}{0.64,0,0.95,0.40}
\definecolor{BrickRed}{cmyk}{0,0.89,0.94,0.28}
\definecolor{DarkOrchid}{cmyk}{0.40,0.80,0.20,0}
\definecolor{Fuchsia}{cmyk}{0.47,0.91,0,0.08}
\definecolor{Mulberry}{cmyk}{0.34,0.90,0,0.02}
\definecolor{Maroon}{cmyk}{0,0.87,0.68,0.32}
\definecolor{Mahogany}{cmyk}{0,0.85,0.87,0.35}
\definecolor{RawSienna}{cmyk}{0,0.72,1,0.45}

\makeatletter
\@addtoreset{equation}{section}

\makeatother

\def\contract{\makebox[1.2em][c]{\mbox{\rule{.6em}
{.01truein}\rule{.01truein}{.6em}}}}

\newtheorem{theorem}{Theorem}

\newtheorem{answer}[theorem]{Answer}

\newtheorem{corollary}[theorem]{Corollary}

\newtheorem{definition}[theorem]{Definition}
\newtheorem{example}[theorem]{Example}
\newtheorem{exercise}[theorem]{Exercise}
\newtheorem{lemma}[theorem]{Lemma}

\newtheorem{problem}[theorem]{Problem}
\newtheorem{proposition}[theorem]{Proposition}
\newtheorem{question}[theorem]{Question}
\newtheorem{remark}[theorem]{Remark}
\newtheorem{solution}[theorem]{Solution}

\numberwithin{theorem}{section}
\newenvironment{proof}[1][Proof]{\textbf{#1.} }{\ \rule{0.5em}{0.5em}}

\def\0{{\bf 0}}

\def\R{\mathbb{R}}

\begin{document}

\title{Applications of Poisson Geometry to Physical Problems\\
$\quad$\\
\large Summer School and Conference on Poisson Geometry
\\ICTP, Trieste, Italy, 4-22 July 2005}

\author{Prof Darryl D. Holm\\
d.holm@imperial.ac.uk\\
dholm@lanl.gov}
\date{August 2, 2007}
\maketitle

\newpage
\tableofcontents

\newpage


\subsection*{Preface}

These being lecture notes for a summer school, one should not seek original material in them. Rather, the most one could hope to find would be the insight arising from incorporating a unified approach (based on reduction by symmetry of Hamilton's principle) with some novel applications. I hope the reader will find insight in the lecture notes, which are meant to be informal, more like stepping stones than a proper path. 

Many excellent encyclopedic texts have already been published on the foundations of
this subject and its links to symplectic and  Poisson geometry. See, for
example, \cite{AbMa1978}, \cite{Ar1979}, \cite{GuSt1984}, \cite{JoSa98},
\cite{LiMa1987}, \cite{MaRa1994}, \cite{McSa1995} and many more. In fact,
the scope encompassed by the modern literature on this subject is a bit
overwhelming. In following the symmetry-reduction theme in geometric
mechanics from the Euler-Poincar\'e viewpoint, I have tried to select only
the material the student will find absolutely necessary for solving the
problems and exercises, at the level of a beginning postgraduate student.
The primary references are \cite{Ma1992}, \cite{MaRa1994},
\cite{Le2003}, \cite{Bl2004}, \cite{RaTuSbSoTe2005}. Other very useful references are \cite{ArKh1998} and \cite{Ol2000}. The reader may see the strong influences of all these references in these lecture notes, but expressed at a considerably lower level of mathematical sophistication than the originals.

The scope of these lectures is quite limited: a list of the topics in geometric mechanics not included in these lectures would fill volumes!
The necessary elements of calculus on smooth manifolds and the basics of Lie group theory are only briefly described here, because these topics were discussed in more depth by other lecturers at the summer school. Occasional handouts are included that add a bit more depth in certain key topics. The main subject of these lecture notes is the use of Lie symmetries in Hamilton's principle to derive symmetry-reduced equations of motion and to analyze their solutions. The Legendre transformation provides the Hamiltonian formulation of these equations in terms of Lie-Poisson brackets.

For example, we consider Lagrangians in Hamilton's principle defined on the tangent space $TG$ of a Lie group $G$. Invariance of such a Lagrangian under the action of $G$ leads to the symmetry-reduced Euler-Lagrange equations called the Euler-Poincar\'e equations. In this case, the invariant Lagrangian is defined on the Lie algebra of the group and its Euler-Poincar\'e equations are defined on the dual Lie algebra, where dual is defined by the operation of taking variational derivative. On the Hamiltonian side, the Euler-Poincar\'e equations are Lie-Poisson and they possess accompanying momentum maps, which encode both their conservation laws and the geometry of their solution space. 

The standard Euler-Poincar\'e examples are treated, including 
particle dynamics, the rigid body, the heavy top and geodesic motion on 
Lie groups. Additional topics deal with Fermat's principle, the
$\mathbb{R}^3$ Poisson bracket, polarized optical traveling waves,
deformable bodies (Riemann ellipsoids) and shallow water waves, including
the integrable shallow water wave system known as the 
Camassa-Holm equation. The lectures end with the semidirect-product
Euler-Poincar\'e reduction theorem for ideal fluid dynamics. This theorem
introduces the Euler--Poincar\'e variational principle for incompressible
and compressible motions of ideal fluids, with applications to geophysical fluids. It also leads to their Lie-Poisson Hamiltonian formulation. 

Some of these lectures were first given at the MASIE (Mechanics and
Symmetry in Europe) summer school in 2000 \cite{Ho2005}. I am grateful to the MASIE participants for their helpful remarks and suggestions  which led to many improvements in those lectures. For their feedback and
comments, I am also grateful to my colleagues at Imperial College London, especially Colin Cotter, Matthew Dixon, J. D. Gibbon, J. Gibbons, G.~Gottwald, J.T. Stuart, J.-L. Thiffeault, Cesare Tronci and the students who attended these lectures in my classes at Imperial College. After each class, the students were requested to turn in a response sheet on which they answered two questions. These questions were,
``What was this class about?'' and ``What question would you like
to see pursued in the class?'' The answers to these questions helped
keep the lectures on track with the interests and understanding of the
students and it enfranchised the students because they themselves selected the material in several of the lectures. 

I am enormously grateful to many friends and colleagues
whose encouragement, advice and support have helped sustained my interest
in this field over the years. I am particularly grateful to J. E. Marsden,
T. S. Ratiu and A. Weinstein for their faithful comraderie in many
research endeavors. 

\newpage

\section{Introduction}
\subsection{Road map for the course}
\begin{itemize}
\item
Spaces -- Smooth Manifolds
\item
Motion -- Flows $\phi_t\circ\phi_s=\phi_{t+s}$ of Lie groups acting on
smooth manifolds 
\item
Laws of Motion and discussion of solutions
\item Newton's Laws
     \begin{itemize}
     \item
     Newton: $dp/dt=F$, for momentum $p$ and prescribed force $F$
(on $\mathbb{R}^n$ historically)
     \item
     Optimal motion 
          \begin{itemize}
          \item
          Euler-Lagrange equations -- optimal ``action'' (Hamilton's
principle)
          \item
          Geodesic motion -- optimal with respect to kinetic energy metric
          \end{itemize}
     \end{itemize}
\end{itemize}

     \begin{itemize}
     \item
     Lagrangian  and Hamiltonian Formalism
          \begin{itemize}
          \item Newton's Law of motion
          \item Euler-Lagrange theorem
          \item Noether theorem
          \item Euler-Poincar\'e theorem
          \item Kelvin-Noether theorem
          \end{itemize}
     \item
     Applications and examples
          \begin{itemize}
          \item Geodesic motion on a Riemannian manifold
          \item Rigid body -- geodesic motion on $SO(3)$
          \item Other geodesic motion, e.g., Riemann ellipsoids on $GL(3,R)$
          \item Heavy top
          \end{itemize}
     \item
     Lagrangian mechanics on Lie groups \& Euler-Poincar\'e (EP) equations
          \begin{itemize}
          \item EP$(G)$, EP equations for geodesics on a Lie group $G$
          \item EPDiff$(\mathbb{R})$ for geodesics on Diff$(\mathbb{R})$
          \item Pulsons, the singular solutions of EPDiff$(\mathbb{R})$) wrt any norm
          \item Peakons, the singular solitons for EPDiff$(\mathbb{R},H^1)$, wrt
                the $H^1$ norm
          \item EPDiff$(\mathbb{R}^n)$ and singular geodesics
          \item Diffeons and momentum maps for EPDiff$(\mathbb{R}^n)$
          \end{itemize}
     \item
     Euler-Poincar\'e (EP) equations for continua
          \begin{itemize}
          \item EP semidirect-product reduction theorem
          \item Kelvin-Noether circulation theorem
          \item EP equations with advected parameters for geophysical fluid dynamics
          \end{itemize}
     \end{itemize}
\newpage
\noindent
{\Large\bfi Hamilton's principle of stationary action:}
\begin{multicols}{2}{
\noindent
{\large\bfi Lagrangians on $T\mathbb{R}^{3N}$:} \\
Euler-Lagrange equations\\
Noether's theorem \\
Symmetry $\implies$ cons. laws \\
Legendre transformation\\
Hamilton's canonical equations\\
Poisson brackets\\
Symplectic manifold\\
Momentum map\\
Reduction by symmetry

\noindent
{\large\bfi G-invariant Lagrangians on $TG$:} \\
Euler-Poincar\'e equations\\
Kelvin-Noether theorem\\
Cons. laws are built-in\\
Legendre transformation\\
Lie-Hamilton equations\\
Lie-Poisson brackets\\
Poisson manifold\\
Momentum map\\
Reduction to coadjoint orbits
}
\end{multicols}

\subsection{Motivation for the geometric approach}

We begin with a series of outline sketches to motivate the geometric approach taken in the course and explain more about its content. \\

\noindent \textbf{\large Why is the geometric approach useful?}
\begin{itemize}
\item Defines problems on manifolds
\begin{itemize}
\item coordinate-free
\begin{itemize}
\item don't have to re-do calculations when changing coordinates
\item more compact
\item unified framework for expressing ideas and using symmetry
\end{itemize}
\end{itemize}
\item ``First principles'' approach
\begin{itemize}
\item variational principles
\item systematic -- unified approach \\
e.g. similarity between tops and fluid dynamics (semi-direct product), and MHD, and \dots
\item POWER \\
Geometric constructions can give useful answers without us having to find and work with
complicated explicit solutions. e.g. stability of rigid body equilibria.
\end{itemize}
\end{itemize}

\noindent \textbf{\large Course Outline}
\begin{itemize}
\item{Geometrical Structure of Classical Mechanics}
\begin{itemize}
\item Smooth manifolds
\begin{itemize}
\item calculus
\item tangent vectors
\item action principles
\end{itemize}
\item Lie groups
\begin{itemize}
\item flow property $\phi_{t+s} = \phi_t\circ\phi_s$
\item symmetries encode conservation laws into geometry
\item richer than vector spaces
\end{itemize}
\item Variational principles with symmetries
\begin{itemize}
\item Euler-Lagrange equations $\rightarrow$ Euler-Poincar\'e equations (more compact)
\item Two main formulations:
\begin{multicols}{2}{
Lagrangian side:\\
Hamilton's principle \\
Noether's theorem \\
symmetry $\implies$ cons. laws \\
momentum maps

Hamiltonian side:\\
Lie-Poisson brackets\\
cons. laws $\iff$ symmetries \\
momentum maps\\
Jacobi identity
}
\end{multicols}
(These two views are mutual beneficial!)
\end{itemize}
\end{itemize}
\item Applications and Modelling 
\begin{itemize}
\item oscillators \& resonance (e.g., LASER)
\item tops -- integrable case
\item fluids
\item waves $\begin{cases} 
\textrm{shallow water waves} \\
\textrm{optical pulses} \\
\textrm{solitons}
\end{cases}$
\end{itemize}
\end{itemize}

\section*{Range of topics}
\subsubsection*{Rigid body}
     \begin{itemize}
     \item
     Euler-Lagrange and Euler-Poincar\'e equations
     \item
     Kelvin-Noether theorem 
     \item
     Lie-Poisson bracket, Casimirs \&  coadjoint orbits
     \item
     Reconstruction and momentum maps
     \item
     The symmetric form of the rigid body
     equations ($\dot{Q}=Q\Omega$, $\dot{P}=P\Omega$)
     \item
     $\mathbb{R}^3$ bracket and intersecting level surfaces
\[
\mathbf{\dot{x}}=\nabla{C}\times\nabla{H}
=\nabla(\alpha C + \beta H)\times\nabla(\gamma C + \epsilon H)
\,,
\quad\hbox{for}\quad
\alpha\epsilon-\beta\gamma=1
\]
Examples: \\
(1) Conversion: rigid body $\Longleftrightarrow$ pendulum,\\
(3) Fermat's principle and ray optics,\\
(2) Self-induced transparency.
     \item
     Nonlinear oscillators: the $n:m$ resonance
     \item
     $SU(2)$ rigid body, Cayley-Klein parameters and Hopf fibration 
     \item
     The Poincar\'e sphere for polarization dynamics 
     \item
     3-wave resonance, Maxwell-Bloch equations, cavity resonators, symmetry reduction
     and the Hopf fibration 
     \item
     4-wave resonance, coupled Hopf fibrations, coupled Poincar\'e
spheres and coupled rigid bodies 
     \item
     Higher dimensional rigid bodies\\ -- Manakov integrable top on $O(n)$
     and its spectral problem
     \item
     Semi-rigid bodies -- geodesic motion on $GL(3)$ and Riemann ellipsoids
     \item
     Reduction with respect to subgroups of $GL(3)$ and Calogero equations
     \end{itemize}
\subsubsection*{Heavy top}
     \begin{itemize}
     \item
     Euler--Poincar\'e variational principle for the heavy top
     \item 
     Kaluza-Klein formulation of the heavy top 
     \end{itemize}
\subsubsection*{Utility}
     Kirchhoff elastica, underwater vehicles, liquid crystals, stratified
     flows, polarization dynamics of telcom optical pulses

\subsubsection*{General theory}
     \begin{itemize}
     \item
     Euler-Poincar\'e semidirect-product reduction theorem
     \item 
     Semidirect-product Lie-Poisson formulation
     \end{itemize}
\subsubsection*{Shallow water waves}
     \begin{itemize}
     \item 
     CH equation -- peakons (geodesics)
     \item 
     EPDiff equation -- (also geodesics)
     \end{itemize}

\subsubsection*{Fluid dynamics}
     \begin{itemize}
     \item
     Euler--Poincar\'e variational principle for incompressible ideal fluids
     \item
     Euler--Poincar\'e variational principle for compressible ideal fluids
     \end{itemize}

\paragraph{Outlook:}The variational principles and the Poisson brackets for the rigid
body and the heavy top provide models of a general construction associated
to Euler--Poincar\'e reduction with respect to any Lie group. The Hamiltonian
counterpart will be the semidirect-product Lie-Poisson formulation. We will often
refer to the rigid body and the heavy top for interpretation and enhanced
understanding of the general results.

\section{Review Newton, Lagrange \& Hamilton}

\begin{itemize}
\item
{\bfi Newton's Law:} $m\ddot{q}=F(q,\dot{q})$, inertial frames, uniform motion,
etc.
\item
{\bfi Lagrange's equations:}
$\frac{d}{dt}\frac{\partial L}{\partial \dot{q}}
=
\frac{\partial L}{\partial q}$ for Lagrangian $L(q,\dot{q},t)$.

Defined on the tangent bundle%
\footnote{
The terms tangent bundle and cotangent bundle are defined in Section
\ref{smooth manifold}. For now, we may think of the tangent bundle as the space of 
positions and velocities. Likewise, the cotangent bundle is the space of positions
and momenta.
}

$TQ$ of the configuration space $Q$ with coordinates
$(q,\dot{q})\in TQ$, the solution is a curve (or trajectory) in $Q$ parameterized
by time $t$. The tangent vector of the curve $q(t)$ through each
point $q\in Q$ is the velocity $\dot{q}$ along the trajectory that passes though
the point $q$ at time $t$. This vector is written $\dot{q}\in T_qQ$.

Lagrange's equations may be expressed compactly in terms of vector fields and
one-forms (differentials). Namely, the Lagrangian vector field
$X_L=\dot{q}\frac{\partial}{\partial q} +F(q,\dot{q})\frac{\partial}{\partial
\dot{q}}$ acts on the one-form
$(\frac{\partial L}{\partial\dot{q}}\,dq)$ just as a time-derivative does, to
yield
\begin{eqnarray*}
\frac{d}{dt}\Big(\frac{\partial L}{\partial \dot{q}}\,dq\Big)
=
\Big(\frac{d}{dt}\frac{\partial L}{\partial \dot{q}}\Big)\,dq
+
\Big(\frac{\partial L}{\partial \dot{q}}\Big)\, d\dot{q}
=
dL
\quad\Longrightarrow\quad
\frac{d}{dt}\frac{\partial L}{\partial \dot{q}}
=
\frac{\partial L}{\partial q}
\end{eqnarray*}
\item
{\bfi Hamiltonian} 
$H(p\cdot q)=p\dot{q}-L$
and {\bfi Hamilton's canonical equations:}
\begin{eqnarray*}
\dot{q}=\,\frac{\partial H}{\partial p}
\,,\quad
\dot{p}=-\,\frac{\partial H}{\partial q}
\end{eqnarray*}

The {\bfi configuration space} $Q$ has coordinates $q\in Q$.
Its phase space, or cotangent bundle $T^*Q$ has coordinates $(q,p)\in T^*Q$.

Hamilton's canonical equations are associated to the {\bfi canonical Poisson
bracket} for functions on phase space, by 
\begin{eqnarray*}
\dot{p}=\{p\,,\,H\}
\,,\quad
\dot{q}=\{q\,,\,H\}
\Longleftrightarrow
\dot{F}(q,p) = \{F\,,\,H\}
= \frac{\partial F}{\partial q}\frac{\partial H}{\partial p}
-
\frac{\partial F}{\partial p}\frac{\partial H}{\partial q}
\end{eqnarray*}
The canonical Poisson bracket has the following familiar properties, which may be
readily verified:
\begin{enumerate}
\item
It is bilinear,
\item
skew symmetric, $\{F\,,\,H\}=-\,\{H\,,\,F\}$,
\item
satisfies the Leibnitz rule (chain rule),
\begin{eqnarray*}
 \{FG\,,\,H\}
= 
\{F\,,\,H\}G+ F\{G\,,\,H\}
\end{eqnarray*}
for the product of any two phase space functions $F$ and $G$,
\item
and satisfies the Jacobi identity
\begin{eqnarray*}
 \{F\,,\,\{G\,,\,H\}\} + \{G\,,\,\{H\,,\,F\}\} + \{H\,,\,\{F\,,\,G\}\} = 0
\end{eqnarray*}
for any three phase space functions $F$, $G$ and $H$.

\end{enumerate}
Its Leibnitz property (chain rule) property means the canonical Poisson bracket
is a type of derivative. This derivation property of the Poisson bracket allows
its use in defining the {\bfi Hamiltonian vector field} $X_H$, by
\begin{eqnarray*}
X_H = \{\cdot\,,\,H\}
= \frac{\partial H}{\partial p}\frac{\partial }{\partial q}
-
\frac{\partial H}{\partial q}\frac{\partial }{\partial p}
\,,
\end{eqnarray*}
for any phase space function $H$. The
action of $X_H$ on phase space functions is given by
\begin{eqnarray*}
\dot{p}=X_H p
\,,\quad
\dot{q}=X_H q
\,,\quad\hbox{and}\quad
X_H (FG) = (X_HF)G + FX_HG
=
\dot{F}G+F\dot{G}
\,.
\end{eqnarray*}
Thus, solutions of Hamilton's canonical equations are the characteristic paths
of the first order linear partial differential operator $X_H$. That is,
$X_H$ corresponds to the time derivative along these characteristic paths,
given by
\begin{eqnarray}\label{HamVecCharEqns}
dt
=\frac{dq}{\partial H/\partial p}
=\frac{dp}{-\partial H/\partial q}
\end{eqnarray}
The union of these paths in phase space is called the {\bfi flow} of
the Hamiltonian vector field $X_H$.

\begin{proposition}[Poisson bracket as commutator of Hamiltonian vector fields]
The Poisson bracket $\{F\,,\,H\}$ is associated to the commutator of the
corresponding Hamiltonian vector fields $X_F$ and $X_H$ by
\begin{eqnarray*}
X_{\{F\,,\,H\}} = X_H X_F - X_F X_H =: -\,[X_F\,,\,X_H]
\end{eqnarray*}
\end{proposition}
\begin{proof}
Verified by direct computation.
\end{proof}

\begin{corollary}
Thus, the Jacobi identity for the canonical Poisson bracket $\{\cdot\,,\,\cdot\}$
is associated to the Jacobi identity for the commutator $[\cdot\,,\,\cdot]$ of the
corresponding Hamiltonian vector fields,
\begin{eqnarray*}
[X_F\,,\,[X_G\,,\,X_H] + [X_G\,,\,[X_H\,,\,X_F] + [X_H\,,\,[X_F\,,\,X_G] = 0
\,.
\end{eqnarray*}
\end{corollary}
\begin{proof}
This is the Lie algebra property of Hamiltonian vector fields, as  
verified by direct computation.
\end{proof}

\subsection{Differential forms}
The {\bfi differential}, or {\bfi exterior derivative} of a function $F$
on phase space is written
\begin{eqnarray*}
dF = F_q dq + F_p dp
\,,
\end{eqnarray*}
in which subscripts denote partial derivatives. For the Hamiltonian itself, the
exterior derivative and the canonical equations yield
\begin{eqnarray*}
dH = H_q dq + H_p dp
= -\,\dot{p} dq + \dot{q} dp
\,.
\end{eqnarray*}
The action of a Hamiltonian vector field $X_H$ on a phase space function $F$ 
commutes with its differential, or exterior derivative. Thus,
\begin{eqnarray*}
d(X_H F)  = X_H ( dF )
\,.
\end{eqnarray*}
This means $X_H$ may also act as a time derivative on differential
forms defined on phase space. For example, it acts on the time-dependent
one-form
$p\,dq(t)$ along solutions of Hamilton's equations as,
\begin{eqnarray*}
X_H\big(p\,dq\big)
=
\frac{d}{dt}\big(p\,dq\big)
&=&
\dot{p}\,dq
+
p\,d\dot{q}
\\
&=&
\dot{p}\,dq-\dot{q}\,dp+d(p\dot{q})
\\
&=&
-\,H_q dq
-\,H_p dp
 +d(p\dot{q})
\\
&=&
d(-H+p\dot{q}) =: dL(q,p)
\end{eqnarray*}
upon substituting Hamilton's canonical equations. \\

The exterior derivative of the one-form $p dq$ yields the canonical, or
symplectic two-form%
\footnote{
The properties of differential forms are summarized in the handouts in sections \ref{ext-calculus-review} and \ref{diff-form-review}.
}
\begin{eqnarray*}
d(p dq) = dp\wedge dq 
\end{eqnarray*}
Here we have used the chain rule for the exterior derivative and its property that
$d^2=0$. (The latter amounts to equality of cross derivatives for continuous
functions.) The result is written in terms of the wedge product
$\wedge$, which combines two one-forms (the line elements $dq$ and $dp$) into a
two-form (the oriented surface element $dp\wedge dq = -\,dq\wedge dp$).  As a
result, the two-form
$\omega=dq\wedge dp$ representing area in phase space is conserved along the
Hamiltonian flows:
\begin{eqnarray*}
X_H\big(dq\wedge dp\big)
=
\frac{d}{dt}\big(dq\wedge dp\big)
=
0
\end{eqnarray*}
This proves
\begin{theorem}[Poincar\'e's theorem]
Hamiltonian flows preserve area in phase space.
\end{theorem}

\begin{definition}[Symplectic two-form]
The phase space area $\omega=dq\wedge dp$ is called the symplectic two-form.
\end{definition}

\begin{definition}[Symplectic flows]
Flows that preserve area in phase space are said to be symplectic.
\end{definition}

\begin{remark}[Poincar\'e's theorem]
Hamiltonian flows are symplectic.
\end{remark}

\end{itemize}

\newpage
\section{Handout on exterior calculus, symplectic forms and Poincar\'e's theorem in higher dimensions}
\label{ext-calculus-review}

Exterior calculus on symplectic manifolds is the geometric language of
Hamiltonian mechanics. As an introduction and motivation for more detailed study, we
begin with a preliminary discussion.

In differential geometry, the operation of {\bfi contraction} denoted as
$\contract$ introduces a pairing between  vector fields and differential forms. 
Contraction is also called {\bfi substitution} of a vector field into a
differential form. For example, there are the dual relations, 
\begin{eqnarray*}
\partial_{q}\contract dq=1
=\partial_{p}\contract dp
\,,\quad\hbox{and}\quad
\partial_{q}\contract dp=0
=\partial_{p}\contract dq
\end{eqnarray*}

A Hamiltonian vector field:
\begin{eqnarray*}
X_H
=
\dot{q}\frac{\partial}{\partial q}
+\dot{p}\frac{\partial}{\partial p}
=
H_p\partial_q-H_q\partial_p
=\{\,\cdot\,,\,H\,\}
\end{eqnarray*}
satisfies 
\begin{eqnarray*}
X_H\contract dq=H_p
\quad\hbox{and}\quad
X_H\contract dp=-\,H_q
\end{eqnarray*}
The rule for contraction or substitution of a vector field into a 
differential form is to sum the substitutions of $X_H$ over the
permutations of the factors in the differential form that bring the corresponding
dual basis element into its leftmost position. For example,
substitution of the Hamiltonian vector field $X_H$ into the symplectic form  
$\omega
=
dq\,\wedge \,dp$
yields
\begin{eqnarray*}
X_H\contract \omega
=
X_H\contract(dq\,\wedge \,dp)
=
(X_H\contract dq)\,dp
-
(X_H\contract dp)\,dq
\end{eqnarray*}
In this example, $X_H\contract dq=H_p$ and $X_H\contract dp=-\,H_q$, so
\begin{eqnarray*}
X_H\contract \omega
=
H_pdp + H_qdq
= dH
\end{eqnarray*}
which follows because $\partial_q\contract dq=1=\partial_p\contract dp$
and $\partial_q\contract dp=0=\partial_p\contract dq$.
This calculation proves
\begin{theorem}[Hamiltonian vector field]
The Hamiltonian vector field $X_H=\{\,\cdot\,,\,H\,\}$ satisfies
\begin{eqnarray}
X_H\contract\,\omega&=&dH
\quad\hbox{with}\quad
\omega
=
dq\,\wedge \,dp
\label{HVF-def}
\end{eqnarray}
\end{theorem}
Relation (\ref{HVF-def}) may be taken as the {\bfi definition} 
of a Hamiltonian vector field. \\

As a consequence of this formula, the flow of $X_H$ 
preserves the closed exact two form $\omega$ for any Hamiltonian
$H$. This preservation may be verified by a formal calculation using
(\ref{HVF-def}). Along
$(dq/dt,dp/dt)=(\dot{q},\dot{p})=(H_p,-H_q)$, we have
\begin{eqnarray*}
\frac{d\omega}{dt}
&=&
d\dot{q}\,\wedge \,dp + dq\,\wedge \,d\dot{p}
=
dH_p\,\wedge \,dp - dq\,\wedge \,dH_q
\\
&=&
d(H_p\,dp +H_q\, dq)
=
d(X_H\contract\omega)
=d(dH)
=
0
\end{eqnarray*}
The first step uses the chain rule for differential forms and the third and last
steps use the property of the exterior derivative $d$ that
$d^2=0$ for continuous forms. The latter is due to equality of cross
derivatives $H_{pq}=H_{qp}$ and antisymmetry of the wedge product:
$dq\,\wedge\,dp = - dp\,\wedge \,dq$. \\

Consequently, the
relation $d(X_H\contract\,\omega)=d^2H=0$ for Hamiltonian vector fields shows

\begin{theorem}[Poincar\'e's theorem for one degree of freedom]\label{Pthm-1dof}
$\quad$\\
The flow of a Hamiltonian vector field is {\bfi symplectic}, which means it
preserves the phase-space area, or two-form, $\omega=dq\,\wedge\,dp$. 
\end{theorem}

\begin{definition}[Cartan's formula for the Lie derivative]
The operation of {\bfi Lie derivative} of a differential form $\omega$ by a vector
field $X_H$ is defined by
\begin{eqnarray}\label{LieDerivGeomDef}
\pounds_{X_H}\omega
:=
d(X_H\contract\omega)
+X_H\contract d\omega
\end{eqnarray}
\end{definition}

\begin{corollary}
Because $d\omega=0$, the symplectic property $d\omega/dt=d(X_H\contract\omega)=0$
in Poincar\'e's Theorem \ref{Pthm-1dof} may be rewritten using Lie derivative
notation as
\begin{eqnarray}\label{symp-thm}
0
=
\frac{d\omega}{dt}
=
\pounds_{X_H}\omega
:=
d(X_H\contract\omega)
+X_H\contract d\omega
=: ({\rm div}X_H)\,\omega
\,.
\end{eqnarray}
The last equality defines the {\bfi divergence} of the vector field $X_H$ in
terms of the Lie derivative.
\end{corollary}

\begin{remark}$\quad$

\begin{itemize}
\item
Relation (\ref{symp-thm}) associates Hamiltonian dynamics with the symplectic
flow in phase space of the Hamiltonian vector field $X_H$, which is
divergenceless with respect to the symplectic form $\omega$.
\item
The Lie derivative operation defined in  (\ref{symp-thm})
is equivalent to the time derivative along the characteristic paths
(flow) of the first order linear partial differential operator $X_H$, which are
obtained from its characteristic equations in (\ref{HamVecCharEqns}). This is
the {\bfi dynamical meaning} of the Lie derivative $\pounds_{X_H}$ in
(\ref{LieDerivGeomDef}) for which invariance $\pounds_{X_H}\omega=0$ gives the
geometric definition of  symplectic flows in phase space.
\end{itemize}
\end{remark}

\begin{theorem}[Poincar\'e's theorem for $N$ degrees of freedom]$\quad$\\
For a system of $N$ particles, or $N$ degrees of freedom,
the flow of a Hamiltonian vector field preserves each subvolume in the
phase space $T^*\mathbb{R}^N$. That is, let 
$\omega_n\equiv dq_n\,\wedge \,dp_n$ be the symplectic form expressed
in terms of the position and momentum of the $n-$th particle. Then
\begin{eqnarray*}
\frac{d\,\omega_M}{dt}
=
0
\,,\quad\hbox{for}\quad
\omega_M = \Pi_{n=1}^M\omega_n
\,,\quad\forall M\le N
\,.
\end{eqnarray*}
\end{theorem}

The proof of the preservation of these {\bfi Poincar\'e invariants}
$\omega_M$ with $M=1,2,\dots,N$ follows the same pattern as the
verification above for a single degree of freedom. Basically, this is because each factor 
$\omega_n = dq_n\,\wedge \,dp_n$ in
the wedge product of symplectic forms is preserved by its
corresponding Hamiltonian flow in the sum
\begin{eqnarray*}
X_H
=
\sum_{n=1}^M\Big(
\dot{q_n}\frac{\partial}{\partial q_n}
+\dot{p_n}\frac{\partial}{\partial p_n}
\Big)
=
\sum_{n=1}^M\big(
H_{p_n}\partial_{q_n}-H_{q_n}\partial_{p_n}
\big)
=
\sum_{n=1}^MX_{H_n}
=\{\,\cdot\,,\,H\,\}
\end{eqnarray*}
That is, $\pounds_{X_{H_n}}\omega_M$ vanishes for each term in the sum
$\pounds_{X_H}\omega_M=\sum_{n=1}^M\pounds_{X_{H_n}}\omega_M$
since $\partial_{q_m}\contract dq_n=\delta_{mn}
=\partial_{p_m}\contract dp_n$
and $\partial_{q_m}\contract dp_n=0
=\partial_{p_m}\contract dq_n$.

\newpage 
\section{Fermat's theorem in geometrical ray optics}
\subsection{Fermat's principle: Rays take paths of least optical length}
In geometrical optics, the ray path is determined by Fermat's principle
of least optical length,
\[
\delta \int n(x,y,z)\,ds = 0
\,.
\]
Here $n(x, y, z)$ is the index of refraction at the
spatial point $(x, y, z)$ and $ds$ is the element of arc
length along the ray path through that point.
Choosing coordinates so that the $z-$axis coincides
with the optical axis (the general direction of propagation), gives
\[
ds
= [(dx)^2 + (dy)^2 + (dz)^2]^{1/2}
= [1+\dot{x}^2 + \dot{y}^2 ]^{1/2}\,dz
\,,
\]
with $\dot{x} = dx/dz$ and $\dot{y}= d y / d z $. Thus, Fermat's
principle can be written in Lagrangian form, with
$z$ playing the role of time, 
\[
\delta \int L(x,y,\dot{x},\dot{y},z)\,dz = 0
\,.
\]
Here, the optical Lagrangian is,
\[
L(x,y,\dot{x},\dot{y},z) 
=  
n(x,y,z)[1+\dot{x}^2 + \dot{y}^2 ]^{1/2}
=:
n/\gamma
\,,
\]
or, equivalently, in two-dimensional vector notation with
$\mathbf{q}=(x,y)$,
\[
L(\mathbf{q},\mathbf{\dot{q}},z) 
=  
n(\mathbf{q},z)[1+|\mathbf{\dot{q}}|^2 ]^{1/2}
=:
n/\gamma
\quad\hbox{with}\quad
\gamma = [1+|\mathbf{\dot{q}}|^2 ]^{-1/2}
\le1
\,.
\]
Consequently, the vector Euler-Lagrange equation
of the light rays is
\[
\frac{d}{ds}\Big(n\frac{d\mathbf{q}}{ds}\Big)
=
\gamma\frac{d}{dz}\Big(n\gamma\frac{d\mathbf{q}}{dz}\Big)
=
\frac{\partial n}{\partial \mathbf{q}}
\,.
\]
The momentum $p$ canonically conjugate to the
ray path position $q$ in an ``image plane'', or on an
``image screen'', at a fixed value of $z$ is given by
\begin{eqnarray*}
\mathbf{p}=\frac{\partial L}{\partial \mathbf{\dot{q}}}
= n\gamma\mathbf{\dot{q}}
&\quad\hbox{which satisfies}\quad&
|\mathbf{p}|^2=n^2(1-\gamma^2)
\,.\\
&\quad\hbox{This implies the velocity}\quad&
\mathbf{\dot{q}}=\mathbf{p}/(n^2-|\mathbf{p}|^2)^{1/2}
\,.
\end{eqnarray*}
Hence, the momentum is real-valued and the Lagrangian is hyperregular,
provided $n^2-|\mathbf{p}|^2>0$. When $n^2=|\mathbf{p}|^2$, the ray
trajectory is vertical and has {\it grazing incidence} with the image
screen.\\

Defining $\sin \theta = d z / d s = \gamma$ leads to
$|\mathbf{p}| = n \cos \theta$, and gives the following geometrical
picture of the ray path. Along the optical axis (the $z-$axis)
each image plane normal to the axis is pierced at
a point $\mathbf{q} = (x, y)$ by a vector of magnitude $n(\mathbf{q}, z)$
tangent to the ray path. This vector makes an angle $\theta$ 
to the plane. The projection of this vector onto
the image plane is the canonical momentum $\mathbf{p}$.
This picture of the ray paths captures all but the
rays of grazing incidence to the image planes.
Such grazing rays are ignored in what follows.\\

Passing now via the usual Legendre transformation from the Lagrangian to
the Hamiltonian description gives
\[
H
=\mathbf{p}\cdot\mathbf{\dot{q}} - L 
=
n\gamma|\mathbf{\dot{q}}|^2 -n/\gamma
= - n \gamma = - \,\big[n(\mathbf{q},z)^2 - |\mathbf{p}|^2\big]^{1/2}
\]
Thus, in the geometrical picture, the component
of the tangent vector of the ray-path along the optical
axis is (minus) the Hamiltonian, i.e. $n(\mathbf{q},z)\sin \theta = - H$.

The phase space description of the ray path now follows from Hamilton's
equations,
\[
\mathbf{\dot{q}}
=\frac{\partial H}{\partial \mathbf{p}}
=\frac{-1}{H}\mathbf{p}
\,,\qquad
\mathbf{\dot{p}}
=
-\,\frac{\partial H}{\partial \mathbf{q}}
=
\frac{-1}{2H}\frac{\partial n^2}{\partial \mathbf{q}}
\,.
\]
\begin{remark} [Translation invariant media] If $n = n(\mathbf{q})$, so
that the medium is translation invariant along the optical axis, $z$, then
$H = - n\sin \theta$ is conserved. (Conservation of $H$ at an interface is
Snell's law.) For translation-invariant media, the vector ray-path equation
simplifies to
\[
\mathbf{\ddot{q}}=-\,
\frac{1}{2H^2}\,
\frac{\partial n^2}{\partial \mathbf{q}}
\,,\quad\hbox{Newtonian dynamics for}\quad
\mathbf{q}\in\mathbb{R}^2
\,.
\]
Thus, in this case geometrical ray tracing reduces
to ``Newtonian dynamics'' in $z$, with potential
$-n^2(\mathbf{q})$ and with ``time'' rescaled along each path
by the value of $\sqrt{2}\,H$ determined from the initial
conditions for each ray.
\end{remark}

\subsection{Axisymmetric, translation invariant materials}
In axisymmetric, translation invariant media, the index of refraction is
a function of the radius alone. Axisymmetry implies an additional constant
of motion and, hence, reduction of the Hamiltonian system for the light
rays to phase plane analysis. For such media, the index of refraction
satisfies
\[
n ( \mathbf{q} , z ) = n ( r )
\,,\qquad
r = |\mathbf{q}|
\,.
\]
Passing to polar coordinates $(r, \phi)$ with $\mathbf{q} = (x, y)
= r(\cos \phi, \sin \phi)$ leads in the usual way to
\[
|\mathbf{p}|^2
=
p_r^2
+  p_\phi^2/r^2
\,.
\]
Consequently, the optical Hamiltonian,
\[
H
= - \,\big[n(r)^2 - p_r^2 -  p_\phi^2/r^2\big]^{1/2}
\]
is independent of the azimuthal angle $\phi$; so its
canonically conjugate ``angular momentum'' $p_\phi$ is
conserved.

Using the relation $\mathbf{q}\cdot\mathbf{p} = rp_r$ leads to an
interpretation of $p_\phi$ in terms of the image-screen
phase space variables $\mathbf{p}$ and $\mathbf{q}$. Namely,
\[
|\mathbf{p}\times\mathbf{q}|^2
=
|\mathbf{p}|^2|\mathbf{q}|^2
-
(\mathbf{p}\cdot\mathbf{q})^2
=
p_\phi^2
\]
The conserved quantity $p_\phi = \mathbf{p}\times\mathbf{q}
= yp_x - xp_y$ is called the skewness function, or the Petzval invariant
for axisymmetric media. Vanishing of $p_\phi$ occurs for {\it meridional
rays}, for which $\mathbf{p}$ and $\mathbf{q}$ are collinear in the image
plane. On the other hand, $p_\phi$ takes its maximum value for {\it
sagittal rays}, for which $\mathbf{p}\cdot\mathbf{q}=0$, so that
$\mathbf{p}$ and
$\mathbf{q}$ are orthogonal in the image plane. 

\begin{exercise}[Axisymmetric, translation invariant materials]
Write Hamilton's canonical equations for axisymmetric, translation
invariant media. Solve these equations for the case of an optical fiber
with radially graded index of refraction in the following form:
\[
n^2(r) = \lambda^2 + (\mu-\nu r^2)^2
\,,\quad\lambda,\,\mu,\,\nu={\rm constants,}
\]
by reducing the problem to phase plane analysis. How does the phase space
portrait differ between $p_\phi=0$ and $p_\phi\ne0$? Show that for
$p_\phi\ne0$ the problem reduces to a Duffing oscillator in a rotating
frame, up to a rescaling of time by the value of the Hamiltonian on each
ray ``orbit.''
\end{exercise}

\subsection{The Petzval invariant and its Poisson bracket relations}
The skewness function 
\[S = p_\phi = \mathbf{p}\times\mathbf{q}= yp_x - xp_y\]
generates rotations of phase space, of $\mathbf{q}$ and $\mathbf{p}$
jointly, each in its plane, around the optical axis.
Its square, $S^2$ (called the Petzval invariant) is conserved for ray
optics in axisymmetric media. That is, $\{S^2,H\}=0$ for optical
Hamiltonians of the form,
\[
H
= - \,\big[n(|\mathbf{q}|^2)^2 - |\mathbf{p}|^2\big]^{1/2}
\,.
\]
We define the axisymmetric invariant coordinates by the map
$T^*\mathbb{R}^2\mapsto\mathbb{R}^3$ $(\mathbf{q},\mathbf{p})\mapsto
(X,Y,Z)$,
\[
X=|\mathbf{q}|^2\ge0
\,,\quad
Y=|\mathbf{p}|^2\ge0
\,,\quad
Z=\mathbf{p}\cdot\mathbf{q}
\,.
\]
The following Poisson bracket relations hold
\[
\{S^2,X\}=0
\,,\quad
\{S^2,Y\}=0
\,,\quad
\{S^2,Z\}=0
\,,
\]
since rotations preserve dot products.
In terms of these invariant coordinates,
the Petzval invariant and optical Hamiltonian satisfy
\[
S^2=XY-Z^2\ge0
\,,\quad\hbox{and}\quad
H^2=n^2(X)-Y\ge0
\,.
\]
The level sets of $S^2$ are hyperboloids of revolution around 
the $X = Y$ axis, extending up through the interior of the $S = 0$
cone, and lying between the $X-$ and $Y-$axes. The level sets of $H^2$
depend on the functional form of the index of refraction, but they are
$Z-$independent.

\subsection{$\mathbb{R}^3$ Poisson bracket for ray optics}
The Poisson brackets among the axisymmetric
variables $X$, $Y$ and $Z$ close among themselves,
\[
\{X,Y\}=4Z
\,,\quad 
\{ Y , Z \} = - 2 Y 
\,,\quad 
\{ Z , X \} = - 2 X 
\,.
\]
These Poisson brackets derive from a single $\mathbb{R}^3$  
Poisson bracket for $\mathbf{X}=(X,Y,Z)$ given by
\[
\{F,H\}=-\nabla{S^2}\cdot\nabla{F}\times\nabla{H}
\]
Consequently, we may re-express the equations of
Hamiltonian ray optics in axisymmetric media
with $H = H(X, Y )$ as
\[
\mathbf{\dot{X}}=\nabla{S^2}\times\nabla{H}
\,.
\]
with Casimir $S^2$, for which $\{ S^2 , H \} = 0$, for every $H$.
Thus,  the flow preserves volume (div$\,\mathbf{\dot{X}}=0$) and the
evolution takes place on intersections of level surfaces of the
axisymmetric media invariants $S^2$ and $H(X,Y)$. 

\subsection{Recognition of the Lie-Poisson bracket for ray optics}

The Casimir invariant $S^2=XY-Z^2$ is quadratic. In such cases, one may
write the $\mathbb{R}^3$ Poisson bracket in the suggestive form, 
\[
\{F,H\}
=
-\,C^k_{ij}X_k
\frac{\partial F}{\partial X_i}
\frac{\partial H}{\partial X_j}
\]
In this particular case, $C^3_{12}=4$, $C^2_{23}=2$ and $C^1_{31}=2$ and
the rest either vanish, or are obtained from antisymmetry of $C^k_{ij}$
under exchange of any pair of its indices. These values are the structure
constants of any of the Lie algebras $sp(2,\mathbb{R})$, $so(2,1)$,
$su(1,1)$, or $sl(2,\mathbb{R})$. Thus, the reduced description of
Hamiltonian ray optics in terms of axisymmetric $\mathbb{R}^3$ variables
is said to be ``Lie-Poisson'' on the dual space of any of these Lie
algebras, say, $sp(2,\mathbb{R})^*$ for definiteness. We will have more to
say about Lie-Poisson brackets later, when we reach the Euler-Poincar\'e
reduction theorem. 

\begin{exercise}
Consider the $\mathbb{R}^3$ Poisson bracket
\begin{equation}\label{r3pb-ex}
\{f, h\} = -\,\nabla{c}\cdot\nabla{f}\times\nabla{h}
\end{equation}
Let  $c=\mathbf{x}^T\cdot\mathbb{C}\mathbf{x}$ be a quadratic form on
$\mathbb{R}^3$, and let $\mathbb{C}$ be the associated symmetric
$3\times3$ matrix. Show that this is the Lie-Poisson bracket for the Lie
algebra structure
\[
[\mathbf{u},\mathbf{v}]_\mathbb{C}
=
\mathbb{C}(\mathbf{u}\times\mathbf{v})
\]
What is the underlying matrix Lie algebra? What are the coadjoint orbits
of this Lie algebra?
\end{exercise}

\begin{remark}[Coadjoint orbits]
As one might expect, the coadjoint orbits of the group $SP(2,\mathbb{R})$
are the hyperboloids corresponding to the level sets of $S^2$.
\end{remark}

\begin{remark}
As we shall see later, the map $T^*\mathbb{R}^2\mapsto sp(2,\mathbb{R})^*$ taking 
$(\mathbf{q},\mathbf{p})\mapsto(X,Y,Z)$ is an example of a {\bfi momentum map}.
\end{remark}


\section{Geometrical Structure of Classical Mechanics}\label{smooth manifold}

\subsection{Manifolds}

Configuration space: coordinates $q\in M,$ where $M$ is a smooth manifold.

$\phi_\beta \circ \phi^{-1}_\alpha$ is a smooth change of variables.

For later, smooth coordinate transformations: $q\to Q$ with
$dQ=\frac{\partial Q}{\partial q} dq$

\begin{definition}
A smooth manifold $M$ is a set of points together with a finite (or
perhaps countable) set of subsets  $U_\alpha\subset M$ and 1-to-1 mappings
$\phi_\alpha:U_\alpha \to \mathbb{R}^n$ such that
\begin{enumerate}
\item \ $\bigcup_\alpha U_\alpha = M$
\item For every nonempty intersection $U_\alpha \cap U_\beta,$ the set 
$\phi_\alpha \left(U_\alpha \cap U_\beta\right)$ is an open subset of
$\mathbb{R}^n$ and the 1-to-1 mapping $\phi_\beta \circ \phi_\alpha^{-1}$
is a smooth function on  $\phi_\alpha \left(U_\alpha \cap U_\beta\right).$ 
\end{enumerate}
\end{definition}

\begin{remark}
The sets $U_\alpha$ in the definition are called {\bfi coordinate charts}.
The mappings $\phi_\alpha$ are called {\bfi coordinate functions} or {\bfi
local coordinates}. A collection of charts satisfying 1 and 2 is called
an {\bfi atlas}. Condition 3 allows the definition of manifold to be made
independently of a choice of atlas. A set of charts satisfying 1 and 2 can
always be extended to a maximal set; so, in practice, conditions 1 and 2
define the manifold.
\end{remark}

\begin{example}
Manifolds often arise as intersections of zero level sets 
\[M=\left\{x \big| f_i(x) = 0, \ i=1,\dots,k\right\},\] for a given
set of functions
$f_i: \mathbb{R}^n\to\mathbb{R},$ $i=1,\dots,k$.\\ If the gradients
$\nabla f_i$ are linearly independent, or more generally if the rank of
$\left\{\nabla f(x)\right\}$ is a constant $r$ for all $x,$  then $M$ is
a smooth manifold of dimension $n-r.$ The proof uses the Implicit Function
Theorem to show that an $(n-r)-$dimensional coordinate chart may be
defined in a neighborhood of each point on $M$. In this situation, the
set $M$ is called a {\bfi submanifold} of $\mathbb{R}^n$ (see
\cite{Le2003}).
\end{example}

\begin{definition} If $r=k,$ then the map $\left\{f_i\right\}$ is called a
{\bfi submersion}. 
\end{definition}

\begin{exercise} 
Prove that all submersions are submanifolds (see \cite{Le2003}).
\end{exercise}

\begin{definition}
[Tangent space to level sets]
Let \[M=\left\{x \big| f_i(x) = 0, \ i=1,\dots,k\right\}\] be a manifold in
$\mathbb{R}^n.$ The {\bfi tangent space} at each $x\in M,\ $ is defined by 
\[T_xM=\left\{ v\in \mathbb{R}^n \ \big| \ \frac{\partial f_i}{\partial
x^a} (x) v^a = 0,\ i=1,\dots,k\right\}.\]  Note: we use the {\bfi summation
convention}, i.e. repeated indices are summed over their range.
\end{definition}

\begin{remark} 
The tangent space is a linear vector space.
\end{remark} 

\begin{example}[Tangent space to the sphere in $\mathbb{R}^3$]
\end{example}

\begin{example}[Tangent space to the sphere in $\mathbb{R}^3$]
The sphere $S^2$ is the set of points $(x,y,z)\in\mathbb{R}^3$ solving
$x^2 + y^2 + z^2 = 1$.  The tangent space to the sphere at such a point
$(x,y,z)$ is the plane containing vectors
$(u,v,w)$ satisfying $xu + yv + zw = 0$.
\end{example}

\begin{definition}\label{tanbun-def}
[Tangent bundle]
The {\bfi tangent bundle} of a manifold $M$, denoted by $TM$, is the
smooth manifold whose underlying set is the disjoint union of the tangent
spaces to $M$ at the points $x\in M$; that is, 
\[TM = \bigcup_{x\in M}\,T_xM\]
Thus, a point of $TM$ is a vector $v$ which is tangent to $M$ at some point
$x\in M$.
\end{definition}

\begin{example}[Tangent bundle $TS^2$ of $S^2$]
The tangent bundle $TS^2$ of $S^2\in \mathbb{R}^3$ is the union of the
tangent spaces of $S^2$:\\
$TS^2 = \left\{ (x,y,z;u,v,w) \in \mathbb{R}^6 \ 
\big| \ x^2 + y^2 + z^2 = 1 
\textrm{ and } xu + yv + zw = 0 \right\}$.\\ 
\end{example}

\begin{remark} [Dimension of tangent bundle $TS^2$]
Defining $TS^2$ requires two independent conditions in
$\mathbb{R}^6$; so {\rm dim}$TS^2=4$.
\end{remark} 

\begin{exercise}
Define the sphere $S^{n-1}$ in $\mathbb{R}^n$. What is the dimension of
its tangent space $TS^{n-1}$?
\end{exercise}

\begin{example}[The two stereographic projections of $S^2\to\mathbb{R}^2$]\index{stereographic projection}\index{Riemann sphere}
\label{2StereoProj-eg}$\quad$\\
The unit sphere 
\[
S^2=\{(x,y,z):\, x^2+y^2+z^2=1\}
\]
is a smooth two-dimensional manifold realized as a submersion in $\mathbb{R}^3$. Let
\[
U_N = S^2\backslash \{0,0,1\}
\,,\quad\hbox{and}\quad
U_S = S^2\backslash \{0,0,-1\}
\]
be the subsets obtained by deleting the North and South poles of $S^2$, respectively. Let 
\[
\chi_N:\,U_N\to (\xi_N,\eta_N)\in \mathbb{R}^2
\,,\quad\hbox{and}\quad
\chi_S:\,U_S\to (\xi_S,\eta_S)\in \mathbb{R}^2
\]
be stereographic projections from the North and South poles onto the equatorial plane, $z=0$. 

\rem{
\begin{figure} [h]

\centering 

\includegraphics[width=10cm]{./Figures/Riemann3D.pdf} 
 \caption{Two Riemann projections from the North pole of the unit sphere onto the $z=0$ plane with coordinates $(\xi_N,\,\eta_N)$. 
} 
 \label{fig:} 
 \end{figure} 

\begin{figure} [h]

\centering 

\includegraphics[width=10cm]{./Figures/Riemann_Sphere2.pdf} 
 \caption{Riemann projection in the $\xi_N-z$ plane at fixed azimuth $\phi=0$. The projection through $\hat{z}=(\sin\theta,\,\cos\theta)$ strikes the $\xi_N-$axis at distance $r=\cot(\theta/2)$.
} 
 \label{fig:Riemanncircle} 
 \end{figure} 
}

\noindent
Thus, one may place two different coordinate patches in $S^2$ intersecting everywhere except at the points along the $z-$axis at $z=1$ (North pole) and $z=-1$ (South
pole).\\

\noindent
In the equatorial plane $z=0$, one may define two sets of
(right-handed) coordinates, 
\[
\phi_\alpha:\,U_\alpha\to\mathbb{R}^2\backslash\{0\}
\,,\quad
\alpha=N,S\,,
\]
obtained by the following two stereographic projections from the North and South poles:
\\  (1)  (valid everywhere except $z=1$)
\[
\phi_N(x,y,z)=(\xi_N,\,\eta_N)
=\left(\frac{x}{1-z},\,\frac{y}{1-z}\right)
\,,\]
(2) (valid everywhere except $z=-1$)
\[
\phi_S(x,y,z)=(\xi_S,\,\eta_S)
=\left(\frac{x}{1+z},\,\frac{-\,y}{1+z}\right)
\,.\]
(The two complex planes are identified differently with the plane $z = 0$. An orientation-reversal is necessary to maintain consistent coordinates on the sphere.) 

One may check directly that on the overlap $U_N\cap U_S$ the map,
\[
\phi_N\circ\phi_S^{-1}:\, \mathbb{R}^2\backslash\{0\}
\to \mathbb{R}^2\backslash\{0\}
\]
is a smooth diffeomorphism, given by the inversion
\[
\phi_N\circ\phi_S^{-1} (x,\,y)
=
\Big(\frac{x}{x^2+y^2},\,\frac{y}{x^2+y^2}\Big)
\,.
\] 
\end{example}

\begin{exercise}
Construct the mapping from
$(\xi_N,\eta_N)\to(\xi_S,\eta_S)$ and verify that it is a diffeomorphism in  $\mathbb{R}^2\backslash\{0\}$.
Hint: $(1+z)(1-z)=1-z^2=x^2+y^2$.
\end{exercise}

\begin{answer}
\[
(\xi_S,-\,\eta_S)=\frac{1-z}{1+z}\,(\xi_N,\eta_N)
=\frac{1}{\xi_N^2+\eta_N^2}\,(\xi_N,\eta_N)
\,.
\]
The map $(\xi_N,\eta_N)\to(\xi_S,\eta_S)$ is smooth and invertible
except at $(\xi_N,\eta_N)=(0,0)$.
\end{answer}

\begin{example}
If we start with two identical circles in the $xz$-plane, of radius $r$ and
centered at $x=\pm 2r$, then rotate them round the $z$ axis in
$\mathbb{R}^3$, we get a torus, written $T^2.$ It's a manifold.
\end{example}

\begin{exercise}
If we begin with a figure eight in the $xz$-plane, along the $x$ axis and
centered at the origin, and spin it round the $z$ axis in $\mathbb{R}^3,$
we get a ``pinched surface'' that looks like a sphere that has been
``pinched'' so that the north and south poles touch.  Is this a manifold?
Prove it.
\end{exercise}

\begin{answer}
The origin has a neighbourhood diffeomorphic to a double cone. This is
not diffeomorphic to $\mathbb{R}^2.$ A proof of this is that, if the
origin of the cone is removed, two components remain; while if the origin
of $\mathbb{R}^2$ is removed, only one component remains.
\end{answer}

\begin{remark}
The sphere will appear in several examples as a reduced space in which
motion takes place after applying a symmetry. Reduction by symmetry is
associated with a classical topic in  celestial mechanics known as normal
form theory. Reduction may be ``singular,'' in which case it leads to
``pointed'' spaces that are smooth manifolds except at one or more points.
For example different resonances of coupled oscillators correspond to the
following reduced spaces: 1:1 resonance -- sphere; 1:2 resonance --
pinched sphere with one cone point;  1:3 resonance -- pinched sphere with
one cusp point; 2:3 resonance -- pinched sphere with one cone point and 
one cusp point.
\end{remark}

\subsection{Motion: Tangent Vectors and Flows}

Envisioning our later considerations of dynamical systems, we shall consider
motion along curves $c(t)$ parameterised by time $t$ on a smooth manifold $M.$
Suppose these curves are trajectories of a flow $\phi_t$ of a vector field. We
anticipate this means 
$\phi_t\left(c(0)\right) = c(t)$ and 
$\phi_t \circ \phi_s = \phi_{t+s}$ (flow property). The flow will be
tangent to $M$ along the curve. To deal with such flows, we will need the
concept of {\bfi tangent vectors}.

Recall from {\bfi Definition \ref{tanbun-def}} that 
the tangent bundle of $M$ is \[TM=\bigcup\limits_{x\in M} T_xM.\]
We will now add a bit more to that definition. The tangent bundle is an
example of a more general structure than a manifold. 

\begin{definition}[Bundle]
A {\bfi bundle} consists of a manifold $B$, another manifold $M$ called
the ``base space'' and a projection between them $\Pi:\,B\to M$. Locally,
in small enough regions of $x$ the inverse images of the projection $\Pi$
exist. These are called the {\bfi fibers} of the bundle. Thus, subsets of
the bundle $B$ locally have the structure of a Cartesian product. An
example is $(B,M,\Pi)$ consisting of
$(\mathbb{R}^2,\mathbb{R}^1,\Pi:\mathbb{R}^2\to\mathbb{R}^1)$. In this
case, $\Pi:(x,y)\in\mathbb{R}^2\to x\in\mathbb{R}^1$. Likewise, the
tangent bundle consists of $M,TM$ and a map $\tau_M:TM\to M$.
\end{definition}

Let $x = \left(x^1,\dots,x^n\right)$ be local coordinates on $M$, and let
$v = \left(v^1,\dots,v^n\right)$ be components of a tangent vector.
\[T_xM = \left\{ v \in \mathbb{R}^n \ \big| \ \frac{\partial f_i}{\partial
x} \cdot v = 0, i=1,\dots, m \right\}\] 
for 
\[M = \left\{ x \in
\mathbb{R}^n \ \big| \ f_i(x) = 0, i=1,\dots, m \right\}\]
These $2n$ numbers $(x,v)$ give local coordinates on $TM$, whose
dimension is  $\dim TM = 2 \dim M$. The {\bfi tangent bundle projection}
is a map $\tau_M : TM \to M$ which takes a tangent vector $v$ to a point
$x\in M$ where the tangent vector $v$ is attached (that is, $v\in T_x
M$). The inverse of this projection $\tau_M^{-1}(x)$ is called the {\bfi
fiber} over $x$ in the tangent bundle.

\subsection*{Vector fields, integral curves and flows}

\begin{definition}\label{VecField-defn}
A {\bfi vector field} on a manifold $M$ is a map $X:M\to TM$ that assigns
a vector $X(x)$ at each point
$x\in M.$ This implies that $\tau_M \circ X = Id.$
\end{definition}

\begin{definition}
An {\bfi integral curve} of $X$ with initial conditions $x_0$ at $t=0$ is a
differentiable map $c:]a,b[ \to M,$ where $]a,b[$ is an open interval
containing $0,$ such that $c(0)=0$ and
$c\,'(t) = X\left(c(t)\right)$ for all $t\in ]a,b[.$
\end{definition}

\rem{
T's note: I believe that Darryl said that the integral curves being
differentiable implies that $X$ is differentiable, but I don't see why.

D's reply: I hope I said that $X$ being Lipschitz implies
the integral curves are differentiable. I recall saying we wouldn't
be discussing solution properties of ODEs.}

\begin{remark}
A standard result from the theory of ordinary differential equations
states that $X$ being Lipschitz implies its integral curves are unique
and $C^1$ \cite{CoLe1984}. The integral curves $c(t)$ are differentiable for
smooth $X$.
\end{remark}

\subsection{Summary}

\begin{definition}
The {\bfi flow} of $X$ is the collection of maps $\phi_t:M\to M$, where
$t\to\phi_t(x)$ is the integral curve  of $X$ with initial condition $x$.
\end{definition}

\begin{remark}$\quad$

\begin{enumerate}
\item Existence and uniqueness results for solutions of $c\,'(t)=X(c(t))$
guarantee that flow $\phi$ of $X$  is smooth in $(x,t),$ for smooth $X$.
\item 
Uniqueness implies the flow property 
\[\phi_{t+s}=\phi_t \circ
\phi_s,
\hspace{2cm} (FP)
\]
for initial condition $\phi_0=Id.$
\item The flow property (FP) generalizes to the nonlinear case the
familiar linear situation where $M$ is a vector space, $X(x)=Ax$ is a
linear vector field for a bounded linear operator $A$, and
$\phi_t(x)=e^{At} x$.
\end{enumerate}
\end{remark}

\subsection*{Differentials of functions and the cotangent bundle}

We are now ready to define differentials of smooth functions and the
cotangent bundle.\\

Let $f:M\to\R$ be a smooth function. We differentiate
$f$ at $x\in M$ to obtain $T_xf:T_xM\to T_{f(x)} \R.$
As is standard, we identify $T_{f(x)}\R$ with $\R$ itself, thereby
obtaining a linear map $df(x):T_xM\to \R.$ The result $df(x)$ is an
element of the cotangent space $T_x^*M$, the dual space of the tangent space
$T_xM$. The natural pairing between elements of the tangent space and the
cotangent space is denoted as $\langle\cdot\,,\,\cdot\rangle:T_x^*M\times
T_xM\mapsto\mathbb{R}$.

In coordinates, the linear map $df(x):T_xM\to \R$ may be written as the
directional derivative,
\[
\langle df(x)\,,\,v\rangle
=
df(x)\cdot v = \frac{\partial f}{\partial x^i}
\cdot v^i
\,,
\] 
for all $v\in T_xM$.
(Reminder: the summation convention is intended over repeated indices.) Hence,
elements $df(x)\in T_x^*M$ are dual to vectors $v\in T_xM$ with respect to the
pairing $\langle\cdot\,,\,\cdot\rangle$.

\begin{definition}
$df$ is the {\bfi differential} of the function $f$.
\end{definition}

\begin{definition}
The dual space of the tangent bundle $TM$ is the {\bfi cotangent bundle}
$T^*M$. That is, 
\[
(T_xM)^*=T_x^*M
\quad\hbox{and}\quad
T^*M=\bigcup_xT_x^*M
\,.
\]
\end{definition}
Thus, replacing $v\in T_xM$ with $df\in T_x^*M,$ for all $x\in M$ and for
all smooth functions $f:M\to \R,$ yields the {\bfi cotangent bundle}
$T^*M.$\\

\paragraph{Differential bases} When the basis of vector fields is denoted
as $\frac{\partial}{\partial x^i}$ for $i=1,\dots,n$, its dual basis is
often denoted as $dx^i$. In this notation, the differential of a function
at a point $x\in M$ is expressed as
\[
df(x)=\frac{\partial f}{\partial x^i} dx^i
\]

The corresponding pairing $\left<\cdot\,,\cdot\right>$ of bases is
written in this notation as
\[\left< dx^j, \frac{\partial}{\partial x^i}\right> = \delta_i^j\] 
Here $\delta_i^j$ is the Kronecker delta, which equals unity for $i=j$ and
vanishes otherwise. That is, defining $T^*M$ requires a pairing
$\left<\cdot\,,\cdot\right>:T^*M\times TM\to \R.$

(Different pairings exist for curvilinear coordinates, Riemannian
manifolds, etc.)

\section{Derivatives of differentiable maps -- the tangent lift}

We next define derivatives of differentiable maps between manifolds
(tangent lifts).

We expect that a smooth map $f:U\to V$ from a chart $U\subset M$ to a chart
$V\subset N,$  will lift to a map between the tangent bundles $TM$ and
$TN$ so as to make sense from the viewpoint of ordinary calculus,
\begin{align*}
U\times \R^m\subset TM 
&\longrightarrow V\times \R^n\subset TN \\
\left(q^1,\dots, q^m; X^1,\dots, X^m\right) 
&\longmapsto \left(Q^1,\dots, Q^n; Y^1,\dots, Y^n\right)
\end{align*}
Namely, the relations between the vector field components should be
obtained from the differential of the map $f:U\to V$. Perhaps not
unexpectedly, these vector field components will be related by
\[
Y^i\frac{\partial }{\partial Q^i}
=
X^j\frac{\partial }{\partial q^j}
\,,\quad\hbox{so}\quad
Y^i = \frac{\partial Q^i}{\partial q^j} X^j
\,,
\] 
in which the quantity called the {\bfi tangent lift} 
\[
Tf=\frac{\partial Q}{\partial q}
\]
of the function $f$
arises from the chain rule and is equal to the Jacobian for the
transformation $Tf:TM\mapsto TN$. \\

The dual of the tangent lift is the cotangent lift, explained later in section
\ref{}. Roughly speaking, the {\bfi cotangent lift} of the function $f$, 
\[
T^*f=\frac{\partial q}{\partial Q}
\]
arises from
\[
 \beta_i dQ^i = \alpha_j dq^j
\,,\quad\hbox{so}\quad
\beta_i=\alpha_j\frac{\partial q^j}{\partial Q^i}
\]
and $T^*f:T^*N\mapsto T^*M$. Note the directions of these maps:
\begin{eqnarray*}
&&
Tf
\,:\quad
q\,,X\in TM \mapsto Q,Y\in TN
\\&&
f
\,:\quad
q\in M \mapsto Q\in N
\\&&
T^*f
\,:\quad
Q\,,\beta\in T^*N \mapsto  q\,,\alpha\in T^*M
\quad\hbox{(map goes the other way, see the \fbox{picture})}
\end{eqnarray*}

\subsection{Summary remarks about derivatives on manifolds}

\begin{definition}[Differentiable map]
A map $f:\,M\to N$ from manifold $M$ to manifold $N$ 
is said to be {\bfi differentiable} (resp. $C^k$) if it is represented in
local coordinates on $M$ and $N$ by differentiable (resp. $C^k$)
functions. 
\end{definition}

\begin{definition}[Derivative of a differentiable map]
The {\bfi derivative} of a differentiable map
\[f:\,M\to N\]
at a point $x\in M$ is defined to be the linear map
\[T_xf:\,T_xM\to T_xN\]
constructed, as follows. For $v\in T_xM$, choose a curve $c(t)$ that maps an open
interval $t\in(-\epsilon,\epsilon\,)$ around the point $t=0$ to the manifold $M$
\begin{eqnarray*}
&&c:\ (-\epsilon,\epsilon\,)\,\longrightarrow M
\\
&&\hbox{with}\quad c(0)=x
\\
&&\hbox{and velocity vector}\quad
c\,'(0):= \frac{dc}{dt}\Big|_{t=0}=v
\,.
\end{eqnarray*}
Then $T_xf\cdot v$ is the velocity vector at $t=0$ of the curve $f\circ
c:\,\mathbb{R}\to N$. That is,
\[
T_xf\cdot v 
= \frac{d}{dt}f(c(t))\Big|_{t=0}
= \frac{\partial f}{\partial c}\frac{d}{dt}c(t)\Big|_{t=0}
\]
\end{definition}
\begin{definition}\label{tang-lift-defn}
The union $Tf=\bigcup_xT_xf$ of the derivatives $T_xf:\,T_xM\to T_xN$
over points $x\in M$ is called the {\bfi tangent lift} of the map
$f:\,M\to N$.
\end{definition}

\begin{remark}
The chain-rule definition of the derivative $T_xf$ of a differentiable
map at a point $x$ depends on the function $f$ and the vector $v$.  Other
degrees of differentiability are possible. For example, if
$M$ and
$N$ are manifolds and
$f:\,M\to N$ is of class $C^{k+1}$, then the tangent lift (Jacobian)
$T_xf:\,T_xM\to T_xN$ is $C^k$.
\end{remark}

\begin{exercise}
Let $\phi_t : S^2 \to S^2$ rotate points on $S^2$ about a fixed axis
through an angle $\psi(t)$. Show that $\phi_t$ is the flow of a certain
vector field on $S^2$.
\end{exercise}

\begin{exercise}
Let $f : S^2 \to\mathbb{R}$ be defined by $f(x, y, z) = z$. Compute $df$
using spherical coordinates $(\theta,\phi)$.
\end{exercise}

\begin{exercise}
Compute the tangent lifts for the two stereographic projections of
$S^2\to\mathbb{R}^2$ in example \ref{2StereoProj-eg}. That is, assuming
$(x,y,z)$ depend smoothly on
$t$, find
\begin{enumerate}
\item
How $(\dot{\xi}_N,\dot{\eta}_N)$ depend on $(\dot{x},\dot{y},\dot{z})$.
Likewise, for $(\dot{\xi}_S,\dot{\eta}_S)$.
\item
How $(\dot{\xi}_N,\dot{\eta}_N)$ depend on $(\dot{\xi}_S,\dot{\eta}_S)$.
\end{enumerate}
Hint: Recall $(1+z)(1-z)=1-z^2=x^2+y^2$ and
use $x\dot{x}+y\dot{y}+z\dot{z}=0$ when
$(\dot{x},\dot{y},\dot{z})$ is tangent to $S^2$ at $(x,y,z)$.
\end{exercise}

\rem{
\begin{exercise}
Show that $(\dot{x},\dot{y},\dot{z})$ being tangent to $S^2$ at $(x,y,z)$
implies that $(\dot{\xi}_N,\dot{\eta}_N)$ is tangent to
$\mathbb{R}^2$ at $({\xi}_N,{\eta}_N)$  and that $(\dot{\xi}_S,\dot{\eta}_S)$ is
tangent to $\mathbb{R}^2$ at $({\xi}_S,{\eta}_S)$. Thus, differentiable
maps preserve the property of tangency. 
\end{exercise}
}

\section{Lie groups and Lie algebras}

\subsection{Matrix Lie groups}

\begin{definition}
A {\bfi group} is a set of elements with 
\begin{enumerate}
\item
A binary product (multiplication), $G \times G \to G$
\begin{description}
\item
The product of $g$ and $h$ is written $gh$
\item
The product is associative, $(gh)k=g(hk)$
\end{description}
\item
Identity element $e:\,eg=g$ and $ge=g$, $\forall g\in G$
\item
Inverse operation $G \to G$, so that 
$gg^{-1}=g^{-1}g=e$
\end{enumerate}
\end{definition}

\begin{definition}
A {\bfi Lie group} is a smooth manifold $G$ which is a group
and for which the group operations of multiplication, 
$(g,h) \to gh$ for $g,h \in G$, and inversion, $g \to g^{-1}$ with
$gg^{-1}=g^{-1}g=e$, are smooth.
\end{definition}

\begin{definition}\label{MatrixLG-def}
A {\bfi matrix Lie group} is a set of invertible $n\times n$ matrices
which is closed under matrix multiplication and which is a submanifold of
$\mathbb{R}^{n\times n}$.
The conditions showing that a matrix Lie group is a Lie group are easily
checked:
\begin{description}
\item
A matrix Lie group is a manifold, because it is a submanifold of
$\mathbb{R}^{n\times n}$
\item
Its group operations are smooth, since they are algebraic operations on
the matrix entries.
\end{description}

\end{definition}

\begin{example}[The general linear group $GL(n,\mathbb{R})$]
The matrix Lie group $GL(n,\mathbb{R})$ is the group of linear
isomorphisms of $\mathbb{R}^n$ to itself.  The dimension of the matrices
in $GL(n,\mathbb{R})$ is $n^2$.
\end{example}

\begin{proposition}\label{Smap-prop}
Let $K \in GL(n,\mathbb{R})$ be a symmetric matrix,
$K^T=K$. Then the subgroup $S$ of $GL(n,\mathbb{R})$ defined by the
mapping 
\[S = \{U \in GL(n,\mathbb{R})|\,U^TKU = K\}\]
is a submanifold of $\mathbb{R}^{n\times n}$ of dimension 
$n(n-1)/2$.
\end{proposition}

\begin{remark}
The subgroup $S$ leaves invariant a certain symmetric quadratic form
under linear transformations,
$S\times\mathbb{R}^n\to\mathbb{R}^n$ given by
$\mathbf{x}\to U\mathbf{x}$, since
\[
\mathbf{x}^T K\mathbf{x}
=
\mathbf{x}^T U^TKU\mathbf{x}
\,.
\]
So the matrices $U \in S$ change the basis for this
quadratic form, but they leave its value unchanged. Thus, $S$ is the {\bfi
isotropy subgroup} of the quadratic form associated with $K$. 
\end{remark}

\paragraph{Proof.} 
\begin{description}
\item
Is $S$ a subgroup? We check the following three defining properties
\begin{enumerate}
\item
Identity:\\
$I\in S$ because $I^TKI=K$
\item
Inverse:\\
$U\in S\Longrightarrow U^{-1}\in S$, because
\\
$K=U^{-T}(U^TKU)U^{-1}=U^{-T}(K)U^{-1}$
\item
Closed under multiplication:\\
$U, V\in S\Longrightarrow UV\in S$, because \\
$(UV)^TKUV=V^T(U^TKU)V=V^T(K)V=K$
\end{enumerate}
\item
Hence, $S$ is a subgroup of $GL(n,\mathbb{R})$. 

\item
Is $S$ is a submanifold of $\mathbb{R}^{n\times n}$  of dimension 
$n(n-1)/2$? 
\begin{description}
\item
Indeed, $S$ is the zero locus of the mapping $UKU^T - K$. This makes
it a submanifold, because it turns out to be a submersion. 
\item
For a submersion, the dimension of the level set is the dimension of the
domain minus the dimension of the range space. In this case, this
dimension is $n^2-n(n+1)/2=n(n-1)/2$
\end{description}
\end{description}

\begin{exercise}
Explain why one can conclude that the zero locus map for $S$ is a
submersion. In particular, pay close attention to establishing the
constant rank condition for the linearization of this map.
\end{exercise}

\paragraph{Solution} Here is why $S$ is a submanifold of $R^{n\times n}$.\\

First, $S$ is the zero locus of the mapping 
\[
U\to U^TKU - K
\,,\qquad\hbox{(locus map)}
\] 
Let $U\in S$, and let $\delta{U}$ be an arbitrary element of $R^{n\times
n}$. Then linearize to find
\[
(U + \delta{U})^TK(U + \delta{U}) - K =
U^TKU - K + \delta{U}^TKU + U^TK\delta{U} + O(\delta{U})^2
,.
\]
We may conclude that $S$ is a submanifold of $R^{n\times n}$ if we can
show that the linearization of the locus map, namely the linear mapping 
defined by
\[
L\equiv \delta{U}\to \delta{U}^TKU + U^TK\delta{U}
\,,\qquad R^{n\times n}\to R^{n\times n}
\] 
has constant rank for all $U\in S$.

\begin{lemma} The linearization map $L$ is onto the space of ${n\times n}$
of symmetric matrices and hence the original map is a submersion.
\end{lemma}

\paragraph{Proof that $L$ is onto.} 
\begin{itemize}
\item
Both the original locus map and the image of $L$ lie in the subspace
of ${n\times n}$ symmetric matrices. 
\item
Indeed, given $U$ and any symmetric
matrix
$S$ we can find $\delta{U}$ such that 
\[
 \delta{U}^TKU + U^TK\delta{U} = S
\,.
\] 
Namely 
\[
 \delta{U}=K^{-1}U^{-T}S/2
\,.
\] 
\item
Thus, the linearization map $L$ is onto the space of ${n\times n}$
of symmetric matrices and the original locus map $U\to UKU^T - K$ to the
space of symmetric matrices is a submersion. 
\end{itemize}

For a submersion,
the dimension of the level set is the dimension of the domain minus the
dimension of the range space. In this case, this dimension is $n^2 - n(n +
1)/2 = n(n - 1)/2$.

\begin{corollary}[$S$ is a matrix Lie group]
$S$ is both a subgroup and a submanifold of the general linear group
$GL(n,\mathbb{R})$. Thus, by definition \ref{MatrixLG-def}, $S$ is
a matrix Lie group.
\end{corollary}

\begin{exercise}
What is the tangent space to $S$ at the identity, $T_IS$? 
\end{exercise}

\begin{exercise}
Show that for any pair of matrices $A,B\in T_IS$, the matrix
commutator $[A,B]\equiv AB-BA \in T_IS$.
\end{exercise}

\begin{proposition} 
The linear space of matrices $A$ satisfying 
\[
A^TK+KA=0
\]
defines $T_IS$, the tangent space at the identity
of the matrix Lie group $S$ defined in Proposition \ref{Smap-prop}.
\end{proposition}

\paragraph{Proof.} Near the identity the defining condition for $S$
expands to 
\[
(I+\epsilon A^T + O(\epsilon^2))K
(I+\epsilon A + O(\epsilon^2)) = K
\,,\quad\hbox{for}\quad \epsilon\ll1
\,.
\]
At linear order $O(\epsilon)$ one finds,
\[
A^TK + KA =0
\,.\]
This relation defines the linear space of matrices $A\in T_IS$.
\\

\noindent
If $A,B\in T_IS$, does it follow that $[A,B]\in T_IS$?\\

\noindent
Using $[A,B]^T=[B^T,A^T]$, we check closure by a direct computation,
\begin{eqnarray*}
[B^T,A^T]K+K[A,B]&=&B^TA^TK-A^TB^TK+KAB-KBA
\\&=&B^TA^TK-A^TB^TK-A^TKB+B^TKA=0
\,.
\end{eqnarray*}
Hence, the tangent space of $S$ at the identity $T_IS$ is closed under the
matrix commutator $[\cdot\,,\,\cdot]$. 

\begin{remark}
In a moment, we will show that the matrix commutator for $T_IS$ also
satisfies the Jacobi identity. This will imply that the condition $A^TK +
KA =0$ defines a matrix Lie algebra.
\end{remark}

\subsection{Defining Matrix Lie Algebras}
We are ready to prove the following, in preparation for defining matrix
Lie algebras. 

\begin{proposition} \label{matrixLG-comprop}
Let $S$ be a matrix Lie group, and let $A,B \in T_IS$
(the tangent space to $S$ at the identity element). Then $AB - BA \in
T_IS$.
\end{proposition}
The proof makes use of a lemma.
\begin{lemma} Let $R$ be an arbitrary element of a matrix Lie group $S$,
and let $B \in T_IS$. Then $RBR^{-1}\in T_IS$.
\end{lemma}

\paragraph{Proof of lemma.} 
Let $R_B(t)$ be a curve in $S$ such that $R_B(0) = I$ and $R^\prime(0)
= B$. Define $S(t) = RR_B(t)R^{-1}\in T_IS$ for all $t$. Then $S(0) = I$
 and $S^\prime(0)=RBR^{-1}$. Hence, $S^\prime(0) \in T_IS$, thereby proving
the lemma. 

\paragraph{Proof of Proposition \ref{matrixLG-comprop}.} Let $R_A(s)$ be a
curve in $S$ such that
$R_A(0) = I$ and $R_A^\prime(0) = A$. 
Define $S(t) = R_A(t)BR_A(t)^{-1}\in T_IS$.
Then the lemma implies  that $S(t)\in T_IS$ for every $t$.
Hence, $S^\prime(t)\in T_IS$, and in particular, $S^\prime(0) = AB - BA\in
T_IS$ . 

\begin{definition} [Matrix commutator]
For any pair of $n \times n$ matrices $A,B$, the {\bfi matrix
commutator} is defined as  $[A,B] = AB - BA$.
\end{definition} 

\begin{proposition} [Properties of the matrix commutator]
The matrix commutator has the following
two properties:
\begin{description}
\item
(i) Any two $n \times n$ matrices $A$ and $B$ satisfy
\[[B,A] = -[A,B]\] (This is
the property of skew-symmetry.)
\item
(ii) Any three $n \times n$ matrices $A$, $B$ and $C$ satisfy
\[[[A,B],C] + [[B,C],A] + [[C,A],B] = 0\]
(This is known as the {\bfi Jacobi identity}.)
\end{description}
\end{proposition}

\begin{definition} [Matrix Lie algebra]
A matrix Lie algebra $\mathfrak{g}$ is a set of $n \times n$ matrices
which is a vector space with respect to the usual operations of matrix
addition and multiplication by real numbers (scalars) and which is closed
under the matrix commutator $[\cdot\,,\,\cdot]$.
\end{definition}

\begin{proposition} For any matrix Lie group $S$, the tangent space at
the identity $T_IS$ is a matrix Lie algebra.
\end{proposition}
\paragraph{Proof.} This follows by proposition \ref{matrixLG-comprop} and
because $T_IS$ is a vector space.

\subsection{Examples of matrix Lie groups}

\begin{example}[The Orthogonal Group $O(n)$]
The mapping condition $U^TKU=K$ in Proposition \ref{Smap-prop}
specializes for $K = I$ to $U^TU=I$, which defines the orthogonal group.
Thus, in this case, $S$ specializes to $O(n)$, the group of ${n\times n}$
orthogonal matrices. The orthogonal group is of special interest in
mechanics.
\end{example}

\begin{corollary}[$O(n)$ is a matrix Lie group]
By Proposition \ref{Smap-prop} the orthogonal group $O(n)$ 
is both a subgroup and a submanifold of the general linear group
$GL(n,\mathbb{R})$. Thus, by definition \ref{MatrixLG-def}, the
orthogonal group $O(n)$ is a matrix Lie group.
\end{corollary}

\begin{example}[The Special Linear Group $SL(n,\mathbb{R})$]
The subgroup of $GL(n,\mathbb{R})$ with $\det(U)=1$ is called
$SL(n,\mathbb{R})$. 
\end{example}

\begin{example}[The Special Orthogonal Group $SO(n)$]
The special case of $S$ with $\det(U)=1$ and $K=I$ is called
$SO(n)$. In this case, the mapping condition $U^TKU=K$ specializes
to $U^TU=I$ with the extra condition $\det(U)=1$. 
\end{example}

\begin{example}[The tangent space of $SO(n)$ at the identity
$T_ISO(n)$] The special case with $K=I$ of $T_ISO(n)$
yields, 
\[
A^T+A=0
\,.
\]
These are antisymmetric matrices. Lying in the tangent space at the
identity of a matrix Lie group, this linear vector space forms a matrix
Lie algebra. 
\end{example}

\begin{example}[The Symplectic Group] Suppose $n = 2l$ (that is, let $n$ be 
even) and consider the nonsingular skew-symmetric matrix
\[
J=
\left[
\footnotesize{
\begin{array}{cc}
0 &I
\\
-I&0
\end{array}}
\right]
\]
where $I$ is the $l\times l$ identity matrix. One may verify that
\[Sp(l) = \{U \in GL(2l,\mathbb{R})|\,U^TJU = J\}\]
is a group. This is called the symplectic group. Reasoning as before, the matrix
algebra $T_ISp(l)$ is defined as the set of $n\times n$ matrices $A$ satisfying
$JA^T+AJ=0$. This algebra is denoted as $sp(l)$.
\end{example}

\begin{example}[The Special Euclidean Group]
Consider the Lie group of $4\times4$ matrices of the form
\[
E(R,v)
=
\left[
\footnotesize{
\begin{array}{cc}
R &v
\\
0&1
\end{array}}
\right]
\]
where $R\in SO(3)$ and $v\in\mathbb{R}^3$. This is the special
Euclidean group, denoted $SE(3)$. The special Euclidean group is of
central interest in mechanics since it describes the set of rigid motions
and coordinate transformations of three-dimensional space. 
\end{example}

\begin{exercise}
A point $P$ in $\mathbb{R}^3$ undergoes a rigid motion associated with
$E(R_1, v_1)$ followed by a rigid motion associated with $E(R_2, v_2)$.
What matrix element of $SE(3)$ is associated with the composition of these
motions in the given order?
\end{exercise}

\begin{exercise}
Multiply the special Euclidean matrices of $SE(3)$. Investigate their
matrix commutators in their tangent space at the identity. (This is an
example of a semidirect product Lie group.)
\end{exercise}

\begin{exercise}[Tripos question]
When does a stone at the equator of the Earth weigh the most? Two hints:
(1) Assume the Earth's orbit is a circle around the Sun and ignore the
declination of the Earth's axis of rotation. (2) This is an exercise in
using $SE(2)$. 
\end{exercise}

\begin{exercise}
Suppose the $n\times n$ matrices $A$ and $M$ satisfy \[AM+MA^T = 0\,.\]
Show that $\exp(At)M \exp(A^T t) = M$ for all $t$. Hint: $A^nM=M(-A^T)^n$. This
direct calculation shows that for $A \in so(n)$ or $A \in sp(l)$, we have $\exp(At)\in
SO(n)$ or
$\exp(At) \in Sp(l)$, respectively.
\end{exercise}

\subsection{Lie group actions} The action of a Lie group $G$ on a
manifold $M$ is a group of transformations of $M$ associated to 
elements of the group $G$, whose composition acting on $M$ is corresponds
to group multiplication in $G$.

\begin{definition} Let $M$ be a manifold and let $G$ be a Lie group. A
{\bfi left action} of a Lie group $G$ on $M$ is a smooth mapping 
$\Phi: G\times M\to M$ such that 
\begin{description}
\item
(i) $\Phi(e, x) = x \hbox{ for all } x \in M$, 
\item
(ii) $\Phi(g, \Phi(h, x)) = \Phi(gh, x)$ for all $g, h \in G$ and $x \in
M$, and 
\item
(iii) $\Phi(g, \cdot)$ is a diffeomorphism on $M$ 
for each $g \in G$. 
\end{description}
We often use the convenient notation $gx$ for $\Phi(g, x)$ and
think of the group element $g$ acting on the point $x \in M$. The
associativity condition (ii) above then simply reads $(gh)x = g(hx)$.
\end{definition}
Similarly, one can define a {\bfi right action}, which is a map $\Psi:
M\times G\to M$  satisfying $\Psi(x,e) = x$ and $\Psi(\Psi(x,g),h) =
\Psi(x,gh)$. The convenient notation for right action is $xg$ for
$\Psi(x, g)$, the right action of a group element $g$ on the point
$x\in M$. Associativity $\Psi(\Psi(x,g),h) = \Psi(x,gh)$ is
then be expressed conveniently as $(xg)h=x(gh)$.

\begin{example}[Properties of group actions]
The action $\Phi: G\times M\to M$ of a group $G$ on a manifold $M$ is said
to be:
\begin{enumerate} 
\item
{\bfi transitive}, if for every $x,y\in M$ there exists a $g\in G$, such
that $gx=y$;
\item
{\bfi free}, if it has no fixed points, that is, $\Phi_g(x)=x$ implies
$g=e$; and 
\item
{\bfi proper}, if whenever a convergent subsequence $\{x_n\}$ in $M$
exists, and the mapping $g_nx_n$ converges in $M$, then $\{g_n\}$ has
a convergent subsequence in $G$.
\end{enumerate}
\end{example}

\paragraph{Orbits.} Given a group action of $G$ on $M$, for a
given point $x \in M$, the subset
\[{\rm Orb}\, x = \{gx |\, g \in G\}\subset M \,,\]
is called the {\bfi group orbit} through $x$. In finite dimensions, it can
be shown that group orbits are always smooth (possibly immersed) manifolds.
Group orbits generalize the notion of orbits of a dynamical system. 

\begin{exercise}
The flow of a vector field on $M$ can be thought of as an action of
$\mathbb{R}$ on $M$. Show that in this case the general notion of group
orbit reduces to the familiar notion of orbit used in dynamical systems.
\end{exercise}

\begin{theorem}
Orbits of proper group actions are embedded submanifolds.
\end{theorem}
This theorem is stated in Chapter 9 of \cite{MaRa1994}, who refer to 
\cite{AbMa1978} for the proof.

\begin{example}[Orbits of $SO(3)$]
A simple example of a group orbit is the action of $SO(3)$ on
$\mathbb{R}^3$ given by matrix multiplication: The action of $A\in SO(3)$
on a point $\mathbf{x}\in\mathbb{R}^3$ is simply the product
$A\mathbf{x}$. In this case, the orbit of the origin is a single point
(the origin itself), while the orbit of any other point is the sphere
through that point.
\end{example}

\begin{example}
[Orbits of a Lie group acting on itself]
The action of a group $G$ on itself from either the left, or the right,
also produces group orbits. This action sets the stage for
discussing the tangent lifted action of a Lie group on its tangent bundle.
\\

{\bfi Left and right translations on the group} are denoted
$L_g$ and $R_g$, respectively. For example, $L_g : G \to G$ is the map
given by $h\to gh$, while $R_g : G \to G$ is the map given
by $h\to hg$, for $g,h\in G$.
\begin{description}
\item[(a)]
{\bfi Left translation} $L_g:G\to G;\,h\to gh$ defines a transitive and
free action of $G$ on itself. Right multiplication 
$R_g:G\to G;\,h\to hg$ defines a right action,  while $h\to hg^{-1}$
defines a left action of $G$ on itself. 
\item[(b)]
$G$ acts on $G$ by conjugation, $g\to I_g=R_{g^{-1}}\circ L_g$. The map
$I_g:G\to G$ given by $h\to ghg^{-1}$ is the {\bfi inner automorphism}
associated with $g$. Orbits of this action are called {\bfi conjugacy
classes}. 
\item[(c)]
Differentiating conjugation at $e$ gives the {\bfi adjoint action}
of $G$ on $\mathfrak{g}$: 
\[
{\rm Ad}_g:=T_eI_g\,:\,T_eG=\mathfrak{g}\to T_eG=\mathfrak{g}.
\]
Explicitly, the {\bfi adjoint action} of $G$ on $\mathfrak{g}$ is given by
\[
{\rm Ad}\,:\,G\times \mathfrak{g}\to \mathfrak{g}
\,,\quad
{\rm Ad}_g(\xi)=T_e(R_{g^{-1}}\circ L_g)\xi
\]
We have already seen an example of adjoint action for matrix Lie groups
acting on matrix Lie algebras, when we defined 
$S(t) = R_A(t)BR_A(t)^{-1}\in T_IS$ as a key step in the proof of
Proposition \ref{matrixLG-comprop}.  
\item[(d)]
The {\bfi coadjoint action} of $G$ on $\mathfrak{g}^*$, the dual of the Lie
algebra $\mathfrak{g}$ of $G$, is defined as follows. Let
${\rm Ad}^*_g:\mathfrak{g}^*\to\mathfrak{g}^*$ be the dual of ${\rm Ad}_g$,
defined by 
\[
\langle{\rm Ad}^*_g\alpha,\xi\rangle
=
\langle\alpha,{\rm Ad}_g\xi\rangle
\]
for $\alpha\in \mathfrak{g}^*$, $\xi\in \mathfrak{g}$ and pairing
$\langle\cdot\,,\,\cdot\rangle:\mathfrak{g}^*\times\mathfrak{g}
\to\mathbb{R}$. Then the map
\[
\Phi^*\,:\,G\times \mathfrak{g}^*\to\mathfrak{g}^*
\quad\hbox{given by}\quad
(g,\alpha)\mapsto {\rm Ad}^*_{g^{-1}}\alpha
\]
is the coadjoint action of $G$ on $\mathfrak{g}^*$.
\end{description}
\end{example}

\subsection{Examples: $SO(3)$, $SE(3)$, etc.}

\subsection*{A basis for the matrix Lie algebra $so(3)$ and a map to
$\mathbb{R}^3$} The Lie algebra of $SO(n)$ is called $so(n)$. A basis 
$(e_1,e_2,e_3)$ for $so(3)$ when
$n=3$ is given by 
\[
\mathbf{\hat{x}}
=
\left[
\footnotesize{
\begin{array}{ccc}
0 &-z&y
\\
z&0&-x
\\
-y&x&0
\end{array}}
\right]
=
xe_1+ye_2+ze_3
\]
\begin{exercise}
Show that $[e_1,e_2]=e_3$ and cyclic permutations, while all other
matrix commutators among the basis elements vanish. 
\end{exercise}

\begin{example}
[The isomorphism between $so(3)$ and $\mathbb{R}^3$]
The previous equation may be written equivalently by defining the
hat-operation $\hat{(\,\cdot\,)}$ as 
\[
\mathbf{\hat{x}}_{ij}
=
\epsilon_{ijk}x^k
\,,\quad\hbox{where}\quad
(x^1,x^2,x^3)=(x,y,z)
\,.
\]
Here $\epsilon_{123}=1$ and $\epsilon_{213}=-1$, with cyclic permutations.
The totally antisymmetric tensor
$\epsilon_{ijk}=-\,\epsilon_{jik}=-\,\epsilon_{ikj}$ also defines the cross
product of vectors in $\mathbb{R}^3$. Consequently, we may write,
\[
(\mathbf{x}\times\mathbf{y})_i
=
\epsilon_{ijk}x^jy^k
=
\mathbf{\hat{x}}_{ij}y^j
\,,\quad\hbox{that is,}\quad
\mathbf{x}\times\mathbf{y}
=
\mathbf{\hat{x}}\mathbf{y}
\]
\end{example}

\begin{exercise}
What is the analog of the hat map $so(3)\mapsto\mathbb{R}^3$ for the three
dimensional Lie algebras $sp(2,\mathbb{R})$, $so(2,1)$,
$su(1,1)$, or $sl(2,\mathbb{R})$?
\end{exercise}

Background reading for this lecture is Chapter 9 of \cite{MaRa1994}.

\subsection*{Compute the Adjoint and adjoint operations by differentiation}

\begin{description}
\rem{\item[1.]
Start from the inner automorphisms by differentiating
\[
I_g\,:\,G\to G
\,,\quad\hbox{where}\quad
I_g(h)=ghg^{-1}
\]
}
\item[1.]
Differentiate $I_g(h)$ wrt $h$ at $h=e$ to produce the {\bfi Adjoint
operation}
\[
{\rm Ad}\,:\,G\times\mathfrak{g}\to\mathfrak{g}
\,:\quad
{\rm Ad}_g\,\eta=T_eI_g\,\eta
\]
\item[2.]
Differentiate ${\rm Ad}_g\,\eta$ wrt $g$ at $g=e$ in the direction $\xi$
to get the {\bfi Lie bracket}
$[\xi,\eta]\,:\,\mathfrak{g}\times\mathfrak{g}\to\mathfrak{g}$ and thereby
to produce the {\bfi adjoint operation}
\[
T_e({\rm Ad}_g\,\eta)\,\xi=[\xi,\eta]={\rm ad}_\xi\,\eta
\]
\end{description}

\subsection*{Compute the co-Adjoint and coadjoint operations by taking
duals}

\begin{description}
\item[1.]
${\rm Ad}^*_g:\mathfrak{g}^*\to\mathfrak{g}^*$, the dual of ${\rm Ad}_g$,
is defined by 
\[
\langle{\rm Ad}^*_g\alpha,\xi\rangle
=
\langle\alpha,{\rm Ad}_g\xi\rangle
\]
for $\alpha\in \mathfrak{g}^*$, $\xi\in \mathfrak{g}$ and pairing
$\langle\cdot\,,\,\cdot\rangle:\mathfrak{g}^*\times\mathfrak{g}
\to\mathbb{R}$. The map
\[
\Phi^*\,:\,G\times \mathfrak{g}^*\to\mathfrak{g}^*
\quad\hbox{given by}\quad
(g,\alpha)\mapsto {\rm Ad}^*_{g^{-1}}\alpha
\]
defines the {\bfi co-Adjoint action} of $G$ on $\mathfrak{g}^*$.
\item[2.]
The pairing 
\[
\langle{\rm ad}^*_\xi\alpha,\eta\rangle
=
\langle\alpha,{\rm ad}_\xi\,\eta\rangle
\]
defines the {\bfi coadjoint action} of $\mathfrak{g}$ on $\mathfrak{g}^*$,
for $\alpha\in \mathfrak{g}^*$ and $\xi,\eta\in \mathfrak{g}$.
\end{description}

See Chapter 9 of \cite{MaRa1994} for more discussion of the Ad and ad
operations.

\subsection*{Example: the rotation group $SO(3)$}
\label{subsection: SO(3)}
\subsection*{The Lie algebra ${so}(3)$ and its dual.} 
The special orthogonal group is defined by
\[
SO(3): = \{A \mid A \text{~a~} 3 \times 3
\text{~orthogonal matrix}, \operatorname{\det}(A) = 1\}
\,.
\]
Its Lie algebra ${so}(3)$ is formed by $3\times 3 $ skew
symmetric matrices, and its dual is denoted ${so}(3)^\ast $. 

\subsection*{The Lie algebra isomorphism $\mathbf{\hat{\,}}\,:\,({so}(3),
[\cdot, \cdot])
\to (\mathbb{R}^3, \times)$}
The Lie algebra $({so}(3), [\cdot, \cdot])$, where $[\cdot,
\cdot]$ is the commutator bracket of matrices, is isomorphic to the Lie
algebra $(\mathbb{R}^3, \times) $, where $\times $ denotes the vector
product in $\mathbb{R}^3$, by the isomorphism
\begin{eqnarray*}
\label{so three isomorphism in coordinates}
\mathbf{u}: =(u^1, u^2, u ^3) \in \mathbb{R}^3 \mapsto 
\mathbf{\hat{u}}: =
\left[
\begin{array}{ccc}
0&-u^3&u^2\\
u^3&0&-u^1\\
-u^2&u^1&0
\end{array}
\right] \in {so}(3)
\,,\hbox{ that is, }
\mathbf{\hat{u}}_{ij}: =-\,\epsilon_{ijk}u^k
\end{eqnarray*}
Equivalently, this isomorphism is given by
\begin{eqnarray*}
\label{so three isomorphism}
\mathbf{\hat{u}} \mathbf{v} = \mathbf{u}\times \mathbf{v} \quad
\text{for all} \quad \mathbf{u}, \mathbf{v}\in \mathbb{R}^3.
\end{eqnarray*}
The following formulas for $\mathbf{u}, \mathbf{v}, \mathbf{w} \in
\mathbb{R}^3$ may be easily verified:
\begin{eqnarray*}
\label{bracket relation}
(\mathbf{u} \times \mathbf{v})\hat{\phantom{u}} 
&=& [\mathbf{\hat{u}},
\mathbf{\hat{v}}]\\
\label{triple product}
[\mathbf{\hat{u}},\mathbf{\hat{v}}]\mathbf{w} 
&=&
(\mathbf{u} \times \mathbf{v}) \times \mathbf{w}\\
\label{dot product}
\mathbf{u}\cdot \mathbf{v} 
&=& - \frac{1}{2}
\operatorname{trace}(\mathbf{\hat{u}} \mathbf{\hat{v}}).
\end{eqnarray*}

\paragraph{Ad action of $SO(3)$ on ${so}(3)$} The corresponding
adjoint action of $SO(3)$ on
$so(3)$ may be obtained as follows. For $SO(3)$ we have $I_A(B)=ABA^{-1}$.
Differentiating $B(t)$ at $B(0)=Id$ gives
\[
{\rm Ad}_A\mathbf{\hat{v}}
=\frac{d}{dt}\Big|_{t=0}AB(t)A^{-1}
=A\mathbf{\hat{v}}A^{-1}
\,,\quad\hbox{with}\quad
\mathbf{\hat{v}}=B^\prime(0)\,.
\] 
One calculates the pairing with a vector $\mathbf{w}\in\mathbb{R}^3$ as
\[
{\rm Ad}_A\mathbf{\hat{v}}(\mathbf{w})
=
A\mathbf{\hat{v}}(A^{-1}\mathbf{w})
=
A(\mathbf{v}\times A^{-1}\mathbf{w})
=
A\mathbf{v}\times\mathbf{w}
=
(A\mathbf{v})\,\hat{\,}\,\mathbf{w}
\]
where we have used a relation
\begin{eqnarray*}
\label{times relation}
A(\mathbf{u} \times \mathbf{v}) = A \mathbf{u} \times A \mathbf{v}
\end{eqnarray*}
which holds for any $\mathbf{u}, \mathbf{v} \in \mathbb{R}^3$ and $A \in
SO(3)$. \\
Consequently,
\[
{\rm Ad}_A\mathbf{\hat{v}}=(A\mathbf{v})\,\hat{}
\]
Identifying $so(3)\simeq\mathbb{R}^3$ then gives 
\[
{\rm Ad}_A\mathbf{v}=A\mathbf{v}.
\]
So (speaking prose all our lives) the adjoint action of $SO(3)$ on $so(3)$
may be identitified with  multiplication of a matrix in $SO(3)$ times a
vector in $\mathbb{R}^3$.

\paragraph{ad-action of $so(3)$ on ${so}(3)$} 
Differentiating again gives the ad-action of the Lie algebra $so(3)$ on
itself:
\[
[\mathbf{\hat{u}}, \mathbf{\hat{v}}] =
\operatorname{ad}_{\mathbf{\hat{u}}}
\mathbf{\hat{v}} =
\left.\frac{d}{dt}\right|_{t = 0} \left(e^{t
\mathbf{\hat{u}}}\mathbf{v}\right)
\!\!\!\hat{\phantom{A}} = (\mathbf{\hat{u}}
\mathbf{v})\!\hat{\phantom{A}} = (\mathbf{u}
\times \mathbf{v})\!\hat{\phantom{A}}.
\]
So in this isomorphism the vector cross product is identified with the
matrix commutator of skew symmetric matrices.

\paragraph{Infinitesimal generator}
Likewise, the {\bfi infinitesimal generator} corresponding to
$\mathbf{u}\in\mathbb{R}^3$ has the expression
\begin{eqnarray*}
\label{inf generator for rotations}
\mathbf{u}_{\mathbb{R}^3} (\mathbf{x}) := \left.\frac{d}{dt}\right|_{t =
0} e^{t \mathbf{\hat{u}}} \mathbf{x} = \mathbf{\hat{u}}\, \mathbf{x} = \mathbf{u}
\times \mathbf{x}.
\end{eqnarray*}

\begin{exercise}
What is the analog of the hat map $so(3)\mapsto\mathbb{R}^3$ for the three
dimensional Lie algebras $sp(2,\mathbb{R})$, $so(2,1)$,
$su(1,1)$, or $sl(2,\mathbb{R})$?
\end{exercise}

\subsection*{The dual Lie algebra isomorphism
$\mathbf{\tilde{\,}}\,:\,{so}(3)^\ast\to
\mathbb{R}^3$}
\paragraph{Coadjoint actions} 
The dual ${so}(3) ^\ast$ is identified with $\mathbb{R}^3$
by the isomorphism 
\[
\mathbf{\Pi} \in \mathbb{R}^3  \mapsto
\mathbf{\tilde{\Pi}} \in {so}(3) ^\ast
\,:
\quad\mathbf{\tilde{\Pi}}(\mathbf{\hat{u}}): = \mathbf{\Pi} \cdot
\mathbf{u}
\quad\hbox{for any}\quad
{\bf u}\in \mathbb{R}^3
\,.\]
In terms of this isomorphism, the co-Adjoint action of $SO(3)$ on
${so}(3)^\ast$ is given by
\begin{eqnarray*}
\label{coadjoint so three action}
\operatorname{Ad}^\ast _{A^{-1}} \mathbf{\tilde{\Pi}} = (A
\boldsymbol{\Pi})\!\tilde{\phantom{A}}
\end{eqnarray*}
and the coadjoint action of ${so}(3) $ on ${so}(3) ^\ast$ is
given by
\begin{eqnarray}
\label{coadjoint so three lie algebra action}
\operatorname{ad}^\ast_{\mathbf{\hat{u}}}
\tilde{\boldsymbol{\Pi}} = (\boldsymbol{\Pi}\times
\mathbf{u})\tilde{\phantom{u}}.
\end{eqnarray}
\paragraph{Computing the co-Adjoint action of $SO(3)$ on ${so}(3)^\ast$}
This is given by
\begin{align*}
\left(\operatorname{Ad}^\ast _{A^{-1}} \mathbf{\tilde{\Pi}}\right)
(\mathbf{\hat{u}}) &= \mathbf{\tilde{\Pi}} \cdot
\operatorname{Ad}_{A^{-1}}\mathbf{\hat{u}} = \mathbf{\tilde{\Pi}} \cdot
(A^{-1}\mathbf{u})\!\hat{\phantom{u}} = \boldsymbol{\Pi}\cdot
A^{T}\mathbf{u} \\ &= A \boldsymbol{\Pi} \cdot \mathbf{u}
=(A \boldsymbol{\Pi})\!\tilde{\phantom{A}}(\mathbf{\hat{u}}),
\end{align*}
that is, the co-Adjoint action of $SO(3) $ on ${so}(3) ^\ast$ has
the expression
\[
\operatorname{Ad}^\ast _{A^{-1}} \mathbf{\tilde{\Pi}} = (A
\boldsymbol{\Pi})\!\tilde{\phantom{A}}
\,,
\]
\rem{
where the isomorphism
$\mathbf{\tilde{\phantom{u}}}: \mathbb{R}^3 \rightarrow {so}(3) ^\ast$
is given by $\mathbf{\tilde{\Pi}}(\mathbf{\hat{u}}): = \mathbf{\Pi} \cdot
\mathbf{u}$ for any ${\bf u}\in \mathbb{R}^3$ and
$\mathbf{\hat{\phantom{u}}}: (\mathbb{R}^3, \times ) \rightarrow ({so}(3),
[\cdot, \cdot])$ is the Lie algebra isomorphism. }

Therefore, the co-Adjoint
orbit $ \mathcal{O} = \left\{A \boldsymbol{\Pi} \mid
 A\in SO(3)\right\} \subset \mathbb{R}^3$ of $SO(3)$ through
$\boldsymbol{\Pi} \in \mathbb{R}^3$ is a $2$-sphere  of
radius $\|\boldsymbol{\Pi}\|$.\\

\paragraph{Computing the coadjoint action of $so(3)$ on ${so}(3)^\ast$} Let
$\mathbf{u}, \mathbf{v} \in \mathbb{R}^3$ and note that
\begin{align*}
\left\langle\operatorname{ad}^\ast_{\mathbf{\hat{u}}}
\tilde{\boldsymbol{\Pi}},
\mathbf{\hat{v}}\right\rangle
&= \left\langle \tilde{\boldsymbol{\Pi}}, \left[ \mathbf{\hat{u}},
\mathbf{\hat{v}} \right] \right\rangle
= \left\langle \tilde{\boldsymbol{\Pi}}, (\mathbf{u} \times \mathbf{v})
\hat{\phantom{u}} \right\rangle
= \boldsymbol{\Pi}\cdot(\mathbf{u}\times \mathbf{v}) \\
&= (\boldsymbol{\Pi}\times \mathbf{u})\cdot \mathbf{v}
= \left\langle \boldsymbol{\Pi}\times \mathbf{u})\tilde{\phantom{u}},
\mathbf{\hat{v}}\right\rangle,
\end{align*}
which shows that $\operatorname{ad}^\ast_{\mathbf{\hat{u}}}
\tilde{\boldsymbol{\Pi}} = (\boldsymbol{\Pi}\times
\mathbf{u})\tilde{\phantom{u}}$, thereby proving \eqref{coadjoint so three
lie algebra action}. Therefore,
$T_{\boldsymbol{\Pi}}\mathcal{O}
= \left\{\boldsymbol{\Pi} \times \mathbf{u} \mid \mathbf{u}\in
\mathbb{R}^3 \right\}$, since the plane perpendicular to
$\boldsymbol{\Pi}$, that is, the tangent space to the sphere centered
at the origin of radius $\| \boldsymbol{\Pi}\|$, is given by
$\left\{\boldsymbol{\Pi} \times \mathbf{u} \mid
\mathbf{u}\in \mathbb{R}^3 \right\}$.\\


\section{Lifted Actions}

\begin{definition} Let  $\Phi:G\times M\to M$ be a left action, and write
$\Phi_g(x)=\Phi(g,x)$ for $x\in M$.  The {\bfi tangent lift action} of $G$ on the tangent bundle
$TM$ is defined by $gv=T_x\Phi_g(v)$ for every $v\in T_xM.$
\rem{
(Similarly for right actions, $vg = T_x\Psi_g(v).$)
Likewise, the {\bfi cotangent lift} action of $G$ on $T^*M$ is defined by
$g\alpha = \left(T_m\Phi_{g^{-1}}\right) (\alpha)$ for every $\alpha\in T_mM.$
}
\end{definition}

\begin{remark} 
\rem{
Let $\Phi(g(t),x)=c(t)\in M$ define a curve in $M$ for an open interval $t\in(-\epsilon, \epsilon)$ with $\Phi(g(0),x)=c(0)=x$. 
}
In standard calculus notation, the expression for tangent lift may be written as 
\[
T_x\Phi\cdot v 
= \frac{d}{dt}\Phi(c(t))\Big|_{t=0}
=
\frac{\partial \Phi}{\partial c}c\,'(t)\Big|_{t=0} 
=: 
D\Phi(x)\cdot v
\,,\quad\hbox{with}\quad
c(0)=x
\,,\quad
c\,'(0)=v
\,.
\]
\end{remark} 

\begin{definition}
If $X$ is a vector field on $M$ and $\phi$ is a differentiable map from $M$ to itself, then the \textbf{push-forward} of $X$ by $\phi$ is the vector field $\phi_*X$ defined by $\left(\phi_*X\right)\left(\phi(x)\right) =
T_x\phi\left(X(x)\right).$ That is, the following diagram commutes:

\hspace{3cm}
\begin{picture}(150,100)(-60,0)
\put(20,75){\vector(1,0){65} }
\put(-10,72){$TM$}
\put(95,72){$TM$}
\put(40,85){$T\phi$}
\put(100,10){\vector(0,1){55} }
\put(110,30){$\phi_*X$}
\put(20,0){\vector(1,0){65} }
\put(45,10){$\phi$}
\put(0,10){\vector(0,1){55} }
\put(-5,-3){$M$}
\put(95,-3){$M$}
\put(-20,30){$X$}
\end{picture}
\vspace{5mm}

\noindent If $\phi$ is a diffeomorphism then the \textbf{pull-back} $\phi^*X$ is also defined:  
$\left(\phi^*X\right)\left(x\right) =
T_{\phi(x)}\phi^{-1}\left(X\left(\phi(x)\right)\right).$
\end{definition}

\begin{definition} \label{defgX}
Let  $\Phi:G\times M\to M$ be a left action, and write $\Phi_g(m)=\Phi(g,m).$ 
Then $G$ has a left action on $X\in\mathfrak{X}(M)$ (the set of vector fields on $M$)
by the push-forward: $gX = \left(\Phi_g\right)_* X.$
\end{definition}

\begin{definition} \label{definvt}
Let $G$ act on $M$ on the left.
A vector field $X$ on $M$ is \textbf{invariant} with respect to this action
(we often say ``$G$-invariant'' if the action is understood) if $gX=X$ for all $g\in G;$
equivalently (using all of the above definitions!) $g\left(X(x)\right) = X(gx)$ for
all $g\in G$ and all $x\in X.$
\end{definition}


\begin{definition}
Consider the left action of $G$ on itself by left multiplication, $\Phi_g(h)=L_g(h) = gh.$
A vector field on $G$ that is invariant with respect to this action is called 
\textbf{left-invariant}. From Definition \ref{definvt}, we see that $X$ is left-invariant
if and only if $g\left(X(h)\right) = X(gh),$ which in less compact notation means
$T_hL_g X(h) = X(gh).$
The set of all such vector fields is written $\mathfrak{X}^L(G).$
\end{definition}

\begin{proposition} Given a $\xi\in T_eG,$ define $X_\xi^L(g)=g\xi$
(recall: $g\xi\equiv T_eL_g \xi$).
Then $X_\xi^L$ is the unique left-invariant vector field such that
$X_\xi^L(e)=\xi.$
\end{proposition}

\vspace{-4pt}
\textbf{Proof} To show that $X_\xi^L$ is left-invariant, we need to show that
$g\left(X_\xi^L(h)\right) = X_\xi^L(gh)$ for every $g,h\in G.$ This follows from the
definition of $X_\xi^L$ and the associativity property of group actions:
\[
g\left(X_\xi^L(h)\right) = g \left(h \xi\right) = \left(gh\right) \xi = X_\xi^L(gh)
\]
We repeat the last line in less compact notation:
\[
T_hL_g\left(X_\xi^L(h)\right) = T_hL_g \left(h \xi\right) = T_eL_{gh} \xi = X_\xi^L(gh)
\]
For uniqueness, suppose $X$ is left-invariant and $X(e)=\xi.$ Then for any $g\in G,$ we have
$X(g)=g\left(X(e)\right)=g\xi=X_\xi^L(g). \quad \blacksquare$

\begin{remark}
Note that the map $\xi\mapsto X_\xi^L$ is an vector space isomorphism from $T_eG$ to 
$\mathfrak{X}^L(G).$
\end{remark}

All of the above definitions have analogues for right actions.
The definitions of \emph{right-invariant}, $\mathfrak{X}^R(G)$ and $X_\xi^R$
use the right action of $G$ on itself defined by $\Phi(g,h)=R_g(h)=hg.$

\begin{exercise} There is a left action of $G$ on itself defined by $\Phi_g (h)=hg^{-1}.$ 
\end{exercise}

We will use the map $\xi\mapsto X_\xi^L$ to relate the Lie bracket on 
$\mathfrak{g},$
defined as $[\xi,\eta] = \mathrm{ad}_\xi \eta,$  with the
Jacobi-Lie bracket on vector fields.

\begin{definition} The \textbf{Jacobi-Lie bracket} on $\mathfrak{X}(M)$ is defined in
local coordinates by
\[
[X,Y]_{J-L} \equiv (DX)\cdot Y- (DY)\cdot X
\]
which, in finite dimensions, is equivalent to
\[
[X,Y]_{J-L} 
\equiv -\,(X\cdot \nabla)Y+ (Y\cdot \nabla) X
\equiv -\,[X,Y]
\]
\end{definition}

\begin{theorem}\label{JL} 
[Properties of the Jacobi-Lie bracket]
$\quad$\vspace{-1mm}
\begin{enumerate}
\item
The Jacobi-Lie bracket satisfies
\[[X,Y]_{J-L}=\mathcal{L}_XY\equiv \left.\frac{d}{dt}\right|_{t=0} \Phi_t^* Y,\]
where $\Phi$ is the flow of $X.$
(This is coordinate-free, and can be used as an alternative definition.)
\item This bracket makes $\mathfrak{X}^L(M)$ a Lie algebra with 
$[X,Y]_{J-L}=-\,[X,Y]$, where $[X,Y]$ is the Lie algebra bracket on
$\mathfrak{X}(M)$.
\item $\phi_*[X,Y] = [\phi_*X,\phi_*Y]$ for any differentiable $\phi:M\to M.$
\end{enumerate}
\end{theorem}

\begin{remark}\label{Lie-bracket}
The first property of the Jacobi-Lie bracket is proved for matrices in
section \ref{sec-JLhandout}. The other two properties are proved below for the case
that $M$ is the Lie group $G$. 
\end{remark}

\begin{theorem}\label{XLsub}
$\mathfrak{X}^L(G)$ is a subalgebra of $\mathfrak{X}(G).$
\end{theorem}
\begin{proof} Let $X,Y\in \mathfrak{X}^L(G).$
Using the last item of the previous theorem, and then the $G$ invariance of $X$ and $Y$,
gives the push-forward relations
\[
\left(L_g\right)_*[X,Y]_{J-L} 
= [\left(L_g \right)_*X,\left(L_g\right)_*Y]_{J-L} 
\]
for all $g\in G.$
Hence $[X,Y]_{J-L}\in \mathfrak{X}^L(G).$
This is the second property in Theorem \ref{JL}.
\end{proof} 

\rem{
\vspace{4pt}\noindent
For the proof of the third property, begin by recalling the definition $[\xi,\eta] =
\mathrm{ad}_\xi\eta.$
}
 
\begin{theorem} \label{Lie-JL}
Set $[X_\xi^L,X_\eta^L]_{J-L}(e)=[\xi,\eta]$ for every
$\xi,\eta\in \mathfrak{g},$ where the bracket on the right is the Jacobi-Lie bracket.
(We say: the Lie bracket on $\mathfrak{g}$ is the pull-back of the Jacobi-Lie bracket by
the map $\xi\mapsto X_\xi^L.$)
\end{theorem}
\begin{proof} The proof of this theorem for matrix Lie algebras is relatively easy: we have already seen that 
$\mathrm{ad}_AB = AB-BA.$ On the other hand, since $X_A^L(C)=CA$ for all $C,$ and this is 
linear in $C,$ we have $DX_B^L(I) \cdot A=AB,$ so
\begin{align*}
[A,B]=[X_A^L,X_B^L]_{J-L}(I) 
&= DX_B^L(I)\cdot X_A^L(I) - DX_A^L(I)\cdot X_B^L(I) \\
&=DX_B^L(I)\cdot A - DX_A^L(I)\cdot B = AB-BA
\end{align*}
This is the third property of the Jacobi-Lie bracket listed in Theorem \ref{JL}.
For the general proof, see Marsden and Ratiu \cite{MaRa1994}, Proposition 9.14.
\end{proof} 

\begin{remark}
This theorem, together with Item 2 in Theorem \ref{JL}, proves that the Jacobi-Lie
bracket makes $\mathfrak{g}$ into a Lie algebra.
\end{remark}

\begin{remark}
By Theorem \ref{XLsub}, the vector field $[X_\xi^L,X_\eta^L]$ is left-invariant.
Since $[X_\xi^L,X_\eta^L]_{J-L}(e)=[\xi,\eta]$, it follows that
\[
[X_\xi^L,X_\eta^L] = X_{[\xi,\eta]}^L.
\]  
\end{remark}

\begin{definition} Let $\Phi:G\times M\to M$ be a left action, and let $\xi\in \mathfrak{g}.$
Let $g(t)$ be a path in $G$ such that $g(0)=e$ and $g'(0)=\xi.$ Then the
\textbf{infinitesimal generator} of the action in the $\xi$ direction is the vector field $\xi_M$ on $M$ 
defined by
\[
\xi_M(x) = \left.\frac{d}{dt}\right|_{t=0} \Phi_{g(t)}(x)
\]
\end{definition}

\begin{remark}
Note: this definition does not depend on the choice of $g(t)$. For example, the choice in 
Marsden and Ratiu \cite{MaRa1994}  is $\exp(t\xi)$, where $\exp$ denotes the exponentiation on Lie groups
(not defined here).
\end{remark}

\begin{exercise} Consider the action of $SO(3)$ on the unit sphere $S^2$ around the
origin, and let $\xi=(0,0,1)\hat{}.$ Sketch the vector field $\xi_M$. (Hint: the
vectors all point ``Eastward.'')
\end{exercise}

\begin{theorem}
For any left action of $G,$ the Jacobi-Lie bracket of infinitesimal generators is related to the Lie bracket 
on $\mathfrak{g}$ as follows (note the minus sign):
\[
[\xi_M,\eta_M] = -\,[\xi,\eta]_M
\]
\end{theorem}

For a proof, see Marsden and Ratiu \cite{MaRa1994}, Proposition 9.3.6.

\begin{exercise}
Express the statements and formulas of this lecture for the case of
$SO(3)$ action on its Lie algebra $so(3)$. (Hint: look at the previous
lecture.) Wherever possible, translate these formulas to $\mathbb{R}^3$ by
using the $\hat{\,}$ map: $so(3)\to\mathbb{R}^3$.\\

\noindent
Write the Lie algebra for $so(3)$ using the Jacobi-Lie bracket in terms of
linear vector fields on $\mathbb{R}^3$. What are the characteristic curves
of these linear vector fields?
\end{exercise}

\newpage
\section{Handout: The Lie Derivative and the Jacobi-Lie Bracket}\label{sec-JLhandout}

Let $X$ and $Y$ be two vector fields on the same manifold $M$.

\begin{definition} The {\bfi Lie derivative} of $Y$ with respect to $X$ is 
$\mathcal{L}_XY\equiv \left.\frac{d}{dt} \Phi_t^* Y\right|_{t=0},$
where $\Phi$ is the flow of $X.$
\end{definition}

The Lie derivative $\mathcal{L}_XY$ is 
``the derivative of $Y$ in the direction given by $X.$''
Its definition is coordinate-independent. By contrast, $DY\cdot X$ (also
written as $X[Y]$) is also ``the derivative of $Y$ in the $X$ direction'',
but the value of $DY\cdot X$ depends on the coordinate system, and in
particular does not usually equal $\mathcal{L}_XY$ in the chosen
coordinate system.

\begin{theorem}\label{JLbrkt-def}
$\mathcal{L}_XY = [X,Y],$ where the bracket on the right is the {\bfi Jacobi-Lie bracket}.
\end{theorem}

\textbf{Proof.} 
In the following calculation, we assume that $M$ is finite-dimensional, and we work in local coordinates. 
Thus we may consider everything as matrices, which allows us to 
use the product rule and
the identities $\left(M^{-1}\right)' = -M^{-1}M'M^{-1}$ and 
$\frac{d}{dt}\left(D\Phi_t(x)\right) = D\left(\frac{d}{dt}\Phi_t\right)(x).$
\begin{align*}
\mathcal{L}_XY(x) &= \left.\frac{d}{dt}\Phi_t^* Y(x)\right|_{t=0}  \\
&=\left.\frac{d}{dt}\left(D\Phi_t(x)\right)^{-1} 
Y\left(\Phi_t(x)\right) \right|_{t=0} \\
&=\left[\left(\frac{d}{dt}\left(D\Phi_t(x)\right)^{-1}\right) Y\left(\Phi_t(x)\right) 
+\left(D\Phi_t(x)\right)^{-1}\frac{d}{dt} Y\left(\Phi_t(x)\right) \right]_{t=0} \\
&=\left[-\left(D\Phi_t(x)\right)^{-1}
\left(\frac{d}{dt}D\Phi_t(x) \right)\left(D\Phi_t(x)\right)^{-1}
Y\left(\Phi_t(x)\right) \right.\\
&\hspace{3cm}+\left.
\left(D\Phi_t(x)\right)^{-1}\frac{d}{dt} Y\left(\Phi_t(x)\right)
\right]_{t=0} \\ &=\left[-\left(\frac{d}{dt}D\Phi_t(x)
\right)Y\left(x\right)  +\frac{d}{dt} Y\left(\Phi_t(x)\right)
\right]_{t=0} \\ &=-D\left(\left.\frac{d}{dt}\Phi_t(x)\right|_{t=0}\right)
Y\left(x\right)  +DY(x)\left(\left.\frac{d}{dt} \Phi_t(x)\right|_{t=0}
\right)\\ &=-DX(x)\cdot Y\left(x\right) + DY(x)\cdot X(x)\\
&=[X,Y]_{J-L}(x)
\end{align*}
Therefore $\mathcal{L}_XY = [X,Y]_{J-L}\,.\quad \blacksquare$\\

\subsection*{Vorticity dynamics}
The same formula applies in infinite dimensions, although the proof is
more elaborate. For example, the equation for the vorticity dynamics of an
Euler fluid with velocity $\mathbf{u}$ (with ${\rm div}\,\mathbf{u}=0$)
and vorticity
${\omega}={\rm curl}\,\mathbf{u}$ may be written as,
\begin{eqnarray*}
\partial_t {\omega}
&=& -\,\mathbf{u}\cdot\nabla{\omega}
+{\omega}\cdot\nabla\mathbf{u}
\\&=&-\,[u,\omega]
\\&=&-\,{\rm ad}_u\omega
\\&=&-\,{\mathcal{L}}_u\omega
\end{eqnarray*}
All of these equations express the invariance of the vorticity vector
field $\omega$ under the flow of its corresponding divergenceless velocity
vector field $u$. This is also encapsulated in the language of fluid
dynamics in {\bfi characteristic form} as
\begin{eqnarray*}
\frac{d}{dt}\Big(\omega\cdot\frac{\partial}{\partial \mathbf{x}}\Big)
=0
\,,\quad\hbox{along}\quad
\frac{d\mathbf{x}}{dt}=\mathbf{u}(\mathbf{x},t)={\rm curl}^{-1}\omega
\,.
\end{eqnarray*}
Here, the curl-inverse operator is defined by the Biot-Savart Law,
\begin{eqnarray*}
\mathbf{u}={\rm curl}^{-1}\omega
={\rm curl}(-\Delta)^{-1}\omega
\,,
\end{eqnarray*}
which follows from the identity
\begin{eqnarray*}
{\rm curl}\,{\rm curl}\mathbf{u}
=-\,\Delta\mathbf{u} + \nabla{\rm div}\,\mathbf{u}
\,,
\end{eqnarray*}
and application of ${\rm div}\,\mathbf{u}=0$.
Thus, in coordinates,
\begin{eqnarray*}
\frac{d\mathbf{x}}{dt}=\mathbf{u}(\mathbf{x},t)
\quad\Longrightarrow\quad
\mathbf{x}(t,\mathbf{x}_0)=\Phi_t\mathbf{x}_0
\quad\hbox{with}\quad
\Phi_0 = Id
\,,
&&\hspace{-6mm}
\hbox{that is}\quad
\mathbf{x}(0,\mathbf{x}_0)=\mathbf{x}_0
\hbox{ at }t=0\,,
\\
\hbox{and}\quad\omega^j(\mathbf{x}(t,\mathbf{x}_0),t)
\frac{\partial}{\partial x^j(t,\mathbf{x}_0)}
&=&
\omega^A(\mathbf{x}_0)\frac{\partial}{\partial x_0^A}
\circ \Phi_t^{-1}
\,.
\end{eqnarray*}
Consequently,
\begin{eqnarray*}
\Phi_{t*}\,\omega^j(\mathbf{x}(t,\mathbf{x}_0),t)
=
\omega^A(\mathbf{x}_0)
\frac{\partial x^j(t,\mathbf{x}_0)}{\partial x_0^A}
=:
D\Phi_t\cdot\omega
\,.
\end{eqnarray*}
This is the Cauchy (1859) solution of {\bfi Euler's equation for
vorticity},
\begin{eqnarray*}
\frac{\partial\omega}{\partial t}
=[\,\omega\,,\,{\rm curl}^{-1}\omega\,]
\,.
\end{eqnarray*}
This type of equation will reappear several more times in the remaining
lectures. In it, the vorticity $\omega$ evolves by the
ad-action of the right-invariant vector field 
$\mathbf{u}={\rm curl}^{-1}\omega$. That is,
\begin{eqnarray*}
\frac{\partial\omega}{\partial t}
=
-\,{\rm ad}_{{\rm curl}^{-1}\omega}\,\omega
\,.
\end{eqnarray*}
The Cauchy solution is the tangent lift of this flow, namely,
\begin{eqnarray*}
\Phi_{t*}\,\omega(\Phi_{t}(\mathbf{x}_0))
=
T_{\mathbf{x}_0}\Phi_{t}(\omega(\mathbf{x}_0))
\,.
\end{eqnarray*}

%
%
%

\newpage
\section{Handout: Summary of Euler's equations for incompressible flow}
\label{sec-FluidsHandout}

\begin{itemize}
\item
{\bfi Euler's equation of incompressible fluid motion}

\begin{align*}\label{}
\underset{\text{(advective time derivative)}}
{\underset{d\mathbf{u}/dt\,\,\text{along}\,\,d\mathbf{x}/dt=\mathbf{u}}
{\underbrace
{\mathbf{u}_{\,t}\,+\,\mathbf{u}\,\cdot\,\nabla\,\mathbf{u}}}}
\hspace{-3mm}
+\,\nabla\,p\,=\,0
\end{align*}
where
$\mathbf{u}:\mathbb{R}^3\times\mathbb{R}\rightarrow\mathbb{R}^3$
satisfies ${\rm div}\,\mathbf{u}=0$.

\item
{\bfi Geometric dynamics of vorticity}
\begin{align*}\label{}
\boldsymbol{\omega}&=\text{curl}\,\mathbf{u}\\
\boldsymbol\omega_t&=
-\mathbf{u}\cdot\nabla\boldsymbol\omega
+\boldsymbol\omega\cdot\nabla\mathbf{u}\\
&=
-[\mathbf{u},\boldsymbol{\omega}]\\
&=
-\text{ad}_{\mathbf{u}}\boldsymbol{\omega} \\
&=-\mathcal{L}_{\mathbf{u}}\boldsymbol\omega
\end{align*}
In these equations, one denotes 
$\frac{d}{dt}=\frac\partial{\partial{t}}+\mathcal{L}_{\mathbf{u}}$ and, hence, may write Euler vorticity dynamics equivalently in any of the following three forms
\[
\frac{d\boldsymbol\omega}{dt}=\boldsymbol\omega\cdot\nabla\mathbf{u}
\,,
\]
as well as
\begin{align*}\label{}
\left(
\frac\partial{\partial{t}}+\mathcal{L}_{\mathbf{u}}
\right)
\left(
\boldsymbol\omega\cdot\frac\partial{\partial\mathbf{x}}
\right)=0
\end{align*}
or
\begin{align*}\label{}
\frac{d}{dt}(\boldsymbol\omega\cdot\nabla)=0
\quad\text{along}\quad
\frac{d\mathbf{x}}{dt}=\mathbf{u}
\end{align*}
The last form is found using the chain rule as
\begin{align*}\label{}
\frac{d}{dt}(\boldsymbol\omega\cdot\nabla)
=
\frac{d\boldsymbol\omega}{dt}\cdot\nabla
+
\boldsymbol\omega\cdot\frac{d}{dt}\nabla
=
\Big(\frac{d}{dt}\boldsymbol\omega
-
\boldsymbol\omega\cdot\nabla\mathbf{u} \Big)\cdot\nabla
=
0
\,.
\end{align*}

\item
{\bfi Ertel's theorem \cite{Er1942}} The operators $d/dt$ and $\boldsymbol{\omega\cdot\nabla}$ commute on solutions of Euler's fluid equations. That is,
\begin{align*}\label{}
\left[\frac{d}{dt},\boldsymbol\omega\cdot\nabla\right]=0
\,,\quad\text{so that}\quad
\frac{d}{dt}(\boldsymbol\omega\cdot\nabla\mathbf{A})
=\boldsymbol\omega\cdot\nabla\frac{d}{dt}\mathbf{A}
\quad\quad\hbox{for all differentiable }
\mathbf{A}
\end{align*}
when $\boldsymbol\omega={\rm curl}\,\mathbf{u}$ and $\mathbf{u}$ is a solution of Euler's equations for incompressible fluid flow. Consequently, one finds the following infinite set of {\bfi conservation laws}:
\begin{align*}\label{}
\text{If}\quad \frac{d\mathbf{A}}{dt}=0,\quad\text{then}\quad
\int\Phi(\boldsymbol\omega\cdot\nabla\mathbf{A})\,d^3\mathbf{x}=\text{const}\quad\hbox{for all differentiable }
\Phi
\end{align*}

\item
{\bfi Ohktani's formula \cite{Ohk1993}}
\begin{align*}\label{}
\frac{d^2\boldsymbol\omega}{dt^2}
&=
\frac{d}{dt}(\boldsymbol\omega\cdot\nabla\mathbf{u})\\
&=
\boldsymbol\omega\cdot\nabla\frac{d\mathbf{u}}{dt}\\
&=
-\boldsymbol\omega\cdot\nabla\nabla{p}\\
&=
-\mathbb{P}\,\boldsymbol\omega
\end{align*}
where
\begin{align*}\label{}
\mathbb{P}_{ij}=\frac{\partial^2 p}{\partial {x^i} \partial {x^j}}
\qquad \text{(``Hessian'' of pressure)}
\end{align*}
In addition, one has the relations
\begin{align*}\label{}
p&=-\Delta^{-1}\,\text{tr}(\nabla\mathbf{u}^T\cdot\nabla\mathbf{u})\\
S&=\frac12 (\nabla\mathbf{u}+\nabla\mathbf{u}^T)\quad\text{(strain rate tensor)}
\end{align*}
so that, the following system of equations results,
\begin{align*}\label{}
\frac{d\boldsymbol\omega}{dt}&=S\boldsymbol\omega\\
\frac{d^2\boldsymbol\omega}{dt^2}&=-\mathbb{P}\boldsymbol\omega
\end{align*}

\item
{\bfi Kelvin (1890's) circulation theorem}
\begin{align*}\label{}
\boldsymbol\omega\cdot\frac{\partial}{\partial\mathbf{x}}
&=
\omega^j\frac{\partial}{\partial{x}^j}\\
\frac{d\mathbf{u}}{dt}+&\nabla{p}=0
\,,
\end{align*}
where div$\,\mathbf{u}=0$, or equivalently $u^j_{\,,j}=0$ in index notation.
The motion equation may be rewritten equivalently as a 1-form relation,
\begin{align*}\label{}
\frac{d{u_i}}{dt}\,dx^i&=-dp=\nabla_i{p}\,dx^i
\quad\text{along}\quad \frac{d\mathbf{x}}{dt}=\mathbf{u}\\
\frac{d}{dt}(u_i\,d&x^i)\,\,\,-
\underset{=u_i du^i=d|\mathbf{u}|^2/2}
{\underbrace{u_i\,\,\,\frac{d}{dt}dx^i}}
=\,\,-dp
\end{align*}
Consequently, 
\begin{align*}\label{}
\frac{d}{dt}(\mathbf{u}\cdot d\mathbf{x})
&=
-\,d\Big(p-\frac{|\mathbf{u}|^2}{2}\Big)
\end{align*}
which becomes
\begin{align*}\label{}
\frac{d}{dt}
\oint_{C(\mathbf{u})}\hspace{-3mm}\mathbf{u}\cdot d\mathbf{x}
&=
-\oint_{C(\mathbf{u})}\hspace{-3mm} d\left(p-\frac{\mathbf{u}^2}{2}\right)
=0
\end{align*}
upon integrating around a closed loop $C(\mathbf{u})$ moving with velocity $\mathbf{u}$
The 1-form relation above may be rewritten as
\begin{align*}\label{}
(\partial_t+\mathcal{L}_{\mathbf{u}})(\mathbf{u}\cdot d\mathbf{x})&=
-d\left(p-\frac{\mathbf{u}^2}{2}\right)
\end{align*}
whose exterior derivative yields using $d^2=0$ 
\begin{align*}\label{}
(\partial_t+\mathcal{L}_{\mathbf{u}})(\boldsymbol\omega\cdot d\mathbf{S})&=0
\end{align*}
where $\boldsymbol\omega\cdot d\mathbf{S}
={\rm curl}\,\mathbf{u}\,\cdot\, d\mathbf{S}
= d\,(\mathbf{u}\cdot d\mathbf{x})$.

For these geometric quantities, one sees that the characteristic, or advective derivative is equivalent to a Lie derivative. Namely, 
\begin{align*}\label{}
\underset{\text{fluids}}
{\underbrace{
\left.\frac{d}{dt}\right|_\text{advect}}}=\quad
\underset{\text{geometry}}
{\underbrace{
\underset{\quad}{
\partial_t+\mathcal{L}_{\mathbf{u}}}}}
\end{align*}

\item
{\bfi Stokes theorem}

The classical theorem due to Stokes
\begin{align*}\label{}
\oint_{\partial{S}} \mathbf{u}\cdot d\mathbf{x}=
\iint_S \text{curl\,}\mathbf{u}\,\cdot\,d\mathbf{S}
\end{align*}
shows that Kelvin's circulation theorem is equivalent to conservation of flux of vorticity
\begin{align*}\label{}
\frac{d}{dt}\iint_S \boldsymbol\omega\,\cdot\,&d\mathbf{S}=0
\,,\quad\hbox{with}\quad
\partial S=\,C(\mathbf{u})
\end{align*}
through any surface comoving with the flow. 

Recall the definition,
\begin{align*}\label{}
\boldsymbol\omega\cdot\frac{\partial}{\partial\mathbf{x}}\,\contract\,d^3 x
=\boldsymbol\omega\cdot d\mathbf{S}
\end{align*}

One may check this formula directly, by computing
\begin{align*}\label{}
&\left(
\omega^1\frac{\partial}{\partial x^1}+
\omega^2\frac{\partial}{\partial x^2}+
\omega^3\frac{\partial}{\partial x^3}
\right)
\,\contract\,
\left(
dx^1\wedge dx^2\wedge dx^3
\right)\\&
=\omega^1 dx^2\wedge dx^3 + \omega^2 dx^3\wedge dx^1 \omega dx^1 \wedge dx^2\\&
=\boldsymbol\omega\cdot d \mathbf{S}
\end{align*}
One may then use the vorticity equation in vector-field form,
\begin{align*}\label{}
(\partial_t+\mathcal{L}_{\mathbf{u}})\,\boldsymbol\omega\cdot\!\frac{\partial}{\partial\mathbf{x}}=0
\end{align*}
to prove that the flux of vorticity through any comoving surface is conserved, as follows. 
\begin{align*}\label{}
&(\partial_t+\mathcal{L}_{\mathbf{u}})\,
\left(
\boldsymbol\omega\cdot\!\frac{\partial}{\partial\mathbf{x}}
\,\contract\,
d^3\!x
\right)\\
&=
\Big(\
\underbrace{\
(\partial_t+\mathcal{L}_{\mathbf{u}})\,\boldsymbol\omega\cdot\!\frac{\partial}{\partial\mathbf{x}}\
}_{\hbox{$=0$}}\
\Big)\,\contract\,d^3\!x
+
\boldsymbol\omega\cdot\!\frac{\partial}{\partial\mathbf{x}}
\,\contract\,
\underbrace{\
(\partial_t\,+\mathcal{L}_{\mathbf{u}})d^3\!x\
}_{\hbox{$=\text{div}\,\mathbf{u}\,d^3\!x=0$}}
=0
\end{align*}
That is, as computed above using the exterior derivative
\[
(\partial_t+\mathcal{L}_{\mathbf{u}})\,
\boldsymbol\omega\cdot d \mathbf{S}
=
0
\,.
\]

\item
{\bfi Momentum conservation}

From Euler's fluid equation $du_i/dt+\nabla_i p=0$ with $u^j_{,j}=0$ one finds,
\begin{eqnarray*}\label{}
\int (\partial_t u_i +u^j \partial_j u_i +\partial_i p)d^3\!x&=&0
\\&=&
\frac{d}{dt}\int\! u_i\,d^3\!x +\int\!\partial_j
(u_i u^j +p\,\delta_i^j)\, d^3\!x
\\&=&
\underbrace{\
\frac{d}{dt}M_i\
}_{\hbox{$=0$}}
+
\underbrace{\
\oint\widehat{n}_j (u_i u^j+p\,\delta_i^j)d\!S\
}_{\hbox{$=0$, if 
$\widehat{\mathbf{n}}\cdot\mathbf{u}=0$}}
\end{eqnarray*}

Local conservation of fluid momentum is expressed using differentiation by parts as
\begin{equation*}\label{}
\partial_t u_i=-\partial_j T_i^{\,j}
\end{equation*}
where $T_i^{\,j}:=u_i\,u^j+p\,\delta_i^j$ is the {\bfi fluid stress tensor}. 

Moreover, each component of the total momentum $M_i=\int\! u_i\,d^3\!x$ for $i=1,2,3,$ is conserved for an incompressible Euler flow, provided the flow is tangential to any fixed boundaries, i.e., $\widehat{\mathbf{n}}\cdot\mathbf{u}=0$.

\item
{\bfi Mass conservation} 

For mass density $D(\mathbf{x},t)$ with total mass
$\int D(\mathbf{x},t)\,d^3\!x$, along $d\mathbf{x}/dt = \mathbf{u}(x,t)$ one finds,
\begin{align*}\label{}
\frac{d}{dt}D\,d^3\!x
=
(\partial_t+\mathcal{L}_{\mathbf{u}})(D d^3\!x)
=
\underbrace{\
(\partial_t D+\text{div} D\mathbf{u})\
}_{\hbox{continuity eqn}}
\,d^3\!x
=
0
\end{align*}
The solution of this equation is written in Lagrangian form as
\begin{align*}\label{}
(D\,d^3\!x)\cdot g^{-1}(t)=D(\mathbf{x}_0)\,d^3\!x
\end{align*}
For incompressible flow, this becomes
\begin{align*}\label{}
\frac1D=\text{det}\,
\frac{\partial \mathbf{x\,\,}}{\partial \mathbf{x}_0}=
\frac{d^3\!x}{d^3\!x_0}=1
\end{align*}
Likewise, in the Eulerian representation one finds the equivalent relations,
\begin{align*}\label{}
\left.
\begin{array}{c}
D=1 \\ 
\partial_t D+ \text{div}(D\mathbf{u})=0
\end{array}
\right\}\Rightarrow
\text{div}\,\mathbf{u}=0
\end{align*}

\item
{\bfi Energy conservation}

Euler's fluid equation for incompressible flow div$\mathbf{u}=0$
\begin{align*}\label{}
\partial_t\mathbf{u}+\mathbf{u}&\cdot\nabla\mathbf{u}+\nabla{p}=0
\end{align*}
conserves the total kinetic energy, defined by
\begin{align*}\label{}
KE&=\int \frac12 |\mathbf{u}|^2 d^3\!x
\end{align*}
The vector calculus identity
\begin{align*}\label{}
\mathbf{u}\cdot\nabla\mathbf{u}&=
-\mathbf{u}\times\text{curl}\,\mathbf{u}+\frac12\nabla|\mathbf{u}|^2
\end{align*}
recasts Euler's equation as 
\begin{align*}\label{}
{\partial_t} \mathbf{u}-\mathbf{u}\times&\text{curl}\,\mathbf{u}+
\nabla\left(p+\frac{\mathbf{u}^2}{2}\right)=0
\end{align*}
So that
\begin{align*}\label{}
\frac{\partial}{\partial t}\frac{|\mathbf{u}|^2}{2}
\
+\
&
\text{div}\!\left(p+\frac{|\mathbf{u}|^2}{2}\right)\!\mathbf{u}=0
\end{align*}
Consequently,
\begin{align*}\label{}
\frac{d}{dt}\int_\Omega \frac{|\mathbf{u}|^2}{2} d^3\!x
=
-&\oint_{\partial\Omega}\!\!\left(p+\frac{|\mathbf{u}|^2}{2}\right)
\!\mathbf{u}\cdot d\mathbf{S}=0
\end{align*}
since $\mathbf{u}\cdot d\mathbf{S}=\mathbf{u}\cdot\widehat{\mathbf{n}}\, dS=0$ on any fixed boundary and one finds
\begin{align*}\label{}
KE=\int\frac12|\mathbf{u}|^2 d^3\!x=\text{const}
\end{align*}
for Euler fluid motion.
\end{itemize}

\newpage

\section{Lie group action on its tangent bundle}
\begin{definition}
A Lie group $G$ acts on its tangent bundle $TG$ by tangent lifts.%
\footnote{Recall definition \ref{tang-lift-defn} of tangent lifts of a
differentiable manifold.} Given $X\in T_hG$ we can consider the action
of $G$ on $X$ by either left or right translations, denoted as 
$T_hL_gX$ or $T_hR_gX$, respectively. These expressions may be
abbreviated as 
\[
T_hL_gX=L^*_gX=gX
\quad\hbox{and}\quad
T_hR_gX=R^*_gX=X g
\,.
\]
Left action of a Lie group $G$ on its tangent bundle $TG$ is
illustrated in the figure below.

\hspace{3cm}
\begin{picture}(150,100)(-60,0)
\put(20,75){\vector(1,0){65} }
\put(-10,72){$TG$}
\put(95,72){$TG$}
\put(40,85){$TL_g$}
\put(100,10){\vector(0,1){55} }
\put(110,30){$gX$}
\put(20,0){\vector(1,0){65} }
\put(45,10){$L_g$}
\put(0,10){\vector(0,1){55} }
\put(-5,-3){$G$}
\put(95,-3){$G$}
\put(-20,30){$X$}
\end{picture}
\vspace{5mm}

\noindent
For matrix Lie groups, this action is just multiplication on the left or
right, respectively. 
\end{definition}

\paragraph{Left- and Right-Invariant Vector Fields.} 
A vector field $X$ on $G$ is called left-invariant, if for every $g\in G$
one has $L^*_gX=X$, that is, if 
\[
(T_hL_g)X(h)=X(gh)
\,,
\]
for every $h\in G$. 
The commutative diagram for a left-invariant vector field is illustrated
in the figure below.

\hspace{3cm}
\begin{picture}(150,100)(-60,0)
\put(20,75){\vector(1,0){65} }
\put(-10,72){$TG$}
\put(95,72){$TG$}
\put(40,85){$TL_g$}
\put(100,10){\vector(0,1){55} }
\put(110,30){$X$}
\put(20,0){\vector(1,0){65} }
\put(45,10){$L_g$}
\put(0,10){\vector(0,1){55} }
\put(-5,-3){$G$}
\put(95,-3){$G$}
\put(-20,30){$X$}
\end{picture}
\vspace{5mm}

\begin{proposition}
The set $\mathfrak{X}_L(G)$ of left invariant vector fields on the
Lie group $G$  is a subalgebra of $\mathfrak{X}(G)$ the set of all
vector fields on $G$.
\end{proposition}
\paragraph{Proof.} If $X,Y\in\mathfrak{X}_L(G)$ and $g\in
G$, then 
\[
L_g^*[X,Y]
=
[L_g^*X,L_g^*Y]
=
[X,Y]
\,,
\]
Consequently, the Lie bracket $[X,Y]\in \mathfrak{X}_L(G)$. Therefore,
$\mathfrak{X}_L(G)$ is a subalgebra of $\mathfrak{X}(G)$, the set of all
vector fields on $G$. 

\begin{proposition}
The linear maps
$\mathfrak{X}_L(G)$ and $T_eG$ are isomorphic as vector spaces.
\end{proposition}
\paragraph{Demonstration of proposition.}
For each $\xi\in T_eG$, define a vector field $X_\xi$ on $G$ by letting 
$
X_\xi(g)=T_eL_g(\xi)
\,.
$
Then 
\begin{eqnarray*}
X_\xi(gh)&=&T_eL_{gh}(\xi)=T_e(L_g\circ L_h)(\xi)
\\
&=& T_hL_g(T_eL_h(\xi))=T_hL_g(X_\xi(h))
\,,
\end{eqnarray*}
which shows that $X_\xi$ is left invariant. 
(This proposition is stated in Chapter 9 of \cite{MaRa1994}, who refer to
\cite{AbMa1978} for the full proof.)

\begin{definition}\label{Jac-Lie-bracket-defn}
[{\bfi Jacobi-Lie bracket of vector fields}] Let $g(t)$ and $h(s)$ be
curves in $G$ with $g(0)=e$, $h(0)=e$ and define vector fields at the
identity of $G$ by the tangent vectors  $g'(0)=\xi$, $h'(0)=\eta$. Compute
the linearization of the Adjoint action of $G$ on $T_eG$ as 
\[
[\xi,\eta]:=
\frac{d}{dt}\frac{d}{ds}g(t)h(s)g(t)^{-1}\Big|_{s=0,t=0}
=
\frac{d}{dt}g(t)\eta g(t)^{-1}\Big|_{t=0}
=
\xi\eta-\eta\xi
\,.
\]
This is the {\bfi Jacobi-Lie bracket} of the vector fields $\xi$ and
$\eta$.
\end{definition}

\begin{definition}
The {\bfi Lie bracket} in $T_eG$ is defined by 
\[
[\xi,\eta]:=[X_\xi,X_\eta](e)
\,,
\]
for $\xi,\,\eta\in T_eG$ and for $[X_\xi,X_\eta]$ the Jacobi-Lie bracket
of vector fields. This makes $T_eG$ into a Lie algebra. Note that 
\[
[X_\xi,X_\eta]=X_{[\xi,\eta]}
\,,
\]
for all $\xi,\,\eta\in T_eG$.
\end{definition}

\begin{definition}
The vector space $T_eG$ with this Lie algebra structure is called the
{\bfi Lie algebra of $G$} and is denoted by $\mathfrak{g}$. 
\end{definition}

If we let $\xi_L(g) = T_eL_g\xi$, then the Jacobi-Lie bracket of two such
left-invariant vector fields in fact gives the Lie algebra bracket: 
\[[\xi_L,\eta_L](g) = [\xi,\eta]_L(g)\]
For the right-invariant case, the right hand side obtains a minus sign,
\[[\xi_R,\eta_R](g) = -[\xi,\eta]_R(g)\,.\]
The relative minus sign arises because of the difference in action
$(xh^{-1})g^{-1}=x(gh)^{-1}$ on the right versus $(gh)x = g(hx)$ on the
left.

\paragraph{Infinitesimal Generator.} In mechanics, group actions 
often appear as symmetry transformations, which arise through their
infinitesimal generators, defined as follows.
\begin{definition} Suppose $\Phi : G \times M \to M$ is an action. For
$\xi\in \mathfrak{g}$, $\Phi^\xi(t,x) : \mathbb{R}\times M \to M$ defined
by $\Phi^\xi(x) = \Phi(\exp{t\xi},x)=\Phi_{\exp{t\xi}}(x)$ is an
$\mathbb{R}-$action on $M$. In other words, $\Phi_{\exp{t\xi}}\to M$ is
a flow on $M$. The vector field on
$M$ defined by%
\footnote{Recall Definition \ref{VecField-defn} of vector fields.}
\[
\xi_M(x) 
= 
\frac{d}{dt}\Big|_{t=0}
\Phi_{\exp{t\xi}}(x)
\]
is called the {\bfi infinitesimal generator} of the action corresponding
to $\xi$.
\end{definition}
The Jacobi-Lie bracket of infinitesimal generators is related to the Lie
algebra bracket as follows:
\[
[\xi_M,\eta_M]=-[\xi,\eta]_M
\,.
\]
See, for example, Chapter 9 of \cite{MaRa1994} for the
proof. 

%
%
%

\section{Lie algebras as vector fields}

\begin{definition} [The ad-operation]
For $A \in \mathfrak{g}$ we define the operator ad$_A$ to be the operator
ad$:\,\mathfrak{g}\times\mathfrak{g}\to\mathfrak{g}$ that maps $B \in
\mathfrak{g}$ to $[A,B]$. We write ad$_A B = [A,B]$.
\end{definition} 

\begin{definition} A {\bfi representation} of a Lie algebra $\mathfrak{g}$
on a vector space $V$ is a mapping $\rho$ from $\mathfrak{g}$ to the
linear transformations of $V$ such that for $A,B \in \mathfrak{g}$ and any
constant scalar $c$,
\begin{eqnarray*}
&(i)&\quad \rho(A + cB) =\rho(A) + c\rho(B)
\\
&(ii)&\quad \rho([A,B]) = \rho(A)\rho(B) - \rho(B)\rho(A).
\end{eqnarray*}
If the map $\rho$ is 1-1 the representation is said to {\bfi faithful}.
\end{definition}

\begin{exercise}
For a Lie algebra $\mathfrak{g}$, show that the map $A \to{\rm ad}A$ is a
representation of the Lie algebra $\mathfrak{g}$, with $\mathfrak{g}$
itself the vector space of the representation. This is called the
{\bfi adjoint representation}.
\end{exercise}

\begin{example}[Vector field representations of Lie algebras] 
The Jacobi-Lie bracket of the vector fields $\xi$ and $\eta$ in
definition \ref{JLbrkt-def} may be represented in coordinate
charts as,
\[
\eta=\frac{dx}{ds}\Big|_{s=0}=v(x)
\,,\quad\hbox{and}\quad
\xi=\frac{dx}{dt}\Big|_{t=0}=u(x)
\,.
\]
The Jacobi-Lie bracket of these two vector fields yields a third vector
field,
\[
\xi\eta-\eta\xi
=
\frac{d\eta}{dt}\Big|_{t=0}-\frac{d\xi}{ds}\Big|_{s=0}
=
\frac{dv}{dx}\frac{dx}{dt}\Big|_{t=0}
-\frac{du}{dx}\frac{dx}{ds}\Big|_{s=0}
=
\frac{dv}{dx}\cdot u
-\frac{du}{dx}\cdot v
=
u\cdot\nabla v - v \cdot\nabla u
\,.
\]
Thus, the Jacobi-Lie bracket of vector fields at the tangent space of the
identity $T_eG$ is closed and may be represented in coordinate charts by
the Lie bracket (commutator of vector fields)
\[
[\xi,\eta]
:=\xi\eta-\eta\xi
=
u\cdot\nabla v - v \cdot\nabla u =: [u,v]
\,.
\]
\end{example}
This example also proves the following 
\begin{proposition}
Let $\mathfrak{X}(\mathbb{R}^n)$ be the set of vector fields defined on
$\mathbb{R}^n$. A Lie algebra $\mathfrak{g}$ may be
represented on coordinate charts by vector fields 
$X_\xi=X_\xi^i\frac{\partial}{\partial x^i}\in
\mathfrak{X}(\mathbb{R}^n)$ for each element $\xi\in\mathfrak{g}$. This
vector field representation satisfies
\[
X_{[\xi,\eta]}=[X_\xi,X_\eta]
\]
where $[\xi,\eta]\in\mathfrak{g}$ is the Lie algebra product and
$[X_\xi,X_\eta]$ is the vector field commutator.
\end{proposition}

\section{Lagrangian  and Hamiltonian Formulations}

\subsection{Newton's equations for particle motion in Euclidean space} 
{\bfi Newton's equations} 
\begin{equation}
m_i\mathbf{\ddot{q}}_i = \mathbf{F}_i\,,\quad
i=1,\ldots, N \,,\quad\hbox{(no sum on }i)
\label{Newton1}
\end{equation}
describe the {\bfi accelerations} $\mathbf{\ddot{q}}_i$ of $N$ particles
with
\begin{eqnarray*}
\hbox{Masses}&&
m_i\,, \quad
i=1,\ldots, N 
\,,\\
\hbox{Euclidean positions}&&
\mathbf{q}: =(\mathbf{q}_1, \dots , \mathbf{q}_N) \in
\mathbb{R}^{3N}
\,,
\end{eqnarray*}
in response to {\bfi prescribed forces},
\begin{eqnarray*}
\mathbf{F} =(\mathbf{F}_1, \dots ,\mathbf{F}_N)
\,,
\end{eqnarray*}
acting on these particles.
Suppose the forces arise from a {\bfi potential}. That is, let
\begin{equation}
{\bf F}_i(\mathbf{q})=
-\,\frac{\partial V(\{\mathbf{q}\})}{\partial \mathbf{q}_i}
\,,\quad V:\mathbb{R}^{3N} \rightarrow \mathbb{R}
\,,
\label{force}
\end{equation}
where $\partial V/ \partial \mathbf{q}_i $ denotes the gradient of
the potential with respect to the variable $\mathbf{q}_i $.  Then Newton's
equations
\eqref{Newton1} become
\begin{equation}
m_i \mathbf{\ddot{q}}_i =
-\,\frac{\partial V}{\partial \mathbf{q}_i}\,,\quad
i=1,\ldots, N \,.
\label{Newton2}
\end{equation}
\begin{remark}
Newton (1620) introduced the gravitational potential for celestial
mechanics, now called the Newtonian potential,
\begin{equation}
 V(\{\mathbf{q}\})
=
\sum_{i,j=1}^N \frac{-\,Gm_im_j}{|\mathbf{q}_i-\mathbf{q}_j|}
\,.
\label{Newton-pot}
\end{equation}
\end{remark}

\subsection{Equivalence Theorem}

\begin{theorem}[Lagrangian  and Hamiltonian formulations] Newton's
equations in potential form,
\begin{equation}
m_i \mathbf{\ddot{q}}_i =
-\,\frac{\partial V}{\partial \mathbf{q}_i}\,,\quad
i=1,\ldots, N \,,
\label{Newton-thm}
\end{equation}
for particle motion in Euclidean space $\mathbb{R}^{3N}$ are
{\bfi equivalent} to the following four statements:
\begin{description}
\item [(i)] 
The Euler-Lagrange equations
\begin{equation}
\frac{d}{dt}\left(\frac{\partial L}{\partial\mathbf{\dot{q}}_i}\right)
-\frac{\partial L}{\partial \mathbf{q}_i}=0
\,, 
\qquad i = 1, \dots, N\,,
\label{EulerLag-eqs}
\end{equation}
hold for the Lagrangian
$
L : \mathbb{R}^{6N} =
\{(\mathbf{q},\mathbf{\dot{q}}) \mid \mathbf{q}, \mathbf{\dot{q}} \in
\mathbb{R}^{3N} \} \rightarrow \mathbb{R}
\,,
$
defined by
\begin{equation}
L(\mathbf{q},\mathbf{\dot{q}}):=\sum_{i=1}^N\frac{m_i}{2}\,
\|\mathbf{\dot{q}}_i\|^2-V(\mathbf{q})
\,,\label{Lagrangian-thm}
\end{equation}
with $\|\mathbf{\dot{q}}_i\|^2=\mathbf{\dot{q}}_i\cdot\mathbf{\dot{q}}_i
=\dot{q}^j_i\dot{q}^k_i\delta_{jk}$ (no sum on $i$).
\item [(ii)] 
{\bfi Hamilton's principle of stationary action}, $\delta\mathcal{S}=0$,
holds  for the {\bfi action functional} (dropping $i$'s)
\begin{equation}
\mathcal{S}[\mathbf{q}(\cdot)]: =\int_a^b
L(\mathbf{q}(t),\mathbf{\dot{q}}(t))\, dt\,.
\label{Action-thm}
\end{equation}
\item [(iii)] 
{\bfi Hamilton's equations of motion},
\begin{equation}
\label{hameq-thm}
\mathbf{\dot{q}} = \frac{\partial H}{\partial \mathbf{p}}\,, \qquad
\mathbf{\dot{p}} =-\,\frac{\partial H}{\partial \mathbf{q}}
\,,
\end{equation}
hold for the {\bfi Hamiltonian} resulting from the {\bfi Legendre
transform}, 
\begin{equation}
H(\mathbf{q},\mathbf{p}) :=
\mathbf{p}\cdot\mathbf{\dot{q}}(\mathbf{q},\mathbf{p}) -
L(\mathbf{q},\mathbf{\dot{q}}(\mathbf{q},\mathbf{p}))
\,,
\end{equation}
where $\mathbf{\dot{q}}(\mathbf{q},\mathbf{p})$ solves for
$\mathbf{\dot{q}}$ from the definition $\mathbf{p}:=\partial
L(\mathbf{q},\mathbf{\dot{q}})/\partial\mathbf{\dot{q}}$.

In the case of Newton's equations in potential form (\ref{Newton-thm}),
the Lagrangian in equation (\ref{Lagrangian-thm}) yields
$\mathbf{p}_i=m_i\mathbf{\dot{q}}_i$ and the resulting Hamiltonian is
(restoring $i$'s)
\[
H
= 
\underbrace{\
\sum_{i=1}^N\frac{1}{2m_i}\|\mathbf{p}_i\|^2\
}_{\hbox{Kinetic energy}}
\
+\
\underbrace{\
V(\mathbf{q})\
}_{\hbox{Potential}}
\]
\item [(iv)] 
Hamilton's equations in their Poisson bracket formulation,
\begin{equation}
\dot{F}=\{F,H\}\quad\mbox{for all $F\in{\cal F}(P)$}\,,
\label{PoissonEqs-thm}
\end{equation}
hold with {\bfi Poisson bracket\/} defined by
\begin{equation}
\{F,G\} := \sum_{i =1}^N\left(\frac{\partial F}{\partial \mathbf{q}_i}\cdot
\frac{\partial G}{\partial \mathbf{p}_i} -\frac{\partial F}{\partial
\mathbf{p}_i} \cdot
\frac{\partial G}{\partial \mathbf{q}_i}\right)\quad \text{for all} \quad
F,G\in{\cal F}(P).
\label{Poisson-thm}
\end{equation}
\end{description}
\end{theorem}
We will prove this theorem by proving a chain of linked
equivalence relations: 
(\ref{Newton-thm}) $\Leftrightarrow$ {\bf (i)} 
$\Leftrightarrow$ {\bf (ii)} $\Leftrightarrow$ {\bf (iii)}
$\Leftrightarrow$ {\bf (iv)} as propositions. (The symbol $\Leftrightarrow$
means ``equivalent to''.)

\rem{
\begin{proposition}[Newton's equations (\ref{Newton-thm}) are equivalent to
(i)] Upon introducing the {\bfi Lagrangian\/} function,
\[
L : \mathbb{R}^{6N} =
\{(\mathbf{q},\mathbf{\dot{q}}) \mid \mathbf{q}, \mathbf{\dot{q}} \in
\mathbb{R}^{3N} \} \rightarrow \mathbb{R}
\,,
\]
defined by
\begin{equation}
L(\mathbf{q},\mathbf{\dot{q}}):=\frac{1}{2}\sum_{i=1}^Nm_i
\|\mathbf{\dot{q}}_i\|^2-V(\mathbf{q})
\,,\label{Lagrangian}
\end{equation}
with $\|\mathbf{\dot{q}}_i\|^2=\mathbf{\dot{q}}_i\cdot\mathbf{\dot{q}}_i$
and assuming that
$\mathbf{\dot{q}}= d\mathbf{q}/dt $, Newton's equations
in potential form
(\ref{Newton-thm}) are equivalent to the  {\bfi Euler-Lagrange equations\/}
\begin{equation}
\frac{d}{dt}\left(\frac{\partial L}{\partial\mathbf{\dot{q}}_i}\right)
-\frac{\partial L}{\partial \mathbf{q}_i}=0
\,, 
\qquad i = 1, \dots, N\,,
\label{Lag-eqs}
\end{equation}
where $\partial L/ \partial \mathbf{\dot{q}}_i , \partial L /\partial
\mathbf{q}_i \in \mathbb{R}^3$ denote the gradients in $\mathbb{R}^3$ of
$L$  with respect to $\mathbf{\dot{q}}_i , \mathbf{q}_i \in \mathbb{R}^3$.
\end{proposition}
}

\paragraph{\bfi Step I. Proof that Newton's equations (\ref{Newton-thm})
are equivalent to (i).\\} Check by direct verification. 

\paragraph{\bfi Step II. Proof that (i) $\Leftrightarrow$ (ii):}
{\bfi The Euler-Lagrange equations (\ref{EulerLag-eqs}) are equivalent
to Hamilton's principle of stationary action}.\\

To simplify notation, we momentarily suppress the particle index $i$. \\

We need to prove the solutions of \eqref{EulerLag-eqs} are critical
points $\delta\mathcal{S}=0$ of the {\bfi action functional}
\begin{equation}
\mathcal{S}[\mathbf{q}(\cdot)]: =\int_a^b
L(\mathbf{q}(t),\mathbf{\dot{q}}(t))\, dt\,,
\label{Action}
\end{equation}
(where $\mathbf{\dot{q}} = d\mathbf{q}(t)/dt$) with respect to
variations on $C^\infty([a,b],\mathbb{R}^{3N})$, the space of
smooth trajectories $\mathbf{q} : [a, b] \rightarrow
\mathbb{R}^{3N}$ with fixed endpoints $\mathbf{q}_a$,
$\mathbf{q}_b$.\\

In $C^\infty([a,b],\mathbb{R}^{3N})$ consider a {\bfi deformation}
$\mathbf{q}(t, s)$, $s \in (-\epsilon, \epsilon)$, $\epsilon>0 $, with
fixed endpoints $\mathbf{q}_a $, $\mathbf{q}_b$, of a curve
$\mathbf{q}_0(t)$, that is, $\mathbf{q}(t, 0)=\mathbf{q}_0(t)$ for all
$t \in [a,b]$ and $\mathbf{q}(a, s) = \mathbf{q}_0(a) = \mathbf{q}_a $,
$\mathbf{q}(b, s) = \mathbf{q}_0(b) = \mathbf{q}_b $ for all $s \in
(-\epsilon, \epsilon)$.\\

Define a {\bfi variation} of the curve $\mathbf{q}_0(\cdot)$ in
$C^\infty([a,b],\mathbb{R}^{3N})$ by
$$
\delta\mathbf{q}(\cdot):=\left.\frac{d}{ds}\right|_{s =
0}\mathbf{q}(\cdot, s) \in
T_{\mathbf{q}_0(\cdot)}C^\infty([a,b],\mathbb{R}^{3N}),
$$
and define the {\bfi first variation\/} of $\mathcal{S}$ at
$\mathbf{q}_0(t)$ to be the following derivative:
\begin{equation}
\delta\mathcal{S}
:=
\mathbf{D}\mathcal{S}[\mathbf{q}_0(\cdot)](\delta\mathbf{q}(\cdot))
:=\left.\frac{d}{ds}\right|_{s = 0} \mathcal{S}[\mathbf{q}(\cdot, s)].
\label{varA}
\end{equation}
Note that $\delta\mathbf{q}(a) = \delta\mathbf{q}(b) = \mathbf{0}$.
With these notations, {\bfi Hamilton's principle of stationary action\/}
states that the curve
$\mathbf{q}_0(t)$ satisfies the Euler-Lagrange equations
\eqref{EulerLag-eqs} if and only if
$\mathbf{q}_0(\cdot)$ is a  critical point of the action functional,
that is,
$\mathbf{D}\mathcal{S}[\mathbf{q}_0(\cdot)] = 0 $. Indeed, using  the
equality of mixed partials, integrating by parts, and taking into
account that $\delta\mathbf{q}(a) = \delta\mathbf{q}(b) = 0 $, one finds
\begin{align*}
\delta\mathcal{S}
:=
\mathbf{D}\mathcal{S}[\mathbf{q}_0(\cdot)](\delta\mathbf{q}(\cdot))&=
\left.\frac{d}{ds}\right|_{s = 0} \mathcal{S}[\mathbf{q}(\cdot, s)]
=\left.\frac{d}{ds}\right|_{s = 0} \int_a^b
L(\mathbf{q}(t, s),\mathbf{\dot{q}}(t,s))\, dt \\
&= \sum_{i=1}^N\int_a^b\left[\frac{\partial L}{\partial \mathbf{q}_i} 
\cdot\delta \mathbf{q}_i(t,s) 
+ \frac{\partial L}{\partial \mathbf{\dot{q}}_i} 
\cdot\delta \mathbf{\dot{q}}_i\right] dt\\
&= -\sum_{i=1}^N\int_a^b\left[
\frac{d}{dt}
\left(\frac{\partial L}{\partial
\mathbf{\dot{q}}_i}\right) 
-\frac{\partial L}{\partial \mathbf{q}_i}
\right]\cdot\delta\mathbf{q}_i \ dt=0
\end{align*}
for all smooth $\delta \mathbf{q}_i(t)$ satisfying 
$\delta\mathbf{q}_i(a)  = \delta\mathbf{q}_i(b) = 0 $.
This proves the equivalence of {\bf (i)} and
{\bf (ii)}, upon restoring particle index $i$ in the last two lines.

\paragraph{\bfi Step III. Proof that (ii) $\Leftrightarrow$ (iii):}
{\bfi Hamilton's principle of stationary action is equivalent to Hamilton's
canonical equations}.

\begin{definition}
The conjugate momenta for the Lagrangian in (\ref{Lagrangian-thm}) are
defined as
\begin{equation}
\mathbf{p}_i:=\frac{\partial L}{\partial \mathbf{\dot{q}}_i} = m_i
\mathbf{\dot{q}}_i\in
\mathbb{R}^3, \quad i=1,\ldots, N 
\,,\quad\hbox{(no sum on }i)
\label{momenta}
\end{equation}
\end{definition}

\begin{definition}
The {\bfi Hamiltonian} is defined via the change of variables
$(\mathbf{q},\mathbf{\dot{q}})\mapsto (\mathbf{q},\mathbf{p})$, called
the {\bfi Legendre transform}, 
\begin{align}
H(\mathbf{q},\mathbf{p}) :&=
\mathbf{p}\cdot\mathbf{\dot{q}}(\mathbf{q},\mathbf{p}) -
L(\mathbf{q},\mathbf{\dot{q}}(\mathbf{q},\mathbf{p}))
\nonumber\\
&=\sum_{i=1}^N\frac{m_i}{2}\|\mathbf{\dot{q}}_i\|^2+V(\mathbf{q})
\nonumber \\ 
&= 
\underbrace{\
\sum_{i=1}^N\frac{1}{2m_i}\|\mathbf{p}_i\|^2\
}_{\hbox{Kinetic energy}}
\
+\
\underbrace{\
V(\mathbf{q})\
}_{\hbox{Potential}}
\label{hamiltonian}
\end{align}
\end{definition}
\begin{remark}
The value of the Hamiltonian coincides with the  total energy of the
system. This value will be shown momentarily to remain constant under the
evolution of Euler-Lagrange equations (\ref{EulerLag-eqs}).
\end{remark}

\begin{remark}
The Hamiltonian $H$ may be obtained from the Legendre transformation as a
function of the variables $(\mathbf{q},\mathbf{p})$, provided one may solve
for $\mathbf{\dot{q}}(\mathbf{q},\mathbf{p})$, which requires the
Lagrangian to be {\bfi regular}, e.g.,
\[
\det\frac{\partial^2 L}
{\partial \mathbf{\dot{q}}_i\partial \mathbf{\dot{q}}_i}\ne0
\quad\hbox{(no sum on }i)
\,.
\] 
\end{remark}
Lagrangian (\ref{Lagrangian-thm}) is regular and the
derivatives of the Hamiltonian may be shown to satisfy,
\[
\frac{\partial H}{\partial \mathbf{p}_i}
= \frac{1}{m_i}\mathbf{p}_i
=\mathbf{\dot{q}}_i=\frac{d \mathbf{q}_i}{dt}
\quad\hbox{and}\quad
\frac{\partial H}{\partial \mathbf{q}_i}
= \frac{\partial V}{\partial \mathbf{q}_i}
= -\,\frac{\partial L}{\partial \mathbf{q}_i}
\,.
\]
Consequently, the Euler-Lagrange equations (\ref{EulerLag-eqs}) imply
\[
\mathbf{\dot{p}}_i = \frac{d \mathbf{p}_i}{dt} = \frac{d}{dt}\left(
\frac{\partial L}{\partial\mathbf{\dot{q}}_i} \right)
= \frac{\partial L}{\partial \mathbf{q}_i} 
= -\,\frac{\partial H}{\partial \mathbf{q}_i}\,.
\]
These calculations show that {\bfi the Euler-Lagrange equations
\eqref{EulerLag-eqs} are equivalent to Hamilton's canonical
equations}
\begin{equation}
\label{hameq}
\mathbf{\dot{q}}_i = \frac{\partial H}{\partial \mathbf{p}_i}\,, \qquad
\mathbf{\dot{p}}_i =-\,\frac{\partial H}{\partial \mathbf{q}_i}\,,
\end{equation}
where $\partial H/ \partial \mathbf{q}_i, \partial H / \partial
\mathbf{p}_i
\in \mathbb{R}^3 $ are the gradients of $H $ with respect to $\mathbf{q}_i,
\mathbf{p}_i \in \mathbb{R}^3$, respectively.
This proves the equivalence of {\bf (ii)} and {\bf (iii)}.

\begin{remark}
The Euler-Lagrange equations are second order and they determine
{\bfi curves in configuration space} $\mathbf{q}_i\in
C^\infty([a,b],\mathbb{R}^{3N})$. In contrast, Hamilton's equations are
first order and they determine {\bfi curves in phase
space} $(\mathbf{q}_i,\mathbf{p}_i)\in
C^\infty([a,b],\mathbb{R}^{6N})$, a space whose dimension is twice the
dimension of the configuration space.
\end{remark}

\paragraph{\bfi Step IV. Proof that (iii) $\Leftrightarrow$ (iv):}
{\bfi Hamilton's canonical equations may be written using a Poisson
bracket.}\\  

By the chain rule and (\ref{hameq}) any $F\in{\cal
F}(P)$ satisfies,
\begin{align*}
\frac{dF}{dt}
&=
\sum_{i =1}^N \left(\frac{\partial F}{\partial \mathbf{q}_i}\cdot
\mathbf{\dot{q}}_i +\frac{\partial F}{\partial \mathbf{p}_i} \cdot 
\mathbf{\dot{p}}_i\right)\\
&
= \sum_{i =1}^N \left(\frac{\partial F}{\partial
\mathbf{q}_i} \cdot \frac{\partial H}{\partial \mathbf{p}_i}
-\frac{\partial F}{\partial \mathbf{p}_i}
\cdot \frac{\partial H}{\partial \mathbf{q}_i}\right)
=\{F,H\}
\,.
\end{align*}
This finishes the proof of the theorem, by proving the
equivalence of {\bf (iii)} and {\bf (iv)}. 

\begin{remark}[Energy conservation]
Since the Poisson bracket is skew symmetric, $\{H,F\}=-\,\{F,H\}$, one
finds that $\dot{H}=\{H,H\}=0$. Consequently, the value of the Hamiltonian
is preserved by the evolution. Thus, the Hamiltonian is said to be a {\bfi
constant of the motion}.
\end{remark}

\begin{exercise}
Show that the Poisson bracket is bilinear, skew symmetric, satisfies the
Jacobi identity and acts as derivation on products of functions in phase
space.
\end{exercise}

\begin{exercise}
Given two constants of motion, what does the Jacobi identity imply about
additional constants of motion?
\end{exercise}

\begin{exercise}
Compute the Poisson brackets among
\[
J_i=\epsilon_{ijk}p_jq_k
\]
in Euclidean space. What Lie algebra do these Poisson brackets recall to
you?
\end{exercise}

\begin{exercise}
Verify that HamiltonÕs equations determined by the function
$\langle J(z), \xi\rangle 
= \xi \cdot (\mathbf{q}\times\mathbf{p})$ give infinitesimal 
rotations about the $\xi-$axis.
\end{exercise}

\section{Hamilton's principle on manifolds}
\label{Ham princ manifolds}

\begin{theorem}
[{\bf Hamilton's Principle of Stationary Action}] \label{teol331}
Let the smooth function $L:TQ\to\mathbb{R}$ be a Lagrangian on $TQ$. A
$C^2$ curve
$c: [a , b]\rightarrow Q$ joining $q_a=c(a)$ to $q_b=c(b)$ satisfies the
 Euler-Lagrange equations if and only if
\[\delta\int_a^b L(c(t),{\dot c}(t)) dt =0.\]
\end{theorem}
\begin{proof}
The meaning of the variational derivative in the statement is the
following. Consider a family of $C^2 $ curves $c(t,s)$ for
$|{s}| < \varepsilon$ satisfying $c_0(t) = c(t)$, $c(a,s) =
q_a $, and $c(b,s) = q_b $ for all ${s} \in (- \varepsilon,
\varepsilon)$. Then
\[
\delta\int_a^b L(c(t),{\dot c}(t)) dt : = \left.
\frac{d}{d{s}}\right|_{{s}= 0 }\int_a^b L(c(t,s),\dot{c}(t,s)) dt.
\]
Differentiating under the integral sign, working in local coordinates
(covering the curve $c(t) $ by a finite number of coordinate charts),
integrating by parts, denoting
\[
v(t): =\left.\frac{d}{d{s}}\right|_{{s}=0 } c(t,s),
\]
and taking into account that $v(a) = v(b) = 0 $,
yields
\[
\int_a^b\left(\frac{\partial L}{ \partial q^i} v^i + \frac{\partial L}{
\partial \dot{q}^i}\dot{v}^i  \right) dt
= \int_a^b\left(\frac{\partial L}{ \partial q^i} - \frac{d}{dt}
\frac{\partial L}{ \partial \dot{q}^i} \right) v^i dt.
\]
This vanishes for any $C^1 $ function $v(t) $ if and only if the
Euler-Lagrange equations hold.
\end{proof}

\begin{remark} 
The integral appearing in this theorem
\[
\mathcal{S}(c(\cdot)): = \int_a^b L(c(t),{\dot c}(t)) dt
\]
is called the {\bfi  action integral\/}. It is defined on $C^2$ curves
$c: [a, b] \rightarrow Q $ with fixed endpoints, $c(a ) = q_a $ and
$c(b) = q_b $.
\end{remark} 

\begin{remark} {\bfi Variational derivatives of functionals vs Lie
derivatives of functions}. The variational derivative of a functional 
$S[u]$ is defined as the linearization
\[
\lim_{\epsilon\to0}\frac{S[u+\epsilon v]-S[u]}{\epsilon}
=
\frac{d}{d\epsilon}\Big|_{\epsilon=0}S[u+\epsilon v]
=
\Big\langle \frac{\delta S}{\delta v}\,,v\Big\rangle
\,.\quad
\]
Compare this to the expression for the Lie derivative of a function.
If $f$ is a real valued function on a manifold $M$ and $X$ is a vector
field on $M$, the Lie derivative of $f$ along $X$ is defined as the
directional derivative  
\[
\mathcal{L}_Xf=X(f):=\mathbf{d}f\cdot X
\,.
\]
If $M$ is finite-dimensional, this is
\[
\mathcal{L}_Xf=X[f]:=\mathbf{d}f\cdot X
=
\frac{\partial f}{\partial x^i}X^i
=
\lim_{\epsilon\to0}\frac{f(x+\epsilon X)-f(x)}{\epsilon}
\,.
\]
The similarity is suggestive: Namely, the Lie derivative of a function and
the variational derivative of a functional are both defined as
linearizations of smooth maps in certain directions.

\end{remark}

The next theorem emphasizes the role of Lagrangian
one-forms and two-forms in the variational principle. The following is a
direct corollary of the previous theorem.

\begin{theorem}\label{EL-manifold}
Given a $C^k$ Lagrangian $L:TQ\rightarrow \mathbb{R}$ for $k\geq 2$,
there exists a unique $C^{k-2}$ map ${\cal EL}(L):\ddot{Q}\rightarrow
T^*Q$, where
\[\ddot{Q}:=\left\{\frac{d^2q}{dt^2}\Big|_{t=0}\in T(TQ) \,\Big| \,
q(t) \text{ is a }\, C^2 \text{ curve in } Q \right\}\]
is a  submanifold of $T(TQ)$, and
a unique $C^{k-1}$ one-form $\Theta_L\in\Lambda^1(TQ)$, such that
for all $C^2$ variations $q(t,s)$ (defined on a fixed
$t$-interval) of $q(t,0)=q_0(t):=q(t)$, we have
\begin{eqnarray}\label{Action-var}
\delta \mathcal{S}
&:=&
\frac{d}{ds}\Big|_{s=0}\mathcal{S}[c(\,\cdot\,,s)]
=
\mathbf{D}\mathcal{S}[q(\cdot)]\cdot \delta q (\cdot)
\nonumber\\&=&
\int_a^b
{\cal EL}(L)\left(q, \dot{q}, \ddot{q}\right)\cdot\delta q\, dt\
+\!\!\underbrace{\
\left.\Theta_L\left(q,\dot{q}\right)\cdot\delta q\Big|_a^b \right.\
}_{\hbox{\bfi cf. Noether Thm}}
\end{eqnarray}
where
\[\delta
q(t)=\left.\frac{d}{d s}\right|_{s = 0}
q(t,s)
\,.
\]
\end{theorem}

\newpage

\section{Summary Handout  for Differential Forms}
\label{diff-form-review}

\noindent \textbf{Vector fields and $1$-forms}

\medskip

Let $M$ be a manifold. In what follows, all maps may be assumed to be $C^\infty,$ although that's
not necessary.

A \textbf{vector field} on $M$ is a map $X:M\to TM$ such that $X(x)\in T_xM$ for every $x\in M.$
The set of all smooth vector fields on $M$ is written $\mathfrak{X}(M).$
(``Smooth'' means differentiable or $C^r$ for some $r\le \infty,$ depending on context.)

A \textbf{(differential) $1$-form} on $M$ is a map $\theta:M\to T^*M$ such that $\theta(x)\in T_x^*M$ for every $x\in M.$

More generally, if $\pi:E\to M$ is a bundle, then a \textbf{section} of the bundle is a map $\varphi:M\to E$
such that $\pi\circ \varphi(x)=x$ for all $x\in M.$ So a vector field is a section of the tangent bundle,
while a $1$-form is section of the cotangent bundle.

Vector fields can added and also multiplied by scalar functions $k:M\to \R,$ as follows:
$\left(X_1 + X_2\right) (x) = X_1(x) + X_2(x), \ 
\left(kX\right)(x) = k(x) X(x).$

Differential forms can added and also multiplied by scalar functions $k:M\to \R,$ as follows:
$\left(\alpha + \beta\right) (x) = \alpha(x) + \beta(x), \
\left(k\theta\right)(x) = k(x) \theta(x).$

We have already defined the push-forward and pull-back of a vector field.
The \textbf{pull-back} of a $1$-form $\theta$ on $N$ by a map $\varphi:M\to N$ is
the $1$-form $\varphi^* \theta$ on $M$ defined by 
\[
\left(\varphi^* \theta \right)(x) \cdot v = \theta\left(\varphi(x) \right) \cdot T\varphi(v)
\]
The \textbf{push-forward} of a $1$-form $\alpha$ on $M$
by a diffeomorphism $\psi:M\to N$ is the pull-back of $\alpha$ by $\psi^{-1}.$

A vector field can be \textbf{contracted} with a differential form, using the pairing between
tangent and cotangent vectors: $\left(X\contract \theta\right)(x) = \theta(x)\cdot X(x).$
Note that $X\contract \theta$ is a map from $M$ to $\R.$
Many books write $i_X\theta$ in place of $X\contract \theta,$ and the contraction
operation is often called \textbf{interior product}.

\medskip

The \textbf{differential} of $f:M\to\R$ is a $1$-form $df$ on $M$ defined by
\[
\displaystyle df(x)\cdot v = \left.\frac{d}{dt} f(c(t)\right|_{t=0}
\]
for any $x\in M,$ any $v\in T_xM$ and any path $c(t)$ in $M$ such that
$c(0)=0$ and $c'(0)=v.$ The left hand side, $df(x)\cdot v,$ means the pairing
between cotangent and tangent vectors, which could also be written $df(x)(v)$ or 
$\left<df(x),v\right>.$

Note:
\[
X\contract df = \pounds_X f = X[f]
\]

\begin{remark} $df$ is very similar to $Tf,$
but $Tf$ is defined for all differentiable $f:M\to N,$ whereas $df$ is only defined when $N=\R$
(in this course, anyway).
In this case, $Tf$ is a map from $TM$ to $T\R,$ and
$Tf(v) = df(x)\cdot v \in T_{f(x)}\R$ for every $v\in T_xM$
(we have identified $T_f(x)\R$ with $\R.$) 
\end{remark}

\noindent \textbf{In coordinates...}
Let $M$ be $n$-dimensional, and 
let $x^1,\dots,x^n$ be differentiable local coordinates for $M.$ This means that there's an open subset $U$ of $M$ and an open subset $V$ of $\R^n$ such that the map $\varphi:U\to V$ defined by $\varphi(x) = \left(x^1(x),\dots,x^n(x)\right)$ is a diffeomorphism. In particular, 
each $x^i$ is a map from $M$ to $\R,$ 
so the differential $dx^i$ is defined.
There is also a vector field $\frac{\partial}{\partial x^i}$ for every $i,$
which is defined by 
$\frac{\partial}{\partial x^i}(x) = \left.\frac{d}{dt} \varphi^{-1} \left(\varphi(x) + t \mathbf{e}_i\right)\right|_{t=0},$
where $\mathbf{e}_i$ is the $i^{\textrm{th}}$ standard basis vector.

\begin{exercise} Verify that 
\[
\frac{\partial}{\partial x^i} \contract dx^j \equiv \delta^i_j 
\]
(where $\equiv$ means the left hand side is a constant function with value $\delta^i_j$)
\end{exercise}

\begin{remark} Of course, given a coordinate system $\varphi=\left(x^1,\dots,x^n\right),$ it is usual to write $x= \left(x^1,\dots,x^n\right),$ which means $x$ is identified with $\left(x^1(x),\dots,x^n(x)\right)=\varphi(x).$
\end{remark}

For every $x\in M,$ the vectors $\frac{\partial}{\partial x^i}(x)$ form a basis for $T_xM,$
so every $v\in T^xM$ can be uniquely expressed as $v=v^i \frac{\partial}{\partial x^i}(x).$
This expression defines the \textbf{tangent-lifted coordinates} $x^1,\dots,x^n,v^1,\dots v^n$
on $TM$ (they are local coordinates, defined on $TU\subset TM$).

For every $x\in M,$ the covectors $dx^i(x)$ form a basis for $T_x^*M,$
so every $\alpha \in T^xM$ can be uniquely expressed as $\alpha=\alpha_i dx^i(x).$
This expression defines the \textbf{cotangent-lifted coordinates} $x^1,\dots,x^n,\alpha_1,\dots \alpha_n$
on $T^*M$ (they are local coordinates, defined on $T^*U\subset T^*M$).

Note that the basis $\left( \frac{\partial}{\partial x^i}\right)$ is
dual to the basis $\left(dx^1,\dots,dx^n\right),$ by the previous exercise.
It follows that, 
\[
\left(\alpha_i dx^i\right) \cdot \left(v^i \frac{\partial}{\partial x^i}\right)
= \alpha_i v^i 
\]
(we have used the summation convention).

In mechanics, the configuration space is often called $Q,$ and the lifted coordinates are written:
$q^1,\dots q^n,\dot{q}^1,\dots,\dot{q}^n$ (on $TQ$) and 
$q^1,\dots q^n,p_1,\dots,p_n$ (on $T^*Q$).

\noindent \textbf{Why the distinction between subscripts and superscripts?} This is to keep track
of how quantities vary if coordinates are changed (see next exercise).
One benefit is that using the summation convention gives
coordinate-independent answers.

\rem{
\begin{exercise}
Consider two sets of local coordinates $q^i$ and $s^i$ on $Q,$ related 
by
$\left(s^1,\dots, s^n\right) = \psi \left(q^1,\dots, q^n\right).$
Verify that the corresponding tangent lifted coordinates $\dot{q}^i$ 
and $\dot{s}^i$
are related by
\[
\dot{s}^i =\left(D\psi\right)_{ij} \dot{q}^j\, .
\]
Note that the last equation can be written as $\mathbf{\dot s} = 
D\psi(q) \mathbf{\dot{q}},$
where $\mathbf{\dot s}$ is the column vector $(\dot s^1,\dots \dot 
s^n).$
}

\begin{exercise}
Consider two sets of local coordinates $q^i$ and $s^i$ on $Q,$ related by
$\left(s^1,\dots, s^n\right) = \psi \left(q^1,\dots, q^n\right).$
Verify that the corresponding tangent lifted coordinates $\dot{q}^i$ and $\dot{s}^i$
are related by
\[
\dot{s}^i =\frac{\partial \psi^i}{\partial q^j} \dot{q}^j\, .
\]
Note that the last equation can be written as $\mathbf{\dot s} = D\psi(q) \mathbf{\dot{q}},$
where $\mathbf{\dot s}$ is the column vector $(\dot s^1,\dots \dot s^n)$,
and similarly for $\mathbf{\dot q}$.\\ 

Do the corresponding calculation on the cotangent bundle side. See
{\bfi Definition \ref{cot-lift}}.
\end{exercise}

\noindent\textbf{The next level: $TTQ, T^*T^*Q,$ et cetera}

Since $TQ$ is a manifold, we can consider vector fields on it, which are sections of $T(TQ).$
In coordinates, every vector field on $TTQ$ has the 
form $X=a^i \frac{\partial}{\partial q^i} + b^i \frac{\partial}{\partial \dot q^i},$
where the $a^i$ and $b^i$ are functions of $q$ and $\dot q.$
Note that the same symbol $q^i$ has two interpretations: as a coordinate on $TQ$ and as a 
coordinate on $Q,$ so $\frac{\partial}{\partial q^i}$ can mean a vector field $TQ$ (as above)
or on $Q$.

The tangent lift of the bundle projection $\tau:TQ\to Q$ is a map $T\tau:TTQ\to TQ.$
If $X$ is written in coordinates as above, then $T\tau\circ X= a^i \frac{\partial}{\partial q^i}.$
A vector field $X$ on $TTQ$ is \textbf{second order} if $T\tau \circ X (v) = v;$ in coordinates,
$a^i = \dot q^i.$ The name comes from the process of reducing of second order equations to
first order ones by introducing new variables $\displaystyle \dot q^i = \frac{d q^i}{dt}.$

One may also consider $T^*TQ,TT^*Q$ and $T^*T^*Q.$ However, the subscript/superscript distinction is problematic here. 

\bigskip

\noindent \textbf{$1$-forms}

The $1$-forms on $T^*Q$ are sections of $T^*T^*Q.$ Given cotangent-lifted
local coordinates
\[\left(q^1,\dots,q^n,p_1,\dots,p_n\right)\] on $T^*Q,$ the general
$1$-form on $T^*Q$ has the form
$a_i dq^i + b_i dp_i,$ where $a_i$ and $b_i$ are functions of $(q,p).$
The \textbf{canonical $1$-form} on $T^*Q$ is 
\[\theta = p_i dq^i,\] 
also written
in the short form $p\ dq.$ 
Pairing $\theta(q,p)$ with an arbitrary tangent vector  
$\displaystyle v= a^i \frac{\partial}{\partial q^i} + b^i
\frac{\partial}{\partial p^i}\in T_{(q,p)}T^*Q$ gives
\begin{align*}
\left<\theta(q,p),v\right> &=
\left<p_i dq^i, a^i \frac{\partial}{\partial q^i} + b^i
\frac{\partial}{\partial p^i}\right> \\ &=p_i a^i \\
&= \left<p_i dq^i, a^i \frac{\partial}{\partial q^i} \right> \\
&= \left< p, T\tau^*(v)\right>,
\end{align*}
where $\tau^*:T^*Q\to Q$ is projection. In the last line we have
interpreted $q^i$ as a coordinate on $Q,$ which implies that $p_i dq^i =
p,$ by definition of the coordinates $p_i$. Note that the last line is
coordinate-free.

\bigskip

\noindent \textbf{$2$-forms}

Recall that a $1$-form on $M,$ evaluated at a point $x \in M,$ is a
linear map from $T_xM$ to $\R.$ \\ A $2$-form on $M,$ evaluated at a 
point
$x \in M,$ is a skew-symmetric bilinear form on $T_xM$; and the
bilinear form has to vary smoothly as $x$ changes.
(Confusingly, bilinear forms can be skew-symmetric, symmetric or 
neither;
\emph{differential} forms are assumed to be skew-symmetric.)

The \textbf{pull-back} of a $2$-form $\omega$ on $N$ by a map
$\varphi:M\to N$ is the $2$-form $\varphi^* \omega$ on $M$ defined by 
\[
\left(\varphi^* \omega \right)(x) \left( v,w\right) =
\theta\left(\varphi(x) \right) 
\left( T\varphi(v), T\varphi(w)\right)
\]
The \textbf{push-forward} of a $2$-form $\omega$ on $M$
by a diffeomorphism $\psi:M\to N$ is the pull-back of $\omega$ by
$\psi^{-1}.$

 A vector field $X$ can be \textbf{contracted} with a $2$-form $\omega$ to
get a $1$-form 
$X\contract \omega$ defined by
\[
\left(X\contract \omega\right) (x) (v) = \omega(x)\left(X(x),v\right)
\]
for any $v\in T_xM.$
A shorthand for this is $\left(X\contract \omega\right) (v)=
\omega(X,v),$ or just
$X\contract \omega = \omega(X,\cdot).$

The \textbf{tensor product} of two $1$-forms $\alpha$ and $\beta$ is the
$2$-form $\alpha \otimes \beta$ defined by
\[
\left(\alpha \otimes \beta\right)(v,w) = \alpha(v) \beta(w)
\]
for all $v,w\in T^*_xM.$

The \textbf{wedge product} of two $1$-forms $\alpha$ and $\beta$ is the skew-symmetric $2$-form 
$\alpha \wedge \beta$ defined by
\[
\left(\alpha \wedge \beta\right)(v,w) = \alpha(v) \beta(w) - \alpha(w) \beta(v)  \ .
\]

\bigskip

\noindent \textbf{Exterior derivative}

The differential $df$ of a real-valued function is also called the exterior derivative of $f.$
In this context, real-valued functions can be called $0$-forms.
The exterior derivative is a linear operation from 
$0$-forms to $1$-forms that satisfies the Leibniz identity, a.k.a. the product rule,
\[
d(fg) = f\ dg + g \ df
\]

The exterior derivative of a $1$-form is an alternating $2$-form, defined as follows:
\[
d\left(a_i dx^i\right) = \frac{\partial a_i}{\partial x^j} dx^j \wedge dx^i .
\]
Exterior derivative is a linear operation from $1$-forms to $2$ forms.
The following identity is easily checked:
\[
d(df) = 0
\] 
for all scalar functions $f.$

\bigskip

\noindent \textbf{$n$-forms}

See \cite{MaRa1994}, or \cite{Le2003}, or \cite{AbMa1978}.
Unless otherwise specified, $n$-forms are assumed to be alternating. 
Wedge products and contractions generalise.

It is a fact that all $n$-forms are linear combinations of wedge products of $1$-forms.
Thus we can define exterior derivative recursively by the properties
\begin{align*}
d\left(\alpha \wedge \beta\right) = d\alpha \wedge \beta + (-1)^k \alpha \wedge d\beta,
\end{align*}
 for all $k$-forms $\alpha$ and all forms $\beta,$
 and 
 \[
d\circ d = 0
 \]
In local coordinates, if $\alpha = \alpha_{i_1\cdots i_k} dx^{i_1} \wedge \cdots \wedge dx^{i_k}$
(sum over all $i_1< \cdots < i_k$), then
\[
d\alpha = \frac{\partial \alpha_{i_1\cdots i_k}}{\partial x^j}  
dx^j \wedge dx^{i_1} \wedge \cdots \wedge dx^{i^k}
\]
 
The \textbf{Lie derivative} of an $n$-form $\theta$ in the direction of the vector field $X$ is defined as
\[
\pounds_X \theta = \left. \frac{d}{dt} \varphi_t^* \theta \right|_{t=0},
\]
where $\varphi$ is the flow of $X.$


Pull-back commutes with the operations $d,\contract,\wedge$ and Lie derivative.

\textbf{Cartan's magic formula}:
\[
\pounds _X \alpha = d\left(X\contract \alpha\right) + X\contract d\alpha
\] 
This looks even more magic when written using the notation $i_X\alpha = X\contract \alpha:$
\[
\pounds_X = d i_X + i_X d
\]

An $n$-form $\alpha$ is \textbf{closed} if $d\alpha=0, $ and \textbf{exact} 
if $\alpha = d\beta$ for some $\beta.$ All exact forms are closed (since
$d\circ d = 0$), but the converse is false. It is true that all closed
forms are \emph{locally} exact; this is the \textbf{Poincar\'e Lemma}.

\begin{remark}
For a survey of the basic definitions, properties, and operations on
differential forms, as well as useful of tables of relations between
differential calculus and vector calculus, see, {\it e.g.}, Chapter 2 of
\cite{Bl2004}.
\end{remark}

\newpage

\section{Euler-Lagrange equations of manifolds}
In Theorem \ref{EL-manifold}, 
\begin{eqnarray}\label{Action-var}
\delta \mathcal{S}
&:=&
\frac{d}{ds}\Big|_{s=0}\mathcal{S}[c(\,\cdot\,,s)]
=
\mathbf{D}\mathcal{S}[q(\cdot)]\cdot \delta q (\cdot)
\nonumber\\&=&
\int_a^b
{\cal EL}(L)\left(q, \dot{q}, \ddot{q}\right)\cdot\delta q\, dt\
+\!\!\underbrace{\
\left.\Theta_L\left(q,\dot{q}\right)\cdot\delta q\Big|_a^b \right.\
}_{\hbox{\bfi cf. Noether Thm}}
\end{eqnarray}
where
\[\delta
q(t)=\left.\frac{d}{d s}\right|_{s = 0}
q(t,s)
\,,
\]
the map $\mathcal{EL}: \ddot{Q}\rightarrow T^*Q$ is called the {\bfi
Euler-Lagrange operator\/} and its expression  in local coordinates
is
\[
\mathcal{EL}(q, \dot{q}, \ddot{q})_{\,i} =
\frac{\partial L}{ \partial q^i} - \frac{d}{dt}
\frac{\partial L}{ \partial \dot{q}^i}
\,.
\]
One understands that the formal time derivative is taken in the
second summand and everything is expressed as a function of 
$(q, \dot{q}, \ddot{q})$.

\begin{theorem}{\bfi Noether (1918) Symmetries and Conservation Laws} 
If the action variation in (\ref{Action-var}) vanishes
$\delta\mathcal{S}=0$ because of a symmetry transformation which does
{\bfi not} preserve the end points and the Euler-Lagrange equations hold,
then the term marked {\bfi cf. Noether Thm} must also vanish. However,
vanishing of this term now is interpreted as a constant of motion. Namely,
the term,
\[
A(v,w):= \langle \mathbb{F} L(v),  w\rangle
\,,\quad\hbox{or, in coordinates}\quad 
A(q,\dot{q},\delta{q})
=
\frac{\partial  L}{\partial \dot{q}^i}\,\delta{q}^i
\,,
\] 
is constant for solutions of the Euler-Lagrange equations. This result
first appeared in Noether \cite{No1918}. In fact, the result in
\cite{No1918} is more general than this. In particular, in the PDE
(Partial Differential Equation) setting one must also include the
transformation of the volume element in the action principle. See {\it
e.g.} \cite{Ol2000} for good discussions of the history, framework and
applications of Noether's theorem.
\end{theorem}

\begin{exercise}
Show that {\bfi conservation of energy results from Noether's theorem} if,
in HamiltonÕs principle, the variations are chosen as
\[\delta
q(t)=\left.\frac{d}{d s}\right|_{s = 0}
q(t,s)
\,,
\]
corresponding to symmetry of the Lagrangian under reparametrizations of
time along the given curve $q(t)\to q(\tau(t,s))$.
\end{exercise}

\paragraph{The canonical Lagrangian one-form and two-form.}The one-form
$\Theta_L$, whose existence and uniqueness is guaranteed by Theorem
\ref{EL-manifold}, appears as the boundary term of the derivative of the
action integral, when the endpoints of the curves on the configuration
manifold are free. In finite dimensions, its local expression is
\[
\Theta_L\left(q,\dot{q}\right)
:=
\frac{\partial L}{\partial \dot{q}^i}\,\mathbf{d}q^i
\Big(=p_i\left(q,\dot{q}\right)\,\mathbf{d}q^i\Big)
\,.
\]
The corresponding closed two-form $\Omega_L=\mathbf{d}\Theta_L$ obtained
by taking its exterior derivative may be expressed as
\[
\Omega_L:=-\mathbf{d}\Theta_L=\frac{\partial^2 L}{\partial\dot{q}^i\partial
q^j}\mathbf{d}q^i\wedge \mathbf{d}q^j
 + \frac{\partial^2
L}{\partial\dot{q}^i\partial\dot{q}^j}\mathbf{d}q^i\wedge
\mathbf{d}\dot{q}^j
\Big(=\mathbf{d}p_i\left(q,\dot{q}\right)\wedge\mathbf{d}q^i\Big)\,.
\]
These coefficients may be written as the $2n \times 2n$
skew-symmetric matrix
\begin{equation}
\Omega_L=\left(\begin{array}{cc}
\mathcal{A} & \frac{\partial^2 L}{\partial\dot{q}^i\partial\dot{q}^j} \\
- \frac{\partial^2 L}{\partial\dot{q}^i\partial \dot{q}^j} &  0
\end{array}\right)
\,,
\label{none}
\end{equation}
where $\mathcal{A}$ is the skew-symmetric
$n \times n $ matrix 
$\left(\frac{\partial^2 L}{\partial\dot{q}^i\partial  q^j}\right)
-\left(\frac{\partial^2 L}{\partial\dot{q}^i\partial  q^j}\right)^T$. 
Non-degeneracy of
$\Omega_L$ is equivalent to the invertibility of the matrix
$\left(\frac{\partial^2 L}{\partial\dot{q}^i\partial\dot{q}^j}\right)$.

\begin{definition}
The {\bfi Legendre transformation\/} $\mathbb{F}
L:TQ\rightarrow T^*Q$  is the smooth map near the
identity defined by
\begin{equation*}
\langle \mathbb{F}L(v_q),  w_q \rangle :=
\left.\frac{d}{ds}\right|_{s=0}L(v_q+s w_q)\,.
\end{equation*}
\end{definition}
In the finite dimensional case, the local expression of $\mathbb{F}L $
is
\[
\mathbb{F} L (q^i,\dot q^i)=\left(q^i,\frac{\partial L}
{\partial \dot{q}^i}\right)=(q^i,p_i\left(q,\dot{q}\right)).
\]

If the skew-symmetric matrix (\ref{none}) is invertible, the Lagrangian
$L$ is said to be {\bfi regular\/}. In this case, by the implicit function 
theorem, $\mathbb{F}L$ is locally invertible. If $\mathbb{F}L$ is a
diffeomorphism,  $L$ is called {\bfi hyperregular}.\\

\begin{definition}
Given a Lagrangian $L$, the {\bfi
action\/} of $L$ is  the map $A : TQ
\rightarrow \mathbb{R}$ given by
\begin{equation}
\label{action of Lagrangian}
A(v):= \langle \mathbb{F} L(v),  v\rangle
\,,\quad\hbox{or, in coordinates}\quad 
A(q,\dot{q})=\frac{\partial  L}{\partial \dot{q}^i}\dot{q}^i
\,,
\end{equation}
and the  {\bfi energy\/} of $L$ is
\begin{equation}
\label{energy of Lagrangian}
E(v):= A(v)-L(v)
\,,\quad\hbox{or, in coordinates}\quad 
E(q,\dot{q})
=\frac{\partial  L}{\partial \dot{q}^i}\dot{q}^i
- L(q,\dot{q})
\,.
\end{equation}
\end{definition}

\subsection{Lagrangian vector fields and conservation laws}

\begin{definition}
A vector field $Z$ on $TQ$ is called a {\bfi Lagrangian vector field\/}
if
\[
\Omega_L(v)(Z(v),w)= \langle \mathbf{d}E(v), w \rangle,
\]
for all $v\in T_qQ$, $w\in T_v(TQ)$.
\end{definition}
\begin{proposition}
The {\bfi energy is conserved} along the flow of a Lagrangian vector field
$Z$.
\end{proposition}

\begin{proof} Let $v(t)\in TQ$ be an integral curve of $Z$. 
Skew-symmetry of $\Omega_L $ implies
\begin{align*}
\frac{d}{dt}E(v(t))&=\langle \mathbf{d}E(v(t)),  \dot v(t) \rangle
=\langle \mathbf{d}E(v(t)), Z(v(t))\rangle \\
&=\Omega_L(v(t)) \left(Z(v(t)),Z(v(t))\right) = 0.
\end{align*}
Thus, $E(v(t))$ is constant in $t$.
\end{proof}

\subsection{Equivalence of dynamics for hyperregular Lagrangians and
Hamiltonians}

Recall that a Lagrangian $L$ is said to be {\bfi hyperregular}
if its Legendre transformation $\mathbb{F}L:TQ\rightarrow T^*Q$  is a
diffeomorphism.

The equivalence between the Lagrangian and Hamiltonian formulations for
hyperregular Lagrangians and Hamiltonians is summarized below, following
\cite{MaRa1994}. 

\begin{enumerate}
\item[(a)]
Let $L$ be a  hyperregular Lagrangian on $TQ$ and $H=E\circ (\mathbb{F}
L)^{-1}$, where $E$ is the energy of $L$ and
$(\mathbb{F}L)^{-1}:T^*Q\rightarrow TQ$ is the inverse of the Legendre
transformation. Then the Lagrangian vector field 
$Z$ on $TQ$ and the Hamiltonian vector field
$X_H$ on $T^*Q$ are related by the identity
$$(\mathbb{F} L)^*X_H=Z .$$
Furthermore, if $c(t)$ is an
integral curve of $Z$ and $d(t)$ an integral curve of $X_H$ with
$\mathbb{F} L(c(0))=d(0)$, then $\mathbb{F} L(c(t))=d(t)$  and their
integral curves coincide on the manifold $Q$. That is,
$\tau_Q(c(t))=\pi_Q(d(t))=\gamma(t)$, where $\tau_Q: TQ \rightarrow Q $
and $\pi_Q: T ^\ast Q \rightarrow Q $ are the canonical bundle projections.

In particular,  the pull back of the inverse 
Legendre transformation $\mathbb{F} L^{-1}$ induces a
one-form
$\Theta$ and a closed two-form $\Omega $ on $T^*Q$ by
\begin{eqnarray*}
\Theta=(\mathbb{F} L^{-1})^*\Theta_L
\,,\qquad
\Omega=-\, \mathbf{d} \Theta=(\mathbb{F}L^{-1})^*\Omega_L
\,.
\end{eqnarray*}
In coordinates, these are the canonical presymplectic and
symplectic forms, respectively,
\begin{eqnarray*}
\Theta=p_i\,\mathbf{d}q^i
\,,\qquad
\Omega=-\, \mathbf{d} \Theta=\mathbf{d}p_i\wedge\mathbf{d}q^i
\,.
\end{eqnarray*}

\item[(b)]
A  Hamiltonian $H:T^*Q \rightarrow \mathbb{R}$ is said to be
{\bfi hyperregular\/} if the smooth map $\mathbb{F} H:T^*Q\rightarrow TQ$,
defined by
\[
\langle  \mathbb{F}H(\alpha_q),  \beta_q\rangle
:= \left.\frac{d}{ds}\right|_{s=0}
H(\alpha_q+s\beta_q),\qquad
\alpha_q, \beta_q \in T^*_qQ,\]
is a diffeomorphism.
Define the {\bfi action\/} of $H$ by $G := \langle \Theta\, ,\,
X_H\rangle$. If $H$ is a hyperregular Hamiltonian then the energies of
$L$ and $H$ and the actions of $L$ and $H$ are related by
 $$E=H\circ (\mathbb{F} H)^{-1},\qquad \quad A=G\circ
 (\mathbb{F} H)^{-1}.$$
Also, the Lagrangian $L = A - E$ is hyperregular and  
$\mathbb{F} L = \mathbb{F} H^{-1}$.

\item[(c)]
These constructions define a bijective
correspondence between hyperregular Lagrangians and Hamiltonians.

\end{enumerate}

\begin{remark}
For thorough discussions of many additional results arising from the
Hamilton's principle for hyperregular Lagrangians see, {\it e.g.} Chapters
7 and 8 of \cite{MaRa1994}.
\end{remark}

\begin{exercise}[Spherical pendulum]
A particle rolling on the interior of a spherical surface under gravity is
called a spherical pendulum. Write down the Lagrangian and the equations
of motion for a spherical pendulum with $S^2$ as its configuration space.
Show explicitly that the Lagrangian is hyperregular. Use the Legendre
transformation to convert the equations to Hamiltonian form. Find the
conservation law corresponding to angular momentum about the axis of
gravity by ``bare hands'' methods.
\end{exercise}

\begin{exercise}[Differentially rotating frames]\label{rot-ex}
The Lagrangian for a free particle of unit mass relative to a moving frame
is obtained by setting 
\[
L(\mathbf{\dot{q}},\mathbf{q},t)
=
\frac{1}{2}\|\mathbf{\dot{q}}\|^2 + \mathbf{\dot{q}}\cdot\mathbf{R}(\mathbf{q},t)
\]
for a function $\mathbf{R}(\mathbf{q},t)$ which prescribes the space and
time dependence of the moving frame velocity. For example, a frame
rotating with time-dependent frequency $\Omega(t)$ about the vertical axis
$\mathbf{\hat{z}}$ is obtained by choosing
$\mathbf{R}(\mathbf{q},t)=\mathbf{q}\times\Omega(t)\mathbf{\hat{z}}$.
Calculate $\Theta_L\left(q,\dot{q}\right),\Omega_L\left(q,\dot{q}\right)$,
the Euler-Lagrange operator ${\cal EL}(L)\left(q,\dot{q},
\ddot{q}\right)$, the Hamiltonian and its corresponding canonical
equations.
\end{exercise}

\begin{exercise}
Calculate the action and the energy for the Lagrangian in Exercise
\ref{rot-ex}.
\end{exercise}

\begin{definition}[Cotangent lift]\label{cot-lift}
Given two manifolds $Q$ and $S$ related by a diffeomorphism
$f:Q\mapsto S$, the {\bfi cotangent lift} $T^*f: T^*S\mapsto T^*Q$
of $f$ is defined by 
\begin{equation}
\langle T^*f(\alpha),v\rangle
=
\langle \alpha,Tf(v)\rangle
\end{equation}
where 
\[
\alpha\in T^*_sS
\,,\quad
v\in T_qQ
\,,\quad\hbox{and}\quad
s=f(q)
\,.
\]
As explained in Chapter 6 of \cite{MaRa1994}, cotangent lifts preserve the
{\bfi action} of the Lagrangian $L$, which we write as
\begin{equation}\label{action-pair}
\langle \mathbf{p}\,,\,\mathbf{\dot{q}}\rangle
=
\langle \alpha\,,\,\mathbf{\dot{s}}\rangle
\,,
\end{equation}
where $\mathbf{p}=T^*f(\alpha)$ is the cotangent lift of $\alpha$ under
the diffeomorphism $f$ and $\mathbf{\dot{s}}=Tf(\mathbf{\dot{q}})$ is the
tangent lift of $\mathbf{\dot{q}}$ under the function $f$, which is
written in Euclidean coordinate components as $q^i\to s^i=f^i(\mathbf{q})$.
Preservation of the action in (\ref{action-pair}) yields the
coordinate relations,
\begin{eqnarray*}
\hbox{\rm (Tangent lift in coordinates)}&&
\dot{s}^j=\frac{\partial f^j}{\partial q^i}\dot{q}^i
\qquad\Longrightarrow\\
p_i=\alpha_k\frac{\partial f^k}{\partial q^i}&&
\quad\hbox{\rm (Cotangent lift in coordinates)}
\end{eqnarray*}
Thus, in coordinates, the cotangent lift is the inverse transpose of the
tangent lift.
\end{definition}

\begin{remark}
The cotangent lift of a function preserves the induced action one-form,
\begin{eqnarray*}\label{action-pair1}
\langle \mathbf{p}\,,\,\mathbf{d}\mathbf{q}\rangle
=
\langle \alpha\,,\,\mathbf{d}\mathbf{s}\rangle
\,,
\end{eqnarray*}
so it is a source of (pre-)symplectic transformations.
\end{remark}

\subsection{The classic Euler-Lagrange example: Geodesic flow}
An important example of a Lagrangian vector field is the geodesic spray
of a Riemannian metric. A {\bfi Riemannian
manifold\/} is a smooth manifold $Q$ endowed with a symmetric
nondegenerate covariant tensor $g$, which is positive definite. Thus, on
each tangent space $T_q Q$ there is a nondegenerate definite inner product
defined by pairing with $g(q)$.

If $(Q, g)$ is a Riemannian manifold, there is a natural
Lagrangian on it given by the {\bfi  kinetic energy\/} $K$ of the
metric $g$, namely,
$$K(v) : = \frac{1}{2} g(q)(v_q,v_q),
$$
for $q \in Q $ and $v_q \in T_q Q $. In
finite dimensions, in a local chart,
$$K(q, \dot{q})= \frac{1}{2} g_{ij}(q) \dot q^i \dot q^j. $$
The Legendre transformation is in this case $\mathbb{F} K
(v_q) = g(q) (v_q, \cdot)$, for $v_q \in T_qQ$. In coordinates, this is 
\[
\mathbb{F} K (q,\dot q)=\left(q^i,\frac{\partial K}
{\partial \dot{q}^i}\right)=(q^i,g_{ij}(q) \dot q^j)=:(q^i,p_i).
\]
The Euler-Lagrange equations become the {\bfi geodesic
equations\/} for  the metric $g$, given (for finite dimensional $Q $ in
a local chart) by
\[\ddot{q}^i+\Gamma_{jk}^i\dot q^j\dot q^k=0,\quad  i
= 1 ,\ldots n,\]
where the three-index quantities
\[
\Gamma_{jk}^h=\frac{1}{2}g^{hl}\left(\frac{\partial g_{jl}}{\partial
q^k}+ \frac{\partial g_{kl}}{\partial q^j}-\frac{\partial
g_{jk}}{\partial q^l}\right)
,\quad\hbox{with}\quad
g_{ih}g^{hl}=\delta_i^l
\,,
\] are the {\bfi Christoffel symbols\/}
 of the Levi-Civita connection on $(Q,g)$.

\begin{exercise}
Explicitly compute the geodesic equation as an
Euler-Lagrange equation for the kinetic energy Lagrangian $K(q, \dot{q})=
\frac{1}{2} g_{ij}(q)
\dot q^i
\dot q^j$.
\end{exercise}

\begin{exercise}
For kinetic energy Lagrangian $K(q, \dot{q})=
\frac{1}{2} g_{ij}(q) \dot q^i\dot q^j$ with $i,j=1,2,\dots,N$:
\begin{itemize}
\item
Compute the momentum $p_i$ canonical to $q^i$ for geodesic motion. 
\item
Perform the Legendre transformation to obtain the Hamiltonian for geodesic
motion.
\item
Write out the geodesic equations in terms of $q^i$ and its canonical
momentum $p_i$.
\item
Check directly that Hamilton's equations are satisfied.
\end{itemize}
\end{exercise}

\begin{remark}
A classic problem is to determine the metric tensors $g_{ij}(q)$ for
which these geodesic equations admit enough additional
conservation laws to be integrable.
\end{remark}

\begin{exercise}
Consider the Lagrangian
\[
L_\epsilon(\mathbf{q},\mathbf{\dot{q}})
= 
\frac{1}{2} \|\mathbf{\dot{q}}\|^2
-
\frac{1}{2\epsilon}
(1 - \|\mathbf{q}\|^2)^2
\]
for a particle in $\mathbb{R}^3$. Let $\gamma_\epsilon(t)$ be the curve in
$\mathbb{R}^3$ obtained by solving the Euler-Lagrange equations for
$L_\epsilon$ with the initial conditions
$\mathbf{q}_0=\gamma_\epsilon(0),\,
\mathbf{\dot{q}}_0 = \dot{\gamma}_\epsilon(0)$. Show that 
\[\lim_{\epsilon\to0}\gamma_\epsilon(t)\]
is a great circle on the two-sphere $S^2$, provided that $\mathbf{q}_0$ has
unit length and the initial conditions satisfy
$\mathbf{q}_0\cdot\mathbf{\dot{q}}_0 = 0$.
\end{exercise}

\begin{remark}
The Lagrangian vector field associated to $K$ is called the {\bfi
geodesic spray\/}. Since the Legendre transformation is a
diffeomorphism (in finite dimensions or in infinite dimensions if  the
metric is assumed to be strong), the geodesic  spray is always a second
order equation.
\end{remark}

\subsection{Covariant derivative}
The variational approach to geodesics recovers the classical
formulation using covariant derivatives, as follows. 
Let $\mathfrak{X}(Q)$ denote the set of vector fields on the manifold $Q$.
The {\bfi covariant derivative}
\[
\nabla :\mathfrak{X}(Q)\times \mathfrak{X}(Q)\rightarrow
\mathfrak{X}(Q)\qquad
(X , Y) \mapsto \nabla_X(Y)
\,,
\]
 of the Levi-Civita connection on $(Q,g)$ is
given in local charts by
\[
\nabla_X(Y)=\Gamma_{ij}^kX^iY^j\frac{\partial}{\partial q^k}+
X^i\frac{\partial Y^k}{\partial q^i}\frac{\partial }{\partial q^k}
\,.
\]
If $c(t)$ is a curve on $Q$ and $Y\in \mathfrak{X}(Q)$, the covariant
derivative of $Y$ along $c(t)$ is defined by
$$\frac{DY}{Dt} := \nabla_{\dot c} Y,$$
or locally,
\[
\left(\frac{DY}{Dt}\right)^{\!k} =\Gamma_{ij}^k (c(t)) \dot c^i(t)
Y^j(c(t)) +
\frac{d}{dt} Y^k(c(t)).
\]
A vector field is said to be {\bfi parallel transported\/} along
$c(t)$ if  \[\frac{DY}{Dt} =0\,.\] Thus $\dot c(t)$ is parallel
transported along $c(t)$ if and only if
$$\ddot{c}\,^i + \Gamma^i_{jk} \dot c^j\, \dot c^k = 0.
$$
In classical differential geometry a {\bfi geodesic\/} is defined to be a 
curve $c(t)$ in $Q$ whose tangent vector $\dot{c}(t)$ is parallel
transported along $c(t)$. As the expression above shows, geodesics are
integral curves of the Lagrangian vector field defined by the kinetic
energy of $g$.

\begin{definition}
A {\bfi classical mechanical system\/} is given by a Lagrangian of
the form $L(v_q) = K(v_q) - V(q) $, for $v_q \in T_qQ $. The smooth
function $V: Q \rightarrow \mathbb{R}$ is called the {\bfi potential
energy\/}. The total energy of this system is given by $E = K + V $
and the Euler-Lagrange equations (which are always second order for a
hyperregular Lagrangian) are
\[
\ddot{q}\,^i+\Gamma_{jk}^i\dot q^j\dot q^k +  g ^{il}
\frac{\partial V}{ \partial q ^l}=0,\quad  i  = 1 ,\ldots n,
\]
where $g^{ij}$ are the entries of the inverse matrix of $(g_{ij})$.
\end{definition}

\begin{definition} 
If $Q = \mathbb{R}^3$ and the metric is given by $g_{ij} = \delta_{ij} $,
these equations are Newton's equations of motion \eqref{Newton2} of a
particle in a potential field  which launched our discussion in
Lecture 10.
\end{definition}

\begin{exercise}\label{gauge-ex}
[Gauge invariance] Show that the Euler-Lagrange equations are unchanged
under 
\begin{equation}
L(\mathbf{q}(t),\mathbf{\dot{q}}(t))
\rightarrow
L^\prime = L + \frac{d}{dt}\gamma(\mathbf{q}(t),\mathbf{\dot{q}}(t))
\,,
\end{equation}
for any function $\gamma: \mathbb{R}^{6N} =
\{(\mathbf{q},\mathbf{\dot{q}}) \mid \mathbf{q}, \mathbf{\dot{q}} \in
\mathbb{R}^{3N} \} \rightarrow \mathbb{R}$.
\end{exercise}

\begin{exercise}\label{covar-ex}
[Generalized coordinate theorem] Show that the Euler-Lagrange equations are
{\bfi unchanged in form} under any smooth invertible mapping
$f:\{\mathbf{q}\mapsto\mathbf{s}\}$. That is, with
\begin{equation}
L(\mathbf{q}(t),\mathbf{\dot{q}}(t))
=
\tilde{L}(\mathbf{s}(t),\mathbf{\dot{s}}(t))
\,,
\end{equation}
show that
\begin{equation}
\frac{d}{dt}
\left(\frac{\partial L}{\partial
\mathbf{\dot{q}}}\right) 
-\frac{\partial L}{\partial \mathbf{q}}
=0
\quad\Longleftrightarrow\quad
\frac{d}{dt}
\left(\frac{\partial \tilde{L}}{\partial
\mathbf{\dot{s}}}\right) 
-\frac{\partial \tilde{L}}{\partial \mathbf{s}}
=0
\,.
\end{equation}
\end{exercise}

\begin{exercise}
How do the Euler-Lagrange equations transform under
$\mathbf{q}(t)=\mathbf{r}(t)+\mathbf{s}(t)$?
\end{exercise}

\begin{exercise}[Other example Lagrangians]
Write the Euler-Lagrange equations, then apply the Legendre
transformation to determine the Hamiltonian and Hamilton's canonical
equations for the following Lagrangians. Determine which of them are
hyperregular.
\begin{itemize}
\item
$L(q, \dot{q})= \Big( g_{ij}(q) \dot q^i \dot q^j\Big)^{1/2}$ (Is it
possible to assume that $L(q,\dot{q})=1$? Why?)
\item
$L(q, \dot{q})
= -\Big(1-\mathbf{\dot{q}}\cdot\mathbf{\dot{q}}\Big)^{1/2}$
\item
$L(q, \dot{q})
= \frac{m}{2}\mathbf{\dot{q}}\cdot\mathbf{\dot{q}}
+
\frac{e}{c} \mathbf{\dot{q}}\cdot\mathbf{A}(\mathbf{q})$, for constants
$m$, $c$ and prescribed function $\mathbf{A}(\mathbf{q})$. How do the
Euler-Lagrange equations for this Lagrangian differ from free motion in a
moving frame with velocity $\frac{e}{mc}\mathbf{A}(\mathbf{q})$?
\end{itemize}
\end{exercise}

\paragraph{Example: Charged particle in a magnetic field.} Consider a
particle of charge
$e$ and mass
$m$ moving in a magnetic field $\mathbf{B}$, where $\mathbf{B} =
\nabla \times\mathbf{A}$ is a given magnetic field on $\mathbb{R}^3$.
The Lagrangian for the motion is given by the ``minimal coupling''
prescription (jay-dot-ay)
\[
L(\mathbf{q}, \dot{\mathbf{q}})
=
\frac{m}{2}\|\dot{\mathbf{q}}\|^2 
+ \frac{e}{c} \mathbf{A}(\mathbf{q})\cdot\mathbf{\dot{q}}
\,,
\]
in which the constant $c$ is the speed of light. The derivatives of this
Lagrangian are
\[
\frac{\partial L}{\partial \mathbf{\dot{q}}}
=
m\mathbf{\dot{q}}+\frac{e}{c}\mathbf{A}
=:
\mathbf{p}
\quad\hbox{and}\quad
\frac{\partial L}{\partial \mathbf{q}}
=
\frac{e}{c}\nabla\mathbf{A}^T\cdot\mathbf{\dot{q}}
\]
Hence, the Euler-Lagrange equations
for this system are 
\begin{eqnarray*}
m\,\mathbf{\ddot{q}} =
 \frac{e}{c}
(\nabla\mathbf{A}^T\cdot\mathbf{\dot{q}}
-\nabla\mathbf{A}\cdot\mathbf{\dot{q}})
=
\frac{e}{c}\,\mathbf{\dot{q}}\times\mathbf{B}
\\
\hbox{(Newton's equations for the Lorentz force)}
\end{eqnarray*}
The Lagrangian $L$ is hyperregular, because
\[
\mathbf{p}= \mathbb{F}L (\mathbf{q}, \mathbf{\dot{q}}) =
m\mathbf{\dot{q}}+\frac{e}{c}\mathbf{A}(\mathbf{q})
\]
has the  inverse
\[
\mathbf{\dot{q}} = \mathbb{F}H (\mathbf{q}, \mathbf{p}) =
\frac{1}{m}\left(\mathbf{p}- \frac{e}{c} \mathbf{A}(\mathbf{q})\right)
\,.
\]
The corresponding Hamiltonian is given by the invertible change of
variables,
\begin{equation}\label{Mag-Ham}
H(\mathbf{q},\mathbf{p}) 
= \mathbf{p} \cdot \dot{\mathbf{q}} - L(\mathbf{q}, \dot{\mathbf{q}})
= \frac{1}{2m}\left\| \mathbf{p} -
\frac{e}{c} \mathbf{A} \right\|^2 
\,.
\end{equation}
The Hamiltonian $H$ is hyperregular since
\[
\dot{\mathbf{q}} = \mathbb{F}H (\mathbf{q}, \mathbf{p}) =
\frac{1}{m}\left(\mathbf{p}- \frac{e}{c} \mathbf{A}\right)
\quad\hbox{has the inverse}\quad
\mathbf{p}= \mathbb{F}L (\mathbf{q}, \mathbf{\dot{q}}) = m
\dot{\mathbf{q}} + \frac{e}{c}\mathbf{A}\,.
\]
The canonical equations for this Hamiltonian recover Newton's equations
for the Lorentz force law. 

\paragraph{Example: Charged particle in a magnetic field by the
Kaluza-Klein construction.} Although the minimal-coupling Lagrangian is
not expressed as the kinetic energy of a metric, Newton's equations for
the Lorentz force law may still be obtained as geodesic equations. This is
accomplished by suspending them in a higher dimensional space via the
{\bfi Kaluza-Klein construction\/}, which proceeds as follows.\\

Let $Q_{KK}$ be the manifold $\mathbb{R}^3\times S ^1$ with variables
$(\mathbf{q},\theta)$. On $Q_{KK} $ introduce the one-form $A +
\mathbf{d}\theta$ (which defines a connection one-form on the
trivial circle bundle $\mathbb{R}^3\times S^1 \rightarrow \mathbb{R}^3$)
and introduce the {\bfi Kaluza-Klein Lagrangian}
$L_{KK}:TQ_{KK}\simeq T\mathbb{R}^3\times TS ^1\mapsto \mathbb{R}$ as
\begin{align}
L_{KK}(\mathbf{q}, \theta, \dot{\mathbf{q}},  \dot \theta) &=
\frac{1}{2} m \|\dot{\mathbf{q}}\|^2 + \frac{1}{2}\left\|\left\langle A
+ \mathbf{d}\theta, (\mathbf{q}, \dot{\mathbf{q}}, \theta, \dot \theta)
\right\rangle \right\|^2 \nonumber \\
&= \frac{1}{2} m \|\dot{\mathbf{q}}\|^2 + \frac{1}{2}\left(\mathbf{A}
\cdot \dot{\mathbf{q}} + \dot \theta
\right) ^2.
\label{KK-metric}
\end{align}
The Lagrangian $L_{KK} $ is positive definite in
$(\dot{\mathbf{q}},\dot \theta) $; so it may be regarded as the kinetic
energy of a metric, the {\bfi Kaluza-Klein metric\/} on $TQ_{KK} $.  (This construction fits the idea of $U(1)$ gauge symmetry for
electromagnetic fields in $\mathbb{R}^3$.  It can be generalized to a
principal bundle with compact structure group endowed with a connection.
The Kaluza-Klein Lagrangian in this generalization leads to Wong's
equations for a color-charged particle moving in a classical Yang-Mills
field.) The Legendre transformation for
$L_{KK}$ gives the momenta
\begin{equation}
\label{KK Legendre}
\mathbf{p} = m \dot{\mathbf{q}} + (\mathbf{A}\cdot \dot{\mathbf{q}} +
\dot \theta) \mathbf{A} \qquad \text{and} \qquad \pi =
\mathbf{A}\cdot \dot{\mathbf{q}} + \dot \theta.
\end{equation}
Since $L_{KK}$ does not depend on $\theta$, the Euler-Lagrange equation
\[
\frac{d}{dt}\frac{\partial L_{KK}}{\partial \dot\theta}  =
\frac{\partial L_{KK}}{\partial \theta} = 0
\,,
\]
shows that $\pi = \partial L_{KK}/\partial\dot \theta $ is conserved.
The {\bfi charge\/} is now defined by $e: = c\pi$. The Hamiltonian
$H_{KK}$ associated to $L_{KK}$ by the Legendre transformation
\eqref{KK Legendre} is
\begin{align}
H_{KK}(\mathbf{q}, \theta, \mathbf{p}, \pi) &= \mathbf{p}\cdot
\dot{\mathbf{q}} + \pi \dot \theta - L_{KK}(\mathbf{q}, \dot
{\mathbf{q}}, \theta, \dot \theta) \nonumber \\
&= \mathbf{p}\cdot \frac{1}{m}\left(\mathbf{p} - \pi\mathbf{A} \right) +
\pi(\pi- \mathbf{A}\cdot \dot{\mathbf{q}}) \nonumber \\
&\qquad - \frac{1}{2} m \|\dot{\mathbf{q}}\|^2 - \frac{1}{2}\pi^2
\nonumber\\
& = \mathbf{p}\cdot \frac{1}{m}\left(\mathbf{p} - \pi\mathbf{A} \right)
+ \frac{1}{2}\pi^2 \nonumber\\
&\qquad - \pi \mathbf{A} \cdot  \frac{1}{m}\left(\mathbf{p} -
\pi\mathbf{A} \right) - \frac{1}{2m}
\| \mathbf{p} - \pi \mathbf{A}\|^2 \nonumber\\
&= \frac{1}{2m} \| \mathbf{p} - \pi \mathbf{A}\|^2 + \frac{1}{2}\pi^2.
\label{KK Hamiltonian}
\end{align}
On the constant level set $\pi = e/c$, the Kaluza-Klein Hamiltonian
$H_{KK}$ is a function of only the variables $(\mathbf{q}, \mathbf{p})$
and is equal to the Hamiltonian \eqref{Mag-Ham} for charged particle
motion under the Lorentz force up to an additive constant. This
example provides an easy but fundamental illustration of the geometry of
(Lagrangian) reduction by symmetry. The canonical equations for the
Kaluza-Klein Hamiltonian $H_{KK}$ now reproduce Newton's equations for the
Lorentz force law.

\section{The rigid body in three dimensions} 

In the absence of external torques, Euler's equations for rigid body
motion are:
\begin{equation}\label{rbe1}
\begin{aligned}
&I_1\dot{\Omega}{_1}  =  (I_2 -I_3)\Omega_2\Omega_3,
\\
&I_2\dot{\Omega}{_2}  =  (I_3 -I_1)\Omega_3\Omega_1,
\\
&I_3\dot{\Omega}{_3}  =  (I_1 -I_2)\Omega_1\Omega_2,
\end{aligned}
\end{equation}
or, equivalently, 
\[
{\mathbb I} \boldsymbol{\dot{\Omega}}
={\mathbb I} \boldsymbol{\Omega} \times \boldsymbol{\Omega}\,,
\]
where $\boldsymbol{\Omega} = (\Omega_1, \Omega_2, \Omega_3)$ is the
body angular velocity vector and  $I_1,  I_2, I_3$  are the moments
of inertia of the rigid body.

\begin{question}
Can these equations -- as they are written -- be cast into Lagrangian or
Hamiltonian form in any sense?  (Since there are an odd number of
equations, they cannot be put into canonical Hamiltonian form.)
\end{question}

\noindent
We could reformulate them as:
\begin{description}
\item
Euler--Lagrange equations on $T{\rm SO(3)}$
or 
\item
Canonical Hamiltonian equations on  $T ^\ast{\rm SO(3)}$, 
\end{description}
by using Euler angles and their velocities, or their conjugate momenta.
However, these reformulations on $T{\rm SO(3)}$ or $T^\ast{\rm SO(3)}$
would answer a different question for a {\it six} dimensional system. We
are interested in these structures for the equations as given above.

\begin{answer}[Lagrangian formulation]
The Lagrangian answer is this: These equations
may be expressed in Euler--Poincar\'e form on the Lie algebra $
\mathbb{R}^3$ using the Lagrangian
\begin{equation}\label{rbl-1}
l (\boldsymbol{\Omega}) = \frac{1}{2} (I _1 \Omega ^2 _1 + I _2
\Omega ^2 _2 + I _3 \Omega ^2 _3 )
=
\frac{1}{2}\boldsymbol{\Omega}^T \cdot \mathbb{I}\boldsymbol{\Omega}
\,,
\end{equation}
which is the (rotational) kinetic energy of the
rigid body.\\

\noindent
The Hamiltonian answer to this question will be discussed later.
\end{answer}

\begin{proposition}
The Euler rigid body equations are equivalent to the {\bfi rigid body
action principle} for a {\bfi reduced action}
\begin{equation}\label{rbvp1}
\delta S_{\rm red}=\delta \int^b_a l (\boldsymbol{\Omega})\, d t = 0,
\end{equation}
where variations of $\boldsymbol{\Omega} $ are restricted to
be of the form
\begin{equation}\label{rpvp2}
\delta \boldsymbol{\Omega} = \boldsymbol{\dot{ \Sigma}} +
\boldsymbol{\Omega} \times \boldsymbol{\Sigma},
\end{equation}
in which $\boldsymbol{\Sigma}(t)$ is a curve in $\mathbb{R}^3$ that
vanishes at the endpoints in time. 
\end{proposition}

\begin{proof}
Since
$ l (\boldsymbol{\Omega}) = \frac{1}{2} \langle  {\mathbb I}
\boldsymbol{\Omega}, \boldsymbol{\Omega} \rangle $, and ${\mathbb I}$
is symmetric, we obtain
\begin{align*}
\delta \int^b_a l(\boldsymbol{\Omega}) \, d t
& =  \int^b_a \langle  {\mathbb I} \boldsymbol{\Omega},\delta
\boldsymbol{\Omega}\rangle \, dt \\
& =  \int^b_a \langle  {\mathbb I} \boldsymbol{\Omega},
\boldsymbol{\dot{ \Sigma}}
      + \boldsymbol{\Omega} \times \boldsymbol{\Sigma}\rangle \, dt\\
& =  \int^b_a \left[ \left\langle
      -\, \frac{d }{d t}{\mathbb I} \boldsymbol{\Omega} ,
      \boldsymbol{\Sigma}\right\rangle + \left\langle {\mathbb I}
\boldsymbol{\Omega},
      \boldsymbol{\Omega} \times
\boldsymbol{\Sigma}\right\rangle\right] {d t}\\
& = \int^b_a \left\langle
-\, \frac{d }{d t} {\mathbb I} \boldsymbol{\Omega}
      + {\mathbb I} \boldsymbol{\Omega}\times \boldsymbol{\Omega},
\boldsymbol{\Sigma} \right\rangle d t,
\end{align*}
upon integrating by parts and using the
endpoint conditions, $ \boldsymbol{\Sigma} (b) =
\boldsymbol{\Sigma} (a) = 0 $. Since $\boldsymbol{\Sigma}$ is
otherwise arbitrary,~(\ref{rbvp1}) is equivalent to
\[ -\,\frac{d }{d t} ({\mathbb I} \boldsymbol{\Omega})
+ {\mathbb I} \boldsymbol{\Omega} \times \boldsymbol{\Omega}
= 0 , \]
which are Euler's equations (\ref{rbe1}).
\end{proof}\bigskip

Let's derive this variational principle from the
{\it standard} Hamilton's principle.

\subsection{Hamilton's principle for rigid body motion on $T{\rm SO}(3)$}

An element $\mathbf{ R} \in {\rm SO}(3)$
gives the configuration of the body as a map of a {\bfi reference
configuration} ${\mathcal B}\subset \mathbb{R}^3$ to the current
configuration $\mathbf{ R}({\mathcal B})$; the map
$\mathbf{ R}$ takes a reference or label point $ X \in
{\mathcal B} $ to a current point $ x = \mathbf{ R}(X) \in
\mathbf{ R}({\mathcal B} )$. 
\\

When the rigid body is in motion,
the matrix $\mathbf{ R}$ is time-dependent. 
Thus, 
\[
x(t) = \mathbf{R}(t)X
\]
with $\mathbf{R}(t)$ a curve parameterized by time in ${\rm SO}(3)$.
The velocity of a point of the body is 
\[\dot{ x}(t) = \mathbf{\dot{R}}(t) X 
= \mathbf{\dot{R}} \mathbf{ R} ^{-1}(t) x(t)\,.\]
Since $\mathbf{ R}$ is an
orthogonal matrix, $\mathbf{ R} ^{-1} \mathbf{\dot{R}}$ and
$\mathbf{\dot{R}}\mathbf{ R}^{-1}$ are skew matrices. Consequently,
we can write (recall the hat map)
\begin{equation}\label{sav}
\dot{x} = \mathbf{\dot{R}}\mathbf{ R} ^{-1} x = \boldsymbol{\omega}
\times x\,.
\end{equation}
This formula defines the {\bfi spatial angular velocity
vector\/} $\boldsymbol{\omega}$.  Thus, $\boldsymbol{\omega}$ is
essentially given by {\em right\/} translation of $\mathbf{\dot{R}}$ to
the identity. That is, the vector 
\[
\boldsymbol{\omega}
=\big(\mathbf{\dot{R}} \mathbf{ R} ^{-1}\big)\boldsymbol{\hat{\,}}
\,.
\]

The corresponding {\bfi body angular velocity} is defined by
\begin{equation}\label{bav}
\boldsymbol{\Omega} = \mathbf{ R} ^{ -1} \boldsymbol{\omega} ,
\end{equation}
so that $\boldsymbol{\Omega}$ is the angular velocity
relative to a body fixed frame. Notice that
\begin{align}\label{lefttranslation}
\mathbf{ R} ^{ -1} \mathbf{\dot{R}} X
&= \mathbf{ R} ^{ -1} \mathbf{\dot{R}}
        \mathbf{ R} ^{ -1} x = \mathbf{ R} ^{ -1}
        (\boldsymbol{\omega} \times x) \nonumber \\
&= \mathbf{ R}^{ -1} \boldsymbol{\omega} \times \mathbf{ R}^{ -1} x 
= \boldsymbol{\Omega} \times X,
\end{align}
so that $\boldsymbol{\Omega}$ is given by {\em left\/} translation
of $\mathbf{\dot{R}}$ to the identity. That is, the vector 
\[
\boldsymbol{\Omega}
=\big(\mathbf{ R} ^{-1}\mathbf{\dot{R}}\big)\boldsymbol{\hat{\,}}
\,.
\]

The {\bfi kinetic energy} is obtained by summing up $m | \dot{x}
| ^2 / 2$ (where $|{\,\cdot\,}|$ denotes the Euclidean norm) over the
body. This yields
\begin{equation}\label{ke}
 K = \frac{1}{2} \int _{\mathcal B} \rho(X) |
\mathbf{\dot{R}} X | ^2 \, d^3 X,
\end{equation}
in which $\rho$ is a given mass density in the reference
configuration. Since
\[
| \mathbf{\dot{R}} X |  = | \boldsymbol{\omega} \times x | = |
\mathbf{ R} ^{-1}(\boldsymbol{\omega} \times x ) | = |
\boldsymbol{\Omega} \times X |,
\]
$K$ is a quadratic function of $\boldsymbol{\Omega}$. Writing
\begin{equation}\label{mit}
K = \frac{1}{2} \boldsymbol{\Omega} ^T \cdot{\mathbb I}\boldsymbol{\Omega}
\end{equation}
defines the {\bfi moment of inertia
tensor\/} ${\mathbb I},$ which, provided the
body does not degenerate to a line, is a positive-definite
$(3\times 3)$ matrix, or better, a quadratic form.  This quadratic
form can be diagonalized by a change of basis; thereby defining the
principal axes and moments of inertia. In this basis, we write
$ {\mathbb I} = {\rm diag}(I _1, I _2, I_3). $ \\

The function $K$ is
taken to be the Lagrangian of the system on
$ T{\rm SO}(3) $ (and by means of the Legendre transformation we obtain
the corresponding Hamiltonian description on
$ T ^{\ast} {\rm SO}(3)$). Notice that $K$ in equation (\ref{ke}) is
{\em left\/} (not right) invariant on $T{\rm SO}(3)$, since
\[
\boldsymbol{\Omega}
=\big(\mathbf{ R} ^{-1}\mathbf{\dot{R}}\big)\boldsymbol{\hat{\,}}
\,.
\]
It follows that the corresponding Hamiltonian will also be {\em left\/}
invariant.\\

In the framework of Hamilton's principle, the relation between
motion in $\mathbf{ R}$ space and motion in body angular velocity
(or $\boldsymbol{\Omega}$) space is as follows.

\begin{proposition}
The curve $ \mathbf{R} (t) \in {\rm SO}(3) $ satisfies the Euler-Lagrange
equations for
\begin{equation}\label{rbl-2}
 L (\mathbf{ R}, \mathbf{\dot{R}} )
= \frac{1}{2} \int_{{\mathcal B}}
\rho (X)
|  \mathbf{\dot{R}} X | ^2 \, d ^3 X,
\end{equation}
if and only if $ \boldsymbol{\Omega} (t) $ defined by
$\mathbf{ R} ^{-1} \mathbf{\dot{R}} \mathbf{v} = \boldsymbol{\Omega}
\times \mathbf{v}$ for all $\mathbf{v} \in \mathbb{R}^3$ satisfies
Euler's equations
\begin{equation}\label{ee}
{\mathbb I} \boldsymbol{\dot{ \Omega}} = {\mathbb
I}\boldsymbol{\Omega}\times \boldsymbol{\Omega}
\,.
\end{equation}
\end{proposition}
The proof of this relation will illustrate how
to reduce variational principles using their symmetry groups.
By Hamilton's principle, $\mathbf{R}(t)$
satisfies the Euler-Lagrange equations, if and only if
\[
\delta\int L (\mathbf{ R}, \mathbf{\dot{R}} ) \, dt = 0 .
\]
 Let $ l (\boldsymbol{\Omega}) = \frac{1}{2} ( {\mathbb I}
\boldsymbol{\Omega})
\cdot \boldsymbol{\Omega} $, so that
$ l (\boldsymbol{\Omega}) = L (\mathbf{ R}, \mathbf{\dot{R}} ) $
where the matrix $\mathbf{ R}$ and the vector $\boldsymbol{\Omega}$ are
related by the hat map, $\boldsymbol{\Omega} =\big(\mathbf{ R}
^{-1}\mathbf{\dot{R}}\big)\boldsymbol{\hat{\,}}$. Thus, the Lagrangian
$L$ is left SO(3)-invariant. That is, 
\[
l (\boldsymbol{\Omega}) = L (\mathbf{ R}, \mathbf{\dot{R}} )
= L (\mathbf{ e}, \mathbf{ R}^{-1}\mathbf{\dot{R}} )
\,.
\]
To see how we should use this left-invariance to transform
Hamilton's principle, define the skew matrix
$\boldsymbol{\hat{\Omega}}$  by
$\boldsymbol{\hat{\Omega}}\mathbf{v} = \boldsymbol{\Omega} \times
\mathbf{v}$ for any $\mathbf{v} \in \mathbb{R}^3$.\\

We differentiate the relation
$\mathbf{ R} ^{-1} \mathbf{\dot{R}} = \hat{\boldsymbol{\Omega} }$
with respect to $\mathbf{ R}$ to get
\begin{equation}\label{rda}
- \mathbf{ R} ^{-1} (\delta \mathbf{ R})\mathbf{ R}^{-1}
\mathbf{\dot{R}}  + \mathbf{ R}^{-1} (\delta
\mathbf{\dot{R}}) = \widehat{\delta \boldsymbol{\Omega}} .
\end{equation}
Let the skew matrix $\boldsymbol{\hat{\Sigma}}$ be defined by
\begin{equation}\label{dsig}
 \boldsymbol{\hat{\Sigma}} = \mathbf{R}^{-1} \delta \mathbf{R},
\end{equation}
and define the vector $\boldsymbol{\Sigma}$ by
\begin{equation}\label{dhat}
\boldsymbol{\hat{\Sigma}} \mathbf{v} = \boldsymbol{\Sigma} \times
\mathbf{v} .
\end{equation}
Note that
\[
\boldsymbol{\dot{ \hat{\Sigma}}}
= -\, \mathbf{ R}^{-1} \mathbf{\dot{R}} \mathbf{ R} ^{-1}
\delta \mathbf{ R} + \mathbf{ R} ^{-1}
\delta\mathbf{\dot{R}},
\]
so
\begin{equation}\label{dela}
\mathbf{ R} ^{-1} \delta \mathbf{\dot{R}}
= \boldsymbol{\dot{ \hat{\Sigma}}} + \mathbf{ R} ^{-1}
\mathbf{\dot{R}} \boldsymbol{\hat{\Sigma}}\,.
\end{equation}
Substituting~(\ref{dela}) and~(\ref{dsig})
into~(\ref{rda}) gives
\[ -\,\boldsymbol{\hat{\Sigma}} \boldsymbol{\hat{\Omega}} 
+ \boldsymbol{\dot{ \hat{\Sigma}}} 
+ \boldsymbol{\hat{\Omega}} \boldsymbol{\hat{\Sigma}}
= \widehat{ \delta \boldsymbol{\Omega}}, \]
that is,
\begin{equation}\label{domh}
\widehat{ \delta \boldsymbol{\Omega}} = 
\boldsymbol{\dot{ \hat{\Sigma}}} 
+ [
\boldsymbol{\hat{\Omega}},
\boldsymbol{\hat{\Sigma}} ] .
\end{equation}
The identity
$ [ \boldsymbol{\hat{\Omega}}, \boldsymbol{\hat{\Sigma}} ] =
(\boldsymbol{\Omega} \times
\boldsymbol{\Sigma}) \boldsymbol{\hat{\ }}$ holds by Jacobi's identity for
the cross product and so
\begin{equation}\label{dom}
\delta \boldsymbol{\Omega} = \boldsymbol{\dot{\Sigma}} +
\boldsymbol{\Omega}
\times \boldsymbol{\Sigma} .
\end{equation}

\noindent
These calculations prove the following:
\begin{theorem}  \label{red/variational/prpl}
For a Lagrangian which is left-invariant under SO(3),
Hamilton's variational principle
\begin{equation}\label{hprin}
\delta S = \delta \int^b_a L (\mathbf{ R}, \mathbf{\dot{R}} )\, d t = 0
\end{equation}   on $T{\rm SO}(3)$ is equivalent to the {\bfi
reduced variational principle\/}
\begin{equation}\label{rvprin}
\delta S_{\rm red}= \delta \int^b_a l(\boldsymbol{\Omega}) \, d t = 0
\end{equation}  
with $\boldsymbol{\Omega} =\big(\mathbf{R}^{-1}
\mathbf{\dot{R}}\big)\boldsymbol{\hat{\,}}$ on $\mathbb{R}^3$ where the
variations $ \delta
\boldsymbol{\Omega} $ are of the form
\[
\delta \boldsymbol{\Omega} = \boldsymbol{\dot{\Sigma}} +
\boldsymbol{\Omega}
\times \boldsymbol{\Sigma}\,,
\]
with
$\boldsymbol{\Sigma}(a) =\boldsymbol{\Sigma}(b) = 0 $.
\end{theorem}

\noindent
Recall {\bf Theorem \ref{red/variational/prpl}}:\\
{\it
For a Lagrangian which is left-invariant under SO(3),
Hamilton's variational principle
\begin{equation}\label{hprin-redux}
\delta S = \delta \int^b_a L (\mathbf{ R}, \mathbf{\dot{R}} )\, d t = 0
\end{equation}   on $T{\rm SO}(3)$ is equivalent to the {\bfi
reduced variational principle\/}
\begin{equation}\label{rvprin-redux}
\delta S_{\rm red}= \delta \int^b_a l(\boldsymbol{\Omega}) \, d t = 0
\end{equation}  
with $\boldsymbol{\Omega} =\big(\mathbf{R}^{-1}
\mathbf{\dot{R}}\big)\boldsymbol{\hat{\,}}$ on $\mathbb{R}^3$ where the
variations $ \delta
\boldsymbol{\Omega} $ are of the form
\[
\delta \boldsymbol{\Omega} = \boldsymbol{\dot{\Sigma}} +
\boldsymbol{\Omega}
\times \boldsymbol{\Sigma}\,,
\]
with
$\boldsymbol{\Sigma}(a) =\boldsymbol{\Sigma}(b) = 0 $.
}

\paragraph{Reconstruction of $\mathbf{ R}(t)\in SO(3)$.} In 
Theorem \ref{red/variational/prpl}, Euler's equations for the rigid body
\[{\mathbb I} \boldsymbol{\dot{ \Omega}} =
{\mathbb I}\boldsymbol{\Omega}\times \boldsymbol{\Omega}
\,,
\]
follow from the reduced variational principle (\ref{rvprin}) for the
Lagrangian 
\begin{equation}\label{rb-redlag}
l(\boldsymbol{\Omega}) =  
 \frac{1}{2} ( {\mathbb I}
\boldsymbol{\Omega})
\cdot \boldsymbol{\Omega} 
\,,
\end{equation}
which is expressed in terms of the left-invariant time-dependent angular
velocity in the body, $\boldsymbol{\Omega}\in \mathfrak{so}(3)$. The
body angular velocity $\boldsymbol{\Omega}(t)$ yields the tangent vector
$\mathbf{\dot{R}}(t)\in T_{\mathbf{ R}(t)}SO(3)$ along the integral curve
in the rotation group $\mathbf{ R}(t)\in SO(3)$ by the relation,
\[
\mathbf{\dot{R}}(t)=\mathbf{ R}(t)\boldsymbol{\Omega}(t)
\,.
\]
This relation provides the {\bfi reconstruction formula}. It's solution as
a linear differential equation with time-dependent coefficients yields the
integral curve $\mathbf{ R}(t)\in SO(3)$ for the orientation of the
rigid body, once the time dependence of $\boldsymbol{\Omega}(t)$ is
determined from the Euler equations.

\subsection{Hamiltonian Form of rigid body motion.}
A dynamical system on a manifold $M$ 
\[
\mathbf{\dot{x}}(t)=\mathbf{F}(\mathbf{x})
\,,\quad \mathbf{x}\in M
\]
is said to be in {\bfi Hamiltonian form}, if it can be expressed as
\[
\mathbf{\dot{x}}(t)=\{\mathbf{x},H\}
\,,\quad\hbox{for}\quad
H:M\mapsto \mathbb{R}
\,,
\]
in terms of a Poisson bracket operation, 
\[
\{\cdot\,,\,\cdot\}:
\mathcal{F}(M)\times\mathcal{F}(M)\mapsto\mathcal{F}(M)
\,,
\]
which is bilinear, skew-symmetric and satisfies the Jacobi identity and
(usually) the Leibnitz rule. 

As we shall explain, reduced equations arising from group-invariant
Hamilton's principles on Lie groups are naturally Hamiltonian. If we
{\bfi Legendre transform} our reduced Lagrangian for the $SO(3)$
left invariant variational principle (\ref{rvprin-redux}) for
rigid body dynamics, then its simple, beautiful and well-known Hamiltonian
formulation emerges.

\begin{definition}
The Legendre transformation $\mathbb{F}l
:\mathfrak{so}(3)\rightarrow\mathfrak{so}(3)^*$ is defined  by
$$
\mathbb{F} l (\Omega) = \frac{\delta l}{\delta \Omega} = \Pi \,.
$$
\end{definition}
The Legendre transformation defines the {\bfi body angular momentum} by the
variations of the rigid-body's reduced Lagrangian with respect to the body
angular velocity. For the Lagrangian in (\ref{rb-redlag}), the
$\mathbb{R}^3$ components of the body angular momentum are
\begin{equation}\label{rbm}
\Pi_i  =  I_i\Omega_i = \frac{\partial l}{\partial \Omega _i } \,,
\quad i = 1, 2, 3.
\end{equation}

\subsection{Lie-Poisson Hamiltonian formulation of rigid body dynamics.}
Let
\[h(\Pi) := \langle \Pi , \Omega\rangle - l(\Omega)\,,\]
where the pairing $\langle\cdot\,,\,\cdot\rangle:
\mathfrak{so}(3)^*\times\mathfrak{so}(3)\to\mathbb{R}$ is understood in
components as the vector dot product on $\mathbb{R}^3$
\[
\langle \Pi , \Omega\rangle := \boldsymbol{\Pi}\cdot\boldsymbol{\Omega}
\,.
\]
Hence, one finds the expected expression for the rigid-body Hamiltonian
\begin{equation}\label{rbh} 
h
=
\frac{1}{2}\,
\boldsymbol{\Pi}\cdot {\mathbb I}^{-1}\boldsymbol{\Pi}
:=
\frac{\Pi^2_1}{2I_1}
+ \frac{\Pi^2_2}{2I_2}
+ \frac{\Pi^2_3}{2I_3} 
\,.
\end{equation}
The Legendre transform $\mathbb{F}l$ for this case
is a diffeomorphism, so we may solve for 
$$\frac{\partial h}{\partial \Pi} = \Omega + \left\langle \Pi\, ,\,
\frac{\partial
\Omega}{\partial
\Pi} \right\rangle -  \left\langle \frac{\partial l}{\partial \Omega}\, ,\,
\frac{\partial \Omega}{\partial
\Pi} \right\rangle  = \Omega .
$$
In $\mathbb{R}^3$ coordinates, this relation expresses the body angular
velocity as the derivative of the reduced Hamiltonian with respect to
the body angular momentum, namely (introducing grad-notation),
\[
\nabla_{\Pi}h
:=
\frac{\partial h}{\partial\boldsymbol{\Pi}}
=
\boldsymbol{\Omega}
\,.
\]
Hence, the reduced Euler-Lagrange equations for $l$ may be expressed
equivalently in angular momentum vector components in $\mathbb{R}^3$ and
Hamiltonian
$h$ as:
\[ \frac{d }{d t} ({\mathbb I} \boldsymbol{\Omega})
= {\mathbb I} \boldsymbol{\Omega} \times \boldsymbol{\Omega}
\Longleftrightarrow 
\boldsymbol{\dot{\Pi}}
= 
\boldsymbol{\Pi}\times\nabla_{\Pi}h
:=
\{\boldsymbol{\Pi},h\}
\,. \]
This expression suggests we introduce the following rigid body Poisson
bracket on functions of the ${\boldsymbol \Pi}$'s:
\begin{equation}\label{rbb}
 \{f,h\}({\boldsymbol  \Pi})
:= -\,{\boldsymbol  \Pi} \cdot (\nabla_{\Pi} f \times \nabla_{\Pi} h)
\,.
 \end{equation} 
For the Hamiltonian (\ref{rbh}), one checks that 
the Euler equations in terms of the rigid-body angular
momenta,
\begin{equation}\label{rbe2}
\begin{aligned}
\dot {\Pi}_1  &=  \frac{I_2 - I_3}{I_2I_3} \Pi_2\Pi_3,\\
		\dot {\Pi}_2  &=  \frac{I_3 - I_1}{I_3I_1} \Pi_3\Pi_1,\\
		\dot {\Pi}_3  &=  \frac{I_1 - I_2}{I_1I_2} \Pi_1\Pi_2,
\end{aligned}
\end{equation}
that is,
\begin{equation}\label{rbe3}
\boldsymbol{\dot{\Pi }} = {\boldsymbol\Pi}
\times {\boldsymbol \Omega}\,.
\end{equation}
are equivalent to
\[
\dot{f} = \{f, h\}
\,,\quad\hbox{with}\quad
f=\boldsymbol{\Pi}
\,.
\]
The Poisson bracket proposed in (\ref{rbb}) is an example of a {\bfi Lie
Poisson bracket}, which we will show separately satisfies the defining
relations to be a Poisson bracket.

\subsection{$\mathbb{R}^3$ Poisson bracket.}
The rigid body Poisson bracket (\ref{rbb}) is a special case of the
Poisson bracket for functions on $\mathbb{R}^3$,
\begin{equation}\label{r3pb}
\{f, h\} = -\,\nabla{c}\cdot\nabla{f}\times\nabla{h}
\end{equation}
This bracket generates the motion 
\begin{equation}
\mathbf{\dot{x}} = \{\mathbf{x}, h\} =
\nabla{c}\times\nabla{h}
\end{equation}
For this bracket the motion takes place along the
intersections of level surfaces of the functions $c$ and $h$ in
$\mathbb{R}^3$. In particular, for the rigid body, the motion takes place
along intersections of angular momentum spheres $c=\|\mathbf{x}\|^2/2$ and
energy ellipsoids $h=\mathbf{x}\cdot \mathbb{I}\mathbf{x}$. (See the cover
illustration of Marsden and Ratiu [2003].)

\begin{exercise}
Consider the $\mathbb{R}^3$ Poisson bracket
\begin{equation}\label{r3pb-ex}
\{f, h\} = -\,\nabla{c}\cdot\nabla{f}\times\nabla{h}
\end{equation}
Let  $c=\mathbf{x}^T\cdot\mathbb{C}\mathbf{x}$ be a quadratic form on
$\mathbb{R}^3$, and let $\mathbb{C}$ be the associated symmetric
$3\times3$ matrix. Determine the conditions on the quadratic function
$c(\mathbf{x})$ so that this Poisson bracket will satisfy the Jacobi
identity.
\end{exercise}

\begin{exercise}
Find the general conditions on the function $\mathbf{c}(\mathbf{x})$ so
that the $\mathbb{R}^3$ bracket
\[
\{f, h\} = -\,\nabla{c}\cdot\nabla{f}\times\nabla{h}
\] 
satisfies the defining properties of a Poisson bracket. Is this
$\mathbb{R}^3$ bracket also a derivation satisfying the Leibnitz relation
for a product of functions on $\mathbb{R}^3$? If so, why?
\end{exercise}

\begin{exercise}
How is the $\mathbb{R}^3$ bracket related to the canonical Poisson
bracket? Hint: restrict to level surfaces of the function
$c(\mathbf{x})$. 
\end{exercise}

\begin{exercise}[Casimirs of the $\mathbb{R}^3$ bracket]
The Casimirs (or distinguished functions, as Lie called them) of a
Poisson bracket satisfy 
\[
\{c, h\}(\mathbf{x}) = 0
\,,\quad\forall h(\mathbf{x})
\]
Suppose the function $\mathbf{c}(\mathbf{x})$ is chosen so that the
$\mathbb{R}^3$ bracket (\ref{r3pb}) satisfies the defining properties of a
Poisson bracket. What are the Casimirs for the
$\mathbb{R}^3$ bracket (\ref{r3pb})? Why?
\end{exercise}

\begin{exercise}
Show that the motion equation
\[
\mathbf{\dot{x}} = \{\mathbf{x}, h\}
\]
for the $\mathbb{R}^3$ bracket (\ref{r3pb}) is invariant under a certain
linear combination of the functions $c$ and $h$. Interpret this invariance
geometrically. 
\end{exercise}

\section{Momentum maps}
\subsection*{The Main Idea}

Symmetries are often associated with conserved quantities.  For example, the flow
of any $SO(3)$-invariant Hamiltonian vector field on $T^*\mathbb{R}^3$
conserves angular momentum,
$\mathbf{q}\times \mathbf{p}.$ More generally, given a Hamiltonian $H$ on
a phase space $P,$ and a group action of $G$  on $P$ that conserves $H,$
there is often an associated ``momentum map'' $J:P\to \mathfrak{g}^*$ that
is conserved by the flow of the Hamiltonian vector field.

\medskip


Note: all group actions in this section will be left actions until
otherwise specified.

\subsection{Hamiltonian systems on Poisson manifolds}

\begin{definition} A \textbf{Poisson bracket} on a manifold $P$ is a 
skew-symmetric bilinear operation on $\mathcal{F}(P) := \mathcal{C}^\infty\left(P,\R\right)$ satisfying
the Jacobi identity and the Leibniz identity,
\[
\{ FG,H\} = F\{G,H\} + \{F,H\} G
\]
The pair $\left(P,\{\cdot, \cdot \}\right)$ is called a \textbf{Poisson manifold}.
\end{definition}

\begin{remark}
The Leibniz identity is sometimes not included in the definition. 
Note that bilinearity, skew-symmetry and the Jacobi identity are the axioms of a Lie algebra.
In what follows, a Poisson bracket is a binary operation that makes
$\mathcal{F}(P)$ into a Lie algebra and also satisfies the Leibniz
identity.
\end{remark}

\begin{exercise}
Show that the \textbf{classical Poisson bracket}, defined in
cotangent-lifted coordinates 
\[\left(q^1,\dots,q^N,p_1,\dots,p_N\right)\]
on an $2N$-dimensional cotangent bundle $T^*Q$ by
\[
\left\{F,G\right\} = \sum_{i=1}^{N} \left( \frac{\partial F}{\partial q^i} \frac{\partial G}{\partial p_i} 
- \frac{\partial F}{\partial p_i} 
\frac{\partial G}{\partial q^i}  \right)
,
\]
satisfies the axioms of a Poisson bracket.
Show also that the definition of this bracket is independent of the choice of local coordinates 
$\left(q^1,\dots,q^N\right).$
\end{exercise}

\begin{definition}
A \textbf{Poisson map} between two Poisson manifolds is a map
$\varphi: \left(P_1,\{\cdot, \cdot\}_1\right) \to \left(P_2,\{\cdot, \cdot\}_2\right)$ that preserves the brackets,
meaning
\[
\{F\circ \varphi, G\circ \varphi\}_1 = \{ F,G\}_2\circ \varphi , \quad \quad \textrm{for all } 
F,G\in \mathcal{F}\left(P_2\right).
\]
\end{definition}

\begin{definition}
An action $\Phi$ of $G$ on a Poisson manifold $\left(  P,\left\{  ,\right\}  \right)
$ is \textbf{canonical} if $\Phi_{g}$ is a Poisson map for every $g,$ i.e.
\[
\left\{  F\circ\Phi_{g},K\circ\Phi_{g}\right\}  =\left\{  F,K\right\}
\circ\Phi_{g}%
\]
for every $F,K\in \mathcal{F}(P).$
\end{definition}

\begin{definition}
Let $\left(P,\{\cdot, \cdot\}\right)$ be a Poisson manifold, and let
$H:P\to \R$ be differentiable. 
The \textbf{Hamiltonian vector field} for $H$ is the vector field $X_H$ defined by
\[
X_H(F) = \{F,H\}, \quad \quad \textrm{for any } F\in \mathcal{F}(P)
\]
\end{definition}

\begin{remark}
$X_H$ is well-defined because of the Leibniz identity and the correspondance between
vector fields and derivations (see \cite{Le2003}). 
\end{remark}

\begin{remark}
$X_H(F)= \pounds_{X_H}F = \dot F,$ the Lie derivative of $F$ along the flow of $X_H.$
The equations \[\dot F = \{F,H\}\,,\] called ``Hamilton's equations'', have
already appeared in Theorem 6.2, and are an equivalent definition of $X_H.$
\end{remark}

\begin{exercise} \label{canHam}
Show that Hamilton's equations for the classical Poisson bracket are the canonical 
Hamilton's equations,
\[
\dot q ^i = \frac{\partial H}{\partial p_i}\,, \qquad
\dot p_i = -\frac{\partial H}{\partial q^i}\,.
\]
\end{exercise}

\subsection{Infinitesimal invariance under Hamiltonian vector fields}

Let $G$ act smoothly on $P,$
and let $\xi \in \mathfrak{g}.$ Recall (from Lecture 9) that the infinitesimal
generator $\xi_{P}$ is the vector field on $P$ defined by 
\[
\xi_{P}\left(  x\right)  =\left. \frac{d}{dt}  {g(t)} x \, \right|  _{t=0},
\]
for some path $g(t)$ in $G$ such that $g(0)=e$ and $g'(0)=\xi.$

\begin{remark} \label{mxinfgen}
For matrix groups, we can take $g(t) = \exp \left(t \xi\right).$
This works in general for the exponential map of an arbitrary Lie group.
For matrix groups,
\[\xi_{P}\left(  \mathbf{x}\right)  
=\left. \frac{d}{dt}  \exp(t \xi)  \mathbf{x} \, \right|  _{t=0} 
= \xi \mathbf{x} 
\quad\hbox{(matrix multiplication).}
\] 
\end{remark}

%




\begin{exercise} \label{infinv}
If $H:P\to \R$ is $G$-invariant, meaning that $H(gx)=H(x)$ for all $g\in G$ and $x\in P,$
then $\pounds_{\xi_P}H=0$ for all $\xi \in\mathfrak{g}.$
This property is called \textbf{infinitesimal invariance}.
\end{exercise}

\begin{example} [The momentum map for the rotation group]
Consider the cotangent bundle of ordinary Euclidean space $\mathbb{R}^3$.
This is the Poisson (symplectic) manifold with coordinates 
$(\mathbf{q},\mathbf{p})\in T^*\mathbb{R}^3\simeq\mathbb{R}^6$, equipped with the canonical
Poisson bracket. An element $g$ of the rotation group $SO(3)$
acts on $T^*\mathbb{R}^3$ according to
\[
g(\mathbf{q},\mathbf{p})=(g\mathbf{q},g\mathbf{p})
\] 
Set $g(t)=\exp(tA)$, so that $\frac{d}{dt}\big|_{t=0}g(t)=A$ and the
corresponding Hamiltonian vector field is
\[
X_A
=(\mathbf{\dot{q}},\mathbf{\dot{p}})
=(A\mathbf{q},A\mathbf{p})
\]
where $A\in so(3)$ is a skew-symmetric matrix. The corresponding
Hamiltonian equations read 
\[
\mathbf{\dot{q}}
=
A\mathbf{q}
=
\frac{\partial J_A}{\partial \mathbf{p}}
\,,\quad
\mathbf{\dot{p}}
=
A\mathbf{p}
=
-\,\frac{\partial J_A}{\partial \mathbf{q}}
\,.
\]
Hence,
\[
J_A \left(  \mathbf{q},\mathbf{p}\right)  
=-A\mathbf{p}\cdot \mathbf{q} = a_i\epsilon_{ijk}p_kq_j
= \mathbf{a}\cdot \mathbf{q}\times \mathbf{p}
\,.\]
for a vector $\mathbf{a}\in\mathbb{R}^3$ with components $a_i$, $i=1,2,3$.
So the momentum map for the rotation group is the
angular momentum $J=\mathbf{q}\times\mathbf{p}$.
\end{example}

\begin{example} Consider angular momentum $J=\mathbf{q}\times \mathbf{p},$
defined on $P=T^*\R^3.$ 
%
For every $\xi\in \R^3,$ define
\[
J_{\xi}\left(  \mathbf{q},\mathbf{p}\right)  
:=\xi\cdot\left(  \mathbf{q}\times \mathbf{p}\right) 
=\mathbf{p}\cdot\left(  \xi\times \mathbf{q}\right) 
\]
Using Exercise \ref{canHam} and Example \ref{mxinfgen},
\begin{align*}
X_{J_\xi} \left(\mathbf{q}, \mathbf{p}\right)
&= \left( \frac{\partial J_\xi}{\partial \mathbf{p}}, 
-\,\frac{\partial J_\xi}{\partial \mathbf{q}}\right) \\
&= \left(\xi \times \mathbf{q}, \xi \times \mathbf{p}\right) \\
&=\hat{\xi}_P\left( \mathbf{q}, \mathbf{p}\right),
\end{align*}
where the last line is the infinitesimal generator corresponding to $\hat{\xi}\in so(3).$
Now suppose $H:P\to \R$ is $SO(3)$-invariant.
From Exercise \ref{infinv}, we have $\pounds_{\hat{\xi}} H = 0.$
It follows that 
\[
\pounds_{X_H} J_\xi
=\left\{  J_{\xi},H\right\}  
=-\left\{  H,J_{\xi}\right\}  
=-\pounds_{X_{J_\xi}} H
=-\pounds_{\xi_P} H
=0.
\]
Since this holds for all $\xi,$ we have shown that $J$ is conserved by the Hamiltonian flow.
\end{example}

\subsection{Defining Momentum Maps}

In order to generalise this example, we recast it using the 
hat map $\hat{}:\R^3\to so(3)$
and the associated map $\tilde{}:\left(\R^3\right)^*\to so(3)^*,$
and the standard identification $\left(\R^3\right)^*\cong \R^3$ via the Euclidean dot product.
We consider $J$ as a function from $P$ to $so(3)^*$ given by
$J\left(  \mathbf{q},\mathbf{p}\right)  = \left(\mathbf{q}\times \mathbf{p}\right)^{\tilde{}}.$
For any $\xi = \hat{\mathbf{v}},$ we define
$J_\xi \left(  \mathbf{q},\mathbf{p}\right)  
= \left<\left(\mathbf{q}\times \mathbf{p}\right)^{\tilde{}}, \hat{\mathbf{v}} \right>
= \left(\mathbf{q}\times \mathbf{p}\right) \cdot {\mathbf{v}}.$
As before, we find that 
$X_{J_\xi}=\xi_P$ for every $\xi,$ and $J$ is conserved by the Hamiltonian flow.
We take the first property, $X_{J_\xi}=\xi_P,$ as the general definition of a momentum map.
The conservation of $J$ follows by the same Poisson bracket calculation as in the example;
the result is Noether's Theorem.


\begin{definition}
A \textbf{momentum map} for a canonical action of $G$ on $P$ is a map $J:P\rightarrow
\mathfrak{g}^{\ast}$ such that,  for every $\xi \in \mathfrak{g},$ the map $J_\xi:P\to \R$
defined by $J_\xi(p) = \left<J(p), \xi\right> $ satisfies
\[
X_{J_{\xi}}=\xi_{P}
\]
\end{definition}

 
\begin{theorem}[Noether's Theorem]
Let $G$ act canonically on $\left(
P,\left\{ \cdot\,,\,\cdot\right\}  \right)  $ with momentum map $J.$ If
$H$ is
$G$-invariant, then $J$ is conserved by the flow of $X_{H}.$
\end{theorem}

\textbf{Proof} 
For every $\xi\in\mathfrak{g},$
\[
\pounds_{X_H} J_\xi
=\left\{  J_{\xi},H\right\}  
=-\left\{  H,J_{\xi}\right\}  
=-\pounds_{X_{J_\xi}} H
=-\pounds_{\xi_P} H
=0.
 \quad \blacksquare
\]

\begin{exercise}
Momentum maps are unique up to a choice of a constant element of
\textrm{g}$^{\ast}$ on every connected component of $M.$
\end{exercise}




\begin{exercise}
Show that the $S^1$ action on the torus $T^2:=S^1\times S^1$ given by 
$\alpha \left(\theta,\phi\right) = \left(\alpha + \theta, \phi\right)$ is canonical
with respect to  the classical bracket (with $\theta,\phi$ in place of $q,p$), but
doesn't have a momentum map.
\end{exercise}

\begin{exercise} \label{Petzval mommap}
Show that the Petzval invariant for Fermat's principle in axisymmetric,
translation-invariant media is a momentum map, $T^*\mathbb{R}^2\mapsto
sp(2,\mathbb{R})^*$ taking $(\mathbf{q},\mathbf{p})\mapsto(X,Y,Z)$. What is its
symmetry? What is its Hamiltonian vector field? 
\end{exercise}

\begin{theorem}\label{Cotlift-mommap}
\textbf{(also due to Noether) }Let $G$ act on $Q,$ and
by cotangent lifts on $T^{\ast}Q.$ 
Then $J:T^{\ast
}Q\rightarrow\mathfrak{g}^{\ast}$ defined by, for every $\xi\in\mathfrak{g},$%
\[
J_{\xi}\left(  \alpha_{q}\right)  =\left\langle \alpha_{q},\xi_{Q}\left(
q\right)  \right\rangle ,\text{ for every }\alpha_{q}\in T_{q}^{\ast}Q,
\]
is a momentum map (the ``standard one'') 
for the $G$ action with respect to the classical Poisson bracket.
\end{theorem}

(A proof using symplectic forms is given in Marsden and Ratiu [2003].)\\

\textbf{Proof} We need to show that $X_{J_\xi} = \xi_{T^*Q},$
for every $\xi\in\mathfrak{g}.$ From the definition of Hamiltonian vector fields, 
this is equivalent to showing that $\xi_{T^*Q}[F] = \{F,J_\xi\}$ for every $F\in \mathcal{F}(T^*Q).$
We verify this for finite-dimensional $Q$ by using
cotangent-lifted local coordinates.
\begin{align*}
\frac{\partial J_\xi}{\partial p}(q,p) & = \xi_Q(q) \\
\frac{\partial J_\xi}{\partial q^i} (q,p)&= 
\left< p, \frac{\partial}{\partial q^i}\left(\xi_Q(q)\right)\right> \\
&=\left< p, \frac{\partial}{\partial q^i}\left( \left.\frac{\partial}{\partial t} \Phi_{\left(\exp (t\xi)\right)}( q) \right|_{t=0}\right) \right> 
&=\left< p, \frac{\partial}{\partial t}\left( \left.\frac{\partial}{\partial q^i} 
\Phi_{\left(\exp (t\xi)\right)} (q) \right)\right|_{t=0} \right> \\
&=\left.\frac{\partial}{\partial t} \left< p, 
T\Phi_{\left(\exp (t\xi)\right)} \frac{\partial}{\partial q^i}(q)  \right> \right|_{t=0}
&=\left.\frac{\partial}{\partial t} \left< T^*\Phi_{\left(\exp (t\xi)\right)} p, 
\frac{\partial}{\partial q^i}(q)  \right> \right|_{t=0} \\
&= \left< -\xi_{T^*Q}(q,p),\frac{\partial}{\partial q^i}(q)\right> \\
\frac{\partial J_\xi}{\partial q} (q,p)& = -\xi_{T^*Q}(q,p)
\end{align*}
So for every $F\in \mathcal{F}(T^*Q),$
\begin{align*}
\xi_{T^*Q}[F] &= \left. \frac{\partial}{\partial t} F\left(\exp (t\xi) q, \exp(t\xi) p\right) \right|_{t=0} \\
&= \frac{\partial F}{\partial q}  \xi_Q(q) + 
\frac{\partial F}{\partial p}  \xi_{T^*Q}(q,p)
&= \frac{\partial F}{\partial q}  \frac{\partial J_\xi}{\partial p} -
\frac{\partial F}{\partial p}  \frac{\partial J_\xi}{\partial q}
= \{F,J_\xi\}
\end{align*}
$\quad \blacksquare$

\begin{example}
Let $G\subset M_{n}\left(  \mathbb{R}\right)  $ be a
matrix group, with cotangent-lifted action on
$(q,p)\in T^{\ast}\mathbb{R}^{n}.$  For every $\mathrm{g}\subset
M_{n}\left(\mathbb{R}\right)$, $q\mapsto g q$. The cotangent-lifted
action is $(q,p)\mapsto(g q,g^{-T}p)$. Thus, writing
$g=\exp(t\xi)$, the linearization of this group action yields the
vector field
\[
X_\xi=(\xi q,-\,\xi^T p)
\]
The corresponding Hamiltonian equations read
\[
\xi q=\frac{\partial J_{\xi}}{\partial p}
\,,\quad
-\,\xi^T p=-\,\frac{\partial J_{\xi}}{\partial q}
\]
This yields the momentum map $J(q,p)$ given by
\[
J_{\xi}\left(q,p\right)  
= \langle J(q,p),\xi\rangle
=p^{T}\xi_{Q}\left(  q\right)  =p^{T}\xi q
\,.
\]
In coordinates, $p^{T}\xi q=p_i\xi^i_jq^j$, so $J(q,p)=q^ip_j$.
\end{example}
\begin{exercise}
Calculate the momentum map of the cotangent lifted
action of the group of translations of $\mathbb{R}^{3}.$
\end{exercise}
\begin{solution}
$\mathbf{x}\in\mathbb{R}^{3}$ acts on
$\mathbf{q}\in\mathbb{R}^{3}$ by addition of vectors,
\[
\mathbf{x}\cdot(\mathbf{q})
=\mathbf{q}+\mathbf{x}
\,.
\]
The infinitesimal generator is 
$\lim_{\mathbf{x}\to0}\frac{d}{d\mathbf{x}}
(\mathbf{q}+\mathbf{x})={\rm Id}$. Thus, $\xi_\mathbf{q}=Id$ and
\[
\langle J_k,\xi\rangle
=\langle (\mathbf{q},\mathbf{p}),\xi_\mathbf{q}\rangle
=\langle \mathbf{p},Id\rangle
=p_i\delta^i_k=p_k
\]
This is also Hamiltonian with $J_\xi=\mathbf{p}$, so that
$\{\mathbf{p},J_\xi\}=0$ and $\{\mathbf{q},J_\xi\}=Id$.
\end{solution}

\begin{example}\label{ex-ad-star}
Let $G$ act on itself by left multiplication, and by
cotangent lifts on $T^{\ast}G.$ We first note that the infinitesimal action on
$G$ is
\[
\xi_{G}\left(  g\right)  =\frac{d}{dt}\left.  \exp\left(  t\xi\right)
g\right|  _{t=0}=TR_{g}\xi.
\]
Let $J_L$ be the  momentum map for this action.
Let $\alpha_{g}\in T_{g}^{\ast}G.$ For every $\xi\in\mathfrak{g},$ we have
\[
\left\langle J_L\left(  \alpha_{g}\right)  ,\xi\right\rangle =\left\langle
\alpha_{g},\xi_{G}\left(  g\right)  \right\rangle =\left\langle \alpha
_{g},TR_{g}\xi\right\rangle =\left\langle TR_{g}^{\ast}\alpha_{g}%
,\xi\right\rangle
\]
so $J_L\left(  \alpha_{g}\right)  =TR_{g}^{\ast}\alpha_{g}.$ Alternatively,
writing $\alpha_{g}=T^{\ast}L_{g^{-1}}\mu$ for some $\mu\in\mathfrak{g}^{\ast
}$ we have
\[
J_L\left(  T^{\ast}L_{g^{-1}}\mu\right)  =TR_{g}^{\ast}T^{\ast}L_{g^{-1}}%
\mu=Ad_{g^{-1}}^{\ast}\mu.
\]
\end{example}

\begin{exercise}
Show that the momentum map for the right multiplication action $R_g (h)=hg$ is 
$J_R\left(  \alpha_{g}\right)  =TL_{g}^{\ast}\alpha_{g}.$
\end{exercise}

For matrix groups, the tangent lift of the left (or right) multiplication action is
again matrix multiplication. Indeed, to compute $TR_G(A)$ for any $A\in T_QSO(3),$ let
$B(t)$ be a path in $SO(3)$ such that $B(0)=Q$ and $B'(0)=A.$ Then
\[
TR_{G}(A) =  \left.\frac{d}{dt} B(t)G\right|_{t=0} = AG.
\]
Similarly, $TL_G(A)=GA.$
To compute the cotangent lift similarly, we need to be able to consider elements of
$T^*G$ as matrices.  This can be done using any nondegenerate bilinear form on each
tangent space $T_QG.$ 
We will use the pairing defined by 
\[
\left<\left<A,B\right>\right> 
:= -\,\frac{1}{2}tr \left(A^TB\right)
=\,-\frac{1}{2}tr \left(AB^T\right).
\]
(The equivalence of the two formulas follows from the properties $tr(CD)=tr(DC)$ and
$tr(C^T)=tr(C)$).

\begin{exercise}
Check that this pairing, restricted to $so(3),$ corresponds to the Euclidean inner
product  via the hat map.
\end{exercise}

\begin{example} \label{mommap-ex}
Consider the previous example for a \textbf{matrix} group $G.$ 
For any $Q\in G,$ the pairing given above allows use to consider any
element $P\in T^*_QG$  as a matrix. The natural pairing of $T^*_QG$ with
$T_QG$ now has the formula,
\[
\left<P,A\right> = -\frac{1}{2}tr\left(P^TA\right), \quad \textrm{for all } A \in T_QG.
\]
We compute the {\bfi cotangent-lifts of the left and right multiplication actions}:
\begin{align*}
\left<T^*L_Q(P),A\right> &= \left<P,TL_Q(A)\right> = \left<P,QA\right> \\ 
&= -\frac{1}{2} tr\left(P^TQA\right)
= -\frac{1}{2} tr\left(\left(Q^TP\right)^TA\right)
=\left<Q^TP,A\right> \\
\left<T^*R_Q(P),A\right> &= \left<P,TR_Q(A)\right> = \left<P,AQ\right> \\ 
&= -\frac{1}{2} tr\left(P(AQ)^T\right)
= -\frac{1}{2} tr\left(PQ^TA^T\right)
=\left<PQ^T,A\right>
\end{align*}
In summary,
\begin{align*}
T^*L_Q(P) = Q^TP \quad \textrm{and} \quad
T^*R_Q(P) = PQ^T
\end{align*}
We thus compute the momentum maps as
\begin{align*}
J_L\left(Q,P\right) &= T^*R_QP =PQ^T \\
J_R\left(Q,P\right) &= T^*L_{Q}P =Q^TP
\end{align*}
In the special case of $G=SO(3),$ these matrices $PQ^T$ and $Q^TP$ are skew-symmetric,
since they are elements of $so(3).$ Therefore,
\begin{align*}
J_L\left(Q,P\right) &= T^*R_QP =\frac{1}{2}\left(PQ^T - QP^T\right)\\
J_R\left(Q,P\right) &= T^*L_{Q}P = \frac{1}{2}\left(Q^TP - P^TQ\right)
\end{align*}
\end{example}

\begin{exercise}
Show that the cotangent lifted action on $SO(n)$ is expressed as
\[
Q\cdot P = Q^TP
\]
as matrix multiplication.
\end{exercise}

\begin{definition}
A momentum map is said to be {\bfi equivariant} when it is equivariant
with respect to the given action on
$P$ and the coadjoint action on $\mathfrak{g}^*$. That is,
\[
J(g \cdot p) = {\rm Ad}_{g^{-1}}^* J(p)
\]
for every $g\in G$, $p \in P$, where $g \cdot p$ denotes the action of $g$
on the point $p$ and where Ad denotes the adjoint action.
\end{definition}

\begin{exercise}
Show that the momentum map derived  from the cotangent lift in {\bf
Theorem \ref{Cotlift-mommap}} is equivariant.
\end{exercise}

\begin{example} [Momentum map for symplectic representations]
\label{symp-mommap}

Let $(V, \Omega) $ be a symplectic vector space and let $G $ be a Lie
group acting linearly and symplectically on $V $. This action admits an
equivariant momentum map $\mathbf{J}: V \rightarrow \mathfrak{g}$ given
by
\[
J^ \xi(v) =
\langle \mathbf{J}(v), \xi \rangle = \frac{1}{2}\Omega(\xi\cdot v , v ),
\]
where $\xi\cdot v $ denotes the Lie algebra representation of the
element $\xi \in  \mathfrak{g}$ on the vector $v \in V $. To verify
this, note that the infinitesimal generator $\xi_V(v) = \xi \cdot v $,
by the definition of the Lie algebra representation induced by the given
Lie group representation, and that $\Omega( \xi \cdot u, v ) = -\Omega(
u, \xi\cdot v ) $ for all $u, v \in V $. Therefore
\[
\mathbf{d}J^ \xi (u) (v) = \frac{1}{2}\Omega(\xi\cdot u , v ) +
\frac{1}{2}\Omega(\xi\cdot v, u ) = \Omega(\xi\cdot u , v ).
\]
Equivariance of $\mathbf{J}$ follows from the obvious relation $g ^{-1}
\cdot \xi \cdot g \cdot v = (\operatorname{Ad}_{g ^{-1}} \xi) \cdot v $
for any $g \in G $, $ \xi\in \mathfrak{g}$, and $v \in V $.
\end{example}

\begin{example} [Cayley-Klein parameters and the Hopf fibration]
\label{CK-param}

Consider the natural action of $SU(2)$ on $\mathbb{C}^2 $. Since this
action is by isometries of the Hermitian metric, it is automatically
symplectic and therefore has a momentum map $\mathbf{J}: \mathbb{C} ^2
\rightarrow \mathfrak{su}(2)^\ast$ given in example \ref{symp-mommap},
that is, 
\[
\langle \mathbf{J}(z, w), \xi \rangle = \frac{1}{2}
\Omega(\xi\cdot (z,w), (z, w)),
\]
where $z, w \in \mathbb{C}$ and $\xi \in \mathfrak{su}(2)$. Now the
symplectic form on $\mathbb{C}^2$ is given by minus the imaginary part of
the Hermitian inner product. That is, $\mathbb{C}^{n}$ has Hermitian
inner product given by ${\bf z} \cdot \mathbf{w}: = \sum_{j=1}^n z_j
\overline{w}_j$, where ${\bf z} = (z_1, \dots, z_n), \mathbf{w}=
(w_1, \dots, w_n) \in \mathbb{C}^n $. The symplectic form is thus
given by $\Omega({\bf z}, \mathbf{w}) : = -\operatorname{Im}({\bf
z}\cdot \mathbf{w})$ and it is identical to the one given before on
$\mathbb{R}^{2n}$ by identifying ${\bf z} = \mathbf{u} + i
\mathbf{v} \in \mathbb{C}^n$ with $(\mathbf{u}, \mathbf{v}) \in
\mathbb{R}^{2n}$ and $\mathbf{w} = \mathbf{u}' + i \mathbf{v}' \in
\mathbb{C}^n$ with $(\mathbf{u}', \mathbf{v}') \in \mathbb{R}^{2n}$.

The Lie algebra $\mathfrak{su}(2)$ of $SU(2)$
consists of $2 \times 2 $ skew Hermitian matrices of trace zero. This Lie
algebra is isomorphic to $\mathfrak{so}(3) $ and therefore to
$(\mathbb{R}^3, \times )$ by the isomorphism given by
\[
\mathbf{x} = (x^1, x^2,x^3) \in \mathbb{R}^3 \mapsto
\widetilde{\mathbf{x}} : =
\frac{1}{2}
\left[
\begin{array}{cc}
-ix^3&-ix^1 - x^2\\
-ix^1 + x^2&ix^3
\end{array}
\right] \in \mathfrak{su}(2).
\]
Thus we have $[\widetilde{\mathbf{x}}, \widetilde{\mathbf{y}}] =
(\mathbf{x}\times \mathbf{y})\widetilde{\phantom{y}}$ for any
$\mathbf{x}, \mathbf{y} \in \mathbb{R}^3$. Other useful relations are
$\operatorname{det}(2\widetilde{\mathbf{x}}) = \| \mathbf{x}\|^2 $ and
$\operatorname{trace}(\widetilde{\mathbf{x}}\widetilde{\mathbf{y}}) =
-\frac{1}{2} \mathbf{x} \cdot \mathbf{y}$. Identify
$\mathfrak{su}(2)^\ast$ with $\mathbb{R}^3$ by the map
$ \mu \in \mathfrak{su}(2)^\ast \mapsto \check{\mu} \in \mathbb{R}^3$
defined by
\[
\check{\mu} \cdot \mathbf{x}: = -2\langle \mu, \widetilde{\mathbf{x}}
\rangle
\]
for any $\mathbf{x} \in \mathbb{R}^3$.
With these notations, the
momentum map $\check{\mathbf{J}}: \mathbb{C}^2 \rightarrow \mathbb{R}^3$
can be explicitly computed in coordinates: for any $\mathbf{x}\in
\mathbb{R}^3$ we have
\begin{align*}
\check{\mathbf{J}}(z, w) \cdot \mathbf{x} &=
-2\langle \mathbf{J}(z, w), \widetilde{\mathbf{x}} \rangle \\
&= \frac{1}{2}\operatorname{Im}\left(
\left[
\begin{array}{cc}
-ix^3&-ix^1 - x^2\\
-ix^1 + x^2&ix^3
\end{array}
\right]
\left[
\begin{array}{c}
z\\
w
\end{array}
\right]
\cdot
\left[
\begin{array}{c}
z\\
w
\end{array}
\right]
\right)\\
&= -\frac{1}{2}(2 \operatorname{Re}(w \overline{z}),
2\operatorname{Im}(w\overline{z}), |z|^2 - |w|^2) \cdot \mathbf{x}.
\end{align*}
Therefore
\[
\check{\mathbf{J}}(z, w) = -\frac{1}{2}(2w \overline{z}, |z|^2 - |w|^2) \in
\mathbb{R}^3.
\]
Thus, $\check{\mathbf{J}}$ is a Poisson map from
$\mathbb{C}^2 $, endowed with the  canonical symplectic structure, to
$\mathbb{R}^3 $, endowed with the $+$ Lie Poisson structure.
Therefore, $- \check{\mathbf{J}}: \mathbb{C}^2 \rightarrow \mathbb{R}^3 $
is a canonical map, if $\mathbb{R}^3$  has the $-$ Lie-Poisson bracket
relative to which the free rigid body equations  are Hamiltonian.
Pulling back the Hamiltonian $H(\boldsymbol{\Pi}) = \boldsymbol{\Pi}
\cdot {\mathbb I}^{-1} \boldsymbol{\Pi}/2$ to $\mathbb{C}^2$ gives a
Hamiltonian function (called collective) on $\mathbb{C}^2$. The
classical Hamilton equations for this function are therefore
projected by $-\check{\mathbf{J}}$ to the rigid body equations $\dot
{\boldsymbol{\Pi}} = \boldsymbol{\Pi} \times {\mathbb I}^{-1}
\boldsymbol{\Pi}$. In this context, the variables $(z, w) $ are
called the {\bfi Cayley-Klein parameters\/}. 
\end{example} 

\begin{exercise}
Show  that $-\check{\mathbf{J}}|_{S^3}: S^3 \rightarrow S^2$
is the {\bfi Hopf fibration\/}. In other words, \textit{the momentum
map of the $SU(2)$-action on $\mathbb{C}^2$, the Cayley-Klein
parameters and the family of
Hopf fibrations on concentric three-spheres in $\mathbb{C}^2$ are all the
same map}.
\end{exercise}

\begin{exercise} {\bf Optical traveling wave pulses}
The equation for the evolution of the complex amplitude of a polarized
optical traveling wave pulse in a material medium is given as
\[
\dot{z}_i 
=
\frac{1}{\sqrt{-1}}\,
\frac{\partial H}{\partial z_i^*}
\]
with Hamiltonian $H:\mathbb{C}^2\to\mathbb{R}$ defined by
\[
H
=
z_i^*\chi^{(1)}_{ij}z_j
+
3z_i^*z_j^*\chi^{(3)}_{ijkl}z_kz_l
\]
and the constant complex tensor coefficients $\chi^{(1)}_{ij}$ and
$\chi^{(1)}_{ijkl}$ have the proper Hermitian and permutation symmetries
for $H$ to be real. Define the Stokes vectors
by the isomorphism,
\[
\mathbf{u} = (u^1, u^2,u^3) \in \mathbb{R}^3 \mapsto
\widetilde{\mathbf{u}} : =
\frac{1}{2}
\left[
\begin{array}{cc}
-iu^3      &  -iu^1 - u^2
\\
-iu^1 + u^2&   iu^3
\end{array}
\right] \in \mathfrak{su}(2).
\]
\begin{enumerate}
\item
Prove that this isomorphism is an equivariant momentum map.
\item
Deduce the equations of motion for the Stokes vectors of this optical
traveling wave and write it as a Lie Poisson Hamiltonian system. 
\item
Determine how this system is related to the equations for an $SO(3)$ rigid
body.
\end{enumerate}
\end{exercise}

\begin{exercise}

The formula determining the momentum map for the cotangent-lifted action
of a Lie group $G$ on a smooth manifold $Q$ may be expressed in terms of
the pairing $\langle\,\cdot\,,\,\cdot\,\rangle: 
\mathfrak{g}^*\times\mathfrak{g}\mapsto\mathbb{R}$
as 
\[
\langle\,J\,,\,\xi\,\rangle
=
\langle\,p\,,\,\pounds_\xi q\,\rangle
\,,
\]
where $(q,p)\in T_q^*Q$ and $\pounds_\xi q$ is the infinitesimal generator of the action of
the Lie algebra element $\xi$ on the coordinate $q$. 

Define appropriate pairings and determine the momentum maps explicitly for the following
actions,

\begin{description}
\item{[a]} 
$\pounds_\xi q=\xi\times q$ for $\mathbb{R}^3\times\mathbb{R}^3\mapsto\mathbb{R}^3$
\item{[b]} 
$\pounds_\xi q={\rm ad}_\xi q$ for ad-action
${\rm ad}:\,\mathfrak{g}\times\mathfrak{g}\mapsto\mathfrak{g}$ in a Lie algebra
$\mathfrak{g}$
\item{[c]} 
$AqA^{-1}$ for $A\in GL(3,R)$ acting on $q\in GL(3,R)$ by matrix conjugation
\item{[d]} 
$Aq$ for left action of $A\in SO(3)$ on $q\in SO(3)$
\item{[e]} 
$AqA^T$ for $A\in GL(3,R)$ acting on $q\in Sym(3)$, that is $q=q^T$.
\end{description} 

\end{exercise}

\begin{answer}$\quad$\\
\begin{description}
\item{[a]} 
$p \cdot \xi\times q = q \times p \cdot \xi \Rightarrow J = q \times p$. (The pairing
is scalar product of vectors.) 
\item{[b]} 
$\langle\,p\,,\,{\rm ad}_\xi q\,\rangle 
=
-\,\langle\,{\rm ad}^*_q\,p\,,\,\xi\,\rangle
\Rightarrow J = {\rm ad}^*_q\,p$ for the pairing $\langle\,\cdot\,,\,\cdot\,\rangle: 
\mathfrak{g}^*\times\mathfrak{g}\mapsto\mathbb{R}$
\item{[c]} 
Compute $T_e(AqA^{-1})=\xi q - q\xi = [\xi,q]$ for $\xi=A^\prime(0)\in gl(3,R)$ acting on
$q\in GL(3,R)$ by matrix Lie bracket $[\cdot\,,\,\cdot]$. For the matrix pairing 
$\langle\,A\,,\,B\,\rangle = {\rm trace}(A^TB)$, we have 
${\rm trace}(p^T[\xi,q])={\rm trace}((pq^T-q^Tp)^T\xi)\Rightarrow J = pq^T-q^Tp$.
\item{[d]} 
Compute $T_e(Aq)=\xi q$ for $\xi=A^\prime(0)\in so(3)$ acting on
$q\in SO(3)$ by left matrix multiplication. For the matrix pairing 
$\langle\,A\,,\,B\,\rangle = {\rm trace}(A^TB)$, we have 
${\rm trace}(p^T\xi q)={\rm trace}((pq^T)^T\xi)\Rightarrow J = \frac{1}{2}(pq^T-q^Tp)$,
where we have used antisymmetry of the matrix $\xi\in so(3)$.
\item{[e]} 
Compute $T_e(AqA^T)=\xi q + q \xi^T$ for $\xi=A^\prime(0)\in gl(3,R)$ acting on $q\in
Sym(3)$. For the matrix pairing 
$\langle\,A\,,\,B\,\rangle = {\rm trace}(A^TB)$, we have 
${\rm trace}(p^T(\xi q + q \xi^T))={\rm trace}(q(p^T+p)\xi)={\rm
trace}(2qp)^T\xi)\Rightarrow J = 2qp$, where we have used symmetry of the matrix $\xi q + q
\xi^T$ to choose $p=p^T$. (The momentum canonical to the symmetric matrix $q=q^T$ should be
symmetric to have the correct number of components!)
\end{description} 

\end{answer}

\bigskip
\noindent \textbf{Equivariance}

\begin{definition}
A momentum map is $Ad^{\ast}$-\emph{equivariant} iff
\[
J\left(  g\cdot x\right)  =Ad_{g^{-1}}^{\ast}J\left(  x\right)
\]
for all $g\in G,x\in P.$
\end{definition}

\begin{proposition}
All cotangent-lifted actions are $Ad^{\ast}$-equivariant.
\end{proposition}

\begin{proposition}
Every $Ad^{\ast}$-equivariant momentum map $J:P\to \mathfrak{g}^*$ is a Poisson map,
with respect to the `+' Lie-Poisson bracket on $\mathfrak{g}^*.$
\end{proposition}

\section{Quick summary for momentum maps}

Let $G$ be a Lie group, $\mathfrak{g}$ its Lie algebra, and
let $\mathfrak{g}^*$ be its dual. Suppose that $G$ acts symplectically on a
symplectic manifold $P$ with symplectic form denoted by $\Omega$. Denote
the infinitesimal generator associated with the Lie algebra element
$\xi$ by $\xi_P$ and let the Hamiltonian vector field associated
to a function $f : P\to\mathbb{R}$ be denoted $X_f$, so that
$df=X_f\contract\Omega$.%
\rem{Some of the discussion below is adapted from Marsden, Ratiu and
Scheurle [2000] and Marsden and Weinstein [2001].}

\subsection{Definition, History and Overview}
A {\bfi momentum map} $J : P\to\mathfrak{g}^*$ is defined by the
condition relating the infinitesimal generator $\xi_P$ of a symmetry to
the vector field of its corresponding conservation law,
$\langle J,\xi\,\rangle$,
\[
\xi_P = X_{\langle J,\xi\,\rangle}
\]
for all $\xi\in \mathfrak{g}$. Here $\langle J,\xi\,\rangle :
P\to\mathbb{R}$ is defined by the natural pointwise pairing.\\

A momentum map is said to be {\bfi equivariant} when it is equivariant
with respect to the given action on
$P$ and the coadjoint action on $\mathfrak{g}^*$. That is,
\[
J(g \cdot p) = {\rm Ad}^*_{g^{-1}} J(p)
\]
for every $g\in G$, $p \in P$, where $g \cdot p$ denotes the action of $g$
on the point $p$ and where Ad denotes the adjoint action. \\

\noindent
According to \cite{We1983a}, \cite{Lie1890} already knew that 
\begin{enumerate}
\item
An action of a Lie group $G$ with Lie algebra
$\mathfrak{g}$ on a symplectic manifold $P$ should be accompanied by such
an equivariant momentum map $J : P\to\mathfrak{g}^*$ and 
\item
The orbits of this action are themselves symplectic manifolds. 
\end{enumerate}
The links with mechanics were developed in the work of Lagrange, Poisson,
Jacobi and, later, Noether. In particular, Noether showed that a momentum
map for the action of a group $G$ that is a symmetry of the Hamiltonian
for a given system is a {\bfi conservation law} for that system. \\ 

In modern form, the momentum map and its equivariance were rediscovered in
\cite{Ko1966} and \cite{So1970} in the general symplectic case,
and in \cite{Sm1970} for the case of the lifted action from a manifold
$Q$ to its cotangent bundle $P = T^*Q$. In this case, the equivariant
momentum map is given explicitly by 
\[
\langle J(\alpha_q),\xi\,\rangle  = \langle \alpha_q,\xi_Q(q)\rangle
\,,
\]
where $\alpha_q \in T^*Q$, $\xi\in\mathfrak{g}$, and where the angular
brackets denote the natural pairing on the appropriate spaces. See
\cite{MaRa1994} and \cite{OrRa2004} for additional history and
description of the momentum map and its properties.

\section{Rigid body equations on $\operatorname{SO}(n)$}

Recall from \cite{Man1976} and \cite{Ra1980} that the left
invariant generalized rigid body equations on
$\operatorname{SO}(n)$ may be written as
\begin{align}
\dot Q&= Q\Omega\,,  \nonumber \\
\dot M& = M\Omega-\Omega M =: [M,\Omega]\,, \tag{RBn}
\label{rbl}
\end{align}
where $Q\in \operatorname{SO}(n)$ denotes the configuration space
variable (the attitude of the body), $\Omega=Q^{-1}\dot{Q} \in
{so}(n)$ is the body angular velocity, and
\[
M:=J(\Omega)=D^2\Omega +\Omega D^2 \in
{so}^*(n)\,,
\]
          is the body angular momentum. Here
$J: {so}(n) \rightarrow  {so}(n)^* $ is the symmetric
(with respect to the above inner product) positive definite operator
defined by
\[
J(\Omega)=D^2\Omega +\Omega D^2 ,
\]
          where $D^2$ is
the square of the constant diagonal matrix $D={\rm diag}\,\{d_1,d_2,d_3\}$
satisfying $d^2_i + d^2_j >0$ for all $i \neq j$. For $n=3$ the elements of
$d^2_i$ are related to the standard diagonal moment of inertia tensor $I$
by 
\[
I={\rm diag}\,\{I_1,I_2,I_3\}
\,,\quad
I_1 = d^2_2 + d^2_3
\,,\quad
I_2 = d^2_3 + d^2_1
\,,\quad
I_3 = d^2_1 + d^2_2
\,.
\]

The Euler equations for the $SO(n)$ rigid body 
$ \dot{ M } =  [ M, \Omega ] $ are readily checked to be the
Euler-Lagrange equations on
${so}(n)$ for the Lagrangian
\[
L(Q,\dot{Q})
=
l ( \Omega ) = \frac{1}{2}  \left\langle  \Omega , J( \Omega )
\right\rangle 
\,,\quad\hbox{with}\quad
\Omega = Q^T\dot{Q}
.
\]
The momentum is found via the Legendre transformation to
be
\[
\frac{\partial l}{\partial \Omega}
=
J( \Omega )
=
M
\,,
\]
and the corresponding Hamiltonian is
\[
H (M)
=
\frac{\partial l}{\partial \Omega}\cdot \Omega - l ( \Omega )
= \frac{1}{2}  \left\langle  M , J^{-1}(M)
\right\rangle \,.
\]
The quantity $M$ is the angular momentum in the body frame. The
corresponding angular momentum in space,
\[
m=QMQ^T
\,,\quad\hbox{is conserved}\quad
\dot{m}=0
\,.
\]
Indeed, conservation of spatial angular momentum $m$ implies Euler's
equations for the body angular momentum $M=Q^TmQ={\rm Ad}^*_Qm$.

\subsection{Implications of left invariance}

This Hamiltonian $H(M)$ is invariant under the action of $SO(n)$ from the
left. The corresponding conserved momentum map under this symmetry is
known from the previous lecture as
\[
J_L:T^*SO(n)\mapsto so(n)^*
\quad\hbox{is}\quad
J_L(Q,P)=PQ^T
\]
On the other hand, we know (from Lectures 18 \& 19) that the momentum map
for right action is 
\[
J_R:T^*SO(n)\mapsto so(n)^*
\,,\quad
J_R(Q,P)=Q^TP
\]
Hence $M=Q^TP=J_R$. Therefore, one computes
\begin{eqnarray*}
H(Q,P)&=&H(Q,Q\cdot M)=H({\rm Id},M)
\quad\hbox{(by left invariance)}
\\
&=& H(M)=\frac{1}{2}\langle M\,,\,J^{-1}(M)\rangle
\\
&=&\frac{1}{2}\langle Q^TP\,,\,J^{-1}(Q^TP)\rangle
\end{eqnarray*}
Hence, we may write the $SO(n)$ rigid body Hamiltonian as 
\[
H(Q,P)
=\frac{1}{2}\langle Q^TP\,,\,\Omega(Q,P)\rangle
\]
\rem{
Left invariance of $H(Q,P)=H(WQ,WP)$ for all $W\in SO(n)$ means that
\[
H(Q,P)=\frac{1}{2}\langle(WPQ^TW^T)\,,\,\mathcal{J}^{-1}(WPQ^TW^T)\rangle
\,,\quad\forall\ W\in SO(n)
\]
Choose $W^T=Q$, then left invariance of $H$ implies
\[
H(Q,P)=\frac{1}{2}\langle(Q^TP)\,,\,\mathcal{J}^{-1}(Q^TP)\rangle
=\frac{1}{2}\langle Q^TP\,,\,\Omega(Q,P)\rangle
\]
That is,
\[
H(Q,P)
=\frac{1}{2}\langle J_L\,,\,\mathcal{J}^{-1}J_L\rangle
=\frac{1}{2}\langle J_R\,,\,\mathcal{J}^{-1}J_R\rangle
\]
The variational derivatives of $H(Q,P)$ may be computed using the pairing
\[
\left\langle A\,,\,B\right\rangle
=
{\rm tr}(A^TB)
=
{\rm tr}(AB^T)
=
{\rm tr}(BA^T)
=
{\rm tr}(B^TA)
\]
}
Consequently, the variational derivatives of $H(Q,P)=\frac{1}{2}\langle
Q^TP\,,\,\Omega(Q,P)\rangle$ are
\begin{eqnarray*}
\delta H 
&=&  \left\langle
Q^T\delta{P}+\delta{Q}^TP\,,\,\Omega(Q,P)
\right\rangle 
\\
&=&  {\rm tr}(\delta{P}^TQ\Omega) + {\rm tr}(P^T\delta{Q}\Omega)
\\
&=&  {\rm tr}(\delta{P}^TQ\Omega) + {\rm tr}(\delta{Q}\Omega P^T)
\\
&=&  {\rm tr}(\delta{P}^TQ\Omega) + {\rm tr}(\delta{Q}^T P\Omega^T)
\\
&=&  \left\langle
\delta{P}\,,\,Q\Omega
\right\rangle 
-
\left\langle
\delta{Q}\,,\,P\Omega
\right\rangle 
\end{eqnarray*}
where skew symmetry of $\Omega$ is used in the last step, i.e.,
$\Omega^T=-\Omega$. Thus, Hamilton's
canonical equations take the form,
\begin{align}\label{HamCan-eqns}
\dot Q&= \frac{\delta H}{\delta P}= Q\Omega 
\,,
\nonumber \\
\dot P&= -\,\frac{\delta H}{\delta Q}= P\Omega
\,.
\end{align}
Equations (\ref{HamCan-eqns}) are the {\bfi symmetric generalized rigid body
equations}, derived earlier in \cite{BlCr1996} and \cite{BlBrCr1997}
from the viewpoint of optimal control. Combining them yields,
\[
Q^{-1}\dot{Q}=\Omega =P^{-1}\dot{P}
\Longleftrightarrow
(PQ^T)\dot{\,}=0
\,,
\]
in agreement with conservation of the momentum map $J_L(Q,P)=PQ^T$
corresponding to symmetry of the Hamiltonian under left action of $SO(n)$.
This momentum map is the angular momentum in space, which is related to
the angular momentum in the body by $PQ^T=m=QMQ^T$. Thus, we recognize the
canonical momentum as $P=QM$ (see exercise \ref{ex-ad-star}), and the
momentum maps for left and right actions as,
\begin{eqnarray*}
J_L&=&m=PQ^T\quad\hbox{(spatial angular momentum)}
\\
J_R&=&M=Q^TP\quad\hbox{(body angular momentum)}
\end{eqnarray*}
Thus, momentum maps $TG^*\mapsto\mathfrak{g}^*$ corresponding to symmetries of the
Hamiltonian produce conservation laws; while momentum maps
$TG^*\mapsto\mathfrak{g}^*$ which do {\it not} correspond to symmetries may
be used to re-express the equations on $\mathfrak{g}^*$, in terms of variables on
$TG^*$. 

\section{Manakov's formulation of the $SO(4)$ rigid body}

The Euler equations on $SO(4)$ are
\[
\frac{dM}{dt}= M\Omega - \Omega M =[M,\Omega]
\,,
\tag{RBn}
\]
where $\Omega$ and $M$ are skew symmetric $4\times4$ matrices. The angular
frequency $\Omega$ is a linear function of the angular momentum, $M$.
\cite{Man1976} ``deformed'' these equations into
\[
\frac{d}{dt}(M+\lambda A)=[(M+\lambda A),(\Omega+\lambda B)]
\,,
\]
where $A$, $B$ are also skew symmetric $4\times4$ matrices and $\lambda$
is a scalar constant parameter. For these equations to hold for any value
of $\lambda$, the coefficent of each power must vanish. 

\begin{itemize}
\item
The coefficent of $\lambda^2$ is
\[
0=[A,B]
\]
So $A$ and $B$ must commute. So, let them be constant and diagonal:
\[
A_{ij}={\rm diag}(a_i)\delta_{ij}
\,,\quad
B_{ij}={\rm diag}(b_i)\delta_{ij}
\tag{no sum}
\]
\item
The coefficent of $\lambda$ is
\[
0=\frac{dA}{dt}=[A,\Omega]+[M,B]
\]
Therefore, by antisymmetry of $M$ and $\Omega$,
\[
(a_i-a_j)\Omega_{ij}=(b_i-b_j)M_{ij}
\qquad\Longleftrightarrow\qquad
\Omega_{ij}
=
\frac{b_i-b_j}{a_i-a_j}M_{ij}
\tag{no sum}
\]
\item
Finally, the coefficent of $\lambda^0$ is the Euler equation,
\[
\frac{dM}{dt}=[M,\Omega]
\,,
\]
but now with the restriction that the moments of inertia are of the form,
\[
\Omega_{ij}
=
\frac{b_i-b_j}{a_i-a_j}M_{ij}
\tag{no sum}
\]
which turns out to possess only 5 free parameters.
\end{itemize}
With these conditions, Manakov's deformation of the $SO(4)$ rigid body
implies for every power $n$ that
\[
\frac{d}{dt}(M+\lambda A)^n=[(M+\lambda A)^n,(\Omega+\lambda B)]
\,,
\]
Since the commutator is antisymmetric, its trace vanishes and one has
\[
\frac{d}{dt}{\rm trace}(M+\lambda A)^n=0
\]
after commuting the trace operation with time derivative. Consequently, 
\[
{\rm trace}(M+\lambda A)^n={\rm constant}
\]
for each power of $\lambda$. That is, all the coefficients of each power
of $\lambda$ are constant in time for the $SO(4)$ rigid body.
\cite{Man1976} proved that these constants of motion are sufficient to
completely determine the solution. 

\begin{remark}
This result generalizes considerably. First, it holds for $SO(n)$. Indeed, as as proven
using the theory of algebraic varieties in \cite{Ha1984}, Manakov's method captures all
the algebraically integrable rigid bodies on $SO(n)$ and the moments of inertia of
these bodies possess only $2n-3$ parameters. (Recall that in Manakov's case for $SO(4)$
the moment of inertia possesses only five parameters.) Moreover,
\cite{MiFo1978} prove that every compact Lie group admits a family of left-invariant metrics with completely integrable geodesic flows.
\end{remark}

\begin{exercise}
Try computing the constants of motion ${\rm trace}(M+\lambda A)^n$ for
the values $n=2,\,3,\,4$. 
How many additional constants of motion are needed for integrability for
these cases? How many for general $n$? Hint: keep in mind that $M$ is a skew symmetric
matrix, $M^T=-M$, so the trace of the product of any diagonal matrix times an odd power
of $M$ vanishes.
\end{exercise}

\begin{answer}
The traces of the powers ${\rm trace}(M+\lambda A)^n$ are given by
\begin{eqnarray*}
\fbox{n=2}\,:&&
{\rm tr}\,M^2 + 2\lambda{\rm tr}\,(AM) + \lambda^2{\rm tr}\,A^2
\\
\fbox{n=3}\,:&&
{\rm tr}\,M^3 + 3\lambda{\rm tr}\,(AM^2) + 3\lambda^2{\rm tr}\,A^2M
+ \lambda^3{\rm tr}\,A^3
\\
\fbox{n=4}\,:&&
{\rm tr}\,M^4 + 4\lambda{\rm tr}\,(AM^3) 
+ \lambda^2(2{\rm tr}\,A^2M^2 + 4{\rm tr}\,AMAM)
+ \lambda^3{\rm tr}\,A^3M + \lambda^4{\rm tr}\,A^4
\end{eqnarray*}
The number of conserved quantities for $n=2,3,4$ are, respectively, 
one ($C_1={\rm tr}\,M^2$), one ($I_1={\rm tr}\,AM^2$) and two ($C_2={\rm tr}\,M^4$ and
$I_2=2{\rm tr}\,A^2M^2 + 4{\rm tr}\,AMAM$). The quantities $C_1$ and $C_2$ are Casimirs
for the Lie-Poisson bracket for the rigid body. Thus, $\{C_1,H\}=0=\{C_2,H\}$ for
any Hamiltonian $H(M)$; so of course $C_1$ and $C_2$ are conserved. However, each
Casimir only reduces the dimension of the system by one. The dimension of the original
phase space is dim$\,T^*SO(n)=n(n-1)$. This is reduced in half by left invariance of the
Hamiltonian to the dimension of the dual Lie algebra dim$\,so(n)^*=n(n-1)/2$. For $n=4$,
dim$\,so(4)^*=6$. One then subtracts the number of Casimirs (two) by passing to
their level surfaces, which leaves four dimensions remaining in this case. The other two
constants of motion $I_1$ and $I_2$ turn out to be sufficient for integrability, because
they are in involution $\{I_1,I_2\}=0$ and because the level surfaces of the Casimirs
are symplectic manifolds, by the Marsden-Weinstein reduction theorem \cite{MaWe1974}. 
For more details, see \cite{Ra1980}.
\end{answer}

\begin{exercise}
How do the Euler equations look on $so(4)^*$ as a matrix equation? Is there
an analog of the hat map for $so(3)^*$? Hint: the Lie algebra $so(4)$ is
locally isomorphic to $so(3)\times so(3)$.
\end{exercise}

\begin{exercise}
Write Manakov's deformation of the rigid body equations in the symmetric form
(\ref{HamCan-eqns}).
\end{exercise}

\section{Free ellipsoidal motion on $\operatorname{GL}(n)$}

Riemann \cite{Ri1860} considered the deformation of a body in $\mathbb{R}^n$ given
by 
\begin{equation}
x(t,x_0)=Q(t)\,x_0
\,,
\label{def-grad}
\end{equation}
with $x\,,x_0\in \mathbb{R}^n$, $Q(t)\in GL_+(n,\mathbb{R})$ and
$x(t_0,x_0)=x_0$, so that $Q(t_0)=Id$. (The subscript $+$ in $GL_+(n,\mathbb{R})$
means $n\times n$ matrices with positive determinant.)
Thus, $x(t,x_0)$ is the current (Eulerian)  position at time $t$ of
a material parcel that was at (Lagrangian) position $x_0$ at time
$t_0$. The ``deformation gradient,'' that is, the Jacobian matrix 
$Q=\partial{x}/\partial{x}_0$ of this ``Lagrange-to-Euler map,'' is a
function of only time, $t$,
\[
\partial{x}/\partial{x}_0=Q(t)
\,,\quad\hbox{with}\quad
\det Q>0
\,.
\] 
The velocity of such a motion is given by
\begin{equation}
\dot{x}(t,x_0)
=\dot{Q}(t)\,x_0
=\dot{Q}(t)Q^{-1}(t)\,x
=u(t,x)\,.
\end{equation}
The kinetic energy for such a body occupying a reference volume
$\mathcal{B}$ defines the quadratic form,
\[
L = \frac{1}{2}\int_{\mathcal{B}}\rho(x_0)|\dot{x}(t,x_0)|^2\,d\,^3x_0
= \frac{1}{2}{\rm tr} \Big(\dot{Q}(t)^T I \dot{Q}(t)\Big)
= \frac{1}{2}\,\dot{Q}^i_A\,I^{AB}\dot{Q}^i_B
\,.
\]
Here $I$ is the constant symmetric tensor,
\[
I^{AB}=\int_{\mathcal{B}} \rho(x_0)x_0^Ax_0^B\,d\,^3x_0
\,,
\]
which we will take as being proportional to the identity
$I^{AB}=c_0^2\delta^{AB}$ for the remainder of these considerations. This
corresponds to taking an initially spherical reference configuration for
the fluid. Hence, we are dealing with the Lagrangian consisting only of
kinetic energy,%
\footnote{\cite{Ri1860} considered the much more difficult problem of a
{\it self-gravitating} ellipsoid deforming according to (\ref{def-grad})
in $\mathbb{R}^3$. See \cite{Ch1969} for the history of this problem.}
\[
L = \frac{1}{2}{\rm tr} \Big(\dot{Q}(t)^T \dot{Q}(t)\Big)
\,.
\]
The Euler-Lagrange equations for this Lagrangian simply represent {\bfi
free motion} on the group $GL_+(n,\mathbb{R})$,
\[
\ddot{Q}(t)=0
\,,
\]
which is immediately integrable as
\[
Q(t)=Q(0)+\dot{Q}(0)t
\,,
\]
where $Q(0)$ and $\dot{Q}(0)$ are the values at the initial time $t=0$.
Legendre transforming this Lagrangian for free motion yields
\[
P=\frac{\partial L}{\partial \dot{Q}^T}=\dot{Q}\,.
\]
The corresponding Hamiltonian is expressed as
\[
H(Q,P)=\frac{1}{2}{\rm tr} \Big(P^T P\Big)=\frac{1}{2}\|P\|^2
\,.
\]
The canonical equations for this Hamiltonian are simply
\[
\dot{Q}=P
\,,\quad\hbox{with}\quad 
\dot{P}=0
\,.
\]

\subsection{Polar decomposition of free motion on $GL_+(n,\mathbb{R})$}
The deformation tensor $Q(t)\in GL_+(n,\mathbb{R})$ for such a  body may be
decomposed as
\begin{equation}
Q(t)=R^{-1}(t)D(t)S(t)
\,.
\label{polar-decomp}
\end{equation}
This is the polar decomposition of a matrix in $GL_+(n,\mathbb{R})$.
The interpretations of the various components of the motion can be
seen from equation (\ref{def-grad}). Namely, 
\begin{itemize}
\item
$R\in SO(n)$ rotates the $x$-coordinates,
\item
$S\in SO(n)$ rotates the $x_0$-coordinates in the reference
configuration\footnote{This is the ``particle relabeling map'' for this
class of motions.} and 
\item
$D$ is a diagonal matrix which represents stretching deformations along
the principal axes of the body. 
\end{itemize}
The two $SO(n)$ rotations lead to their corresponding angular frequencies,
defined by
\begin{equation}
\Omega = \dot{R}R^{-1}
\,,\quad
\Lambda= \dot{S}S^{-1}
\,.
\end{equation}
Rigid body motion will result, when $S$ restricts to the identity
matrix and $D$ is a constant diagonal matrix.  

\begin{remark}
The combined motion of a set of fluid parcels governed by
(\ref{def-grad}) along the curve $Q(t)\in GL_+(n,\mathbb{R})$ is called
``ellipsoidal,'' because it can be envisioned in three dimensions as a
fluid ellipsoid whose orientation in space is governed by $R\in SO(n)$,
whose shape is determined by $D$ consisting of its instantaneous principle
axes lengths and whose internal circulation of material is described by
$S\in SO(n)$. In addition, fluid parcels initially arranged along a
straight line within the ellipse will remain on a straight line. 
\end{remark}

\subsection{Euler-Poincar\'e dynamics of free Riemann ellipsoids}

In Hamilton's principle, $\delta\int L
\,dt = 0$, we chose a Lagrangian 
$L:TGL_+(n,\mathbb{R})\to \mathbb{R}$ in the form 
\begin{equation}
L(Q,\dot{Q}) = T(\Omega,\Lambda,D,\dot{D})
\,,
\label{Lag-RE}
\end{equation}
in which the kinetic energy $T$ is given by using the polar
decomposition
$Q(t)=R^{-1}(t)D(t)S(t)$ in (\ref{polar-decomp}), as follows.
\begin{eqnarray}
\dot{Q}=R^{-1}(-\Omega D+\dot{D}+D\Lambda )S
\,.
\end{eqnarray}
Consequently, the kinetic energy for ellipsoidal motion becomes
\begin{eqnarray}
T &=&\frac{1}{2}\,{\rm trace}
\left[- \Omega D^2\Omega - \Omega D\dot{D} 
+ \Omega D\Lambda D + \dot{D}D\Omega 
+ \dot{D}^2 
\right.
\nonumber\\
&&
\left. \hspace{.6in}
- D\Lambda^2D 
- \dot{D}\Lambda D + D\Lambda D\Omega + D\Lambda \dot{D}
\right]
\nonumber\\
&=&
\frac{1}{2}\,{\rm trace} \,
\Big[
- 
\Omega^2D^2 - \Lambda^2D^2 
+ 
\hspace{-.25in}
\underbrace{\
2\Omega D\Lambda D\,
}_{\hbox{Coriolis coupling}}
\hspace{-.2in}
+\,
 \dot{D}^2
\Big]
\,.
\label{KE-RE}
\end{eqnarray}
\begin{remark}
Note the discrete exchange symmetry of the kinetic energy: $T$ is
invariant under $\Omega\leftrightarrow\Lambda$.%
\footnote{According to \cite{Ch1969} this discrete symmetry was
first noticed by Riemann's friend, \cite{De1860}.}
\end{remark}

For $\Lambda=0$ and $D$ constant expression (\ref{KE-RE}) for $T$ reduces
to the usual kinetic energy for the rigid-body,
\begin{equation}
T\Big|_{\Lambda=0,\,D=const}
=
-\,\frac{1}{4}\,{\rm trace}
\Big[\Omega (D\Omega+\Omega D)\Big]
\,.
\label{ke-rb}
\end{equation}

This Lagrangian (\ref{Lag-RE}) is invariant under the right action, $R\to
Rg$ and $S\to Sg$, for $g\in SO(n)$. In taking variations we shall use
the formulas\footnote{These variational formulas are obtained
directly from the definitions of $\Omega$ and $\Lambda$.}
\begin{eqnarray}
\delta\Omega &=& \dot{\Sigma} + [\Sigma, \Omega]
\equiv
\dot{\Sigma} - {\rm ad}_\Omega\Sigma
\,,\quad
\Sigma \equiv \delta{R}\,R^{-1}
\,,\\
\delta\Lambda &=& \dot{\Xi} + [\Xi, \Lambda]
\equiv
\dot{\Xi} - {\rm ad}_\Lambda\Xi
\,,\quad
\Xi \equiv \delta{S}\,S^{-1}
\,,
\end{eqnarray}
in which the ad-operation is defined in terms of the
Lie-algebra (matrix) commutator $[\cdot,\cdot]$ as, e.g.,
${\rm ad}_\Omega\Sigma \equiv [\Omega,\Sigma\,]$.
Substituting these formulas into Hamilton's principle gives
\begin{eqnarray}
0 &=& \delta\!\!\int \!\!L\,dt 
=  \int \!\!dt\
\frac{\partial L}{\partial\Omega}\cdot\delta\Omega
+
\frac{\partial L}{\partial\Lambda}\cdot\delta\Lambda
+
\frac{\partial L}{\partial D}\delta D
+
\frac{\partial L}{\partial \dot{D}}\delta \dot{D}
\,,\nonumber\\
&=& \int \!\!dt\
\frac{\partial L}{\partial\Omega}\cdot
\Big[\dot{\Sigma} - {\rm ad}_\Omega\Sigma\Big]
+
\frac{\partial L}{\partial\Lambda}\cdot
\Big[\dot{\Xi} - {\rm ad}_\Lambda\Xi\Big]
+
\Big[
\frac{\partial L}{\partial D}
-
\frac{d}{dt}\frac{\partial L}{\partial \dot{D}}
\Big]\delta D
\,,\nonumber\\
&=& -\int \!\!dt
\Big[\frac{d}{dt}\frac{\partial L}{\partial\Omega}
- 
{\rm ad}^*_\Omega\frac{\delta L}
{\delta\Omega}\Big]\cdot\Sigma 
+
\Big[\frac{d}{dt}\frac{\partial L}{\partial\Lambda}
- 
{\rm ad}^*_\Lambda\frac{\partial L}{\partial\Lambda}\Big]\cdot\Xi
\\&&
+
\Big[
\frac{d}{dt}\frac{\partial L}{\partial \dot{D}}
-\frac{\partial L}{\partial D}
\Big]\delta D
\,,
\nonumber
\end{eqnarray}
where, the operation ${\rm ad}^*_\Omega$, for example, is defined
by
\begin{equation}
{\rm ad}^*_\Omega
\frac{\partial L}{\partial\Omega}\cdot\Sigma
=
-\,\frac{\partial L}{\partial\Omega}\cdot{\rm ad}_\Omega\Sigma
=
-\,\frac{\partial L}{\partial\Omega}\cdot[\Omega,\Sigma\,]
\,,
\end{equation}
and the dot `$\cdot$' denotes pairing between the Lie algebra and
its dual. This could also have been written in the notation using 
$\langle\cdot\,,\,\cdot\rangle:\,
\mathfrak{g}^*\times\mathfrak{g}\to\mathbb{R}$ as,
\begin{equation}
\left\langle
{\rm ad}^*_\Omega
\frac{\partial L}{\partial\Omega}\,,\,\Sigma
\right\rangle
=
-\,\left\langle
\frac{\partial L}{\partial\Omega}\,,\,{\rm ad}_\Omega\Sigma
\right\rangle
=
-\,\left\langle
\frac{\partial L}{\partial\Omega}\,,\,[\Omega,\Sigma\,]
\right\rangle
.
\end{equation}
The Euler-Poincar\'e dynamics is given by the stationarity
conditions for Hamilton's principle,
\begin{eqnarray}
\Sigma:&&
\frac{d}{dt}\frac{\partial L}{\partial\Omega}
- 
{\rm ad}^*_\Omega\frac{\partial L}{\partial\Omega} = 0
\,,\label{EP-eqn1}\\
\Xi:&&
\frac{d}{dt}\frac{\partial L}{\partial\Lambda}
- 
{\rm ad}^*_\Lambda\frac{\partial L}{\partial\Lambda} = 0
\,,\label{EP-eqn2}\\
\delta D:&&
\frac{d}{dt}\frac{\partial L}{\partial \dot{D}}
-\frac{\partial L}{\partial D}
= 0\,.
\label{EP-eqn3}
\end{eqnarray}
These are the {\bfi Euler-Poincar\'e equations} for the ellipsoidal
motions generated by Lagrangians of the form given in equation
(\ref{Lag-RE}). For example, such Lagrangians determine the
dynamics of the Riemann ellipsoids -- circulating, rotating,
self-gravitating fluid flows at constant density within an
ellipsoidal boundary. 

\subsection{Left and right momentum maps: Angular momentum versus
circulation}

The Euler-Poincar\'e equations (\ref{EP-eqn1}-\ref{EP-eqn3}) involve {\bfi
angular momenta} defined in terms of the angular velocities $\Omega$,
$\Lambda$ and the shape $D$ by
\begin{eqnarray}
M
&=&
\frac{\partial T}{\partial \Omega}
=
-\Omega D^2-D^2\Omega +2D\Lambda D
\label{leg1}
\,,\\
N
&=&
\frac{\partial T}{\partial \Lambda}
=
-\Lambda D^2-D^2\Lambda +2D\Omega D\,.
\label{leg2}
\end{eqnarray}
These angular momenta are related to the original deformation
gradient $Q=R^{-1}DS$ in equation (\ref{def-grad}) by the two momentum
maps from Example \ref{mommap-ex}
\begin{eqnarray}
PQ^T-QP^T=
\dot{Q}Q^{\,T}- Q\dot{Q}^{\,T}
&=&
R^{-1}MR
\label{leg1-F}
\,,\\
P^TQ-Q^TP=
\dot{Q}^{\,T}Q - Q^{\,T}\dot{Q}
&=&
S^{-1}NS
\,.
\label{leg2-F}
\end{eqnarray}
To see that $N$ is related to the {\bfi vorticity}, we consider the
exterior derivative of the circulation one-form
$\mathbf{u}\cdot d\mathbf{x}$ defined as
\begin{equation}
d(\mathbf{u}\cdot d\mathbf{x})
=
{\rm curl}\,\mathbf{u}\cdot d\mathbf{S}
=
\frac{1}{2}\,(\dot{Q}^{\,T}Q - Q^{\,T}\dot{Q})_{jk}\,
dx_0^j\wedge dx_0^k
=
(S^{-1}NS)_{jk}
dx_0^j\wedge dx_0^k
\,.
\label{circ}
\end{equation}
Thus, $S^{-1}NS$ is the fluid vorticity, referred to the Lagrangian
coordinate frame. For Euler's fluid equations, Kelvin's
circulation theorem implies $(S^{-1}NS)\dot{\,}=0$. 

Likewise, $M$ is related to the {\bfi angular momentum} by considering
\begin{equation}
u_ix_j-u_jx_i
=
\dot{Q}_{ik}x_0^kx_0^l Q^{\,T}_{lj}
- 
Q_{ik}x_0^kx_0^l \dot{Q}^{\,T}_{lj}
\,.
\label{ang-mom1}
\end{equation}
For spherical symmetry, we may choose $x_0^kx_0^l = \delta^{kl}$
and, in this case, the previous expression becomes
\begin{equation}
u_ix_j-u_jx_i
=
[\dot{Q}Q^{\,T}
- 
Q\dot{Q}^{\,T}]_{ij}
=
[R^{-1}MR]_{ij}
\,.
\label{ang-mom2}
\end{equation}
Thus, $R^{-1}MR$ is the angular momentum of the motion, referred to
the Lagrangian coordinate frame for spherical symmetry. In this
case, the angular momentum is conserved, so that
$(R^{-1}MR)\dot{\,}=0$.

In terms of these angular momenta, the Euler-Poincar\'e-Lagrange
equations (\ref{EP-eqn1}-\ref{EP-eqn3}) are expressed as
\begin{eqnarray}
\dot{M}&=&[\Omega, M]
\,,\label{epl1}\\
\dot{N}&=&[\Lambda ,N]
\,,\label{epl2}
\\
\frac{d}{dt}\left(\frac{\partial L}{\partial \dot{D}}\right)
&=&
\frac{\partial L}{\partial D}
\,.
\label{epl3}
\end{eqnarray}
Perhaps not unexpectedly, because of the combined symmetries of the
kinetic-energy Lagrangian (\ref{Lag-RE}) under both left and right
actions of $SO(n)$, the first two equations are consistent with the
conservation laws,
\[(R^{-1}MR)\dot{\,} = 0
\quad\hbox{and}\quad
(S^{-1}NS)\dot{\,} = 0
\,,
\]
respectively.
Thus, equation (\ref{epl1}) is the angular momentum equation while
(\ref{epl2}) is the vorticity equation. (Fluids have both types of
circulatory motions.) The remaining equation (\ref{epl3}) for the diagonal
matrix $D$ determines the shape of the ellipsoid undergoing free motion
on $GL(n,\mathbb{R})$.

\subsection{Vector representation of free Riemann ellipsoids in 3D}

In three dimensions these expressions may be written in vector form
by using the {\bfi hat map}, written now using upper and lower case Greek
letters as,
\[
\Omega_{ij}=\epsilon_{ijk}\omega_k
\,,\quad
\Lambda_{ij}=\epsilon_{ijk}\lambda_k
\,,
\]
with $\epsilon_{123}=1$, and $D={\rm diag}\,\{d_1,d_2,d_3\}$. 
\begin{exercise}
What is the analog of the hat map in four dimensions? Hint: locally the
Lie algebra $so(4)$ is isomorphic to $so(3)\times so(3)$.
\end{exercise}
Hence, the angular-motion terms in the kinetic energy may be rewritten as
\begin{equation}
-\,\frac{1}{2}{\rm trace}\,( \Omega^2D^2 )
=
\frac{1}{2}\Big[(d_1^2+d_2^2)\omega_3^2
+ (d_2^2+d_3^2)\omega_1^2
+ (d_3^2+d_1^2)\omega_2^2\Big]
\,,
\label{ke-rb-vec}
\end{equation}
\begin{equation}
-\,\frac{1}{2}{\rm trace}\,( \Lambda^2D^2 )
=
 \frac{1}{2}\Big[(d_1^2+d_2^2)\lambda_3^2
+ (d_2^2+d_3^2)\lambda_1^2
+ (d_3^2+d_1^2)\lambda_2^2\Big]
\,,
\end{equation}
and
\begin{equation}
-\,\frac{1}{2}{\rm trace}\,( \Omega D\Lambda D )
=
\Big[d_1d_2(\omega_3\lambda_3)
+ d_2d_3(\omega_1\lambda_1)
+ d_3d_1(\omega_2\lambda_2)\Big]
\,.
\end{equation}
On comparing equations (\ref{ke-rb}) and (\ref{ke-rb-vec}) for the
kinetic energy of the rigid body part of the motion, we identify
the usual moments of inertia as 
\[
I_k=d_i^2+d_j^2\,,
\quad\hbox{with}\quad
i,\,j,\,k\,
\hbox{ cyclic.}
\] 
The antisymmetric matrices $M$ and $N$ have vector representations
in 3D given by
\begin{eqnarray}
M_k
&=&
\frac{\partial T}{\partial \omega_k} 
=
(d_i^2+d_j^2)\omega_k
- 2d_id_j\lambda_k
\,,\label{leg1-vec}
\\
N_k
&=&
\frac{\partial T}{\partial \lambda_k}
=
(d_i^2+d_j^2)\lambda_k
- 2d_id_j\omega_k
\,,
\label{leg2-vec}
\end{eqnarray}
again with $i$, $j$, $k$ cyclic permutations of $\{1,2,3\}$.

\paragraph{Vector representation in 3D}
In terms of their 3D vector representations of the angular momenta in
equations (\ref{leg1-vec}) and (\ref{leg2-vec}), the two equations
(\ref{epl1}) and (\ref{epl2}) become
\begin{equation}
\dot{\mathbf{M}}
=
(\dot{R}R^{-1})\mathbf{M}
=
\Omega\mathbf{M}
=
 \boldsymbol\omega\times\mathbf{M}
\,,\quad
\dot{\mathbf{N}}
=
(\dot{S}S^{-1})\mathbf{N}
=
\Lambda\mathbf{N}
=
 \boldsymbol\lambda\times\mathbf{N}
\,.
\end{equation}
Relative to the Lagrangian fluid frame of reference, these
equations become
\begin{eqnarray}
(R^{-1}\mathbf{M})\dot{\,}
&=&
R^{-1}(\dot{\mathbf{M}}
-
 \boldsymbol\omega\times\mathbf{M}) = 0
\,,\\
(S^{-1}\mathbf{N})\dot{\,}
&=&
S^{-1}(\dot{\mathbf{N}}
-
 \boldsymbol\lambda\times\mathbf{N}) = 0
\,.
\end{eqnarray}
So each of these degrees of freedom represents a rotating, deforming body,
whose ellipsoidal shape is governed by the Euler-Lagrange equations
(\ref{epl3}) for the lengths of its three principal axes.

\begin{exercise}[Elliptical motions with potential energy on
$GL(2,\mathbb{R})$]
Compute equations (\ref{epl1}-\ref{epl3}) for elliptical motion in the
plane. Find what potentials $V(D)$ are solvable for
$L=T(\Omega,\Lambda,D,\dot{D})-V(D)$ by reducing these equations to the 
separated Newtonian forms,
\[
\frac{d^2r^2}{dt^2}=-\,\frac{dV(r)}{dr^2}
\,,\qquad
\frac{d^2\alpha}{dt^2}=-\,\frac{dW(\alpha)}{d\alpha}
\,,
\]
for $r^2=d_1^2+d_2^2$ and $\alpha=\tan^{-1}(d_2/d_1)$ with $d_1(t)$ and
$d_2(t)$ in two dimensions. Hint: consider the potential energy,
\[
V(D)=V\Big({\rm tr}D^2,\det(D)\Big)
\,,
\] 
for which the equations become homogeneous in $r^2(t)$.
\end{exercise}

\begin{exercise}[Ellipsoidal motions with potential energy on
$GL(3,\mathbb{R})$] Choose the Lagrangian in 3D,
\[
L = \frac{1}{2}{\rm tr} \Big(\dot{Q}^T \dot{Q}\Big)
- V\Big({\rm tr}\,(Q^TQ),\det(Q)\Big)
\,,
\]
where $Q(t)\in GL(3,\mathbb{R})$ is a $3\times3$ matrix function of time
and the potential energy $V$ is an arbitrary function of 
${\rm tr}\,(Q^TQ)$ and $\det(Q)$.
\begin{enumerate}
\item
Legendre transform this Lagrangian. That is, find the 
momenta $P_{ij}$ canonically conjugate to $Q_{ij}$, construct the
Hamiltonian $H(Q,P)$ and write Hamilton's canonical equations
of motion for this problem.
\item
Show that the Hamiltonian is invariant under $Q\to OQ$ where $O\in SO(3)$.
Construct the cotangent lift of this action on $P$. Hence, construct the
momentum map of this action. 
\item
Construct another distinct action of $SO(3)$ on this system which also
leaves its Hamiltonian $H(Q,P)$ invariant. Construct its
momentum map. Do the two momentum maps Poisson commute? Why? 
\item
How are these two momentum maps related to the angular momentum and
circulation in equations (\ref{leg1}) and (\ref{leg2})?
\item
How does the 2D restriction of this problem inform the previous one?
\end{enumerate}
\end{exercise}

\begin{exercise}[$GL(n,\mathbb{R})-$invariant motions]\label{GLn-ex}
Begin with the Lagrangian
\[
L=\frac{1}{2}\,{\rm tr}\Big(\dot{S}S^{-1}\dot{S}S^{-1}\Big)
+
\frac{1}{2}\,\mathbf{\dot{q}}^T S^{-1}\mathbf{\dot{q}}
\]
where $S$ is an $n\times n$ symmetric matrix and
$\mathbf{q}\in\mathbb{R}^n$ is an
$n-$component column vector.
\begin{enumerate}
\item
Legendre transform to construct the corresponding Hamiltonian and canonical
equations.
\item
Show that the system is invariant under the group action 
\[
\mathbf{q}\to A\mathbf{q}
\quad\hbox{and}\quad
S\to ASA^T
\]
for any constant invertible $n\times n$ matrix, $A$.
\item
Compute the infinitesimal generator for this group action and construct
its corresponding momentum map. Is this momentum map equivariant?
\item
Verify directly that this momentum map is a conserved $n\times n$
matrix quantity by using the equations of motion.
\item
Is this system completely integrable for any value of $n>2$?
\end{enumerate}
\end{exercise}

\section{Heavy top equations}
\subsection{Introduction and definitions}
A top is a rigid body of mass $m$ rotating with a fixed point of support
in a constant gravitational field of acceleration $-g\mathbf{\hat{z}}$
pointing vertically downward. The orientation of  the body relative to the
vertical axis $\mathbf{\hat{z}}$ is defined by the unit vector
$\boldsymbol{\Gamma}=\mathbf{ R} ^{-1}(t)\mathbf{\hat{z}}$ for a curve
$\mathbf{R}(t)\in SO(3)$. According to its definition, the unit vector
$\boldsymbol{\Gamma}$ represents the motion of the vertical direction as
seen from the rotating body. Consequently, it satisfies the auxiliary
motion equation,
\[
\boldsymbol{\dot{\Gamma}}
=
-\,\mathbf{ R} ^{-1} \mathbf{\dot{R}}(t)\boldsymbol{\Gamma}
=
\boldsymbol{\Gamma}\times\boldsymbol{\Omega}
\,.
\]
Here the rotation matrix $\mathbf{R}(t)\in SO(3)$, the skew matrix
$\boldsymbol{\hat{\Omega}}=\mathbf{ R} ^{-1} \mathbf{\dot{R}}\in so(3)$
and the body angular frequency vector $\boldsymbol{\Omega}\in\mathbb{R}^3$
are related by the hat map,  $\boldsymbol{\Omega} =\big(\mathbf{ R}
^{-1}\mathbf{\dot{R}}\big)\boldsymbol{\hat{\,}}$, where 
$\mathbf{\hat{\,}}\,:\,({so}(3), [\cdot, \cdot])
\to (\mathbb{R}^3, \times)$ with 
$\boldsymbol{\hat{\Omega}}\mathbf{v} = \boldsymbol{\Omega} \times
\mathbf{v}$ for any $\mathbf{v} \in \mathbb{R}^3$.

The motion of a top is determined from Euler's equations in vector form,
\begin{eqnarray}\label{top-eqns-vector}
{\mathbb I} \boldsymbol{\dot{\Omega}}
&=&{\mathbb I} \boldsymbol{\Omega} \times \boldsymbol{\Omega}
+ mg\, \boldsymbol{\Gamma}\times\boldsymbol{\chi}
\,,\\
\boldsymbol{\dot{\Gamma}}
&=&
\boldsymbol{\Gamma}\times\boldsymbol{\Omega}
\,,
\label{chi-eqns-aux}
\end{eqnarray}
where $\boldsymbol{\Omega},\,
\boldsymbol{\Gamma},\,\boldsymbol{\chi}\in\mathbb{R}^3$ are vectors in the
rotating body frame. Here 
\begin{itemize}
\item
$\boldsymbol{\Omega} = (\Omega_1,\Omega_2, \Omega_3)$ is the body angular
velocity vector.
\item
${\mathbb I}={\rm diag}(I_1,  I_2, I_3)$  is the
moment of inertia tensor, diagonalized in the body principle axes. 
\item
$\boldsymbol{\Gamma}=R^{-1}(t)\mathbf{\hat{z}}$ represents the motion of
the unit vector along the vertical axis, as seen from the body. 
\item
$\boldsymbol{\chi}$ is the constant vector in the body from the point of
support to the body's center of mass. 
\item
$m$ is the total mass of the body
and $g$ is the constant acceleration of gravity.
\end{itemize}

\subsection{Heavy top action principle}

\begin{proposition}\label{ht-actprinc}
The heavy top equations are equivalent to the {\bfi heavy top
action principle} for a {\bfi reduced action}
\begin{equation}\label{rbvp1}
\delta S_{\rm red}=0
\,,\quad\hbox{with}\quad
S_{\rm red}= \int^b_a l(\boldsymbol{\Omega},\boldsymbol{\Gamma}) \, d t
= \int^b_a 
\frac{1}{2} \langle\,  {\mathbb I} \boldsymbol{\Omega}
\,,\,
\boldsymbol{\Omega}\rangle 
-
\langle  mg\,\boldsymbol{\chi}
\,,\,
\boldsymbol{\Gamma}\rangle 
\, d t,
\end{equation}
where variations of $\boldsymbol{\Omega}$ and $\boldsymbol{\Gamma}$ are
restricted to be of the form
\begin{equation}\label{rpvp2}
\delta \boldsymbol{\Omega} = \boldsymbol{\dot{ \Sigma}} +
\boldsymbol{\Omega} \times \boldsymbol{\Sigma}
\quad\hbox{and}\quad
\delta \boldsymbol{\Gamma} = 
\boldsymbol{\Gamma} \times \boldsymbol{\Sigma}
\,,
\end{equation}
arising from variations of the definitions $\boldsymbol{\Omega}
=\big(\mathbf{ R} ^{-1}\mathbf{\dot{R}}\big)\boldsymbol{\hat{\,}}$ and
$\boldsymbol{\Gamma}=\mathbf{ R} ^{-1}(t)\mathbf{\hat{z}}$  in which
$\boldsymbol{\Sigma}(t) =\big(\mathbf{R}^{-1}
\delta\mathbf{R}\big)\boldsymbol{\hat{\,}}$ is a curve in
$\mathbb{R}^3$ that vanishes at the endpoints in time. 
\end{proposition}

\begin{proof}
Since ${\mathbb I}$ is symmetric and $\boldsymbol{\chi}$ is constant, we
obtain the variation,
\begin{align*}
\delta \int^b_a l(\boldsymbol{\Omega},\boldsymbol{\Gamma}) \, d t
& =  \int^b_a 
\langle\,  {\mathbb I} \boldsymbol{\Omega}
\,,\,
\delta\boldsymbol{\Omega}\rangle 
-
\langle  mg\,\boldsymbol{\chi}
\,,\,
\delta\boldsymbol{\Gamma}\rangle \, dt \\
& =  \int^b_a 
\langle\,  {\mathbb I} \boldsymbol{\Omega}
\,,\,
\boldsymbol{\dot{ \Sigma}}
      + \boldsymbol{\Omega} \times \boldsymbol{\Sigma}\rangle
-
\langle  mg\,\boldsymbol{\chi}\,,\,
\boldsymbol{\Gamma} \times \boldsymbol{\Sigma}
\rangle 
\, dt\\
& =  \int^b_a  \left\langle
      -\, \frac{d }{d t}{\mathbb I} \boldsymbol{\Omega}
\,,\,
      \boldsymbol{\Sigma}\right\rangle + \left\langle\, {\mathbb I}
\boldsymbol{\Omega}
\,,\,
      \boldsymbol{\Omega} \times
\boldsymbol{\Sigma}\right\rangle 
-
\langle  mg\,\boldsymbol{\chi}\,,\,
\boldsymbol{\Gamma} \times \boldsymbol{\Sigma}
\rangle 
{d t}\\
& = \int^b_a \left\langle
-\, \frac{d }{d t} {\mathbb I} \boldsymbol{\Omega}
      + {\mathbb I} \boldsymbol{\Omega}\times \boldsymbol{\Omega}
+
mg\, \boldsymbol{\Gamma}\times\boldsymbol{\chi}
\,,\,
\boldsymbol{\Sigma} \right\rangle d t,
\end{align*}
upon integrating by parts and using the
endpoint conditions, $ \boldsymbol{\Sigma} (b) =
\boldsymbol{\Sigma} (a) = 0 $. Since $\boldsymbol{\Sigma}$ is
otherwise arbitrary,~(\ref{rbvp1}) is equivalent to
\[-\, \frac{d }{d t} {\mathbb I} \boldsymbol{\Omega}
      + {\mathbb I} \boldsymbol{\Omega}\times \boldsymbol{\Omega}
+
mg\, \boldsymbol{\Gamma}\times\boldsymbol{\chi}
= 0 , \]
which is Euler's motion equation for the heavy top
(\ref{top-eqns-vector}). This motion equation is completed by the
auxiliary equation
$\boldsymbol{\dot{\Gamma}}=\boldsymbol{\Gamma}\times\boldsymbol{\Omega}$
in (\ref{chi-eqns-aux}) arising
from the definition of
$\boldsymbol{\Gamma}$. 
\end{proof}\bigskip

The Legendre transformation for
$l(\boldsymbol{\Omega},\boldsymbol{\Gamma})$ gives the body angular
momentum
\begin{eqnarray*}
\label{ht Legendre}
\boldsymbol{\Pi} 
= \frac{\partial l}{\partial\boldsymbol{\Omega}}
= {\mathbb I} \boldsymbol{\Omega}
\,.
\end{eqnarray*}
The well known energy Hamiltonian for the heavy top then emerges as 
\begin{eqnarray}
\label{ht Hamiltonian}
h(\boldsymbol{\Pi},\boldsymbol{\Gamma})
=
\boldsymbol{\Pi}\cdot\boldsymbol{\Omega}
-
l(\boldsymbol{\Omega},\boldsymbol{\Gamma})
=
\frac{1}{2}\langle\boldsymbol{\Pi}\,,\,
\mathbb{I}^{-1}\boldsymbol{\Pi}\rangle 
+
\langle  mg\,\boldsymbol{\chi}
\,,\,
\boldsymbol{\Gamma}\,\rangle 
\,,
\end{eqnarray}
which is the sum of the kinetic and potential energies of the top.

\paragraph{The Lie-Poisson Equations.}
Let
$f,h:\mathfrak{g} ^{\ast}\to\mathbb{R}$ be real-valued functions on the
dual space 
$\mathfrak{g} ^{\ast}$.  Denoting elements of $\mathfrak{g} ^{\ast}$ by
$\mu$, the functional derivative of $f$   at $\mu$ is defined as the
unique element $\delta  f/ \delta \mu$ of $\mathfrak{g}$ defined by
\begin{equation}\label{fd}
\lim_{\varepsilon \rightarrow 0}
\frac{1}{\varepsilon }[f(\mu
+ \varepsilon \delta \mu) - f(\mu )]
=  \left\langle \delta  \mu , \frac{
\delta  f}{\delta  \mu } \right\rangle ,
\end{equation}
for all $\delta\mu \in \mathfrak{g}^*$, where 
$\left\langle\cdot \,,\,\cdot\right\rangle$ denotes the pairing between
$\mathfrak{g} ^{\ast} $ and  $\mathfrak{g} $. \\

\begin{definition}[Lie-Poisson brackets \& Lie-Poisson equations]$\quad$\\
The {\bfi $(\pm)$  Lie-Poisson brackets} are defined by
\end{definition}

\begin{equation}\label{lpb}
\{f, h\}_\pm (\mu )  
=  \pm \left\langle \mu , \left[ \frac{ \delta f}{\delta  \mu},
\frac{\delta  h}{\delta \mu } \right] \right\rangle 
=  \mp \left\langle \mu , 
\operatorname{ad}_{\delta h/\delta \mu}
\frac{ \delta f}{\delta  \mu}\right\rangle
\,.
\end{equation}
The corresponding {\bfi Lie-Poisson equations}, determined by $\dot f =
\{ f,h\}
$ read
\begin{equation}\label{alpe}
\dot \mu  = \{ \mu,h\}
= \mp \operatorname{ad} ^{* }_{\delta h/\delta \mu} \mu
\,,
\end{equation}
where one defines the ad$^*$ operation in terms of the pairing
$\left\langle\cdot \,,\,\cdot\right\rangle$, by
\[
\{ f,h\}
=
\left\langle \mu , 
\operatorname{ad}_{\delta h/\delta \mu}
\frac{ \delta f}{\delta  \mu}\right\rangle
=
\left\langle \operatorname{ad}^*_{\delta h/\delta \mu}\mu , 
\frac{ \delta f}{\delta  \mu}\right\rangle
\,.
\]
The Lie-Poisson setting of mechanics is a special case of the general
theory of systems on Poisson manifolds, for which there is now
an extensive theoretical development. (See Marsden and Ratiu [2003] for a
start on this literature.) 

\subsection{Lie-Poisson brackets and momentum maps.}  An important 
feature of the rigid body bracket carries  over to general Lie algebras.
Namely, {\em Lie-Poisson brackets on $\mathfrak{g}^*$ arise from canonical
brackets on the cotangent bundle\/} (phase space) $T^\ast G$ associated
with a Lie group $G$ which has $\mathfrak{g}$ as its associated Lie
algebra. Thus, the process by which the Lie-Poisson brackets arise is the
momentum map
\[
T^\ast G \mapsto \mathfrak{g}^\ast
\,.
\]

For example, a rigid body is free to rotate about its center of
mass and $G$ is the (proper) rotation group ${\rm SO}(3)$. The choice of
$T^\ast G$ as the primitive phase space is made according to
the classical procedures of mechanics described earlier.  For
the description using Lagrangian mechanics, one forms the
velocity phase space $TG$.  The Hamiltonian description on 
$T^\ast G$ is then obtained
by standard procedures: Legendre transforms, etc.

The passage from $T^\ast G$ to the space of ${\boldsymbol
\Pi}$'s (body angular momentum space) is determined by {\em
left\/} translation on the group.  This mapping is  an example
of a {\em momentum map\/}; that is, a
mapping whose components are the ``Noether  quantities''
associated with a symmetry group.  The map from $T^\ast G$ to
$\mathfrak{g}^\ast$ being a Poisson map {\em is a
general fact about momentum maps\/}.  The Hamiltonian point of
view of all this is a standard subject.

\begin{remark}[Lie-Poisson description of the heavy top]

As it turns out, the underlying Lie algebra for the Lie-Poisson
description of the heavy top consists of the Lie algebra $se(3,\mathbb{R})$
of infinitesimal Euclidean motions in $\mathbb{R}^3$. This is a bit
surprising, because heavy top motion itself does {\em not\/} actually
arise through actions of the Euclidean group of rotations and translations
on the body, since the body has a fixed point! Instead, the Lie algebra
$se(3,\mathbb{R})$ arises for another reason associated with the 
{\bfi breaking} of the ${\rm SO}(3)$ isotropy by the presence of the
gravitational field. This symmetry breaking introduces a semidirect-product
Lie-Poisson structure which happens to coincide with the dual of the Lie algebra
$se(3,\mathbb{R})$ in the case of the heavy top. As we shall see later, a close
parallel exists between this case and the Lie-Poisson structure for compressible
fluids.
\end{remark}

\subsection{The heavy top Lie-Poisson brackets}
The {\bfi Lie algebra of the special Euclidean group in 3D} is
${se}(3) = \mathbb{R}^3  \times \mathbb{R}^3$ with the  Lie bracket
\begin{equation}\label{htla} [(\boldsymbol{\xi}, \mathbf{u} ),
(\boldsymbol{\eta}, \mathbf{ v} )]  =  (\boldsymbol{\xi}
\times \boldsymbol{\eta},
\boldsymbol{\xi} \times \mathbf{ v} - \boldsymbol{\eta}
\times
\mathbf{u})  
\,.
\end{equation}
%
We identify the dual space with pairs $({\boldsymbol
\Pi},\boldsymbol{\Gamma} )$; the corresponding
$  (-)  $ Lie-Poisson bracket called the {\bfi heavy top
bracket\/} is
\begin{eqnarray}\label{ht-LPB}
\{f\,,\,h\}(\boldsymbol{\Pi},\boldsymbol{\Gamma})
=
-\,\boldsymbol{\Pi}\cdot\nabla_{\Pi}f\times\nabla_{\Pi}h
\,-\,\boldsymbol{\Gamma}\cdot\big(\nabla_{\Pi}f\times\nabla_{\Gamma}h
-\nabla_{\Pi}h\times\nabla_{\Gamma}f\big)
.
\end{eqnarray}
This Lie-Poisson bracket and the Hamiltonian (\ref{ht
Hamiltonian}) recover the equations (\ref{top-eqns-vector}) and
(\ref{chi-eqns-aux}) for the heavy top, as
\begin{eqnarray*}
\boldsymbol{\dot{\Pi}}
=
\{\boldsymbol{\Pi}\,,\,h\}
&=&
\boldsymbol{\Pi}\times\nabla_{\Pi}h
+\boldsymbol{\Gamma}\times\nabla_{\Gamma}h
=
\boldsymbol{\Pi}\times\mathbb{I}^{-1}\boldsymbol{\Pi}
+\boldsymbol{\Gamma}\times mg\,\boldsymbol{\chi}
\,,\\
\boldsymbol{\dot{\Gamma}}
=
\{\boldsymbol{\Gamma}\,,\,h\}
&=&
\boldsymbol{\Gamma}\times\nabla_{\Pi}h
=
\boldsymbol{\Gamma}\times\mathbb{I}^{-1}\boldsymbol{\Pi}
\,.
\end{eqnarray*}
\begin{remark}[Semidirect products and symmetry breaking]
The Lie algebra of the Euclidean group has a structure which is
a special case of what is called a {\bfi semidirect product\/}.
Here, it is the semidirect product action $so(3)\,\circledS\,\mathbb{R}^3$ of the Lie
algebra of rotations $so(3)$ acting on the infinitesimal translations $\mathbb{R}^3$,
which happens to coincide with $se(3,\mathbb{R})$.  In general, the Lie bracket for
semidirect product action $\mathfrak{g}\,\circledS\,V$ of a Lie algebra $\mathfrak{g}$
on a vector space $V$ is given by 
\[
\Big[(X,a), (\overline{X},\overline{a})\Big]
=
\Big([X,\overline{X}\,],\overline{X}(a)-X(\overline{a})\Big)
\]
in which $X,\overline{X}\in\mathfrak{g}$ and $a,\overline{a}\in V$.
Here, the action of the Lie algebra on the vector space is denoted, e.g.,
$X(\overline{a})$. Usually, this action would be the Lie derivative.  

Lie-Poisson brackets defined on the dual spaces of semidirect product Lie algebras
tend to occur under rather general circumstances when the symmetry in $T^\ast G $ is
broken, e.g., reduced to an isotropy subgroup of a set of parameters.  In
particular, there are similarities in structure between the Poisson bracket for
compressible flow and that for the heavy top. In the latter case, the vertical
direction of gravity breaks isotropy of $\mathbf{R}^3$ from
$SO(3)$ to
$SO(2)$. The general theory for semidirect products is reviewed in a
variety of places, including \cite{MaRaWe1984a,MaRaWe1984b}. Many
interesting examples of Lie-Poisson brackets on semidirect products
exist for fluid dynamics. These semidirect-product Lie-Poisson 
Hamiltonian theories range from
simple fluids, to charged fluid plasmas, to magnetized fluids, to
multiphase fluids, to super fluids, to Yang-Mills fluids, relativistic, or
not, and to liquid crystals. See, for example, \cite{GiHoKu1982},
\cite{HoKu1982}, \cite{HoKu1983}, \cite{HoKu1988}.
For discussions of many of these theories from the Euler-Poincar\'e viewpoint, see
\cite{HoMaRa1998a} and \cite{Ho2002a}.
\end{remark}

\subsection{The heavy top formulation by the Kaluza-Klein construction} 
The Lagrangian in the heavy top action principle (\ref{rbvp1}) may be transformed into
a quadratic form. This is accomplished by suspending the system in a higher dimensional
space via the {\bfi Kaluza-Klein construction\/}. This construction proceeds for the
heavy top as a slight modification of the well-known Kaluza-Klein construction for a
charged particle in a prescribed magnetic field.\\

Let $Q_{KK}$ be the manifold $SO(3)\times\mathbb{R}^3$ with
variables
$(\mathbf{R},\mathbf{q})$. On $Q_{KK} $ introduce the {\bfi Kaluza-Klein
Lagrangian} $L_{KK}:TQ_{KK}\simeq TSO(3)\times T\mathbb{R}^3\mapsto
\mathbb{R}$ as
\begin{equation}
L_{KK}(\mathbf{R},\mathbf{q}, \mathbf{\dot{R}},\mathbf{\dot{q}}
;\mathbf{\hat{z}}) 
=
L_{KK}(\boldsymbol{\Omega},\boldsymbol{\Gamma},\mathbf{q},\mathbf{\dot{q}}) 
=
\frac{1}{2} \langle\,  \mathbb{I} \boldsymbol{\Omega}
\,,\,
\boldsymbol{\Omega}\,\rangle 
+
\frac{1}{2}|\boldsymbol{\Gamma}+\mathbf{\dot{q}}|^2
\,,
\label{KK-metric}
\end{equation}
with $\boldsymbol{\Omega} =\big(\mathbf{ R}^{-1}
\mathbf{\dot{R}}\big)\boldsymbol{\hat{\,}}$
and $\boldsymbol{\Gamma}=\mathbf{ R}^{-1}\mathbf{\hat{z}}$.
The Lagrangian $L_{KK} $ is positive definite in
$(\boldsymbol{\Omega},\boldsymbol{\Gamma},\mathbf{\dot{q}})$;
so it may be regarded as the kinetic energy of a metric, the {\bfi
Kaluza-Klein metric\/} on $TQ_{KK} $.

The Legendre transformation for $L_{KK}$ gives the momenta
\begin{equation}
\label{KK Legendre}
\boldsymbol{\Pi} = {\mathbb I} \boldsymbol{\Omega}
\qquad \text{and} \qquad 
\mathbf{p} =
\boldsymbol{\Gamma}+\mathbf{\dot{q}}
\,.
\end{equation}
Since $L_{KK}$ does not depend on $\mathbf{q}$, the Euler-Lagrange
equation
\[
\frac{d}{dt}\frac{\partial L_{KK}}{\partial \mathbf{\dot{q}}}  
=
\frac{\partial L_{KK}}{\partial \mathbf{q}} = 0
\,,
\]
shows that $\mathbf{p} = \partial L_{KK}/\partial\mathbf{\dot{q}}$
is conserved. The {\bfi constant vector\/} $\mathbf{p}$ is now
identified as the vector in the body,
\[
\mathbf{p} = \boldsymbol{\Gamma}+\mathbf{\dot{q}}
=
-\,mg\,\boldsymbol{\chi}
\,.\]
After this identification, the heavy top action principle in Proposition
\ref{ht-actprinc} with the Kaluza-Klein Lagrangian returns Euler's
motion equation for the heavy top (\ref{top-eqns-vector}). 

The Hamiltonian
$H_{KK}$ associated to $L_{KK}$ by the Legendre transformation
\eqref{KK Legendre} is
\begin{eqnarray*}
H_{KK}(\boldsymbol{\Pi},\boldsymbol{\Gamma},\mathbf{q}, \mathbf{p}) 
&=& 
\boldsymbol{\Pi}\cdot\boldsymbol{\Omega}
 + 
\mathbf{p}\cdot\mathbf{\dot{q}}
- 
L_{KK}(\boldsymbol{\Omega},\boldsymbol{\Gamma},\mathbf{q},\mathbf{\dot{q}})  
\nonumber \\
&=& 
\frac{1}{2}\boldsymbol{\Pi}\cdot\mathbb{I}^{-1}\boldsymbol{\Pi}
-
\mathbf{p}\cdot\boldsymbol{\Gamma}
+
\frac{1}{2}|\mathbf{p}|^2
\nonumber \\
&=& 
\frac{1}{2}\boldsymbol{\Pi}\cdot\mathbb{I}^{-1}\boldsymbol{\Pi}
+
\frac{1}{2}|\mathbf{p}-\boldsymbol{\Gamma}|^2
-
\frac{1}{2}|\boldsymbol{\Gamma}|^2
\,.
\label{KK Hamiltonian}
\end{eqnarray*}
Recall that $\boldsymbol{\Gamma}$ is a unit vector. 
On the constant level set $|\boldsymbol{\Gamma}|^2=1$, the
Kaluza-Klein Hamiltonian $H_{KK}$ is a positive quadratic function,
shifted by a constant. Likewise, on the constant level set $\mathbf{p} =
-\,mg\,\boldsymbol{\chi}$, the Kaluza-Klein Hamiltonian $H_{KK}$ is a
function of only the variables
$(\boldsymbol{\Pi},\boldsymbol{\Gamma})$ and is equal to the Hamiltonian
\eqref{ht Hamiltonian} for the heavy top up to an additive constant.
Consequently, the Lie-Poisson equations for the Kaluza-Klein Hamiltonian
$H_{KK}$ now reproduce Euler's motion equation for the heavy top
(\ref{top-eqns-vector}).

\begin{exercise}
Write the Kaluza-Klein construction on $SE(3)=SO(3)\circledS\mathbb{R}^3$.
\end{exercise}


\section{Euler-Poincar\'e (EP) reduction theorem}

\begin{remark}[Geodesic motion] As emphasized by \cite{Ar1966}, in
many interesting cases, the Euler--Poincar\'e equations on the dual of
a Lie algebra $ \mathfrak{g}^* $ correspond to {\it geodesic motion}
on the corresponding group $G$. The relationship between the equations on
$\mathfrak{g}^*$ and on $G$ is the content of the basic
Euler-Poincar\'e theorem discussed later. Similarly, on the
Hamiltonian side, the preceding paragraphs described the relation between
the Hamiltonian equations on $ T ^\ast G $ and the Lie--Poisson equations
on $
\mathfrak{g}^\ast $.  The issue of geodesic motion is especially simple:
if either the Lagrangian on $\mathfrak{g}$ or the Hamiltonian on
$\mathfrak{g}^\ast$ is purely quadratic, then the corresponding motion on
the group is geodesic motion.
\end{remark}

\subsection{We were already speaking prose (EP)}
Many of our previous considerations may be recast immediately as Euler-Poincar\'{e} equations.
\begin{itemize}
\item
Rigid bodies $\simeq$
\Big(EP$SO(n)$\Big), 
\item
Deforming bodies $\simeq$\Big(EP$GL_{+}(n,\mathbb{R})$\Big),
\item
Heavy tops $\simeq$ \Big(EP$SO(3)\times \mathbb{R}^3$\Big), 
\item
EPDiff
\end{itemize}

\subsection{Euler-Poincar\'{e} Reduction}
\label{Elpor}

This lecture applies reduction by symmetry
to Hamilton's principle. For a $G-$invariant Lagrangian defined on $TG$,
this reduction takes Hamilton's principle from $TG$ to
$TG/G\simeq\mathfrak{g}$. Stationarity of the symmetry-reduced Hamilton's
principle yields the Euler-Poincar\'e equations on $\mathfrak{g}^*$. The corresponding
reduced Legendre transformation yields the Lie-Poisson Hamiltonian formulation of
these equations. 

\bigbreak

{\bfi Euler-Poincar\'{e} Reduction\/} starts with a {right}
(respectively, {left}) invariant Lagrangian $L:TG\rightarrow
\mathbb{R}$ on the tangent bundle of a Lie group $G $. This means that $L
(T_g R_h(v)) = L(v)$, respectively $L (T_g L_h(v)) = L(v)$, for all $g, h
\in G $ and all $v\in T_gG$. In shorter notation, right invariance of
the Lagrangian may be written as 
\[L(g(t) ,\dot g(t))=L(g(t)h ,\dot g(t)h)
\,,
\]
for all $h\in G$.

\begin{theorem}[Euler-Poincar\'{e} Reduction]
\label{EP theorem}
 Let $G$ be a Lie group, $L : TG \rightarrow \mathbb{R}$
a {right}-invariant Lagrangian, and $l: = L|_{\mathfrak{g}} : \mathfrak{g}
\rightarrow \mathbb{R}$ be its restriction to $\mathfrak{g}$. For a
curve $g(t)\in G$, let 
\[
\xi(t) =  \dot g(t)\cdot g(t)^{-1} := T_{g(t)}
R_{g(t)^{-1}} \dot g(t) \in \mathfrak{g}
\,.
\]

\noindent
Then the following four statements are equivalent:
\begin{enumerate}
\item[{\rm \textbf{(i)}}] $g(t)$ satisfies the Euler-Lagrange
equations for Lagrangian $L$ defined on $G$.
\item[{\rm \textbf{(ii)}}] The variational principle
$$
\delta \int_a^b L(g(t) ,\dot g(t)) dt = 0
$$
\noindent holds, for variations with fixed endpoints.
\item[{\rm \textbf{(iii)}}] The ({right} invariant) {\bfi
Euler-Poincar\'{e} equations\/} hold:
$$\frac{d}{dt} \frac{\delta l}{\delta \xi} 
=-\, {\rm ad}^*_\xi
\frac{\delta l}{\delta \xi}\, .
$$
\item[{\rm \textbf{(iv)}}] The variational principle
$$\delta \int_a^b l (\xi(t)) dt = 0
$$
\noindent holds on $\mathfrak{g}$, using variations of the form
 $\delta \xi = \dot \eta - [\xi , \eta ]$, where $\eta(t)$ is an
arbitrary path in $\mathfrak{g}$ which vanishes  at
 the endpoints, i.e., $\eta (a) = \eta (b) = 0$.
\end{enumerate}
\end{theorem}

\begin{proof}  

\paragraph{\bfi Step I. Proof that $\textbf{(i)}
\Longleftrightarrow \textbf{(ii)}$:} This is Hamilton's principle: the
Euler-Lagrange equations follow from stationary action for variations
$\delta g$ which vanish at the endpoints. (See Lecture 10.) 

\paragraph{\bfi Step II. Proof that $\textbf{(ii)}
\Longleftrightarrow\textbf{(iv)}$:}
Proving equivalence of the variational principles $\textbf{(ii)}$ on $TG$
and $\textbf{(iv)}$ on $\mathfrak{g}$ for a {right}-invariant Lagrangian
requires calculation of the variations $\delta\xi$ of
$\xi = \dot{g}g^{-1}$ induced by $\delta g$. To simplify the exposition,
the calculation will be done first for matrix Lie groups, then generalized
to arbitrary Lie groups.

\paragraph{Step IIA. Proof that $\textbf{(ii)}
\Longleftrightarrow\textbf{(iv)}$ for a matrix Lie group.} For $\xi =
\dot{g}g^{-1}$, define $g_\epsilon(t)$ to be a family of curves in $G$
such that $g_0(t) = g(t)$ and  denote 
\[\delta g := \frac{dg_{\epsilon}(t)}{d\epsilon}\Big|_{\epsilon =0}
\,.
\] 
The variation of $\xi$ is computed in terms of $\delta g$ as
\begin{equation}\label{ep1}
\delta \xi =\left.\frac{d}{d\epsilon}\right|_{\epsilon = 0}
(\dot g_\epsilon g_\epsilon^{-1}) 
= 
\left.\frac{d^2 g}{dt d\epsilon}\right|_{\epsilon = 0}g^{-1}
-\dot{g} g^{-1} (\delta g) g^{-1} 
\,.
\end{equation}
Set $\eta := g^{-1} \delta g$. That is, $\eta(t)$ is an arbitrary
curve in $\mathfrak{g}$ which vanishes at the endpoints. The time
derivative of $\eta$ is computed as
\begin{equation}\label{ep2}
\dot\eta
=
\frac{d \eta}{dt} 
=
\frac{d}{d t}\left(\Big(
\left.\frac{d}{d\epsilon}\right|_{\epsilon = 0} g_\epsilon 
\Big)g^{-1}\right) 
=
\left.\frac{d^2 g}{dt d\epsilon}\right|_{\epsilon = 0}g^{-1}
- (\delta g) g^{-1} \dot{g} g^{-1} 
\,.
\end{equation}
Taking the difference of (\ref{ep1}) and (\ref{ep2}) implies
$$
\delta \xi - \dot\eta 
= 
-\dot{g} g^{-1} (\delta g) g^{-1}
+ (\delta g) g^{-1} \dot{g} g^{-1}  
= -\,\xi \eta + \eta \xi = -\,[\xi, \eta ]
\,.
$$
That is, for matrix Lie algebras,
\[\delta \xi = \dot \eta - [\xi ,\eta ]\,,\]
where $[\xi ,\eta ]$ is the matrix commutator.
Next, we notice that {right} invariance of $L$ allows one to change
variables in the Lagrangian by applying $g^{-1}(t)$ from the {right}, as
\[
L(g(t) ,\dot g(t))
=
L(e , \dot g(t)g^{-1}(t))
=: 
l (\xi(t))
\,.
\]
Combining this definition of the symmetry-reduced Lagrangian
$l:\mathfrak{g}\to\mathbb{R}$ together with the formula for variations 
$\delta \xi$ just deduced proves the equivalence of \textbf{(ii)} and
\textbf{(iv)} for matrix Lie groups.

\paragraph{Step IIB. Proof that $\textbf{(ii)}
\Longleftrightarrow\textbf{(iv)}$ for an arbitrary Lie group.} The same
proof extends to any Lie group $G$ by using the following lemma.
\begin{lemma}
\label{mixed partials lemma}
Let $g : U \subset \mathbb{R}^2 \rightarrow G $ be a smooth map and
denote its partial derivatives by
\begin{equation}
\label{partial derivatives of g}
\xi(t, \varepsilon): = T_{g(t, \varepsilon)} 
R_{g(t,\varepsilon)^{-1}}
\frac{\partial g(t,\varepsilon)}{\partial t}\,, \qquad
\eta(t, \varepsilon): = T_{g(t, \varepsilon)} 
R_{g(t,\varepsilon)^{-1}}
\frac{\partial g(t,\varepsilon)}{\partial \varepsilon}\,.
\end{equation}
Then
\begin{equation}
\label{mixed partials}
\frac{\partial \xi}{\partial \varepsilon} - \frac{\partial
\eta}{\partial t} = -\,[ \xi, \eta]\,,
\end{equation}
where $[ \xi, \eta]$ is the Lie algebra bracket on $\mathfrak{g}$. 
Conversely, if $U \subset \mathbb{R}^2 $ is simply connected and
$\xi, \eta:U \rightarrow \mathfrak{g}$ are smooth functions
satisfying \eqref{mixed partials}, then there exists a smooth
function $g : U \rightarrow G $ such that  \eqref{partial
derivatives of g} holds.
\end{lemma}

\begin{proof}
Write $\xi=\dot{g}g^{-1}$ and $\eta=g^\prime g^{-1}$ in natural notation
and express the partial derivatives $\dot{g}=\partial g/\partial t$ and
$g^\prime=\partial g/\partial \epsilon$ using the right translations as 
\[
\dot{g}=\xi\circ{g}
\quad\hbox{and}\quad
{g}^\prime=\eta\circ{g}
\,.
\]
By the chain rule, these definitions have mixed partial derivatives
\[
\dot{g}^\prime=\xi^\prime=\nabla\xi\cdot\eta
\quad\hbox{and}\quad
\dot{g}^\prime=\dot{\eta}=\nabla\eta\cdot\xi
\,.
\]
The difference of the mixed partial derivatives implies the desired
formula (\ref{mixed partials}),
\[
\xi^\prime-\dot{\eta}=\nabla\xi\cdot\eta-\nabla\eta\cdot\xi
= -\,[ \xi, \eta]
= -\,{\rm ad}_\xi \eta
\,.
\]
(Note the minus sign in the last two terms.)
\end{proof}

\paragraph{\bfi Step III. Proof of equivalence 
$\textbf{(iii)} \Longleftrightarrow \textbf{(iv)}$}
Let us show that the reduced variational principle produces the
Euler-Poincar\'e equations. We write the functional derivative of the
reduced action $S_{red}=\int_a^b l(\xi)\,dt$ with Lagrangian $l(\xi)$ in
terms of the natural pairing 
$\left\langle \cdot\, , \,\cdot\right\rangle$ between $\mathfrak{g}^\ast $
and
$\mathfrak{g}$ as
\begin{align*}
\delta \int_a^b l(\xi(t)) dt &=  \int_a^b \left\langle
\frac{\delta l}{\delta \xi}, \delta\xi\right\rangle dt
=  \int_a^b \left\langle \frac{\delta l}{\delta \xi} , \dot \eta -
{\rm ad}_\xi \eta\right\rangle dt\\
&=  \int_a^b \left\langle
\frac{\delta l}{\delta \xi}, \dot \eta \right\rangle dt -  \int_a^b
\left\langle \frac{\delta l}{\delta \xi} , {\rm ad}_\xi
\eta\right\rangle dt\\
& = -\int_a^b\, \left\langle 
 \frac{d}{dt} \frac{\delta l}{\delta \xi}
+ {\rm ad}^*_\xi \frac{\delta l}{\delta \xi}\, ,\, \eta \right\rangle\,
dt \,.
\end{align*}
The last equality follows from integration by parts and
vanishing of the variation $\eta(t)$  at the endpoints. Thus,
stationarity $\delta \int_a^b l(\xi(t)) dt = 0$ for any $\eta(t) $ that
vanishes at  the endpoints is equivalent to
$$\frac{d}{dt} \frac{\delta l}{\delta \xi} 
= -\,{\rm ad}^*_\xi\,\frac{\delta l}{\delta \xi}
\,,
$$
which are the Euler-Poincar\'e equations.
\end{proof}

\begin{remark}[{Left}-invariant Euler-Poincar\'e equations]
The same theorem holds for {left} invariant Lagrangians on $TG$, except
for a sign in the Euler-Poincar\'e equations,
$$
\frac{d}{dt} \frac{\delta l}{\delta \xi} = +\, {\rm ad}^*_\xi
\frac{\delta l}{\delta \xi}
\,,
$$
which arises because {left}-invariant variations satisfy 
$\delta \xi = \dot\eta + [\xi ,\ \eta ]$ (with the opposite sign).
\end{remark}

\begin{exercise}
Write out the corresponding proof of the Euler-Poincar\'e reduction
theorem for left-invariant Lagrangians defined on the tangent space $TG$
of a group $G$.
\end{exercise}

\paragraph{Reconstruction.} The procedure for reconstructing the solution
$v(t)\in T_{g(t)}G$ of the  Euler-Lagrange equations with initial
conditions $g(0) = g_0$ and $\dot{g}(0) = v_0 $ starting from the solution
of  the Euler-Poincar\'{e} equations is as follows. First, solve the
initial value problem for the {right}-invariant Euler-Poincar\'{e}
equations:
\begin{eqnarray*}
\frac{d}{dt} \frac{\delta l}{\delta \xi} 
= -\,{\rm ad}^*_\xi \frac{\delta l}{\delta \xi}
\quad\hbox{with}\quad
\xi(0) = \xi_0 : =  v_0 g_0^{-1}
\,.
\end{eqnarray*}
Then from the solution for $\xi(t)$ reconstruct the curve $g(t)$ on the
group by solving the ``linear differential equation with time-dependent
coefficients"
\[
\dot g(t) = \xi(t) g(t) 
\quad\hbox{with}\quad
g(0) = g_0\,.
\]
The Euler-Poincar\'{e} reduction theorem guarantees then that $v(t) =
\dot g(t) =  \xi(t) \cdot g(t) $ is a solution of the Euler-Lagrange
equations with initial condition $v_0 = \xi_0 g_0$.

\begin{remark}
A similar statement holds, with obvious changes for {left}-invariant
Lagrangian systems on $TG $.
\end{remark}

\subsection{Reduced Legendre transformation}
As in the equivalence relation between the Lagrangian and Hamiltonian
formulations discussed earlier, the relationship between
symmetry-reduced Euler-Poincar\'{e} and Lie-Poisson formulations is
determined by the Legendre transformation. 

\begin{definition}
The Legendre
transformation $\mathbb{F}l :\mathfrak{g}\rightarrow\mathfrak{g}^*$ is
defined  by
$$
\mathbb{F} l (\xi) = \frac{\delta l}{\delta \xi} = \mu \,.
$$
\end{definition}

\paragraph{Lie-Poisson Hamiltonian formulation.}
Let $h(\mu) := \langle \mu , \xi\rangle - l(\xi)$.
Assuming that $\mathbb{F} l$ is a diffeomorphism yields
$$\frac{\delta h}{\delta \mu} = \xi + \left\langle \mu\, ,\,
\frac{\delta
\xi}{\delta
\mu} \right\rangle -  \left\langle \frac{\delta l}{\delta \xi}\, ,\,
\frac{\delta \xi}{\delta
\mu} \right\rangle  = \xi .
$$
So the Euler-Poincar\'{e} equations for $l$ are equivalent to the
Lie-Poisson equations for $h$:
$$\frac{d}{dt} \left( \frac{\delta l}{\delta \xi}\right) = -\,{\rm
ad}^*_\xi
\frac{\delta l}{\delta \xi} \Longleftrightarrow \dot \mu = -\,{\rm
ad}^*_{\delta h/\delta \mu} \mu.
$$
The Lie-Poisson equations may be written in the Poisson bracket
form
\begin{equation} \label{lppoisson}
\dot{f} = \left\{ f, h \right\},
\end{equation}
where $f: \mathfrak{g}^\ast \rightarrow \mathbb{R}$ is
an arbitrary smooth function and the bracket is the
(right) Lie-Poisson bracket given by
\begin{equation}
\label{lpb} \{f, h\}(\mu )  
=  \left\langle \mu , \left[ \frac{ \delta f}{\delta  \mu},
\frac{\delta  h}{\delta \mu } \right] \right\rangle 
=
-
\left\langle \mu , 
\operatorname{ad}_{\delta h/\delta \mu}
\frac{ \delta f}{\delta  \mu}\right\rangle
=
-
\left\langle \operatorname{ad}^*_{\delta h/\delta \mu}\mu , 
\frac{ \delta f}{\delta  \mu}\right\rangle
\,.
\end{equation} 

In the important case when $\ell$ is quadratic, the Lagrangian $L $
is  the quadratic form associated to a right invariant Riemannian
metric on $G$. In this case, the Euler--Lagrange equations for $L$
on $G$ describe geodesic motion relative to this metric and these
geodesics are then equivalently described by either the
Euler-Poincar\'e, or the Lie-Poisson equations.

\begin{exercise}
Exercise \ref{GLn-ex} requires an extension of the pure EP reduction
theorem for a Lagrangian $L:\,(TG\times TQ)\to \mathbb{R}$. Following the
proof of the EP reduction theorem, make this extension.
\end{exercise}

\begin{exercise}
Compute the pure EP equations for geodesic motion on $SE(3)$. These equations turn out to be applicable to the motion of an ellipsoidal body through a fluid. 
\end{exercise}

\section{EPDiff: the Euler-Poincar\'e equation on the diffeomorphisms}

\subsection{The $n-$dimensional EPDiff equation and its properties}
Eulerian geodesic motion of a fluid in $n-$dimensions is
generated as an EP equation via Hamilton's principle, when the 
Lagrangian is given by the kinetic energy. The kinetic energy defines a
norm $\|\mathbf{u}\|^2$ for the Eulerian fluid velocity, taken as 
$\mathbf{u}(\mathbf{x},t):\,R^n\times R^1 \to R^n$. The
choice of the kinetic energy as a positive functional of fluid velocity
$\mathbf{u}$ is a modeling step that depends upon the physics of the
problem being studied. We shall choose the Lagrangian,
\begin{equation}\label{Lag-ansatz-1}
\|\mathbf{u}\|^2 
= \int \mathbf{u}\cdot Q_{op}\mathbf{u}\,d\,^nx
= \int \mathbf{u}\cdot\mathbf{m}\,d\,^nx
\,,\end{equation}
so that the positive-definite, symmetric, operator $Q_{op}$ defines the
norm $\|\mathbf{u}\|$, for appropriate (homogeneous, say, or periodic)
boundary conditions.  The EPDiff equation is the Euler-Poincar\'e
equation for this Eulerian geodesic motion of a fluid. Namely, 
\begin{equation}\label{EP-eqn-abstract}
\frac{d}{dt}\frac{\delta \ell}{\delta \mathbf{u}}
+
{\rm ad}^*_\mathbf{u}
\frac{\delta \ell}{\delta \mathbf{u}}
=
0
\,,
\quad\hbox{with}\quad
\ell[\mathbf{u}]=\frac{1}{2}\|\mathbf{u}\|^2
\,.
\end{equation}
Here ${\rm ad}^*$ is the dual of the vector-field ad-operation (the
commutator) under the natural $L^2$ pairing
$\langle\cdot\,\,,\,\cdot\rangle$ induced by the variational derivative
$\delta\ell[\mathbf{u}]=\langle\delta\ell/\delta\mathbf{u}\,
\,,\,\delta\mathbf{u}\rangle$. This pairing provides the definition of
${\rm ad}^*$,
\begin{equation}\label{ad*-eqn}
\langle{\rm ad}^*_\mathbf{u}\,\mathbf{m}\,,\,\mathbf{v}\rangle
=
-\,\langle\mathbf{m}\,,\,{\rm ad}_\mathbf{u}\mathbf{v}\rangle
\,,
\end{equation}
where $\mathbf{u}$ and $\mathbf{v}$ are vector fields, ${\rm
ad}_\mathbf{u}\mathbf{v}=[\mathbf{u},\mathbf{v}]$ is the commutator, i.e.,
the {\bfi  Lie bracket\/} given in
components by (summing on repeated indices)
\begin{equation}\label{jlb-ad} [\mathbf{u}, \mathbf{v}]^i
= 
u^j\frac{\partial v^i}{\partial x^j} 
 - 
v^j\frac{\partial u^i}{\partial x^j}
\,,\quad\hbox{or}\quad
[\mathbf{u}, \mathbf{v}]
=
\mathbf{u}\cdot\nabla \mathbf{v}
-
\mathbf{v}\cdot\nabla \mathbf{u}
\,.
\end{equation}
The notation $
\operatorname{ad}_\mathbf{u} \mathbf{v} := [\mathbf{u},\,\mathbf{v}]$
formally denotes the adjoint action of the {\it right\/} Lie algebra of
$\operatorname{Diff}(\mathcal{D})$ on itself,
and  $\mathbf{m}=\delta\ell/\delta\mathbf{u}$ is the fluid momentum, a
{\bfi one-form density} whose co-vector components are also denoted as
$\mathbf{m}$. 

If $\mathbf{u} = u^j \partial/\partial x^j,\,
\mathbf{m}  = m_i dx^i\otimes  dV$, then the preceding formula for $
\operatorname{ad}^\ast_\mathbf{u}(\mathbf{m}\otimes  dV)$ has the {\bfi
coordinate expression} in $\mathbb{R} ^n$,
\begin{eqnarray}\label{continuumcoadjoint-coord-form}
\Big(\operatorname{ad}^\ast_\mathbf{u}\mathbf{m}\Big)_{\!i}
dx^i\otimes  dV
&=&
\left (\frac{\partial}{\partial x^j}(u^jm_i) +
m_j \frac{\partial u^j}{\partial x^i}\right ) dx^i\otimes  dV
\,.
\end{eqnarray}

In this notation,  the abstract EPDiff equation (\ref{EP-eqn-abstract})
may be written explicitly in Euclidean coordinates as a partial
differential equation for a co-vector function
$\mathbf{m}(\mathbf{x},t):\,R^{\,n}\times R^1 \to R^{\,n}$. Namely, 
\begin{equation}\label{ep-eqn-coord-form}
\frac{\partial }{\partial t}\mathbf{m}
\ 
+ 
\underbrace{\
\mathbf{u}\cdot\nabla \mathbf{m}\
}_{\hbox{Convection}}
+\ 
\underbrace{\
\nabla \mathbf{u}^T\cdot\mathbf{m}\
}_{\hbox{Stretching}}\
+ 
\underbrace{\
\mathbf{m}({\rm div\,}\mathbf{u})\
}_{\hbox{Expansion}}
=0
\,,
\quad\hbox{with}\quad
\mathbf{m}= \frac{\delta \ell}{\delta \mathbf{u}}
=Q_{op}\mathbf{u}
\,.
\end{equation}
To explain the terms in underbraces, we rewrite EPDiff as preservation of the one-form
density of momentum along the characteristic curves of the velocity. Namely, 
\begin{equation}\label{EPDiff-char-form}
\frac{d}{dt}\Big(\mathbf{m}\cdot d\mathbf{x}\otimes dV\Big)=0
\quad\hbox{along}\quad
\frac{d\mathbf{x}}{dt}=\mathbf{u}=G*\mathbf{m}
\,.
\end{equation} 
This form of the EPDiff equation also emphasizes its nonlocality, since the velocity is
obtained from the momentum density by convolution against the Green's function $G$ of the
operator $Q_{op}$. Thus, $\mathbf{u}=G*\mathbf{m}$ with $Q_{op}G=\delta(\mathbf{x})$, the
Dirac measure. We may check that this ``characteristic form'' of EPDiff recovers its
Eulerian form by computing directly,
\begin{eqnarray*}
\frac{d}{dt}\Big(\mathbf{m}\cdot d\mathbf{x}\otimes dV\Big)
&=&
\frac{d\mathbf{m}}{dt}\cdot d\mathbf{x}\otimes dV
+
\mathbf{m}\cdot d\frac{d\mathbf{x}}{dt}\otimes dV
+
\mathbf{m}\cdot d\mathbf{x}\otimes \Big(\frac{d}{dt}dV\Big)
\hbox{ along }
\frac{d\mathbf{x}}{dt}=\mathbf{u}=G*\mathbf{m}
\\
&=&
\Big(
\frac{\partial }{\partial t}\mathbf{m}
+ 
\mathbf{u}\cdot\nabla \mathbf{m}\
+
\nabla \mathbf{u}^T\cdot\mathbf{m}
+ 
\mathbf{m}({\rm div\,}\mathbf{u})
\Big)\cdot d\mathbf{x}\otimes dV
=0
\,.
\end{eqnarray*}
\begin{exercise}
Show that EPDiff may be written as 
\begin{equation}\label{EPDiff-Lie-form}
\Big(\frac{\partial}{\partial t}+\mathcal{L}_\mathbf{u}\Big)
\Big(\mathbf{m}\cdot d\mathbf{x}\otimes dV\Big)=0
\,,
\end{equation} 
where $\mathcal{L}_\mathbf{u}$ is the Lie derivative with respect to the vector field 
with components $\mathbf{u}=G*\mathbf{m}$. Hint: How does this Lie-derivative form of
EPDiff in (\ref{EPDiff-Lie-form}) differ from its characteristic form
(\ref{EPDiff-char-form})?
\end{exercise}

EPDiff may also be written equivalently in terms of the operators div, grad and curl in
2D and 3D as,
\begin{equation}\label{div-grad-curl-eqn}
\frac{\partial }{\partial t}\mathbf{m} 
- \mathbf{u}\times{\rm\,curl\,}\mathbf{m} 
+ \nabla (\mathbf{u}\cdot\mathbf{m}) 
+ \mathbf{m}({\rm div\,}\mathbf{u})
=0
\,.
\end{equation}
Thus, for example, its numerical solution would require an algorithm which has the
capability to deal with the distinctions and relationships among the operators div, grad
and curl.\\

\subsection{Derivation of the $n-$dimensional EPDiff equation as geodesic flow}
Let's derive the EPDiff equation (\ref{ep-eqn-coord-form}) by following the proof of the
EP reduction theorem leading to the Euler-Poincar\'e equations for right invariance in
the form (\ref{EP-eqn-abstract}). Following this calculation for the present
case yields
\begin{eqnarray*}
\delta \int_a^b l(\mathbf{u}) dt 
&=&  \int_a^b \left\langle
\frac{\delta l}{\delta \mathbf{u}}, \delta\mathbf{u}\right\rangle dt
\
=  \int_a^b \left\langle 
  \frac{\delta l}{\delta \mathbf{u}} , \mathbf{\dot{v}} 
-
{\rm ad}_\mathbf{u} \mathbf{v}\right\rangle dt
\\
&&\hspace{-2cm}=  
\int_a^b \left\langle
\frac{\delta l}{\delta \mathbf{u}}, \mathbf{\dot{v}} \right\rangle dt 
-  
\int_a^b
\left\langle \frac{\delta l}{\delta \mathbf{u}} , {\rm ad}_\mathbf{u}
\mathbf{v}\right\rangle dt
\
=\
 -\int_a^b\, \left\langle 
  \frac{d}{dt} \frac{\delta l}{\delta \mathbf{u}} 
+ 
{\rm ad}^*_\mathbf{u} \frac{\delta l}{\delta \mathbf{u}}\, ,\, \mathbf{v} \right\rangle\,
dt \,,
\end{eqnarray*}
where $\langle\cdot\,,\,\cdot\rangle$ is the pairing between elements of
the Lie algebra and its dual. In our case, this is the $L^2$ pairing, e.g.,
\[
\left\langle
\frac{\delta l}{\delta \mathbf{u}}, \delta\mathbf{u}\right\rangle
=
\int \frac{\delta l}{\delta {u}^i}\,\delta {u}^i\,d\,^nx
\]
This pairing allows us to compute the coordinate form of the EPDiff equation explicitly,
as
\begin{eqnarray*}
\int_a^b \left\langle
\frac{\delta l}{\delta \mathbf{u}}, \delta\mathbf{u}\right\rangle dt
&=&
\int_a^b dt 
\int \frac{\delta l}{\delta {u}^i}
\Big(
\frac{\partial v^i}{\partial t}
+
u^j\frac{\partial v^i}{\partial x^j} 
 - 
v^j\frac{\partial u^i}{\partial x^j}
\Big)
\,d\,^nx
\nonumber\\&=&
-
\int_a^b dt 
\int \bigg\{
\frac{\partial }{\partial t}\frac{\delta l}{\delta {u}^i}
+
\frac{\partial }{\partial x^j}\Big( \frac{\delta l}{\delta {u}^i}u^j \Big)
+
\frac{\delta l}{\delta {u}^j}\frac{\partial u^j}{\partial x^i} 
\bigg\} v^i
\,d\,^nx
\end{eqnarray*}
Substituting $\mathbf{m}=\delta l/\delta \mathbf{u}$ now recovers the coordinate forms for
the coadjoint action of vector fields in (\ref{continuumcoadjoint-coord-form}) and the
EPDiff equation itself in (\ref{ep-eqn-coord-form}). When
$\ell[\mathbf{u}]=\frac{1}{2}\|\mathbf{u}\|^2$, EPDiff describes geodesic motion on the
diffeomorphisms with respect to the norm $\|\mathbf{u}\|$.

\begin{lemma} 
In Step IIB of the proof of the Euler-Poincar\'e reduction theorem that $\textbf{(ii)}
\Longleftrightarrow\textbf{(iv)}$ for an arbitrary Lie group, a certain formula for the
variations for time-dependent vector fields was employed. That formula was employed again
in the calculation above as,
\begin{equation}\label{mixed-partials}
\delta \mathbf{u}
=
\mathbf{\dot{v}} 
-
{\rm ad}_\mathbf{u} \mathbf{v}
\,.
\end{equation}
This formula may be rederived as follows in the present context.
We write $\mathbf{u}=\dot{g}g^{-1}$ and $\mathbf{v}=g^\prime g^{-1}$ in
natural notation and express the partial derivatives $\dot{g}=\partial
g/\partial t$ and $g^\prime=\partial g/\partial \epsilon$ using the right translations as 
\[
\dot{g}=\mathbf{u}\circ{g}
\quad\hbox{and}\quad
{g}^\prime=\mathbf{v}\circ{g}
\,.
\]
To compute the mixed partials, consider the chain rule for 
say $\mathbf{u}(g(t,\epsilon)\mathbf{x}_0)$ and set 
$\mathbf{x}(t,\epsilon)=g(t,\epsilon)\cdot\mathbf{x}_0$. 
Then,
\[
\mathbf{u}^\prime
=
\frac{\partial \mathbf{u}}{\partial \mathbf{x}}
\cdot
\frac{\partial \mathbf{x}}{\partial \epsilon} 
=
\frac{\partial \mathbf{u}}{\partial \mathbf{x}}
\cdot
g^\prime(t,\epsilon)\mathbf{x}_0
=
\frac{\partial \mathbf{u}}{\partial \mathbf{x}}
\cdot
g^\prime g^{-1}\mathbf{x}
=
\frac{\partial \mathbf{u}}{\partial \mathbf{x}}
\cdot\mathbf{v}(\mathbf{x})
\,.
\]
The chain rule for $\mathbf{\dot{v}}$  gives a similar formula
with $\mathbf{u}$ and $\mathbf{v}$ exchanged. Thus, the chain rule gives
two expressions for the mixed partial derivative
$\dot{g}^\prime$ as 
\[
\dot{g}^\prime=\mathbf{u}^\prime=\nabla\mathbf{u}\cdot\mathbf{v}
\quad\hbox{and}\quad
\dot{g}^\prime=\mathbf{\dot{v}}=\nabla\mathbf{v}\cdot\mathbf{u}
\,.
\]
The difference of the mixed partial derivatives then implies the desired
formula (\ref{mixed-partials}), since
\[
\mathbf{u}^\prime-\mathbf{\dot{v}}
=\nabla\mathbf{u}\cdot\mathbf{v}
-\nabla\mathbf{v}\cdot\mathbf{u}
= -\,[ \mathbf{u}, \mathbf{v}]
= -\,{\rm ad}_\mathbf{u}\mathbf{v}
\,.
\]
\end{lemma}

\section{EPDiff: the Euler-Poincar\'e equation on the diffeomorphisms}

In this lecture, we shall discuss the solutions of EPDiff for pressureless
compressible geodesic motion in one spatial dimension.  This is the {\bfi
EPDiff equation} in 1D,%
\footnote{
A one-form density in 1D takes the form
$m\,(dx)^2$ and the EP equation is given by
\[
\frac{d}{dt}\big(m\,(dx)^2\big)
=
\frac{dm}{dt}(dx)^2 + 2m\,(du)(dx)
=0
\quad\hbox{with}\quad
\frac{d}{dt}dx=du=u_xdx
\quad\hbox{and}\quad
u=G*m
\,,
\]
where $G*m$ denotes convolution with a function $G$ on the real
line. 
}
\begin{eqnarray}\label{1dEPDiff}
&&\partial_tm + {\rm ad}^*_um=0
\,,\quad\hbox{or, equivalently,}\quad
\\&&
\partial_tm+um_x+2u_xm=0
\,,\quad\hbox{with}\quad
m= Q_{op}u
\,.
\label{1dEPDiff-pde}
\end{eqnarray}
\begin{itemize}
\item
The EPDiff equation describes geodesic motion on the diffeomorphism group
with respect to a family of metrics for the fluid velocity $u(t,x)$, with
notation,
\begin{eqnarray}\label{mom-def-1D}
m
&=& \frac{\delta\ell}{\delta u} = Q_{op}u
\quad\hbox{for a kinetic-energy Lagrangian}\quad
\\&&
\ell(u) 
= \frac{1}{2}\int u\,Q_{op}u\,dx
= \frac{1}{2}\|u\|^2
\,.
\end{eqnarray}
\item
In one dimension, $Q_{op}$ in equation (\ref{mom-def-1D}) is a positive,
symmetric operator that defines the kinetic energy metric for the
velocity.  
\item
The EPDiff equation (\ref{1dEPDiff-pde}) is written in terms of
the variable $m=\delta\ell/\delta u$. It is appropriate to call this variational
derivative $m$, because it is the momentum density associated with the
fluid velocity $u$. 
\item
Physically, the first nonlinear term in the EPDiff equation
(\ref{1dEPDiff-pde}) is fluid  transport. 
\item
The coefficient 2 arises in the
second nonlinear term, because, in one dimension, two of the  summands in
ad$^*_u\,m=um_x+2u_xm$ are the {\it same}, cf. equation
(\ref{continuumcoadjoint-coord-form}).  
\item
The momentum is expressed in terms of the velocity by $m=\delta\ell/\delta{u}=Q_{op}u$.
Equivalently, for solutions that vanish at spatial infinity, one may think of the velocity
as being obtained from the convolution,
\begin{equation}\label{filter-reln-1D}
u(x)=G*m(x)=\int G(x-y)m(y)\,dy
\,,
\end{equation}
where $G$ is the Green's function for the operator $Q_{op}$ on the real
line. 
\item
The operator $Q_{op}$ and its Green's function $G$ are chosen to be
even under reflection, $G(-x)=G(x)$, so that $u$ and $m$ have the same parity. Moreover, the
EPDiff equation (\ref{1dEPDiff-pde}) conserves the total momentum $M=\int  m(y)\,dy$,
for any even Green's function.
\begin{exercise}
Show that equation (\ref{1dEPDiff-pde}) conserves $M=\int  m(y)\,dy$ for any even Green's
function $G(-x)=G(x)$, for either periodic, or homogeneous boundary conditions.
\end{exercise}
\item
The traveling wave solutions of 1D EPDiff when the Green's function $G$ is chosen to be
even under reflection are the ``pulsons,'' 
\[
u(x,t)=c\,G(x-ct)
\,.
\] 
\begin{exercise}
Prove this statement, that the traveling wave solutions of 1D EPDiff are pulsons when the
Green's function is even. What role is played in the solution by the Green's function being
even? Hint: Evaluate the derivative of an even function at $x=0$.
\end{exercise}
\item
See Fringer and Holm {FrHo2001} and references therein for further discussions and
numerical simulations of the pulson solutions of the 1D EPDiff equation.
\end{itemize}

\subsection{Pulsons} 
The EPDiff equation (\ref{1dEPDiff-pde}) on the real line has the
remarkable property that its solutions {\bfi collectivize}%
\footnote{See \cite{GuSt1984} for discussions of the concept of collective variables
for Hamiltonian theories. We will discuss the collectivization for the EPDiff
equation later from the viewpoint of momentum maps.}
into the finite dimensional solutions of the ``$N-$pulson'' form that was discovered
for a special form of $G$ in Camassa and Holm {CaHo1993}, then was
extended for {\it any} even $G$ in Fringer and Holm {FrHo2001},
\begin{equation}\label{u-pulson}
u(x,t) = \sum_{i=1}^N p_i(t)\,G(x-q_i(t))
\,.
\end{equation}
Since $G(x)$ is the Green's function for the operator $Q_{op}$, the
corresponding solution for the momentum $m=Q_{op}u$ is given by a sum of
delta functions,
\begin{equation}\label{m-delta1}
m(x,t) = \sum_{i=1}^N p_i(t)\,\delta(x-q_i(t))
\,.
\end{equation}
Thus, the time-dependent ``collective coordinates'' $q_i(t)$ and $p_i(t)$
are the positions and velocities of the $N$ pulses in this solution. These
parameters satisfy the finite dimensional geodesic motion equations
obtained as canonical Hamiltonian equations 
\begin{eqnarray} \label{q-eqn} 
\dot{q}_i & = & 
\frac{\partial H_N}{\partial p_i}
=
\sum_{j=1}^{N} p_j\, G(q_i-q_j)
\,, 
\\
\dot{p}_i & = & 
-\,\frac{\partial H_N}{\partial q_i}
=
-\,p_i\sum_{j=1}^{N} p_j\, G\,'(q_i-q_j)
\,,
\label{p-eqn}
\end{eqnarray} 
in which the Hamiltonian is given by the quadratic form,
\begin{equation}\label{Ham-pulson}
H_N = \frac{1}{2} \sum_{i,j=1}^Np_i\,p_j\,
G(q_i-q_j)
\,.
\end{equation}
\begin{remark}In a certain sense, equations (\ref{q-eqn}-\ref{p-eqn})
comprise the analog for the peakon momentum relation (\ref{m-delta1}) of
the ``symmetric generalized rigid body equations'' in (\ref{HamCan-eqns}).
\end{remark}

Thus, the canonical equations for the Hamiltonian $H_N$ describe the
nonlinear collective interactions of the $N-$pulson solutions of the
EPDiff equation (\ref{1dEPDiff-pde}) as finite-dimensional geodesic
motion of a particle on an
$N-$dimensional surface whose co-metric is 
\begin{equation}\label{metric-pulson}
G^{ij}(q)=G(q_i-q_j)
\,.
\end{equation}
Fringer and Holm {FrHo2001} showed numerically that the $N-$pulson
solutions describe the emergent patterns in the solution of the
initial value problem for EPDiff equation (\ref{1dEPDiff-pde}) with 
spatially confined initial conditions. 

\begin{exercise}
Equations (\ref{q-eqn}-\ref{p-eqn}) describe geodesic motion. 
\begin{enumerate}
\item
Write the Lagrangian and Euler-Lagrange equations for this motion. 
\item
Solve equations (\ref{q-eqn}-\ref{p-eqn}) for $N=2$ when
$\lim_{|x|\to\infty}G(x)=0$. 
\begin{enumerate}
\item
Why should the solution be described as exchange of momentum in elastic collisions?
\item
Consider both head-on and overtaking collisions.  
\item
Consider the antisymmetric case, when the total momentum
vanishes.
\end{enumerate}
\end{enumerate}
\end{exercise}

\paragraph{Integrability} Calogero and Francoise
\cite{Ca1995}, \cite{CaFr1996} found that for any finite number 
$N$  the Hamiltonian equations for $H_{N}$ in (\ref{Ham-pulson}) are 
completely integrable in the Liouville  sense%
\footnote{A Hamiltonian system is integrable in the Liouville  sense, if
the number of independent constants of motion in involution is the same 
as the number of its degrees of freedom.}
  for $G \equiv
G_{1}(x)=\lambda + \mu \cos(\nu x) +
\mu_{1} \sin(\nu |x|)$ and 
$G \equiv G_{2}(x)=\alpha + \beta |x| + \gamma x^{2}$, 
with $\lambda$, $\mu$, $\mu_{1}$, 
$\nu$, and $\alpha$, $\beta$, $\gamma$ 
being arbitrary constants, such that $\lambda$ and $\mu$ 
are real and $\mu_{1}$  and $\nu$ both real or both imaginary.%
\footnote{
This choice of the constants keeps $H_{N}$ real in (\ref{Ham-pulson}).} 
Particular cases of $G_{1}$ and $G_{2}$  are the peakons
$G_1(x)=e^{-|x|/\alpha}$ of  \cite{CaHo1993} and  the compactons 
$G_2(x)={\rm max}(1-|x|,0)$ of the Hunter-Saxton equation,
\cite{HuZh1994}.  The latter is the EPDiff equation (\ref{1dEPDiff-pde}),
with $\ell(u)=\frac{1}{2}\int u_x^2\,dx$ and thus
$m=-u_{xx}$.

\paragraph{Lie-Poisson Hamiltonian form of EPDiff} 
In terms of $m$, the conserved energy Hamiltonian for the EPDiff equation
(\ref{1dEPDiff-pde}) is obtained by Legendre transforming the kinetic energy
Lagrangian, as
\[
h = \Big\langle \frac{\delta\ell}{\delta u}\,,u\Big\rangle- \ell(u)
\,.\]
Thus, the Hamiltonian depends on $m$, as
\[
h(m)=\frac{1}{2}\!\int\! m(x)G(x-y)m(y)\,dxdy
\,,\]
which also reveals the geodesic nature of the EPDiff equation
(\ref{1dEPDiff-pde}) and the role of $G(x)$ in the kinetic energy metric on
the Hamiltonian side. 

The corresponding {\bfi Lie-Poisson bracket} for EPDiff as a Hamiltonian evolution
equation is given by,
\[
\partial_tm = \big\{m,h\big\} = -\,(\partial{m}+m\partial)\frac{\delta h}{\delta m}
\quad\hbox{and}\quad
\frac{\delta h}{\delta m}=u
\,,
\]
which recovers the starting equation and indicates some of its connections with fluid
equations on the Hamiltonian side. For any two smooth functionals $f,h$ of $m$ in the space
for which the solutions of EPDiff exist, this Lie-Poisson bracket may be expressed as,
\[
\big\{f,h\big\} 
= -\int \frac{\delta f}{\delta m}
(\partial{m}+m\partial)\frac{\delta h}{\delta m} dx
= -\int m \bigg[\frac{\delta f}{\delta m}
\,,\,\frac{\delta h}{\delta m}\bigg] dx
\]
where $[\cdot\,,\,\cdot]$ denotes the Lie algebra bracket of vector fields. That is,
\[
\bigg[\frac{\delta f}{\delta m}
\,,\,\frac{\delta h}{\delta m}\bigg]
=
\frac{\delta f}{\delta m}\partial\frac{\delta h}{\delta m}
-
\frac{\delta h}{\delta m}\partial\frac{\delta f}{\delta m}
\,.
\]
\begin{exercise}
What is the Casimir for this Lie Poisson bracket? What does it mean from the viewpoint of
coadjoint orbits? 
\end{exercise}

\subsection{Peakons} 
The case $G(x)=e^{-|x|/\alpha}$ with a constant lengthscale $\alpha$ is the
Green's function for which the operator in the kinetic
energy Lagrangian (\ref{mom-def-1D}) is
$ Q_{op}=1-\alpha^2\partial_x^2$. For this (Helmholtz)
operator $ Q_{op}$, the Lagrangian and corresponding kinetic energy
norm are given by, 
\[
\ell\,[u]
=\frac{1}{2}\|u\|^2
= \frac{1}{2}\int u\, Q_{op}u\,dx
=\frac{1}{2}\int u^2+\alpha^2u_x^2\
dx
\,,\quad\hbox{for} \quad
\lim_{|x|\to\infty}u=0
\,.\] 
This Lagrangian is the $H\,^1$ norm of the velocity in one
dimension. In this case, the EPDiff equation (\ref{1dEPDiff-pde}) is also
the zero-dispersion limit of the completely integrable CH 
equation for unidirectional shallow water waves first derived in 
Camassa and Holm {CaHo1993},
\begin{equation} \label{CH-eqn}
m_t + um_x + 2 mu_x = -c_0u_x+\gamma{u}_{xxx}
\,,\qquad
m=u-\alpha^2u_{xx}
\,.
\end{equation}
This equation describes shallow water dynamics as completely integrable
soliton motion at quadratic order in the asymptotic expansion for
unidirectional shallow water waves on a free surface under gravity. 
See Dullin, Gottwald and Holm {DGH[2001,2003,2004]} for more details and explanations of
this asymptotic expansion for unidirectional shallow water waves to quadratic order.

Because of the relation $m=u-\alpha^2u_{xx}$, equation (\ref{CH-eqn}) is 
nonlocal. In other words, it is an integral-partial differential equation. 
In fact, after writing equation (\ref{CH-eqn}) in the equivalent form,
\begin{equation} \label{CH-eqn-nonloc}
(1-\alpha^2\partial^2)(u_t+uu_x)
=
-\,\partial\Big(u^2+\frac{\alpha^2}{2}u_x^2\Big)
-c_0u_x+\gamma{u}_{xxx}
\,,
\end{equation}
one sees the interplay between local and nonlocal linear dispersion in its phase
velocity relation,
\begin{equation} \label{disp-reln}
\frac{\omega}{k}
=
\frac{c_0 - {\gamma}\, k^2}{1 + \alpha^{\,2} k^2}
\,,
\end{equation}
for waves with frequency $\omega$ and wave number $k$ linearized
around ${u}=0$. For $\gamma/c_0<0$, short waves and long
waves travel in the same direction. Long waves travel faster than short ones
(as required in shallow water) provided $\gamma/c_0 > -\,  \alpha^2$.
Then the phase velocity lies in the interval $\omega/k\in(-\,\gamma/\alpha^{\,2}, c_0]$.

The famous Korteweg-de Vries (KdV) soliton equation, 
\begin{equation} \label{KdV-eqn}
u_t+3uu_x
=
-c_0u_x+\gamma{u}_{xxx}
\,,
\end{equation}
emerges at {\it linear} order in the asymptotic expansion for shallow water waves, in which
one takes $\alpha^2\to0$ in (\ref{CH-eqn-nonloc}) and (\ref{disp-reln}). In KdV, the  
parameters $c_0$ and $\gamma$ are seen as deformations of the {\it Riemann equation}, 
\[
u_t+3uu_x=0
\,.
\]
The parameters $c_0$ and $\gamma$ represent linear wave dispersion, which modifies and
eventually balances the tendency for nonlinear waves to steepen and break. The parameter
$\alpha$, which introduces nonlocality, also regularizes this nonlinear tendency, even in
the absence of $c_0$ and $\gamma$.

\section{Diffeons -- singular momentum solutions of the EPDiff equation for
geodesic motion in higher dimensions}
\label{sec-strings}
As an example of the EP theory  in higher dimensions, we shall 
generalize the one-dimensional pulson solutions of the previous section to
$n-$dimensions. The corresponding singular momentum solutions of the EPDiff
equation in higher dimensions are called ``diffeons.''


\subsection{$n-$dimensional EPDiff equation}
Eulerian geodesic motion of a fluid in $n-$dimensions is
generated as an EP equation via Hamilton's principle, when the 
Lagrangian is given by the kinetic energy. The kinetic energy defines a norm
$\|\mathbf{u}\|^2$ for the Eulerian fluid velocity,
$\mathbf{u}(\mathbf{x},t):\,R^n\times R^1 \to R^n$. As mentioned earlier, the
choice of the kinetic energy as a positive functional of fluid velocity
$\mathbf{u}$ is a modeling step that depends upon the physics of the
problem being studied. Following our earlier procedure, as in equations
(\ref{Lag-ansatz-1}) and (\ref{EP-eqn-abstract}), we shall choose
the Lagrangian,
\begin{equation}\label{Lag-ansatz-2}
\|\mathbf{u}\|^2 
= \int \mathbf{u}\cdot Q_{op}\mathbf{u}\,d\,^nx
= \int \mathbf{u}\cdot\mathbf{m}\,d\,^nx
\,,\end{equation}
so that the positive-definite, symmetric, operator $Q_{op}$ defines the
norm $\|\mathbf{u}\|$, for appropriate boundary conditions and the EPDiff
equation for Eulerian geodesic motion of a fluid emerges,
\begin{equation}\label{EP-eqn}
\frac{d}{dt}\frac{\delta \ell}{\delta \mathbf{u}}
+
{\rm ad}^*_\mathbf{u}
\frac{\delta \ell}{\delta \mathbf{u}}
=
0
\,,
\quad\hbox{with}\quad
\ell[\mathbf{u}]=\frac{1}{2}\|\mathbf{u}\|^2
\,.
\end{equation}
\paragraph{Legendre transforming to the Hamiltonian side} 
The corresponding Legendre transform yields the following 
invertible relations between momentum and velocity,
\begin{equation}\label{Legendre-dual-rel}
\mathbf{m} =  Q_{op}\mathbf{u}
\quad\hbox{and}\quad\
\mathbf{u} = G*\mathbf{m}
\,,
\end{equation}
where $G$ is the {\bf Green's function} for the operator $ Q_{op}$, assuming
appropriate boundary conditions (on $\mathbf{u}$) that allow inversion of
the operator $ Q_{op}$ to determine $\mathbf{u}$ from $\mathbf{m}$.

The corresponding {\bf Hamiltonian} is,
\begin{equation}\label{Ham-ansatz}
h[\mathbf{m}]
= \langle\mathbf{m}\,,\,\mathbf{u}\rangle
- \frac{1}{2}\|\mathbf{u}\|^2
= \frac{1}{2}\!\int\! \mathbf{m}\cdot G*\mathbf{m}\ d\,^nx
\equiv \frac{1}{2}\|\mathbf{m}\|^2 
\,,
\end{equation}
which also defines a norm $\|\mathbf{m}\|$ via a convolution kernel $G$
that is symmetric and positive, when the Lagrangian $\ell[\mathbf{u}]$ is a
norm. As expected, the norm $\|\mathbf{m}\|$ given by the
Hamiltonian $h[\mathbf{m}]$ specifies the velocity
$\mathbf{u}$ in terms of its Legendre-dual momentum $\mathbf{m}$ by the variational 
operation,
\begin{equation}\label{Ham-u-def}
\mathbf{u} 
= \frac{\delta h}{\delta \mathbf{m}}
= G*\mathbf{m}
\equiv \int G(\mathbf{x}-\mathbf{y})\,\mathbf{m}(\mathbf{y})\,d\,^ny
\,.
\end{equation}
We shall choose the kernel $G(\mathbf{x}-\mathbf{y})$ to be
translation-invariant (so Noether's theorem implies that total momentum
$\mathbf{M}=\int \mathbf{m}\,d\,^nx$ is conserved) and symmetric under
spatial reflections (so that $\mathbf{u}$ and $\mathbf{m}$ have the same
parity). 
 
After the Legendre transformation (\ref{Ham-ansatz}), the EPDiff equation
(\ref{EP-eqn}) appears in its equivalent {\bf Lie-Poisson Hamiltonian
form},
\begin{equation}\label{LP-eqn}
\frac{\partial}{\partial t}\mathbf{m}
=
\{\mathbf{m},h\}
=
-\,
{\rm ad}^*_{{\delta h}/{\delta \mathbf{m}}}\mathbf{m}
\,.
\end{equation}
Here the operation $\{\cdot\,,\,\cdot\,\}$ denotes the Lie-Poisson
bracket dual to the (right) action of vector fields amongst themselves by
vector-field commutation
\[
\{f\,,\,h\,\}
=
-
\left\langle
\mathbf{m}\,,\,\left[
\frac{\delta f}{\delta \mathbf{m}}
\,,\,
\frac{\delta h}{\delta \mathbf{m}}
\right]
\right\rangle
\]
For more details and additional background
concerning the relation of classical EP theory to Lie-Poisson
Hamiltonian equations, see \cite{MaRa1994,HoMaRa1998a}.

In a moment we will also consider the momentum maps for EPDiff.

\subsection{Diffeons: $n-$dimensional Analogs of Pulsons for the EPDiff 
equation} 

The momentum for the one-dimensional pulson solutions (\ref{m-delta}) on
the real line is  supported at points via the Dirac delta measures in its
solution ansatz,
\begin{equation}\label{pulson-ansatz}
m(x,t)=\sum_{i=1}^Np_i(t)\,\delta\big(x-q_i(t)\big)
\,,\quad
m\in{R^1}
\,.
\end{equation} 
We shall develop $n-$dimensional  analogs of these one-dimensional pulson
solutions for the Euler-Poincar\'e equation (\ref{div-grad-curl-eqn}) by
generalizing this solution ansatz to allow measure-valued $n-$dimensional
vector solutions $\mathbf{m}\in{R^n}$ for which the Euler-Poincar\'e
momentum is supported  on co-dimension$-k$ {\it subspaces}
$R^{n-k}$ with integer $k\in[1,n]$. For example, one may consider the 
two-dimensional vector momentum $\mathbf{m}\in{R^2}$ in the plane that is
supported on one-dimensional curves (momentum fronts). Likewise, in three
dimensions, one could consider two-dimensional momentum surfaces (sheets),
one-dimensional momentum filaments, etc. The corresponding vector momentum
ansatz that we shall use is the following, cf. the pulson solutions
(\ref{pulson-ansatz}),
\begin{equation}\label{m-ansatz}
\mathbf{m}(\mathbf{x},t)
=
\sum_{i=1}^N\int\mathbf{P}_i(s,t)\,
\delta\big(\,\mathbf{x}-\mathbf{Q}\,_i(s,t)\,\big)ds
\,,\quad
\mathbf{m}\in{R^n}
\,.
\end{equation} 
Here, $\mathbf{P}_i,\mathbf{Q}_i\in{R^n}$ for $i=1,2,\dots,N$.
For example, when $n-k=1$, so that $s\in R^1$ is one-dimensional, the delta
function in solution (\ref{m-ansatz}) supports an evolving family of
vector-valued curves, called {\bf momentum filaments}. (For
simplicity of notation, we suppress the implied subscript
$i$ in the arclength $s$ for each $\mathbf{P}_i$ and $\mathbf{Q}_i$.) The
Legendre-dual relations (\ref{Legendre-dual-rel}) imply that the velocity
corresponding to the momentum filament ansatz (\ref{m-ansatz}) is,
\begin{equation} \label{u-ansatz}
\mathbf{u}(\mathbf{x},t)
= G*\mathbf{m}
=
\sum_{j=1}^N\int\mathbf{P}_j(s^{\prime},t)\,
G\big(\,\mathbf{x}-\mathbf{Q}\,_j(s^{\prime},t)\,\big)ds^{\prime}
\,.
\end{equation} 
Just as for the 1D case of the pulsons, we shall show that substitution of
the $n-$D solution ansatz (\ref{m-ansatz}) and (\ref{u-ansatz}) into the
EPDiff equation (\ref{ep-eqn-coord-form}) produces canonical geodesic
Hamiltonian equations for the $n-$dimensional vector parameters
$\mathbf{Q}_i(s,t)$ and $\mathbf{P}_i(s,t)$, $i=1,2,\dots,N$.


\subsubsection{Canonical Hamiltonian dynamics of diffeon momentum
filaments in ${R^n}$}

For definiteness in what follows, we shall consider the
example of momentum filaments $\mathbf{m}\in{R^n}$ supported on
one-dimensional space curves in ${R^n}$, so $s\in{R^1}$ is the arclength
parameter of one of these curves. This solution ansatz is reminiscent of the
Biot-Savart Law for vortex filaments, although the flow is not
incompressible. The dynamics of momentum surfaces, for $s\in{R^k}$ with
$k<n$, follow a similar analysis.

Substituting the momentum filament ansatz (\ref{m-ansatz}) for
$s\in{R^1}$ and its corresponding velocity (\ref{u-ansatz}) into the
Euler-Poincar\'e equation (\ref{ep-eqn-coord-form}), then integrating against a
smooth test function $\phi(\mathbf{x})$ implies the following
canonical equations (denoting explicit summation on $i,j\in1,2,\dots N$),
\begin{eqnarray} 
\hspace{-3mm}
\frac{\partial }{\partial t}\mathbf{{Q}}_i (s,t)
\!\!&=&\!\!
\!\!\sum_{j=1}^{N} \int\mathbf{P}_j(s^{\prime},t)\,
G(\mathbf{Q}_i(s,t)-\mathbf{Q}_j(s^{\prime},t))\,\big)ds^{\prime}
\nonumber\\
\!\!&=&\!\! 
\frac{\delta H_N}{\delta \mathbf{P}_i}
\,,\label{IntDiffEqn-Q}\\
\hspace{-3mm}
\frac{\partial }{\partial t}\mathbf{{P}}_i (s,t)
\!\!&=&\!\! 
-\,\!\!\sum_{j=1}^{N} \int 
\big(\mathbf{P}_i(s,t)\!\cdot\!\mathbf{P}_j(s^{\prime},t)\big)
\, \frac{\partial }{\partial \mathbf{Q}_i(s,t)}
G\big(\mathbf{Q}_i(s,t)-\mathbf{Q}_j(s^{\prime},t)\big)\,ds^{\prime}
\nonumber\\
\!\!&=&\!\! 
-\,\frac{\delta H_N}{\delta \mathbf{Q}_i}
\,,\quad\hbox{(sum on $j$, no sum on $i$)}
\,.
\label{IntDiffEqn-P}
\end{eqnarray}
The dot product $\mathbf{P}_i\cdot\mathbf{P}_j$ denotes the inner, or
scalar, product of the two vectors $\mathbf{P}_i$ and $\mathbf{P}_j$ in
$R^n$. Thus, the solution ansatz (\ref{m-ansatz}) yields a closed set of
{\bfi integro-partial-differential equations (IPDEs)} given by 
(\ref{IntDiffEqn-Q}) and (\ref{IntDiffEqn-P}) for the vector parameters
$\mathbf{Q}_i(s,t)$ and $\mathbf{P}_i(s,t)$ with $i=1,2\dots N$. These
equations are generated canonically by the following Hamiltonian function
$H_N:(R^n\times R^n)^{\otimes N}\to R$,
\begin{equation} \label{H_N-def}
H_N = \frac{1}{2}\!\int\!\!\!\!\int\!\!\sum_{i\,,\,j=1}^{N} 
\big(\mathbf{P}_i(s,t)\cdot\mathbf{P}_j(s^{\prime},t)\big) 
\,G\big(\mathbf{Q}_i(s,t)-\mathbf{Q}_{\,j}(s^{\prime},t)\big)
\,ds\,ds^{\prime}
\,.
\end{equation}
This Hamiltonian arises by substituting the momentum ansatz
(\ref{m-ansatz}) into the Hamiltonian (\ref{Ham-ansatz}) obtained from the
Legendre transformation of the Lagrangian corresponding to the kinetic
energy norm of the fluid velocity. Thus, the evolutionary IPDE system
(\ref{IntDiffEqn-Q}) and (\ref{IntDiffEqn-P}) represents canonically
Hamiltonian geodesic motion on the space of curves in $R^n$ with respect to
the co-metric given on these curves in (\ref{H_N-def}). The Hamiltonian
$H_N=\frac{1}{2}\|\mathbf{P}\|^2$ in (\ref{H_N-def}) defines the norm
$\|\mathbf{P}\|$ in terms of this co-metric that combines
convolution using the Green's function $G$ and sum over filaments with the
scalar product of momentum vectors in $R^n$. 

\begin{remark} Note the Lagrangian property of the $s$ coordinate, since
\[
\frac{\partial }{\partial t}\mathbf{{Q}}_i (s,t)
=
\mathbf{u}(\mathbf{Q}_i (s,t),t) 
\,.
\]

\end{remark}

\section{Singular solution momentum map $\mathbf{J}_{\rm Sing}$ for diffeons}
\label{sing-mom-map-sec}
The diffeon momentum filament ansatz (\ref{m-ansatz}) reduces, and {\bfi
collectivizes} the solution of the geodesic EP PDE  (\ref{ep-eqn-coord-form})
in $n+1$ dimensions into the system (\ref{IntDiffEqn-Q}) and
(\ref{IntDiffEqn-P}) of $2N$ canonical evolutionary IPDEs. One can
summarize the mechanism by which this process occurs, by saying that the
map that implements the canonical $(\mathbf{Q},\mathbf{P})$ variables in
terms of singular solutions is a (cotangent bundle) momentum  map. Such
momentum maps are Poisson maps; so the canonical Hamiltonian
nature of the dynamical equations for
$(\mathbf{Q}, \mathbf{P})$ fits into a general theory which also
provides a framework for suggesting other avenues of
investigation.

\begin{theorem}\label{mom-map}
The momentum ansatz \textup{(\ref{m-ansatz})} for
measure-valued solutions of the \textup{EPDiff} equation
\textup{(\ref{ep-eqn-coord-form})}, defines an equivariant
momentum map
\[
\mathbf{J}_{\rm Sing}: T ^{\ast} \operatorname{Emb}(S, \mathbb{R}^n)
\rightarrow
\mathfrak{X}(\mathbb{R}^n)^{\ast}
\]
that is called the {\bfi singular solution momentum map} in
\cite{HoMa2004}.
\end{theorem}

We shall explain the notation used in the theorem's statement in the
course of its proof. Right away, however, we note that the sense
of ``defines'' is that the momentum solution ansatz
(\ref{m-ansatz}) expressing $\mathbf{m}$ (a vector function of
spatial position $\mathbf{x}$) in terms of
$\mathbf{Q}, \mathbf{P}$ (which are functions of $s$)
can be regarded as a map from the space of $(\mathbf{Q}(s),
\mathbf{P} (s))$ to the space of $\mathbf{m}$'s. This will turn
out to be the Lagrange-to-Euler map for the fluid description of
the singular solutions. 
\medskip

Following \cite{HoMa2004}, we shall give two proofs of this result
from two rather different viewpoints. The first proof below uses the
formula for a momentum map for a cotangent lifted action, while
the second proof focuses on a Poisson bracket computation. Each
proof also explains the context in which one has a momentum map.
(See \cite{MaRa1994} for general background on momentum maps.)
\medskip

\noindent{\bf First Proof.} For simplicity and without loss of
generality, let us take $N = 1$ and so suppress the index $a$.
That is, we shall take the case of an isolated singular
solution. As the proof will show, this is not a real
restriction. 
\medskip

To set the notation, fix a
$k$-dimensional manifold
$S$ with a given volume element and whose points are denoted $s
\in S$. Let $\operatorname{Emb}(S, \mathbb{R}^n)$ denote the
set of smooth embeddings $\mathbf{Q}: S \rightarrow \mathbb{R}^n$.
(If the EPDiff equations are taken on a manifold $M$, replace
$\mathbb{R}^n $ with $M$.) Under appropriate technical
conditions, which we shall just treat formally here,
$\operatorname{Emb}(S, \mathbb{R}^n)$ is a smooth manifold. (See,
for example, \cite{EbMa1970} and \cite{MaHu1983} for a discussion
and references.)
\medskip

The tangent space $T _{\mathbf{Q}} \operatorname{Emb}(S,
\mathbb{R}^n)$ to $\operatorname{Emb}(S, \mathbb{R}^n)$ at the
point $\mathbf{Q} \in \operatorname{Emb}(S, \mathbb{R}^n)$ is
given by the space of {\bfi material velocity fields}, namely the
linear space of maps
$\mathbf{V}: S \rightarrow \mathbb{R}^n$ that are vector fields
over the map $\mathbf{Q}$.  The dual space to this space will be
identified with the space of one-form densities over
$\mathbf{Q}$, which we shall regard as maps $\mathbf{P}: S
\rightarrow \left(\mathbb{R}^n\right) ^{\ast}$. In summary, the
cotangent bundle $T ^{\ast} \operatorname{Emb}(S, \mathbb{R}^n)$
is identified with the space of pairs of maps $\left( \mathbf{Q},
\mathbf{P} \right)$.
\medskip

These give us the domain space for the singular solution momentum
map. Now we consider the action of the symmetry group.  Consider
the group $\mathfrak{G} = \operatorname{Diff}$ of diffeomorphisms
of the space $\mathfrak{S}$ in which the EPDiff equations are
operating, concretely in our case $\mathbb{R}^n$. Let it act on
$\mathfrak{S}$ by composition on the {\it left}. Namely for $\eta
\in \operatorname{Diff} (\mathbb{R}^n)$, we let 
\begin{equation} \label{action}
\eta \cdot \mathbf{Q} = \eta \circ \mathbf{Q}.
\end{equation} 
Now lift this action to the cotangent bundle $T ^{\ast} \operatorname{Emb}(S, \mathbb{R}^n)$
in the standard way (see, for instance, {MaRa1994} for this
construction). This lifted action is a symplectic (and hence Poisson)
action and has an equivariant momentum map. {\it We claim that
this momentum map is precisely given by the ansatz
\textup{(\ref{m-ansatz})}.}

To see this, one only needs to recall and then apply the general formula
for the momentum map associated with an action of a general Lie group
$\mathfrak{G}$ on a configuration manifold $Q$ and cotangent lifted to
$T^{\ast}Q$.

First let us recall the general formula. Namely, the momentum map  is
the map $\mathbf{J}: T^{\ast}Q \rightarrow
\mathfrak{g}^\ast$  ($\mathfrak{g}^\ast$ denotes the dual
of the Lie algebra $\mathfrak{g}$ of $\mathfrak{G}$) defined by
\begin{equation} \label{momentummap}
\mathbf{J} (\alpha _q) \cdot \xi = \left\langle \alpha _q, \xi_Q (q)
\right\rangle,
\end{equation}
where $\alpha_q \in T ^{\ast} _q Q $ and $\xi \in \mathfrak{g}$, 
where $\xi _Q $ is the infinitesimal generator of the action of
$\mathfrak{G}$ on $Q$ associated to the Lie algebra element $\xi$,
and where
$\left\langle \alpha _q, \xi_Q (q)
\right\rangle$ is the natural pairing of an element of $T ^{\ast}_q
Q $ with an element of $T _q Q $.  

Now we apply this formula to the special case in which 
the group $\mathfrak{G}$ is the diffeomorphism group
$\operatorname{Diff} (\mathbb{R}^n)$, the manifold $Q$ is
$\operatorname{Emb}(S, \mathbb{R}^n)$ and where the action of
the group on
$\operatorname{Emb}(S, \mathbb{R}^n)$ is given by
(\ref{action}). The sense in which the Lie algebra of
$\mathfrak{G} = \operatorname{Diff}$ is the space
$\mathfrak{g} = \mathfrak{X}$ of vector fields is
well-understood. Hence, its dual is naturally regarded as
the space of one-form densities. The momentum map is thus a map
$\mathbf{J}: T ^{\ast}
\operatorname{Emb}(S, \mathbb{R}^n)
\rightarrow \mathfrak{X}^{\ast}$.
\medskip

With $\mathbf{J}$ given by (\ref{momentummap}), we only need to
work out this formula. First, we shall work out the infinitesimal
generators. Let $X \in \mathfrak{X}$ be a Lie algebra element. By
differentiating the action (\ref{action}) with respect to $\eta$
in the direction of $X$ at the identity element we find that the
infinitesimal generator is given by
\[
X _{\operatorname{Emb}(S, \mathbb{R}^n)} (\mathbf{Q}) 
= X \circ \mathbf{Q}.
\]
Thus, taking $\alpha _q$ to be the cotangent vector
$(\mathbf{Q}, \mathbf{P})$, equation
(\ref{momentummap}) gives
\begin{align*}
\left\langle \mathbf{J} (\mathbf{Q}, \mathbf{P} ),
X \right\rangle & 
= 
\left\langle (\mathbf{Q}, \mathbf{P}),
X \circ \mathbf{Q} \right\rangle \\
& = \int_{S} P _i(s) X ^i (\mathbf{Q}(s)) d^k s .
\end{align*}
On the other hand, note that the right hand side of
(\ref{m-ansatz}) (again with the index $a$ suppressed,
and with $t$ suppressed as well), when paired with the
Lie algebra element $X$ is
\begin{align*}
\left\langle \int _S \mathbf{P}(s)\,
\delta \left( \mathbf{x}-\mathbf{Q}(s) \right) d^k s,
X \right\rangle & 
= 
\int _{\mathbb{R}^n} 
\int _S  \left( P _i (s)\,
\delta \left( \mathbf{x}-\mathbf{Q}(s) \right) d ^k s \right)
X ^i (\mathbf{x}) d ^n x \\
& = \int _S P _i (s) X ^i (\mathbf{Q} (s) d ^k s.
\end{align*}
This shows that the expression given by
(\ref{m-ansatz}) is equal to $\mathbf{J}$ and 
so the result is proved. \quad $\blacksquare$

\medskip
\noindent{\bf Second Proof.} As is standard (see, for example,
{MaRa1994}), one can characterize momentum maps
by means of the following relation, required to hold for all functions
$F$ on $ T^{\ast} \operatorname{Emb}(S, \mathbb{R}^n)$; that is,
functions of $\mathbf{Q}$ and $\mathbf{P}$:
\begin{equation} \label{mom_map_Poisson}
\left\{ F, \left\langle \mathbf{J}, \xi \right\rangle \right\}
=  \xi_P [F]
\,.
\end{equation}
In  our case, we shall take $\mathbf{J}$ to be given by the
solution ansatz and verify that it satisfies this relation.
To do so, let $\xi\in\mathfrak{X}$ so that the left side
of (\ref{mom_map_Poisson}) becomes
\[
\left\{ F,
\int _S P _i (s) \xi\,^i (\mathbf{Q} (s)) d\, ^k s
\right\}
=
\int _S \left[
\frac{\delta F}{\delta Q^i }\xi\,^i (\mathbf{Q} (s))
-
P _i (s)\frac{\delta F}{\delta P_j }
\frac{\delta}{\delta Q^j }
\xi\,^i (\mathbf{Q} (s))
\right]d\, ^k s
\,.
\]
On the other hand, one can directly compute from the definitions
that the infinitesimal generator of the action on the
space
$T ^{\ast}  \operatorname{Emb}(S, \mathbb{R}^n)$ corresponding to the
vector field $\xi^i(\mathbf{x})\frac{\partial }{\partial Q ^i}$ (a
Lie algebra element), is given by (see {MaRa1994}, formula
(12.1.14)):
\[
\delta\mathbf{Q}
= \xi \circ \mathbf{Q}
\,,\quad
\delta\mathbf{P}
=  -\,P _i (s)\frac{\partial}{\partial \mathbf{Q} }
\xi\,^i(\mathbf{Q} (s)),
\]
which verifies that (\ref{mom_map_Poisson}) holds.
\medskip

An important element left out in this proof so far is that it
does not make clear that the momentum map is {\it equivariant}, a condition
needed for the momentum map to be Poisson. The first proof took
care of this automatically since {\bfi momentum maps for cotangent
lifted  actions are always equivariant} and hence are Poisson.
\medskip

Thus, to complete the second proof, we need to check  directly
that the momentum map is equivariant. Actually, we shall only
check that it is infinitesimally invariant by showing that it is
a Poisson map from $T ^{\ast} \operatorname{Emb}(S,
\mathbb{R}^n)$ to the space of $\mathbf{m}$'s (the dual of the
Lie algebra of
$\mathfrak{X}$) with its Lie-Poisson bracket. This sort of
approach to characterize equivariant momentum maps is discussed
in an interesting way in \cite{We2002}.
\medskip

The following direct computation shows that the singular solution
momentum map (\ref{m-ansatz}) is Poisson. This is accomplished by using
the canonical Poisson brackets for $\{\mathbf{P}\},\,\{\mathbf{Q}\}$ and
applying the chain rule to compute
$\big\{m_i(\mathbf{x}),m_j(\mathbf{y})\big\}$, with notation
$\delta\,^\prime_k(\mathbf{y})
\equiv\partial\delta(\mathbf{y})/\partial{y^k}$. We get

\begin{eqnarray}\label{m-m-bracket1}
\big\{m_i(\mathbf{x}),m_j(\mathbf{y})\big\}\hspace{-1in}&&
\nonumber\\
&=&
\bigg\{\sum_{a=1}^N\!\int\!\!ds\, P_i^a(s,t)
\, \delta (\mathbf{x}-\mathbf{Q}^a(s,t))
\,,\,
\sum_{b=1}^N \!\int\!\!ds^\prime P_j^b(s^\prime,t)\,
\delta (\mathbf{y}-\mathbf{Q}^b(s^\prime,t))\bigg\}
\nonumber\\
&=&
\sum_{a,b=1}^N\int\!\!\!\int\!\!dsds^\prime
\bigg[ \{ P_i^a(s),P_j^b(s^\prime) \}
\,\delta (\mathbf{x}-\mathbf{Q}^a(s))
\,\delta (\mathbf{y}-\mathbf{Q}^b(s^\prime))
\nonumber\\
&&-\,\{ P_i^a(s),Q_k^b(s^\prime) \} P_j^b(s^\prime)
\,\delta (\mathbf{x}-\mathbf{Q}^a(s))
\,\delta\,^\prime_k(\mathbf{y}-\mathbf{Q}^b(s^\prime))
\nonumber\\
&&-\,\{ Q_k^a(s),P_j^b(s^\prime) \} P_i^a(s)
\,\delta\,^\prime_k (\mathbf{x}-\mathbf{Q}^a(s))
\,\delta (\mathbf{y}-\mathbf{Q}^b(s^\prime))
\nonumber\\
&&+\,\{ Q_k^a(s),Q_\ell^b(s^\prime) \} P_i^a(s)P_j^b(s^\prime)
\,\delta\,^\prime_k (\mathbf{x}-\mathbf{Q}^a(s))
\,\delta\,^\prime_\ell (\mathbf{y}-\mathbf{Q}^b(s^\prime))
\bigg]
\,.
\nonumber
\end{eqnarray}
Substituting the canonical Poisson bracket relations
\begin{align*}
\{ P_i^a(s),P_j^b(s^\prime) \} & = 0 \\
\{ Q_k^a(s),Q_\ell^b(s^\prime) \} & = 0, \quad \mbox{and} \; \\
\{ Q_k^a(s),P_j^b(s^\prime) \} & =
\delta^{ab}\delta_{kj}\delta(s-s^\prime)
\end{align*}
into the preceding computation yields,
\begin{eqnarray}\label{m-m-bracket2}
\big\{m_i(\mathbf{x}),m_j(\mathbf{y})\big\}\hspace{-1in}&&
\nonumber\\
&=&
\bigg\{\sum_{a=1}^N\!\int\!\!ds P_i^a(s,t)
\, \delta (\mathbf{x}-\mathbf{Q}^a(s,t))
\,,\,
\sum_{b=1}^N \!\int\!\!ds^\prime P_j^b(s^\prime,t)\,
\delta (\mathbf{y}-\mathbf{Q}^b(s^\prime,t))\bigg\}
\nonumber\\
&=&
\sum_{a=1}^N\int\!\!ds P_j^a(s)
\,\delta (\mathbf{x}-\mathbf{Q}^a(s))
\,\delta\,^\prime_i (\mathbf{y}-\mathbf{Q}^a(s))
\nonumber\\
&&-\sum_{a=1}^N\int\!\!ds P_i^a(s)
\,\delta\,^\prime_j (\mathbf{x}-\mathbf{Q}^a(s))
\,\delta (\mathbf{y}-\mathbf{Q}^a(s))
\nonumber\\
&=&
-\,
\Big(
m_j(\mathbf{x})\frac{\partial}{\partial x^i}
+
\frac{\partial}{\partial x\,^j}\,m_i(\mathbf{x})\Big)
\delta(\mathbf{x}-\mathbf{y})
\,.\nonumber
\end{eqnarray}

\noindent
Thus,
\begin{eqnarray}\label{momentum-map-bracket1}
\big\{m_i(\mathbf{x})\,,\,m_j(\mathbf{y})\big\}
=
-\,
\Big(
m_j(\mathbf{x})\frac{\partial}{\partial x^i}
+
\frac{\partial}{\partial x\,^j}\,m_i(\mathbf{x})\Big)
\delta(\mathbf{x}-\mathbf{y})
\,,\label{LP-bracket}
\end{eqnarray}
which is readily checked to be the Lie-Poisson bracket on the space
of $\mathbf{m}$'s, restricted to their singular support. This completes
the second proof of theorem.
\quad
$\blacksquare$
\medskip

Each of these proofs has shown the following basic fact.

\begin{corollary} \label{Poisson_mom-map}
The singular solution momentum map defined by the
singular solution ansatz \textup{(\ref{m-ansatz})}, namely,
\[
\mathbf{J}_{\rm Sing}: 
T ^{\ast} \operatorname{Emb}(S, \mathbb{R}^n)
\rightarrow
\mathfrak{X}(\mathbb{R}^n)^{\ast}
\]
is a Poisson map from the canonical Poisson structure on $T
^{\ast} \operatorname{Emb}(S, \mathbb{R}^n)$ to the Lie-Poisson
structure on $\mathfrak{X}(\mathbb{R}^n)^{\ast}$.
\end{corollary}

This is perhaps the most basic property of the singular solution
momentum map. Some of its more sophisticated properties are
outlined in \cite{HoMa2004}.

\paragraph{Pulling Back the Equations.}
Since the solution ansatz (\ref{m-ansatz}) has been shown
in the preceding Corollary to be a Poisson map, the pull back of
the Hamiltonian from $\mathfrak{X}^{\ast}$ to $T ^{\ast}
\operatorname{Emb}(S, \mathbb{R}^n)$ gives equations of motion on
the latter space that project to the equations on $\mathfrak{X}
^{\ast}$. 
\begin{quotation}
{\it Thus, the basic fact that the momentum map $\mathbf{J}_{\rm Sing}$ 
is Poisson explains why the functions $\mathbf{Q}^a(s,t)$ and
$\mathbf{P}^a(s,t)$  satisfy canonical Hamiltonian equations.}
\end{quotation}
Note that the coordinate $s\in{\mathbb{R}}^{k}$ that labels these
functions is a ``Lagrangian coordinate'' in the sense that it does
not evolve in time but rather labels the solution. 

In terms of the pairing
\begin{equation}
\langle\cdot\,,\,\cdot\rangle:
\,\mathfrak{g}^*\times\mathfrak{g}\to{\mathbb{R}}
\,,
\end{equation}
between the Lie algebra $\mathfrak{g}$ (vector fields in $\mathbb{R}^n$)
and its dual $\mathfrak{g}^*$  (one-form densities in $\mathbb{R}^n$), the
following relation holds for measure-valued solutions
under the momentum map (\ref{m-ansatz}),
\begin{align}
\label{momentum-map-relation}
\langle\mathbf{m}\,,\,\mathbf{u}\rangle
&= 
\int \mathbf{m}\,\cdot\,\mathbf{u}\,d\,^n\mathbf{x}
\,,\quad\hbox{$L^2$ pairing for }
\mathbf{m}\,\&\,\mathbf{u}\in{\mathbb{R}^n},
\nonumber\\
&= 
\!\int\!\!\!\!\int\!\!\sum_{a\,,\,b=1}^{N}
\big(\mathbf{P}^a(s,t)\cdot\mathbf{P}^b(s^{\prime},t)\big)
\,G\big(\mathbf{Q}^a(s,t)-\mathbf{Q}^{\,b}(s^{\prime},t)\big)
\,ds\,ds^{\prime}
\nonumber\\
&= 
\!\int\!\!\sum_{a=1}^{N}
\mathbf{P}^a(s,t)\cdot\frac{\partial\mathbf{Q}^a(s,t)}{\partial t}
\,ds
\nonumber\\
&\equiv 
\langle\!\langle\mathbf{P}\,,\,\mathbf{\dot{Q}}\rangle\!\rangle,
\end{align}
which is the natural pairing between the points $(\mathbf{Q},
\mathbf{P})
\in  T^{\ast} \operatorname{Emb}(S, \mathbb{R}^n)$ and $(\mathbf{Q},
\dot{\mathbf{Q}}) \in T \operatorname{Emb}(S, \mathbb{R}^n)$.  This 
corresponds to preservation of the action of the Lagrangian $\ell[\mathbf{u}]$ under
cotangent lift of $\operatorname{Diff} (\mathbb{R}^n)$.

The pull-back of the Hamiltonian  $H[\mathbf{m}]$ defined on the dual
of the Lie algebra $\mathfrak{g}^*$, to $ T^{\ast}
\operatorname{Emb}(S, \mathbb{R}^n)$ is easily seen to be consistent
with what we had before:
\begin{equation}
H[\mathbf{m}]
\equiv
\frac{1}{2}\langle\mathbf{m}\,,\,G*\mathbf{m}\rangle
=
\frac{1}{2}\langle\!\langle\mathbf{P}\,,\,G*\mathbf{P}\rangle\!\rangle
\equiv
H_N[\mathbf{P},\mathbf{Q}]
\,.
\label{geodesic-ham}
\end{equation}

In summary, in concert with the Poisson nature of the singular
solution momentum map, we see that the singular solutions in
terms of $\mathbf{Q}$ and $\mathbf{P}$ satisfy Hamiltonian
equations and also define an invariant solution set for the
EPDiff equations. In fact, 
\begin{quotation}\noindent
{\it This invariant solution set is a
special coadjoint orbit for the diffeomorphism group, as we
shall discuss in the next section.}
\end{quotation}

\section{The Geometry of the Momentum Map} \label{geom_mommap}

In this section we explore the geometry of the singular
solution momentum map discussed earlier in
a little more detail. The treatment is formal, in the sense that 
there are a number of
technical issues in the infinite dimensional case that will be left
open. We will mention a few of these as we proceed.

\subsection{Coadjoint Orbits.} We claim  that {\it the
image of the singular solution momentum map is a coadjoint orbit
in $\mathfrak{X} ^{\ast}$.}  This means that (modulo some issues of
connectedness and smoothness, which we do not consider here) the
solution ansatz given by  (\ref{m-ansatz}) defines a
coadjoint orbit in the space of all one-form densities, regarded as
the dual of the Lie algebra of the diffeomorphism group. 
These coadjoint orbits should be thought of as singular
orbits---that is, due to their special nature, they are not
generic. 

Recognizing them as coadjoint orbits is one way of
gaining further insight into why the singular solutions form
dynamically invariant sets---it is a general fact that coadjoint
orbits in $\mathfrak{g}^\ast$ are {\it symplectic submanifolds} of
the Lie-Poisson manifold $\mathfrak{g}^\ast$ (in our case
$\mathfrak{X}(\mathbb{R}^n)^{\ast}$) and, correspondingly, are dynamically invariant
for any Hamiltonian system on $\mathfrak{g}^\ast$.

The idea of the proof of our claim is simply this: whenever
one has an equivariant momentum map $\mathbf{J}: P \rightarrow
\mathfrak{g}^\ast$ for the action of a group $G$ on a symplectic
or Poisson manifold $P$, and that action is transitive, then the
image of $\mathbf{J}$ is an orbit (or at least a piece of an
orbit). This general result, due to Kostant, is stated more
precisely in \cite{MaRa1994}, Theorem 14.4.5. Roughly speaking, the
reason that transitivity holds in our case is because one can ``move
the images of the manifolds $S$ around at will with arbitrary velocity fields'' using
diffeomorphisms of $\mathbb{R}^n$. 

\subsection{The Momentum map $\mathbf{J}_S$ and the Kelvin
circulation theorem.}
The momentum map $\mathbf{J}_{\rm Sing}$  involves
$\operatorname{Diff}(\mathbb{R}^n)$, the left action of the
diffeomorphism group on the space of embeddings
$\operatorname{Emb}(S,\mathbb{R}^n)$ by smooth maps of the target
space $\mathbb{R}^n$, namely,
\begin{equation}\label{leftDiff}
\operatorname{Diff}(\mathbb{R}^n):\
\mathbf{Q}\cdot\eta=\eta\circ\mathbf{Q},
\end{equation}
where, recall, $\mathbf{Q}:S\to \mathbb{R}^n$.
As above, the cotangent bundle $T ^{\ast}
\operatorname{Emb}(S,\mathbb{R}^n)$ is identified with the space of
pairs of maps $( \mathbf{Q}, \mathbf{P})$, with
$\mathbf{Q}:S\to \mathbb{R}^n$ and $\mathbf{P}:S\to T^*\mathbb{R}^n$.  

However, there is another momentum map $\mathbf{J}_S$  associated
with the {\it right action} of the diffeomorphism group of $S$ on
the embeddings
$\operatorname{Emb}(S,\mathbb{R}^n)$ by smooth maps of the
``Lagrangian labels'' $S$ (fluid particle relabeling by $\eta:S\to
S$). This action is given by
\begin{equation}\label{rightDiff}
\operatorname{Diff}(S):\
\mathbf{Q}\cdot\eta=\mathbf{Q}\circ\eta
\,.
\end{equation}
The infinitesimal generator of this right action is 
\begin{equation}\label{infgen-right}
X_{\operatorname{Emb}(S,\mathbb{R}^n)}(\mathbf{Q})
=
\frac{d}{dt}\Big|_{t=0}\mathbf{Q}\circ\eta_t
=
T\mathbf{Q} \circ X.
\end{equation}
where $X \in \mathfrak{X}$ is tangent to the curve $\eta _t$ at $t
 = 0$. Thus, again taking $N = 1$ (so we suppress the index $a$) and
also letting $\alpha _q$ in the momentum map formula
(\ref{momentummap}) be the cotangent vector
$(\mathbf{Q}, \mathbf{P})$, one computes $\mathbf{J}_S$:
\begin{align*}
\left\langle \mathbf{J}_S (\mathbf{Q}, \mathbf{P} ),
X \right\rangle 
& =  
\left\langle (\mathbf{Q}, \mathbf{P}),
T\mathbf{Q}\cdot X \right\rangle
\\
& =  
\int_S
P _i(s) 
\frac{\partial Q^i(s)}{\partial s^m}
X ^m (s)
\,d\,^ks
\\
& =  
\int_S
X\Big(
\mathbf{P}(s) 
\cdot d\mathbf{Q}(s)\Big)
\,d\,^ks
\\
& =  
\left(
\int_S
\mathbf{P}(s) 
\cdot d\mathbf{Q}(s)
\otimes\,d\,^ks
\,,
X (s) \right)
\\
& = 
\langle\,
\mathbf{P}
\cdot d\mathbf{Q}
\,,
X \,\rangle
\,.
\end{align*}
Consequently, the momentum map formula (\ref{momentummap}) yields
\begin{equation} \label{momentummap-Js}
\mathbf{J}_S (\mathbf{Q}, \mathbf{P} ) 
= 
\mathbf{P}\cdot d\mathbf{Q}
\,,
\end{equation}
with the indicated pairing of the one-form density $\mathbf{P}
\cdot d\mathbf{Q}$ with the vector
field $X$. 
\medskip
We have set things up so that the following is true.
\begin{proposition} The momentum map $\mathbf{J}_S$ is preserved
by the evolution equations (\ref{IntDiffEqn-Q}-\ref{IntDiffEqn-P}) for
$\mathbf{Q}$ and
$\mathbf{P}$.
\end{proposition}

\begin{proof} It is enough to notice that the Hamiltonian $H _N$
in equation \eqref{H_N-def}  is invariant under the cotangent lift
of the action of $\operatorname{Diff} (S)$; it merely amounts to the
invariance of the integral over $S$ under reparametrization;
that is, the change of variables formula; keep in mind that
$\mathbf{P}$ includes a density factor.
\end{proof}

\begin{remark}$\quad$

\begin{itemize}
\item
This result is similar to the Kelvin-Noether theorem for circulation
$\Gamma$ of an ideal fluid, which may be written as 
$\Gamma=\oint_{c(s)}D(s)^{-1}\mathbf{P}(s)\cdot d\mathbf{Q}(s)$
for each Lagrangian circuit $c(s)$, where $D$ is the mass
density and $\mathbf{P}$ is again the canonical momentum density. 
This similarity should come as no surprise, because the
Kelvin-Noether theorem for ideal fluids arises from invariance of
Hamilton's principle under fluid parcel relabeling
by the {\it same} right action of the diffeomorphism group, as
in (\ref{rightDiff}).  

\item
Note that, being an equivariant momentum map, the map $\mathbf{J}_S$, as with
$\mathbf{J}_{\rm Sing}$,  is also a Poisson map. That is, substituting the canonical
Poisson bracket into relation \eqref{momentummap-Js}; that is, the relation
$\mathbf{M}(\mathbf{x}) = \sum_iP_i(\mathbf{x})\nabla Q^i(\mathbf{x})$ yields the
Lie-Poisson bracket on the space of $\mathbf{M}$'s. We use the different notations
$\mathbf{m}$ and $\mathbf{M}$ because these quantities are analogous to the body and
spatial angular momentum for rigid body mechanics. In fact, the quantity $\mathbf{m}$ given
by the solution Ansatz; specifically,
$\mathbf{m} = \mathbf{J}_{\rm Sing}(\mathbf{Q}, \mathbf{P})$ gives the singular solutions
of the EPDiff equations, while $\mathbf{M}(\mathbf{x}) = \mathbf{J}_S(\mathbf{Q},
\mathbf{P}) = \sum_iP_i(\mathbf{x})\nabla Q^i(\mathbf{x})$ is a conserved quantity.

\item
In the language of fluid mechanics, the expression of $\mathbf{m}$ in
terms of $(\mathbf{Q}, \mathbf{P})$ is an example of a   {\bfi Clebsch
representation}, which expresses the solution of the EPDiff equations in
terms of canonical variables that evolve by standard canonical Hamilton
equations.  This has been known in the case of fluid mechanics for more
than 100 years.  For modern discussions of the Clebsch representation for
ideal fluids, see, for example, \cite{HoKu1983,MaWe1983}. 

\item
One more remark is in order; namely the special case in which $S = M $ is
of course allowed. In this case, $ \mathbf{Q}$ corresponds to the map $\eta$ itself and
$\mathbf{P}$ just corresponds to its conjugate momentum. The quantity $\mathbf{m}$
corresponds to the spatial (dynamic) momentum density (that is, right translation of
$\mathbf{P}$ to the identity), while $\mathbf{M} $ corresponds to the conserved ``body''
momentum density (that is, left translation of $\mathbf{P}$ to the identity).
\end{itemize}
\end{remark}

\subsection{Brief summary}

$\operatorname{Emb}(S,\mathbb{R}^n)$ admits two group actions. These are: the group 
$\operatorname{Diff} (S)$ of diffeomorphisms of $S$, which acts by composition on
the {\it right}; and the group $\operatorname{Diff}(\mathbb{R}^n)$ which acts 
by composition on the {\it left}. The group
$\operatorname{Diff}(\mathbb{R}^n)$ acting from the left produces the singular solution
momentum map, $\mathbf{J}_{\rm Sing}$. The action of $\operatorname{Diff} (S)$ from the
right produces the conserved momentum map $\mathbf{J}_S: T ^{\ast} 
\operatorname{Emb}(S, \mathbb{R}^n) \rightarrow
\mathfrak{X}(S)^{\ast}$. We now assemble both momentum maps into
one figure as follows: 

\begin{picture}(150,100)(-50,0)%
\put(100,75){$T^{\ast} \operatorname{Emb}(S,M)$} 

\put(78,50){$\mathbf{J}_{\rm Sing}$}        

\put(160,50){$\mathbf{J}_S$}   

\put(72,15){$\mathfrak{X} (M)^{\ast}$}       

\put(170,15){$\mathfrak{X}(S)^{\ast}$}       

\put(130,70){\vector(-1, -1){40}}  

\put(135,70){\vector(1,-1){40}}  

\end{picture}

\section{The Euler-Poincar\'e framework: fluids \`a la
\cite{HoMaRa1998a}}
\label{sec-EPframe}

Almost all fluid models of interest admit the following
general assumptions. These assumptions form the basis of the
Euler-Poincar\'e theorem for Continuua that we shall state later in this
section, after introducing the notation necessary for dealing
geometrically with the reduction of Hamilton's Principle from the material
(or Lagrangian) picture of fluid dynamics, to the spatial (or Eulerian)
picture. This theorem was first stated and proved in \cite{HoMaRa1998a}, to
which we refer for additional details, as well as for abstract definitions
and proofs.

\paragraph{Basic assumptions underlying the Euler-Poincar\'e theorem for
continua}
\begin{itemize}
\item There is a {\it right\/} representation of a Lie group $G$ on
the vector space $V$ and $G$ acts in the natural way on the {\it
right\/} on $TG \times V^\ast$: $(U_g, a)h = (U_gh, ah)$.
\item The Lagrangian function $ L : T G \times V ^\ast
\rightarrow \mathbb{R}$ is right $G$--invariant under the isotropy group of $a_0\in V ^\ast$.%
\footnote{For fluid dynamics, right $G$--invariance of the Lagrangian
function $L$ is traditionally called ``particle relabeling symmetry.''}
\item In particular, if $a_0 \in V^\ast$, define the
Lagrangian $L_{a_0} : TG \rightarrow \mathbb{R}$ by
$L_{a_0}(U_g) = L(U_g, a_0)$. Then $L_{a_0}$ is right
invariant under the lift to $TG$ of the right action of
$G_{a_0}$ on $G$, where $G_{a_0}$ is the isotropy group of $a_0$.
\item  Right $G$--invariance of $L$ permits one to define the Lagrangian
on the Lie algebra $\mathfrak{g}$ of the group $G$. Namely,
$\ell: {\mathfrak{g}} \times V^\ast \rightarrow \mathbb{R}$
is defined by, 
\[
\ell({u},a)
=
L\big(U_gg^{-1}(t), a_0g^{-1}(t)\big) 
= L(U_g, a_0)
\,,
\]
where
$
{u}=U_gg^{-1}(t)
$ and
$
a = a_0g^{-1}(t)
\,.
$
Conversely,  this relation defines for any
$\ell: {\mathfrak{g}} \times V^\ast \rightarrow
\mathbb{R} $ a right $G$--invariant function
$ L : T G \times V ^\ast
\rightarrow \mathbb{R} $.
\item For a curve $g(t) \in G, $ let
$u  (t) := \dot{g}(t) g(t)^{-1}$ and define the curve
$a(t)$ as the unique solution of the linear differential equation
with time dependent coefficients $\dot a(t) = -a(t)u (t)$, where the action
of an element of the Lie algebra $u\in\mathfrak{g}$ on an advected quantity
$a\in V^*$ is denoted by concatenation from the right.
The solution with initial condition $a(0)=a_0\in V^*$ can be
written as $a(t) = a_0g(t)^{-1}$.
\end{itemize}

\subsubsection*{Notation for reduction of Hamilton's Principle by
symmetries}
\begin{itemize}
\item
Let $\mathfrak{g}(\mathcal{D})$ denote the space of
vector fields on $\mathcal{D}$ of some fixed differentiability class.
These vector fields are endowed with the {\bfi  Lie bracket\/} given in
components by (summing on repeated indices)
\begin{equation}\label{jlb} [\mathbf{u}, \mathbf{v}]^i
= 
u^j\frac{\partial v^i}{\partial x^j} 
 - 
v^j\frac{\partial u^i}{\partial x^j}
\,.
\end{equation}
The notation $
\operatorname{ad}_\mathbf{u} \mathbf{v} := [\mathbf{u},\,\mathbf{v}]$
formally denotes the adjoint action of the {\it right\/} Lie algebra of
$\operatorname{Diff}(\mathcal{D})$ on itself.

\item
Identify the Lie algebra of vector fields $\mathfrak{g}$ with its 
dual $\mathfrak{g}^\ast$ by using the $L^2$ pairing
\begin{equation}\label{l2p}
\left\langle \mathbf{ u},\mathbf{v}\right\rangle
=\int_\mathcal{D} \mathbf{ u}
\cdot
\mathbf{v}
\,  dV\,.
\end{equation}
%
\item
Let $\mathfrak{g}(\mathcal{D})^\ast$ denote the geometric dual
space of $\mathfrak{g}(\mathcal{D})$, that is,
$\mathfrak{g}(\mathcal{D})^\ast := \Lambda^1({\mathcal{D}}) \otimes
{\rm Den}({\mathcal{D}})$. This is the space of one--form densities on
$\mathcal{D}$. If $\mathbf{m}\, \otimes  dV \in
\Lambda^1({\mathcal{D}}) \otimes {\rm Den}({\mathcal{D}})$, then the
pairing of  $\mathbf{m} \otimes  dV$ with $\mathbf{u} \in
\mathfrak{g}(\mathcal{D})$ is given by the $L^2$ pairing,
\begin{equation}\label{continuumpairing}
\langle \mathbf{m} \otimes  dV, \mathbf{u} \rangle
= \int_{\mathcal{D}} \mathbf{m}\cdot \mathbf{u}\,  dV
\end{equation}
where $\mathbf{m}\cdot \mathbf{u}$ is the standard contraction of a
one--form {\bfi m} with a vector field $\mathbf{u}$. 
\item
For $\mathbf{u} \in \mathfrak{g}(\mathcal{D})\,$ and $\mathbf{m}\,\otimes 
dV \in \mathfrak{g}(\mathcal{D})^\ast$, the dual of the adjoint
representation is defined by
\begin{equation}\label{ad-star-def}
 \langle  \operatorname{ad}^\ast_\mathbf{u}(\mathbf{m}\otimes  dV),
\mathbf{v}
\rangle
= -\int_{\mathcal{D}}\mathbf{m} \cdot
\operatorname{ad}_\mathbf{u}\!\mathbf{v}\,dV 
= -\int_{\mathcal{D}}\mathbf{m}
\cdot [\mathbf{u}, \mathbf{v}]\,dV
\end{equation}
and its expression is
\begin{equation}\label{continuumcoadjoint}
 \operatorname{ad}^\ast_\mathbf{u}(\mathbf{m}\otimes  dV) 
= (\pounds_\mathbf{u}
\mathbf{m} + ( \operatorname{div}_{dV}\mathbf{u})\mathbf{m})\otimes  dV
= \pounds_\mathbf{u}(\mathbf{m}\otimes  dV)\,,
\end{equation}
where  ${\rm div}_{dV}\mathbf{u}$ is the divergence of $\mathbf{u}$ relative
to the measure $dV$, that is, $\pounds_\mathbf{u}dV = ({\rm
div}_{dV}\mathbf{u})dV$. Hence,
$\operatorname{ad}^\ast_\mathbf{u}$ coincides with the Lie-derivative
$\pounds_\mathbf{u}$ for one-form densities. 
\item
If $\mathbf{u} = u^j \partial/\partial x^j,\,
\mathbf{m}  = m_i dx^i$, then the one--form factor in the preceding
formula for $ \operatorname{ad}^\ast_\mathbf{u}(\mathbf{m}\otimes  dV)$ has the
{\bfi coordinate expression}
\begin{eqnarray}\label{continuumcoadjoint-coord}
\Big(\operatorname{ad}^\ast_\mathbf{u}\mathbf{m}\Big)_{\!i} dx^i
&=&
 \left ( u^j \frac{\partial
m_i}{\partial
x^j} + m_j \frac{\partial u^j}{\partial x^i} +
( \operatorname{div}_{dV}\mathbf{u})m_i \right )dx^i
\\
&=&
\left (\frac{\partial}{\partial x^j}(u^jm_i) +
m_j \frac{\partial u^j}{\partial x^i}\right ) dx^i
\,.
\end{eqnarray}
The last equality assumes that the divergence is taken
relative to the standard measure $dV = d^n\mathbf{x}$ in $\mathbb{R} ^n$.
(On a Riemannian manifold the metric divergence needs to be used.)
\end{itemize}

\subsubsection*{Conventions and terminology in continuum mechanics}
Throughout the rest of the lecture notes, we shall follow
\cite{HoMaRa1998a}  in using the conventions and terminology for the
standard quantities in continuum mechanics. 

\begin{definition}
Elements of $\mathcal{D}$ representing the material particles of the system
are denoted by $X$; their coordinates $X^A, A=1,...,n$ may thus be
regarded as the {\bfi particle labels}. 
\begin{itemize}
\item
A {\bfi configuration}, which we typically denote by
$\eta$, or $g$, is an element of $\operatorname{Diff}(\mathcal{D})$.
\item
A {\bfi motion\/}, denoted as $\eta_t$ or alternatively as $g(t)$, is a time
dependent curve in $\operatorname{Diff}(\mathcal{D})$. 
\end{itemize}

\end{definition}

\begin{definition}
The {\bfi Lagrangian}, or {\bfi material velocity\/} ${\bf U}(X,t)$ of the
continuum along the motion
$\eta_t$ or $g(t)$ is defined by taking the time derivative of the motion
keeping the particle labels $X$ fixed:
\[
{\bf U}(X, t) := \frac{d\eta_t(X)}{dt}:=
\left.\frac{\partial}{\partial t}\right|_{X}\eta_t(X)
:= \dot{g}(t)\cdot X
\,.
\]
These are convenient shorthand notations for the
time derivative at fixed Lagrangian coordinate $X$.

Consistent with this definition of material velocity, the tangent space to
$\operatorname{Diff}(\mathcal{D})$ at $\eta \in
\operatorname{Diff}(\mathcal{D})$ is given by
\[
T_\eta \operatorname{Diff}(\mathcal{D})
= \{ {\bf U}_\eta: {\mathcal{D}} \rightarrow T {\mathcal{D}}  \mid {\bf U}_\eta(X)
\in T_{\eta(X)}\mathcal{D} \}.
\]
Elements of $T_\eta \operatorname{Diff}(\mathcal{D})$ are usually
thought of as vector fields on $\mathcal{D}$ covering $\eta$. The
tangent lift of right translations on
$T\operatorname{Diff}(\mathcal{D})$ by $\varphi \in
\operatorname{Diff}(\mathcal{D})$
is given by
\[
{\bf U}_\eta\varphi := T_\eta R_\varphi ({\bf U}_\eta)
= {\bf U}_\eta \circ \varphi\,.
\]
\end{definition}

\begin{definition}
During a motion $\eta_t$ or $g(t)$, the particle labeled by $X$ describes a
path in $\mathcal{D}$, whose points 
\[x(X, t):= \eta_t(X):=g(t)\cdot X\,,\]
are called the {\bfi Eulerian} or {\bfi spatial points\/} of this path,
which is also called the {\bfi Lagrangian trajectory\/}, because a
Lagrangian fluid parcel follows this path in space. The derivative
$\mathbf{u}(x, t)$ of this path, evaluated at fixed Eulerian point $x$, is
called the {\bfi Eulerian} or {\bfi spatial velocity\/} of the system:
\[
\mathbf{u}(x, t)
:= \mathbf{u}(\eta_t(X), t)
:= {\bf U}(X, t)
:= \left.\frac{\partial}{\partial t}\right|_X\eta_t(X)
:= \dot{g}(t) \cdot X
:= \dot{g}(t)g^{-1}(t)\cdot x
\,.
\]
Thus the Eulerian velocity $\mathbf{u}$ is a time dependent vector
field on $\mathcal{D}$, denoted as $\mathbf{u}_t \in
\mathfrak{g}(\mathcal{D})$, where $\mathbf{u}_t(x) := \mathbf{u}(x, t)$.
We also have the fundamental relationships
\[
{\bf U}_t = \mathbf{u}_t \circ \eta_t
\quad\hbox{and}\quad
\mathbf{u}_t=\dot{g}(t)g^{-1}(t)
\,,
\]
where we denote ${\bf U}_t(X):= {\bf U}(X, t)$.
\end{definition}

\begin{definition}
The {\bfi representation space} $V^\ast$ of
$\operatorname{Diff}(\mathcal{D})$ in continuum mechanics is often some
subspace of the tensor field densities  on $\mathcal{D}$, denoted as
$\mathfrak{T} (\mathcal{D})\otimes {\rm Den}({\mathcal{D}})$, and the
representation is given by pull back. It is thus a {\it right\/}
representation of
$\operatorname{Diff}(\mathcal{D})$ on
$\mathfrak{T}(\mathcal{D})\otimes {\rm Den}({\mathcal{D}})$. The right
action of the Lie algebra $\mathfrak{g}({\mathcal{D}})$ on $V^\ast$ is
denoted as {\bfi concatenation from the right}. That is, we denote 
\[
a\mathbf{u} :=
\pounds_\mathbf{u} a
\,,
\]
which is the Lie derivative of the tensor field density
$a$ along the vector field $\mathbf{u}$.
\end{definition}

\begin{definition}
The {\bfi Lagrangian of a continuum mechanical system} is a function 
\[
L:
T\operatorname{Diff}(\mathcal{D}) \times V^\ast \rightarrow \mathbb{R} 
\,,
\]
which is right invariant relative to the tangent lift of right
translation of $\operatorname{Diff}(\mathcal{D})$ on itself and pull
back on the tensor field densities.
Invariance of the Lagrangian $L$ induces
a function $\ell: \mathfrak{g}(\mathcal{D}) \times V^\ast
\rightarrow \mathbb{R} $ given by
\[
\ell(\mathbf{u}, a) 
= L(\mathbf{u}\circ \eta, \eta^\ast a)
= L(\mathbf{U}, a_0)
\,,
\]
where $  \mathbf{u} \in \mathfrak{g}({\mathcal{D}})$ and $  a \in V^\ast
\subset {\mathfrak{T}}({\mathcal{D}})\otimes {\rm Den}({\mathcal{D}})$,
and where $\eta^\ast a$ denotes the pull back of $a$ by the
diffeomorphism $\eta$ and $\mathbf{u}$ is the Eulerian velocity. 
That is, 
\begin{equation}\label{fund-relns}
\mathbf{U}=\mathbf{u}\circ \eta
\quad\hbox{and}\quad
a_0 = \eta^\ast a
\,.
\end{equation}
The evolution of $a$ is by right action, given by the equation
\begin{equation}\label{eqn-adv-qty}
\dot a = -\,{\pounds}_\mathbf{u}\, a = -\, a\mathbf{u}.
\end{equation}
The solution of this equation, for the
initial condition $a_0$, is 
\begin{equation}\label{soln-adv-qty}
a(t) = \eta_{t\ast} a_0=a_0g^{-1}(t)
\,,
\end{equation}
where the lower star denotes the push forward operation and
$\eta_t$ is the flow of $\mathbf{u}=\dot{g}g^{-1}(t)$.
\end{definition}

\begin{definition}
{\bfi Advected Eulerian quantities} are defined in continuum mechanics to
be those variables which are Lie transported by the flow of the Eulerian
velocity field. Using this standard terminology, 
equation (\ref{eqn-adv-qty}), or its solution (\ref{soln-adv-qty}) states
that the tensor field density $a(t)$ (which may include mass density and
other Eulerian quantities) is advected.
\end{definition}

\begin{remark}[Dual tensors]
As we mentioned, typically $V^\ast \subset
{\mathfrak{T}}({\mathcal{D}})\otimes {\rm Den}({\mathcal{D}})$ for
continuum mechanics.  On a general manifold, tensors of a given type have
natural duals. For example, symmetric covariant tensors are dual to
symmetric contravariant tensor densities, the pairing being given by the
integration of the natural contraction of these tensors. Likewise,
$k$--forms are naturally dual to
$(n-k)$--forms, the pairing being given by taking the integral of
their wedge product.
\end{remark}

\begin{definition}
The {\bfi diamond operation} $\diamond$ between elements of $V$ and
$V^\ast$  produces an element of the dual Lie algebra
$\mathfrak{g}({\mathcal{D}})^\ast$ and is defined as 
\begin{equation}\label{continuumdiamond}
\langle b \diamond a, \mathbf{w}\rangle
= -\int_{\mathcal{D}} b \cdot \pounds_\mathbf{w}\,a\;,
\end{equation}
where $b\cdot \pounds_\mathbf{w}\,a$ denotes the contraction, as
described above, of elements of $V$ and elements of $V^\ast$ and
$\mathbf{w}\in\mathfrak{g}({\mathcal{D}})$. (These operations do {\it not}
depend on a Riemannian structure.)
\end{definition}

For a path $\eta_t \in \operatorname{Diff}(\mathcal{D})$, let 
$\mathbf{u}(x, t)$
be its Eulerian velocity and consider the curve $a(t)$ with initial
condition $a_0$ given by the equation
\begin{equation}
\dot a + \pounds_\mathbf{u} a = 0.
\label{continuityequation}
\end{equation}
Let the Lagrangian $L_{a_0}({\bf U}) := L({\bf U}, a_0)$ be right-invariant
under $\operatorname{Diff}(\mathcal{D})$. We can now state the
Euler--Poincar\'{e} Theorem for Continua of \cite{HoMaRa1998a}.

\begin{theorem}[Euler--Poincar\'{e} Theorem for Continua.]
\label{EPforcontinua}
Consider a path $\eta_t$ in  $\operatorname{Diff}(\mathcal{D})$ with
Lagrangian velocity ${\bf U}$ and Eulerian velocity $\mathbf{u}$. The
following are equivalent:

\begin{enumerate}
\item [{\bf i}] Hamilton's variational principle
\begin{equation} \label{continuumVP}
\delta \int_{t_1}^{t_2} L\left(X, {\bf U}_t (X),
a_0(X)\right)\,dt=0
\end{equation}
holds, for variations $\delta\eta_t$ vanishing at the endpoints.
\item [{\bf ii}] $\eta_t$ satisfies the Euler--Lagrange
equations for $L_{a_0}$ on
$\operatorname{Diff}(\mathcal{D})$.
\item [{\bf iii}] The constrained variational principle in
Eulerian coordinates
\begin{equation}\label{continuumconstrainedVP}
   \delta \int_{t_1}^{t_2} \ell(\mathbf{u},a)\ dt=0
\end{equation}
holds on $\mathfrak{g}(\mathcal{D}) \times V^\ast$, using
variations of the form
\begin{equation}\label{continuumvariations}
   \delta \mathbf{u} = \frac{\partial \mathbf{w}}{\partial t}
                   +[\mathbf{u},\mathbf{w}]
= \frac{\partial \mathbf{w}}{\partial t}
                   +{\rm\,ad\,}_\mathbf{u}\mathbf{w}
\,, \qquad
   \delta a = - \pounds_\mathbf{w}\,a,
\end{equation}
where $\mathbf{w}_t = \delta\eta_t \circ \eta_t^{-1}$ vanishes at
the endpoints.
\item [{\bf iv}] The Euler--Poincar\'{e} equations for continua
\begin{equation}\label{continuumEP}
   \frac{\partial }{\partial t}\frac{\delta \ell}{\delta \mathbf{u}}
   = -\, \operatorname{ad}^{\ast}_\mathbf{u}\frac{\delta \ell}
        {\delta \mathbf{u}}
   +\frac{\delta \ell}{\delta a}\diamond a
   =-\pounds_\mathbf{u} \frac{\delta \ell}{\delta \mathbf{u}}
    +\frac{\delta \ell}{\delta a}\diamond a\,,
\end{equation}
hold, with auxiliary equations $(\partial_t + \pounds_\mathbf{u})a = 0$
for each advected quantity $a(t)$. The
$\diamond$ operation defined in (\ref{continuumdiamond}) needs to be
determined on a case by case basis, depending on the nature of the tensor
$a(t)$. The variation $\mathbf{m}=\delta \ell/\delta \mathbf{u}$ is a
one--form density and we have used relation (\ref{continuumcoadjoint}) in
the last step of equation (\ref{continuumEP}).
\end{enumerate}
\end{theorem}

We refer to \cite{HoMaRa1998a} for the proof of this theorem in
the abstract setting. We shall see some of the features of this result
in the concrete setting of continuum mechanics shortly.

\subsection*{Discussion of the Euler-Poincar\'e equations}

The following string of equalities
shows {\it directly} that {\bf iii} is equivalent to {\bf iv}:
\begin{eqnarray}\label{continuumEPderivation}
0
&=&\delta \int_{t_1}^{t_2} l(\mathbf{u}, a) dt
=\int_{t_1}^{t_2}\left(\frac{\delta l}{\delta \mathbf{u}}\cdot
\delta\mathbf{u} +\frac{\delta l}{\delta a}\cdot \delta a\right)dt
\nonumber \\
&=&\int_{t_1}^{t_2} \left[\frac{\delta l}{\delta \mathbf{u}}
\cdot \left(\frac{\partial\mathbf{w}}
{\partial t}-{\rm ad}_\mathbf{u}\,\mathbf{w}\right) -\frac{\delta
l}{\delta a}\cdot \pounds_\mathbf{w}\, a \right]dt
\nonumber \\
&=&\int_{t_1}^{t_2} \mathbf{w}\cdot
\left[-\,\frac{\partial}{\partial t}
\frac{\delta l}{\delta\mathbf{u}} -{\rm ad}^*_\mathbf{u}\frac
{\delta l}{\delta\mathbf{u}} +\frac{\delta l}{\delta a} \diamond
a\right]dt\,.
\end{eqnarray}
The rest of the proof follows essentially the same track as the proof of the pure
Euler-Poincar\'e theorem, modulo slight changes to accomodate the advected quantities.

In the absence of dissipation, most Eulerian fluid equations%
\footnote{Exceptions to this statement are certain multiphase fluids, and
complex fluids with active internal degrees of freedom such as liquid
crystals. These require a further extension, not discussed here.} 
can be written in the EP form in equation (\ref{continuumEP}), 
\begin{equation}\label{EPeqn}
\frac{\partial}{\partial t} \frac{\delta\ell}{\delta\mathbf{u}}
+
{\rm ad}_\mathbf{u}^*\frac{\delta\ell}{\delta\mathbf{u}}
=
\frac{\delta\ell}{\delta{a} }\diamond{a}
\,,\quad\hbox{with}\quad
\big(\partial_t + \pounds_\mathbf{u}\big)a = 0
\,.
\end{equation}
Equation (\ref{EPeqn}) is {\bfi Newton's Law}: The Eulerian time derivative
of the momentum density $\mathbf{m}=\delta\ell/\delta\mathbf{u}$ (a
one-form density dual to the velocity $\mathbf{u}$) is equal to the force
density $(\delta\ell/\delta a)\diamond a$, with the $\diamond$ operation
defined in (\ref{continuumdiamond}).
Thus, Newton's Law is written in the Eulerian fluid representation as,%
\footnote{
In coordinates, a one-form density takes the form
$\mathbf{m}\cdot d\mathbf{x}\otimes{dV}$ and the EP equation
(\ref{continuumEP}) is given  neumonically by
\[
\frac{d}{dt}\Big|_{Lag}\!\!\big(\mathbf{m}\cdot d\mathbf{x}\otimes{dV}\big)
=
\underbrace{\
\frac{d\mathbf{m}}{dt}\Big|_{Lag}\hspace{-3mm}\cdot\,
d\mathbf{x}\otimes{dV}\
}_{\hbox{Advection}}
\
+\
\underbrace{\
\mathbf{m}\cdot d\mathbf{u}\otimes{dV}\
}_{\hbox{Stretching}}
\
+\
\underbrace{\
\mathbf{m}\cdot d\mathbf{x}\otimes(\nabla\cdot\mathbf{u}){dV}\
}_{\hbox{Expansion}}
=
\frac{\delta\ell}{\delta{a} }\diamond{a}
\]
with
$
\frac{d}{dt}\Big|_{Lag}\!\!\!\!d\mathbf{x}
:=
\big(\partial_t + \pounds_\mathbf{u}\big)d\mathbf{x}
=d\mathbf{u}=\mathbf{u}_{,j}dx^j
\,,
$ upon using commutation of Lie derivative and exterior derivative.
Compare this formula with the definition of 
$ \operatorname{ad}^\ast_\mathbf{u}(\mathbf{m}\otimes  dV)$ in
equation (\ref{continuumcoadjoint-coord}).}
\begin{equation}\label{EPeqn-m}
\hspace{-3mm}\frac{d}{dt}\Big|_{Lag}\!\! \mathbf{m}
:=
\big(\partial_t + \pounds_\mathbf{u}\big)\mathbf{m}
=
\frac{\delta\ell}{\delta{a} }\diamond{a}
\,,\quad\hbox{with}\quad
\frac{d}{dt}\Big|_{Lag}\!\!a
:=
\big(\partial_t + \pounds_\mathbf{u}\big)a = 0
\,.
\end{equation}
\begin{itemize}
\item
The left side of the EP equation in (\ref{EPeqn-m}) describes the fluid's
dynamics due to its kinetic energy. A fluid's kinetic energy typically
defines a norm for the Eulerian fluid velocity,
$KE=\frac{1}{2}\|\mathbf{u}\|^2$. The left side of the EP equation is the
{\bfi geodesic} part of its evolution, with respect to this norm.
See \cite{ArKh1998} for discussions of this interpretation of ideal
incompressible flow and references to the literature. However, in a
gravitational field, for example, there will also be dynamics due to
potential energy.  And this dynamics will by governed by the right side of
the EP equation. 

\item
The right side of the EP equation in (\ref{EPeqn-m}) modifies the geodesic
motion. Naturally, the right side of the EP equation is also a geometrical
quantity. The diamond operation $\diamond$ represents the dual of the Lie
algebra action of vectors fields on the tensor $a$.  Here
$\delta\ell/\delta{a}$ is the dual tensor, under the natural pairing
(usually, $L^2$ pairing) $\langle\,\cdot\,,\cdot\,\rangle$ that is induced
by the variational derivative of the Lagrangian $\ell(\mathbf{u},a)$. The
diamond operation $\diamond$ is defined in terms of this pairing in
(\ref{continuumdiamond}). For the $L^2$ pairing, this is integration
by parts of (minus) the Lie derivative in (\ref{continuumdiamond}). 

\item
The quantity $a$ is typically a tensor (e.g., a density, a scalar, or a
differential form) and we shall sum over the various types of tensors
$a$ that are involved in the fluid description.  The second equation in
(\ref{EPeqn-m}) states that each tensor
$a$ is carried along by the Eulerian  fluid velocity $\mathbf{u}$. Thus,
$a$ is for fluid ``attribute,'' and its Eulerian evolution is given by 
minus its Lie derivative, $-\,\pounds_\mathbf{u}a$. That is, $a$ stands for
the set of fluid attributes that each Lagrangian fluid  parcel carries
around (advects), such as its buoyancy, which is determined by its
individual salt, or heat content, in ocean circulation.  

\item
Many examples of how equation (\ref{EPeqn-m})
arises in the dynamics of continuous media are given in \cite{HoMaRa1998a}.
The EP form of the Eulerian fluid description in (\ref{EPeqn-m}) is
analogous to the classical dynamics of rigid bodies (and tops, under
gravity) in body coordinates. Rigid bodies and tops are also governed by
Euler-Poincar\'e equations, as Poincar\'e showed in a
two-page paper with no references, over a century ago
\cite{Po1901}. For modern discussions of the EP theory, see, e.g., 
\cite{MaRa1994}, or \cite{HoMaRa1998a}. 
\end{itemize}

\begin{exercise}
For what types of tensors $a_0$ can one recast the EP equations for continua
(\ref{continuumEP}) as geodesic motion, by using a version of the Kaluza-Klein
construction?
\end{exercise}

\subsection{Corollary of the EP theorem: the Kelvin-Noether circulation
theorem}

\begin{corollary}[Kelvin-Noether Circulation Theorem.]\label{KelThmforcontinua}
Assume $\mathbf{u}(x, t)$
satisfies the Euler--Poincar\'e equations for continua:
\[
\frac{\partial }{\partial t}\left(\frac{\delta \ell}{\delta \mathbf{u}}\right)
   = -\pounds_\mathbf{u} \left(\frac{\delta \ell}{\delta \mathbf{u}}\right)
    +\frac{\delta \ell}{\delta a}\diamond a
\]
and the quantity $a$ satisfies the {\bfi advection relation}
\begin{equation}\label{advect-def}
\frac{\partial a}{\partial t} +
\pounds_\mathbf{u} a = 0.
\end{equation}
Let $\eta_t$ be the flow of the Eulerian
velocity field $\mathbf{u}$, that is, $\mathbf{u}  =
(d\eta_t/dt)\circ \eta_t^{-1}$. Define the advected fluid loop 
$\gamma_t  := \eta_t\circ \gamma_0$ and the circulation map $I(t)$ by
\begin{equation}\label{circ-def}
I(t) = \oint_{\gamma_t }\frac{1}{D }\frac{\delta \ell}{\delta \mathbf{u}}
\,.\end{equation}
In the circulation map $I(t)$ the advected mass density $D_t$ satisfies
the push forward relation $D_t=\eta_*D_0$. This implies the 
advection relation (\ref{advect-def}) with $a=D$, namely, the continuity
equation,
\[
\partial_tD+{\rm div}\,D\mathbf{u}=0
\,.
\]
Then the map $I(t)$
satisfies the {\bfi Kelvin circulation relation},
\begin{equation}\label{KN-theorem}
\frac{d}{dt}I(t) = \oint_{\gamma_t }
\frac{1}{D}\frac{\delta \ell}{\delta a}\diamond a\;.
\end{equation}
\end{corollary}
Both an abstract proof of the Kelvin-Noether Circulation Theorem and a
proof tailored for the case of continuum mechanical systems are given in
\cite{HoMaRa1998a}. We provide a version of the latter below.\\

\begin{proof} First we change variables in the expression
for $I(t)$:
\[
I(t) = \oint_{\gamma_t}\frac{1}{D_t}
\frac{\delta l}{\delta \mathbf{u}}
=\oint_{\gamma_0} \eta_t^\ast\left[\frac{1}{D_t}\frac{\delta l}
{\delta \mathbf{u}}\right] = \oint_{\gamma_0}
\frac{ 1 }{ D _0 } \eta_t^\ast\left[\frac{\delta
l} {\delta \mathbf{u}}\right].
\]
Next, we use the Lie derivative formula, namely
\[
\frac{d}{dt}\left(\eta_t^*\alpha_t\right) =
\eta_t^*\left(\frac{\partial}{\partial t}\alpha_t
+ \pounds_\mathbf{u} \alpha_t \right)\;,
\]
applied to a one--form density $\alpha_t$.
This formula gives
\begin{eqnarray*}
      \frac{d}{dt} I(t)
  & = & \frac{d}{dt}  \oint_{\gamma_0}
\frac{ 1 }{ D _0 } \eta_t^\ast\left[\frac{\delta
l} {\delta \mathbf{u}}\right] \\
  & = & \oint_{\gamma_0} \frac{1}{D _0} \frac{d}{dt}
\left( \eta_t^\ast\left[
\frac{\delta l} {\delta \mathbf{u}}\right]\right)  \\
  & = & \oint_{\gamma_0} \frac{1}{D _0} \eta_t^*\left[
\frac{\partial}{\partial t}
\left(\frac{\delta l} {\delta \mathbf{u}}\right) +
\pounds_\mathbf{u}
\left(\frac{\delta l} {\delta \mathbf{u}} \right)\right].
\end{eqnarray*}
By the Euler--Poincar\'e equations (\ref{continuumEP}), this becomes
\[
      \frac{d}{dt} I(t)
  =   \oint_{\gamma_0} \frac{1}{D _0} \eta_t^*\left[
\frac{\delta  l}{\delta  a}
\diamond a \right] = \oint_{\gamma_t} \frac{1}{D _t} \left[
\frac{\delta  l}{\delta  a}
\diamond a \right],
\]
again by the change of variables formula.
\end{proof}

\begin{corollary}
Since the last expression holds for every loop $\gamma_t$, we may
write it as
\begin{equation}
\left(\frac{\partial}{\partial t} + \pounds_\mathbf{u} \right)
\frac{1}{D} \frac{\delta l}{\delta \mathbf{u}}
= \frac{1}{D} \frac{\delta  l}{\delta  a} \diamond a\,.
\label{KThfm}
\end{equation}
\end{corollary}

\begin{remark}
The Kelvin-Noether theorem is called so here because its derivation relies
on the invariance of the Lagrangian $L$ under the particle relabeling
symmetry, and  Noether's theorem is associated with this symmetry. However,
the result (\ref{KN-theorem}) is the {\bfi Kelvin circulation theorem}: the
circulation integral $I(t)$ around any fluid loop ($\gamma_t$, moving with
the velocity of the fluid parcels $\mathbf{u}$) is invariant under the
fluid motion.  These two statements are equivalent.  We note that {\bfi two
velocities} appear in the integrand $I(t)$: the fluid velocity $\mathbf{u}$
and $D^{-1}\delta\ell/\delta\mathbf{u}$. The latter velocity is the
momentum density $\mathbf{m}=\delta\ell/\delta\mathbf{u}$ divided by the
mass density $D$.  These two velocities are the basic ingredients for
performing  modeling and analysis in any ideal fluid problem. One
simply needs to put these ingredients together in the Euler-Poincar\'e
theorem and its corollary, the Kelvin-Noether theorem.
\end{remark}

\section{Euler--Poincar\'{e} theorem \& GFD (geophysical
fluid dynamics)}
\label{sec-GFD-applic}

\subsection{Variational Formulae in Three Dimensions} 
We compute explicit formulae for the variations $\delta a$ in the cases
that the set of tensors $a$ is drawn from a set of scalar fields and
densities on ${\mathbb{R}}^3$. We shall denote this symbolically by
writing
\begin{equation}
a\in\{b,D\,d^3x\}\,.
\label{Eul-ad-qts}
\end{equation}
We have seen that invariance of the
set $a$ in the Lagrangian picture under the dynamics of
$\mathbf{u}$ implies in the Eulerian picture that
\[
     \left( \frac{\partial}{\partial t}
         + \pounds_\mathbf{u} \right) \,a=0 \, ,
\]
where $\pounds_\mathbf{u}$ denotes Lie derivative with respect to the
velocity vector field $\mathbf{u}$. Hence, for a fluid dynamical
Eulerian action $\mathfrak{S}=\int\,dt\ \ell(\mathbf{u};b,D)$, the advected
variables
$b$ and $D$ satisfy the following Lie-derivative relations,
\begin{eqnarray}
\left(\frac{\partial}{\partial t}+ \pounds_\mathbf{u}\right) b=0,
&{\rm or}&
\frac{\partial b}{\partial t}
= -\ \mathbf{u}\cdot\nabla\,b\,,
\label{eqn-b} \\
\left(\frac{\partial}{\partial t}+ \pounds_\mathbf{u}\right)D\,d^3x=0,
&{\rm or}&
\frac{\partial D}{\partial t}
= -\ \nabla\cdot(D\mathbf{u})\,.
\label{eqn-D}
\end{eqnarray}
In fluid dynamical applications, the advected Eulerian variables $b$
and $D\,d^3x$ represent the buoyancy $b$ (or specific entropy, for the
compressible case) and volume element (or mass density) $D\,d^3x$,
respectively. According to Theorem \ref{EPforcontinua}, equation
(\ref{continuumconstrainedVP}), the variations of the tensor
functions $a$ at fixed $\mathbf{x}$ and $t$ are also given by Lie
derivatives, namely
$\delta a = -\,\pounds_\mathbf{w}\,a$, or
\begin{eqnarray}
\delta b
&=& -\,\pounds_\mathbf{w}\ b = -\,\mathbf{w}\cdot\nabla\,b\,,
\nonumber \\
\delta D\ d^3x&=&-\,\pounds_\mathbf{w}\,(D\,d^3x)
= -\,\nabla\cdot(D\mathbf{w})\ d^3x
\,.
\end{eqnarray}
Hence, Hamilton's principle (\ref{continuumconstrainedVP}) with this
dependence yields
\begin{eqnarray}
0 &=&\delta \int dt\ \ell(\mathbf{u}; b,D)
\nonumber \\
&=&\int dt\ \bigg[\frac{\delta \ell}{\delta
\mathbf{u}}\cdot \delta\mathbf{u}
+\frac{\delta \ell}{\delta b}\ \delta b
+\frac{\delta \ell}{\delta D}\ \delta D
\bigg]
\nonumber \\
&=&\int dt\ \bigg[\frac{\delta \ell}{\delta \mathbf{u}}
\cdot \Big(\frac{\partial \mathbf{w}}{\partial t}
-{\rm ad}_\mathbf{u}\,\mathbf{w}\Big)
-\frac{\delta \ell}{\delta b}\ \mathbf{w}\cdot\nabla\,b
-\frac{\delta \ell}{\delta D}\ \Big(\nabla
\cdot(D\mathbf{w})\Big)\bigg]
\nonumber \\
&=&\int dt\ \mathbf{w}\cdot
\bigg[-\frac{\partial }{\partial t}
\frac{\delta \ell}{\delta \mathbf{u}}
-{\rm ad}^*_\mathbf{u}\ \frac{\delta \ell}{\delta \mathbf{u}}
-\frac{\delta \ell}{\delta b}\ \nabla\,b
+D\ \nabla\frac{\delta \ell}{\delta D}\bigg]
\nonumber \\
&=&-\int dt\ \mathbf{w}\cdot
\bigg[\Big(\frac{\partial }{\partial t}
 +\pounds_\mathbf{u} \Big)\frac{\delta \ell}{\delta \mathbf{u}}
+\frac{\delta \ell}{\delta b}\ \nabla\,b
-D\ \nabla\frac{\delta \ell}{\delta D}\bigg]\,,
\label{eq-EP-Eul}
\end{eqnarray}
where we have consistently dropped boundary terms arising from
integrations by parts, by invoking natural boundary conditions.
Specifically, we may impose $\mathbf{\hat{n}}\cdot\mathbf{w}=0$ on the
boundary, where $\mathbf{\hat{n}}$ is the boundary's outward unit normal
vector and $\mathbf{w} = \delta\eta_t \circ \eta_t^{-1}$ vanishes at
the endpoints.

\subsection{Euler--Poincar\'e framework for GFD}
The Euler--Poincar\'e equations for continua (\ref{continuumEP}) may
now be summarized in vector form for advected Eulerian variables $a$ in the
set (\ref{Eul-ad-qts}). We adopt the notational convention of the
circulation map $I$ in equations (\ref{circ-def}) and (\ref{KN-theorem})
that a one form density can be made into a one form (no longer a density)
by dividing it by the mass density $D$ and we use the Lie-derivative
relation for the continuity equation
$({\partial}/{\partial t}+\pounds_\mathbf{u})Dd^3x = 0$. Then, the
Euclidean components of the Euler--Poincar\'e equations for continua
in equation (\ref{eq-EP-Eul}) are expressed in Kelvin theorem form
(\ref{KThfm}) with a slight abuse of notation as
\begin{equation}
\Big(\frac{\partial }{\partial t} + \pounds_\mathbf{u}\Big)
\Big(\frac{1}{D}\frac{\delta \ell}{\delta
\mathbf{u}}\cdot d\mathbf{x}\Big)
\,+\,\frac{1}{D}\frac{\delta \ell}{\delta b}\nabla b \cdot d\mathbf{x}
\,-\,\nabla\Big(\frac{\delta \ell}{\delta D}\Big)\cdot d\mathbf{x}
 = 0\,,
\label{EP-Kthm}
\end{equation}
in which the variational derivatives of the Lagrangian $\ell$ are to be
computed according to the usual physical conventions, i.e., as
Fr\'echet derivatives. Formula (\ref{EP-Kthm}) is the
Kelvin--Noether form of the equation of motion for ideal continua.
Hence, we have the explicit Kelvin theorem expression, cf. equations
(\ref{circ-def}) and (\ref{KN-theorem}),
\begin{equation} \label{KN-theorem-bD}
\frac{d}{dt}\oint_{\gamma_t(\mathbf{u})} \frac{1}{D}\frac{\delta \ell}{\delta \mathbf{u}}\cdot d\mathbf{x}
= -\oint_{\gamma_t(\mathbf{u})}
\frac{1}{D}\frac{\delta \ell}{\delta b}\nabla b \cdot d\mathbf{x}\;,
\end{equation}
where the curve $\gamma_t(\mathbf{u})$ moves with the fluid velocity
$\mathbf{u}$. Then, by Stokes' theorem, the Euler equations generate
circulation of
$\mathbf{v}:=(D^{-1}\delta{l}/\delta\mathbf{u})$
whenever the gradients~$\nabla b$ and
$\nabla(D^{-1}\delta{l}/\delta{b})$ are not collinear. The
corresponding {\bfi conservation of potential vorticity} $q$ on fluid
parcels is given by
\begin{equation} \label{pv-cons-EP}
\frac{\partial{q}}{\partial{t}}+\mathbf{u}\cdot\nabla{q} =
0\,,
\quad \hbox{where}\quad
{q}=\frac{1}{D}\nabla{b}\cdot{\rm curl}
\left(\frac{1}{D}\frac{\delta \ell}{\delta \mathbf{u}}\right).
\end{equation}
This is also called {\bfi PV convection}.
Equations (\ref{EP-Kthm}-\ref{pv-cons-EP}) embody most of the panoply
of equations for GFD.  The vector form of equation (\ref{EP-Kthm}) is,
\begin{eqnarray}\label{vec-EP-eqn}
\underbrace{
\Big({\partial\over \partial t} + {\bf u}\cdot\nabla\Big)
\Big({1\over D}{\delta l\over \delta {\bf u}}\Big)
+{1\over D}{\delta l\over \delta u^j}\nabla u^j}
_{\hbox{\bf Geodesic Nonlinearity: Kinetic energy}}
= 
\underbrace{\
\nabla {\delta l\over \delta D}
-{1\over D} {\delta l\over \delta b}\nabla b\
}_{\hbox{\bf Potential energy}}
\end{eqnarray}
In geophysical applications, the Eulerian variable $D$ represents
the frozen-in volume element and $b$ is the buoyancy. In this case,
{\bf Kelvin's theorem} is
\[
{dI\over dt}=\int\int_{S(t)}
\nabla \left({1\over D} {\delta l\over \delta b}\right)
\times\nabla b \cdot d {\bf S} 
\,,
\]
with circulation integral
\[
I=\oint_{\gamma(t)} {1\over D}{\delta l\over \delta {\bf u}}\cdot
d{\bf x}
\,.
\]
\subsection{Euler's Equations for a Rotating Stratified Ideal
Incompressible Fluid}
\label{sec-Euler}

\paragraph{The Lagrangian.} In the Eulerian velocity representation, we
consider Hamilton's principle for fluid motion in a three dimensional
domain with action functional ${S}=\int\, l\,dt$ and
Lagrangian
$l(\mathbf{u},b,D)$ given by
\begin{equation}
l(\mathbf{u},b,D) = \int 
\rho_0 D (1+b)
\bigg(\frac{1}{2} |\mathbf{u}|^2
+ \mathbf{u}\cdot\mathbf{R}(\mathbf{x}) - gz\bigg)
 - p(D-1)\,d^{\,3}x\,,
\label{lag-v}
\end{equation}
where $\rho_{tot}=\rho_0 D (1+b)$ is the total mass density, $\rho_0$
is a dimensional constant and $\mathbf{R}$ is a given function of
$\mathbf{x}$. This variations at fixed $\mathbf{x}$ and $t$ of this
Lagrangian are the following, 
\begin{eqnarray}
&&\frac{1}{D}\frac{{\delta} l}{{\delta} \mathbf{u}}
= \rho_0(1+b)(\mathbf{u}+ \mathbf{R})\,,
\quad
\frac{{\delta} l}{{\delta} b}
= \rho_0 D\Big(\frac{1}{2}|\mathbf{u}|^2
+ \mathbf{u}\cdot\mathbf{R} - gz\Big)\,,
\nonumber \\
&&\frac{{\delta} l}{{\delta} D}
= \rho_0(1+b)\Big(\frac{1}{2}|\mathbf{u}|^2
+ \mathbf{u}\cdot\mathbf{R} -
gz\Big) - p\,,
\quad
\frac{{\delta} l}{{\delta} p} = -\,(D-1)\,.
\hspace{.25in}
\label{vds-1}
\end{eqnarray}
Hence, from the Euclidean component formula (\ref{vec-EP-eqn}) for
Hamilton principles of this type and the fundamental vector identity, 
\begin{equation}
(\mathbf{b}\cdot\nabla)\mathbf{a} + a_j\nabla
b^j =-\ \mathbf{b}\times(\nabla\times \mathbf{a})
+ \nabla(\mathbf{b}\cdot\mathbf{a})\,,
\label{fvid}
\end{equation}
we find the motion equation for an Euler fluid in three dimensions,
\begin{equation}
\frac{d\mathbf{u}}{dt} - \mathbf{u} \times {\rm curl}\, \mathbf{R}
+g\hat{\bf z}
+ \frac{1}{\rho_0(1+b)}\nabla p = 0\,,
\label{Eul-mot}
\end{equation}
where ${\rm curl}\,\mathbf{R}=2\boldsymbol{\Omega}(\mathbf{x})$ is the
Coriolis parameter (i.e., twice the local angular rotation
frequency). In writing this equation, we have used advection of
buoyancy,
$$
\frac{\partial{b}}{\partial{t}}+\mathbf{u}\cdot\nabla{b} = 0,
$$
from equation (\ref{eqn-b}). The pressure $p$ is determined by requiring
preservation of the constraint $D=1$, for which the continuity equation
(\ref{eqn-D}) implies ${\rm div}\,\mathbf{u}=0$. The Euler motion equation
(\ref{Eul-mot}) is Newton's Law for the acceleration of a fluid due to
three forces: Coriolis, gravity and pressure gradient. The dynamic
balances among these three forces produce the many circulatory flows of
geophysical fluid dynamics. The {\bfi conservation of potential vorticity}
$q$ on fluid parcels for these Euler GFD flows is given by
\begin{equation} \label{pv-cons-EulGFD}
\frac{\partial{q}}{\partial{t}}+\mathbf{u}\cdot\nabla{q} =
0\,,
\quad \hbox{where, on using }D=1\,,\quad
{q}=\nabla{b}\cdot{\rm curl}
\big(\mathbf{u}+ \mathbf{R}\big).
\end{equation}

\begin{exercise}[Semidirect-product Lie-Poisson bracket for
compressible ideal fluids]$\quad$
\begin{enumerate}
\item
Compute the Legendre transform for the Lagrangian,
\[
l(\mathbf{u},b,D):\,
\mathfrak{X}\times\Lambda^0\times\Lambda^3\mapsto\mathbb{R}\]
whose advected variables satisfy the auxiliary equations,
\[
\frac{\partial b}{\partial t}
= -\ \mathbf{u}\cdot\nabla\,b
\,,\qquad
\frac{\partial D}{\partial t}
= -\ \nabla\cdot(D\mathbf{u})\,.
\]
\item
Compute the Hamiltonian, assuming the Legendre transform is a linear
invertible operator on the velocity $\mathbf{u}$. For definiteness
in computing the Hamiltonian, assume the Lagrangian is given by
\begin{equation}
l(\mathbf{u},b,D) = \int 
D \Big(\frac{1}{2} |\mathbf{u}|^2
+ \mathbf{u}\cdot\mathbf{R}(\mathbf{x})
 - e(D,b)\Big)
\,d^{\,3}x\,,
\end{equation}
with prescribed function $\mathbf{R}(\mathbf{x})$ and specific internal
energy $e(D,b)$ satisfying the First Law of Thermodynamics,
\[
de=\frac{p}{D^2}dD + T db
\,,
\]
where $p$ is pressure, $T$ temperature.
\item
Find the semidirect-product Lie-Poisson bracket for the Hamiltonian
formulation of these equations. 
\item
Does this Lie-Poisson bracket have Casimirs? If so, what are the
corresponding symmetries and momentum maps?
\end{enumerate}
\end{exercise}

\section{Hamilton-Poincar\'e reduction and Lie-Poisson equations}

In the Euler-Poincar\'e framework one starts with a Lagrangian defined on
the tangent bundle of a Lie group G
\[
L:TG\rightarrow\mathbb{R}
\]
and the dynamics is given by Euler-Lagrange equations arising from the variational
principle
\[
\delta\int_{t_0}^{t_1} L(g,\dot{g}) dt =0
\]
The Lagrangian L is taken left/right invariant and because of this property
one can \emph{reduce} the problem obtaining a new system which is defined
on the Lie algebra $\mathfrak{g}$ of G, obtaining a new set of equations,
the Euler-Poincar\'e equations, arising from a reduced variational principle
\[
\delta\int_{t_0}^{t_1} l(\xi) dt =0
\]
where $l(\xi)$ is the reduced lagrangian and $\xi\in\mathfrak{g}$.

\begin{problem}
Is there a similar procedure for Hamiltonian systems? More precisely: given
a Hamiltonian function defined on the cotangent bundle $T^*G$
\[
H:T^*G\rightarrow\mathbb{R}
\]
one wants to perform a similar procedure of reduction and derive the equations of motion on the dual of the Lie algebra $\mathfrak{g}^*$, provided the Hamiltonian
is again left/right invariant.
\end{problem}
Hamilton-Poincar\'e reduction gives a positive answer to this problem, in
the context of variational principles as it is done in the Euler-Poincar\'e
framework: we are going to explain how this procedure is performed.\\
More in general, we will also consider advected quantities belonging to a
vector space $V$ on which $G$ acts, so that the Hamiltonian is written in this case as \cite{HoMaRa98}
[HoMaRa98]
\[
H:T^*G\times V^*\rightarrow \mathbb{R}
\]
The space $V$ is regarded here exactly the same as in the Euler-Poincar\'e theory.\\
The equations of motion, i.e. Hamilton's equations, may be derived from the
following variational principle
\[
\delta\int_{t_0}^{t_1} \{\langle p(t),\dot{g}(t) \rangle - H_{a_0} (g(t),p(t))\}\, dt =0
\]
as it is well know from ordinary classical mechanics ($\dot{g}(t)$ has
to be considered as the tangent vector to the curve $g(t)$, so that $\dot{g}(t)\in
T_{\!g(t)} G$).

\begin{problem}
What happens if $H_{a_0}$ is left/right invariant?
\end{problem}
It turns out that in this case the whole function
\[
F(g,\dot{g},p)=\langle p,\dot{g} \rangle - H_{a_0} (g,p)
\]
is also invariant. The proof is straightforward once the action is
specified (from now on we consider only left invariance):
\[
h\,(g,\dot{g},p)=(hg,\,T_g L_h\, \dot{g},\,T^*_{hg} L_{h^{-1}}\, p)
\]
where $T_g L_h : T_g G \rightarrow T_{hg}G$ is the tangent of the left translation map $L_h\,g=hg\in{G}$ at the point $g$ and 
$T^*_{hg}L_{h^{-1}}:T^*_g G \rightarrow T^*_{hg}G$ is the dual of the map $T_{hg}L_{h^{-1}} : T_{hg}G\rightarrow T_g G$.\\
We now check that
\begin{align*}
\langle h\, p,\, h\, \dot{g} \rangle
&=\langle T^*_{hg} L_{h^{-1}}\, p,\, T_g L_h\, \dot{g} \rangle
\\
&=\langle p,\,T_{hg} L_{h^{-1}}\circ T_g L_h\, \dot{g} \rangle
\\
&=\langle p,\,T_g (L_{h^{-1}}\circ L_h)\, \dot{g} \rangle
=\langle  p, \dot{g} \rangle
\end{align*}
where the chain rule for the tangent map has been used. The same
result holds for the right action.\\
Due to this invariance property, one can write the variational principle
as
\[
\delta\int_{t_0}^{t_1} \{\langle \mu,\xi \rangle - h (\mu,a)\}\, dt =0
\]
with
\[
\mu(t)=g^{-1}(t)\,p(t)\in\mathfrak{g}^*,
\qquad
\xi(t)=g^{-1}(t)\,\dot{g}(t)\in\mathfrak{g},
\qquad
a(t)=g^{-1}(t)\,a_0\in V^*
\]
In particular $a(t)$ is the solution of
\[
\dot{a}(t)=-\xi(t)\, a_0.
\]
where a Lie algebra action of $\mathfrak{g}$ on $V^*$ is implicitly defined.
In order to find the equations of motion one calculates the variations
\[
\delta\int_{t_0}^{t_1} \{\langle \mu,\xi \rangle - h (\mu,a)\}\, dt 
=
\int_{t_0}^{t_1} \left\{\langle \delta\mu,\xi \rangle + 
\langle \mu,\delta\xi \rangle  - 
\left\langle \delta\mu,\frac{\delta h}{\delta \mu}\right\rangle -
\left\langle \delta a,\frac{\delta h}{\delta a}\right\rangle \right\}\, dt
\]
As in the Euler-Poincar\'e theorem, we use the following expressions for the
variations
\begin{align*}
\delta\xi=\dot\eta + [\xi,\eta],
\qquad
\delta a=-\eta a
\end{align*}
and using the definition of the diamond operator we find
\begin{align*}
&\int_{t_0}^{t_1} \left\{\langle \delta\mu,\xi \rangle + 
\langle \mu,\delta\xi \rangle  - 
\left\langle \delta\mu,\frac{\delta h}{\delta \mu}\right\rangle -
\left\langle \delta a,\frac{\delta h}{\delta a}\right\rangle \right\}\, dt\\
&=
\int_{t_0}^{t_1} \left\{\left\langle 
\delta\mu,\,\xi - \frac{\delta h}{\delta \mu} \right\rangle + 
\langle \mu,\,\dot\eta + \text{ad}_\xi \eta \rangle  + 
\left\langle \eta a,\frac{\delta h}{\delta a}\right\rangle \right\}\, dt\\
&=
\int_{t_0}^{t_1} \left\{\left\langle 
\delta\mu,\,\xi - \frac{\delta h}{\delta \mu} \right\rangle + 
\langle -\dot\mu + \text{ad}^*_\xi \mu,\,\eta \rangle  - 
\left\langle \frac{\delta h}{\delta a}\diamond a,\, \eta\right\rangle \right\}\, dt\\
\end{align*}
so that 
\[
\xi = \frac{\delta h}{\delta \mu}
\]
and the equations of motion are
\[
\dot\mu = \text{ad}^*_\xi \mu - \frac{\delta h}{\delta a}\diamond a
\]
together with
\[
\dot{a}=-\frac{\delta h}{\delta \mu}\, a.
\]
This equations of motion written on the dual Lie algebra $\mathfrak{g}$ are
called \emph{Lie-Poisson} equations. We have now proven the following
\begin{theorem}{\textbf{[Hamilton-Poincar\'e reduction theorem]\\}}
With the preceding notation, the following statements are equivalent:\\
\begin{enumerate}
\item
With $a_0$ held fixed, the variational principle
\[
\delta\int_{t_0}^{t_1} \{\langle p(t),\dot{g}(t) \rangle - H_{a_0} (g(t),p(t))\}\, dt =0
\]
holds, for variations $\delta{g(t)}$ of $g(t)$ vanishing at the endpoints.
\item
$(g(t),p(t))$ satisfies Hamilton's equations for $H_{a_0}$ on G.
\item
The constrained variational principle
\[
\delta\int_{t_0}^{t_1} \{\langle \mu(t),\xi(t) \rangle - h (\mu(t),a(t))\}\, dt =0
\]
holds for $\mathfrak{g}\times V^*$, using variations of $\xi$ and $a$ of the form
\begin{align*}
\delta\xi=\dot\eta + [\xi,\eta],
\qquad
\delta a=-\eta a
\end{align*}
where $\eta(t)\in\mathfrak{g}$ vanishes at the endpoints
\item
The \emph{Lie-Poisson} equations hold on $\mathfrak{g}\times V^*$
\[
(\dot\mu,\dot{a}) = \left(\textnormal{ad}^*_\xi \mu - \frac{\delta h}{\delta a}\diamond a,\,-\frac{\delta h}{\delta \mu}\, a\right)
\]
\end{enumerate}
\end{theorem}

\begin{remark}
More exactly one should start with an invariant Hamiltonian defined on
\[
T^*(G\times V)=T^*G\times V\times V^*
\]
However, as mentioned in \cite{}[HoMaRe98], such an approach turns out to
be equivalent to the treatment presented here.
\end{remark}

\begin{remark}{\textnormal{[Legendre transform]}\\}
Lie-Poisson equations may arise from the Euler-Poincar\'e setting by Legendre
transform
\[
\mu=\frac{\delta l}{\delta \xi}.
\]
If this is a diffeomorphism, then the Hamilton-Poincar\'e
theorem is equivalent to the Euler-Poincar\'e theorem.
\end{remark}

\begin{remark}{\textnormal{[Lie-Poisson structure]}\\}
One shows that $\mathfrak{g}^*\times V^*$ is a Poisson manifold:
\begin{align*}
\dot{F}(\mu,a)  & =
\left\langle \dot{\mu},\frac{\delta F}{\delta\mu}\right\rangle 
+\left\langle \dot{a},\frac{\delta F}{\delta a}\right\rangle =
\\
& =\left\langle 
\text{ad}_{\delta H/\delta\mu}^{\ast}\mu
-\frac{\delta H}{\delta a}\diamond a,
\frac{\delta F}{\delta\mu}
\right\rangle -
\left\langle
\frac{\delta H}{\delta\mu}a,\frac{\delta F}{\delta a}
\right\rangle =
\\
& =\left\langle \mu,
\left[  \frac{\delta H}{\delta\mu},\frac{\delta F}{\delta\mu}\right]  
\right\rangle -
\left\langle 
\frac{\delta H}{\delta a}\diamond a,
\frac{\delta F}{\delta\mu}
\right\rangle -
\left\langle
\frac{\delta H}{\delta\mu}a,\frac{\delta F}{\delta a}
\right\rangle =
\\
& =-
\left\langle 
\mu,\left[  \frac{\delta F}{\delta\mu},\frac{\delta H}{\delta\mu}\right] \right\rangle -
\left\langle 
\frac{\delta H}{\delta a}\diamond a,
\right\rangle -
\left\langle
\frac{\delta H}{\delta\mu}a,\frac{\delta F}{\delta a}
\right\rangle =
\\
& =-\left\langle 
\mu,\left[  \frac{\delta F}{\delta\mu},\frac{\delta H}{\delta\mu}\right]
\right\rangle -
\left\langle a,\frac{\delta F}{\delta\mu}\frac{\delta H}{\delta a}-
\frac{\delta H}{\delta\mu}\frac{\delta F}{\delta a}\right\rangle
\end{align*}
In fact it can be easily shown that this structure
\[
\{F,H\}(\mu,a)=
-\left\langle 
\mu,\left[  \frac{\delta F}{\delta\mu},\frac{\delta H}{\delta\mu}\right]
\right\rangle -
\left\langle a,\frac{\delta F}{\delta\mu}\frac{\delta H}{\delta a}-
\frac{\delta H}{\delta\mu}\frac{\delta F}{\delta a}\right\rangle
\]
 satisfies the definition of a Poisson
structure. In particular one finds that \textbf{any dual Lie algebra $\mathfrak{g}$
is a Poisson manifold}.
\\\\\textbf{Please note}: this structure has been found during lectures for the simpler case without advected quantities.
\end{remark}

\begin{remark}{\textnormal{[right invariance]}\\}
It can be shown that for a right invariant Hamiltonian one has
\begin{align*}
\{F,H\}(\mu,a)&=+
\left\langle 
\mu,\left[  \frac{\delta F}{\delta\mu},\frac{\delta H}{\delta\mu}\right]
\right\rangle +
\left\langle a,\frac{\delta F}{\delta\mu}\frac{\delta H}{\delta a}-
\frac{\delta H}{\delta\mu}\frac{\delta F}{\delta a}\right\rangle
\\
(\dot\mu,\dot{a})&=-
\left(\textnormal{ad}^*_\xi \mu - \frac{\delta h}{\delta a}\diamond a,\,
-\frac{\delta h}{\delta \mu}\, a\right)
\end{align*}
with all signs changed respect to the case of left invariance.
\end{remark}

\bigskip

\bigskip

\section{Two applications}

\bigskip

\subsection{The Vlasov equation}
In plasma physics a main topic is collisionless particle dynamics, whose
main equation, the Vlasov equation, will be heuristically derived here. In this context a central role is held by the distribution function on phase space $f(\mathbf{q,p},t)$, basically expressing the particle density on phase space.
Intended as a density one defines $F:=f(\mathbf{q,p},t)d\mathbf{q}d\mathbf{p}$:
because of the conservation of particles, one writes the continuity equation
just as one does as in the context of fluid dynamics
\[
\dot{F}+\nabla \cdot (\mathbf{u}\, F)=0
\] 
where $\mathbf{u}$ is a ``velocity'' vector field on phase space, which is given by the
single particle motion
\[
\mathbf{u}=(\dot{\mathbf{q}},\dot{\mathbf{p}})\in\mathfrak{X}\,(T^*\mathbb{R}^N)
\]
if we now assume that the generic single particle undergoes a Hamiltonian
motion, the Hamiltonian function $h(\mathbf{q,p},t)$ can be introduced directly by means of the single particle Hamilton's equations
\[
(\dot{\mathbf{q}},\dot{\mathbf{p}})=
\left( \frac{\partial h}{\partial \mathbf{p}}, 
-\frac{\partial h}{\partial \mathbf{q}}\right)
\]
which shows that $\mathbf{u}$ has zero divergence, assuming the Hessian of
$h$ is symmetric. Therefore, the Vlasov equation written in terms of the distribution function $f(\mathbf{q,p},t)$ is 
\[
\dot{f}+\mathbf{u}\cdot\nabla f=0
\]
Expanding now the Hamiltonian $h$ as the total single particle energy
\[
h(\mathbf{q,p},t)=\frac1{2m}\mathbf{p}^2+V(\mathbf{q},\mathbf{p},t)
\]
one obtains the more common form
\[
\frac{\partial f}{\partial t}+
\frac{\mathbf{p}}{m}\cdot\frac{\partial f}{\partial \mathbf{q}}-
\frac{\partial V}{\partial \mathbf{q}}\cdot\frac{\partial f}{\partial \mathbf{p}}=0
\]
\begin{problem}
Can Vlasov equation be cast in Lie-Poisson form?
\end{problem}
We show here why the answer is yes. First we write the Vlasov equation in
terms of a generic single particle Hamiltonian $h$ as
\[
\dot{f}+\{f,h\}=0
\]
where we recall the canonical Poisson bracket
\[
\{f,h\}=
\frac{\partial f}{\partial \mathbf{q}}\cdot\frac{\partial h}{\partial \mathbf{p}}-
\frac{\partial f}{\partial \mathbf{p}}\cdot\frac{\partial h}{\partial \mathbf{q}}
\]
The main point of this discussion is that the canonical Poisson bracket provides
the set $\mathcal{F}(T^*\mathbb{R}^N)$ of the functions on the phase space
with a Lie algebra structure 
\[
[k,h]=\{k,h\}
\]
At this point, in order to look for a Lie-Poisson equation, one calculates
the coadjoint operator such that
\[
\langle f,\{h,k\} \rangle=
\langle f,\text{ad}_h k \rangle=
\langle \text{ad}^*_h f,k \rangle=
\langle -\{h,f\},k \rangle
\]
where the last equality is justified by the Leibniz property of the Poisson
bracket, with the pairing defined as
\[
\langle f,g \rangle=\int f\,g\, d\mathbf{q}d\mathbf{p}.
\]
In conclusion, the argument above shows that the Vlasov equation can in fact be written in the Lie Poisson form
\[
\dot{f}+\text{ad}^*_h\,f=0
\]

\bigskip

\subsection{Ideal barotropic compressible fluids}

The reduced Lagrangian for ideal compressible fluids is written as
\[
l(\mathbf{u},D)
=\int \frac{D}{2}|\mathbf{u}|^2 
- De(D)
\,  d\mathbf{x}
\]
where $\mathbf{u}\in\mathfrak{X}(M\!\!\subset\!\mathbb{R}^3)$ is tangential on the boundary $\partial M$ and $D$ is the advected density, which satisfies the {\bfi continuity equation}
\[
\partial_t{D}+\mathcal{L}_\mathbf{u} D=0.
\]
Moreover, the internal energy satisfies the barotropic First Law of Thermodynamics
\[
de=-p(D)d(D^{-1})=\frac{p(D)}{D^2}\,dD
\]
for the pressure $p(D)$. 
The ``reduced'' Legendre transform on this Lie algebra $\mathfrak{X}(\mathbb{R}^3)$
is given by
\[
\mathbf{m}=D\mathbf{u}
\]
and the Hamiltonian is then written as
\[
h(\mathbf{m},D)=
\langle \mathbf{m},\mathbf{u}\rangle-l(\mathbf{u},D)
\]
that is
\[
h(\mathbf{m},D)= \int \frac{1}{2D}\,|\mathbf{m}|^2\, 
+ De(D) \,d\mathbf{x}
\]
The Lie Poisson equations in this case are as from the general theory
\begin{align*}
\partial_t{\mathbf{m}}
&=
-\text{ad}^*_{\delta h/\delta \mathbf{m}}\, \mathbf{m}
-
\frac{\delta h}{\delta D}\diamond D
\\
\partial_t{D}&=-\mathcal{L}_{\delta h/\delta \mathbf{m}}\, D
\end{align*}
Earlier we found that the coadjoint action is given by the Lie derivative. On the other hand we may calculate the expression of the diamond operation from
its definition
\[
\left\langle 
\frac{\delta h}{\delta D}\,,\,-\pounds_\eta D
\right\rangle
=
\left\langle 
\frac{\delta h}{\delta D} \diamond D \,,\,\eta 
\right\rangle
\]
to be
\[
\left\langle 
\frac{\delta h}{\delta D}\,,\,-{\rm div}\, D\eta
\right\rangle
=
\left\langle 
D\nabla\frac{\delta h}{\delta D}  \,,\,\eta 
\right\rangle
\]
Therefore, we have
\[
\frac{\delta h}{\delta D}\diamond D =  D \nabla\frac{\delta h}{\delta D}
\]
where
\[
\delta h/\delta D= -\,\frac{|\mathbf{m}|^2}{2D^2} + \Big(e + \frac{p}{D}\Big)
\]
Substituting into the momentum equation and using the First Law to find $d(e+p/D)=(1/D)dp$ yields
\[
\partial_t{\mathbf{m}}=-\mathcal{L}_{\mathbf{u}}\, \mathbf{m}-\nabla p
\]
Upon expanding the Lie derivative for the momentum density $\mathbf{m}$ and using the continuity equation for the density, this quickly becomes
\[
\partial_t{\mathbf{u}}=-\,\mathbf{u}\!\cdot\!\nabla{\mathbf{u}}
\, -\,
\frac{1}{D}\nabla p
\]
which is Euler's equation for a barotropic fluid.

\subsection{Euler's equations for ideal incompressible fluid motion} 
The barotropic equations recover Euler's equations for ideal incompressible fluid motion when the internal energy in the reduced Lagrangian for ideal compressible fluids is replaced by the constraint $D=1$, as
\[
l(\mathbf{u},D)
=\int \frac{D}{2}|\mathbf{u}|^2 
- p(D-1)\,  d\mathbf{x}
\]
where again $\mathbf{u}\in\mathfrak{X}(M\!\!\subset\!\mathbb{R}^3)$ is tangential on the boundary $\partial M$ and the advected density $D$ satisfies the continuity equation,
\[
\partial_t{D}+{\rm div}\,D\mathbf{u} =0
\,.
\]
This equation enforces incompressibility ${\rm div}\,\mathbf{u}=0$ when evaluated on the constraint $D=1$. The pressure $p$ is now a Lagrange multiplier, which is determined by the condition that incompressibility be preserved by the dynamics.


\begin{thebibliography}{999}

\bibitem[AbMa1978]{AbMa1978}
 Abraham, R. and Marsden, J.E. [1978], {\em Foundations of Mechanics}.
\newblock Addison-Wesley, second edition.

\bibitem[AcHoKoTi1997]{AcHoKoTi1997}
Aceves, A., D. D. Holm, G. Kovacic and I. Timofeyev [1997]
Homoclinic Orbits and Chaos in a Second-Harmonic
Generating Optical Cavity.
{\it Phys. Lett. A} {\bf 233} 203-208.

\bibitem[AbMaRa1988]{AbMaRa1988}
 Abraham, R.,  Marsden, J.E. and  Ratiu, T.S. [1988], {\em Manifolds,
 Tensor  Analysis, and Applications}. Volume~75 of {\em Applied
 Mathematical  Sciences}.
\newblock Springer-Verlag, second edition.

\bibitem[AlHo1996]{AlHo1996}
     \textrm{Allen, J. S., and Holm, D. D.} [1996]
     {Extended-geostrophic Hamiltonian models for
     rotating shallow water motion.}
     {\it Physica D}, {\bf 98} 229-248.

\bibitem[AlHoNe2002]{AlHoNe2002}
     \textrm{Allen, J. S., Holm, D. D., and Newberger, P. A.} [2002]
     {Toward an extended-geostrophic Euler--Poincar\'e
     model for mesoscale oceanographic flow.} 
     In {\it Large-Scale Atmosphere-Ocean Dynamics 1:
     Analytical Methods and Numerical Models}. Edited by J. Norbury \& I.
     Roulstone, Cambridge University Press: Cambridge, pp. 101--125.

\bibitem[AmCZ96]{AmbCo}
Ambrosetti, A. and  Coti Zelati, V.  [1996], {\em Periodic
Solutions of Singular Lagrangian Systems}.
\newblock Birkh\"auser

\bibitem[AnMc1978]{AnMc1978}
     \textrm{Andrews, D. G. and McIntyre, M. E.} [1978]
     {An exact theory of nonlinear waves on a Lagrangian-mean flow.}
     \textit{J. Fluid Mech.} \textbf{89}, 609--646.

\bibitem[Ar1966]{Ar1966}
Arnold V.I. [1966], Sur la g\'eometrie differentialle des groupes de Lie
de dimiension infinie et ses applications \`a l'hydrodynamique des fluids
parfaits. {\it Ann. Inst. Fourier, Grenoble} {\bf16}, 319-361.

\bibitem[Ar1979]{Ar1979}
Arnold V.I. [1979], {\em Mathematical Methods of Classical
Mechanics}. Volume 60 of \textit{Graduate Texts in Mathematics}.
\newblock Springer-Verlag.

\bibitem[ArKh1998]{ArKh1998}
     \textrm{Arnold, V. I and Khesin, B. A.} [1998]
     \textit{Topological Methods in Hydrodynamics.} 
     Springer: New York.

\bibitem[Be1986]{Benci}
Benci,  V. [1986], Periodic solutions of Lagrangian Systems on a
compact manifold,
  {\em Journal of Differential Equations} \textbf{63}, 135--161.

\bibitem[BlBr92]{BlBr92}
Blanchard,  P. and Bruning, E. [1992], {\em Variational Methods in
Mathematical Physics}.
\newblock Springer Verlag.

\bibitem[Bl2004]{Bl2004}
Bloch, A. M.  (with the collaboration of:
J. Baillieul, P. E. Crouch, and J. E. Marsden)
[2004]
\textit{Nonholonomic Mechanics and Control},
Springer-Verlag NY.

\bibitem[BlBrCr1997]{BlBrCr1997}
 Bloch, A.M., R. W. Brockett, and P.E. Crouch [1997]
Double bracket equations and geodesic flows on symmetric
spaces. {\it Comm. Math Phys} {\bf 187},357-373.

\bibitem[BlCr1997]{BlCr1996}
Bloch, A. M. and P. E. Crouch [1996]
Optimal control and geodesic flows.
{\it Systems and Control Letters}
{\bf 28}, 65-72.

\bibitem[BlCrHoMa2001]{BlCrHoMa2001} 
Bloch, A.M., P.E. Crouch, D. D. Holm and J. E. Marsden [2001]
Proc. of the 39th IEEEE Conference on Decision and Control,
{\it Proc. CDC} {\bf39}, 1273Ð1279.

\bibitem [Bo71]{Bo1971} 
Bourbaki, N. [1971],
\textit{Vari\'et\'es diff\'erentielles et analytiques. Fascicule de
r\'esultats}. Hermann.

\bibitem[Bo1989]{Bo1989}  Lie {1-3} Bourbaki, N. [1989],
{\it Lie Groups and Lie Algebras. Chapters 1--3.} Springer Verlag.

\bibitem[Ca1995]{Ca1995} 
     \textrm{Calogero, F.} [1995]
     {An integrable hamiltonian system.} 
     \textit{Phys. Lett. A} {\bf 201}, 306-310. 

\bibitem[CaFr1996]{CaFr1996}  
     \textrm{Calogero, F. and Francoise, J.-P.}  [1996]
     {An integrable hamiltonian system.}
     \textit{J. Math. Phys.} {\bf 37}, 2863-2871. 

\bibitem[CaHo1993]{CaHo1993} 
     \textrm{Camassa, R. and Holm, D. D.} [1993]
     {An integrable shallow water equation with peaked solitons.}  
     \textit{Phys. Rev. Lett.} {\bf 71}, 1661-64.

\bibitem[CeMaPeRa2003]{CeMaPeRa2003}
Cendra, H, Marsden, J.E., Pekarsky, S., and Ratiu, T.S. [2003],
Variational principles for Lie-Poisson and Hamilton-Poincar\'e equations,
\textit{Moskow Math. Journ.}, \textbf{3}(3), 833--867.

\bibitem[Ch1969]{Ch1969} 
Chandrasekhar, S. [1969] {\it Ellipsoidal Figures of Equilibrium}. Yale University
Press.

\bibitem[ChMa1974]{ChMa1974}
Chernoff, P. R. and Marsden, J.E. [1974], \textit{Properties of
Infinite Dimensional Hamiltonian Systems}, Lecture Notes in
Mathematics, {\bf 425}, Springer-Verlag.

\bibitem[CoLe1984]{CoLe1984} 
Coddington, E. A. and Levinson, N. \textit{Theory of Ordinary Differential
Equations}, Melbourne, FL: Robert E. Krieger, 1984.

\bibitem[CuBa1997]{CuBa1997}
Cushman, R.H. and Bates, L.M. [1997], {\em Global Aspects of
Integrable Systems}.
\newblock Birkh\"auser.

\bibitem[DaHo1992]{DaHo1992}
David, D. and D. D. Holm, [1992]
Multiple Lie-Poisson Structures, Reductions, and 
Geometric Phases for the Maxwell-Bloch Traveling-Wave Equations.  
{\it J. Nonlin. Sci.} {\bf 2}  241--262.

\bibitem[De1860]{De1860}
Dedekind, R. [1860] 
Zusatz zu der vorstehenden Abhandlung. {\it J. Riene Angew. Math.}
{\bf58}, 217-228.

\bibitem[DeMe1993]{DeMe1993} 
Dellnitz M. and Melbourne I. [1993], The equivariant Darboux Theorem,
\textit{Lectures in Appl. Math.}, \textbf{29}, 163--169.

\bibitem[DrJo1989]{DrJo1989}
Drazin,  P.G. and Johnson, R.S. [1989], {\em Solitons: An
Introduction}.
\newblock Cambridge University Press.

\bibitem[DuNoFo1995]{DuNoFo1995}
Dubrovin,  B., Novikov, S.P., and Fomenko, A.T. [1995], {\em Modern
Geometry I,  II, III}, Volumes 93, 104, 124 of \textit{Graduate
Texts in Mathematics}.
\newblock Springer-Verlag.

\bibitem[DuGoHo2001]{DuGoHo2001}
     \textrm{Dullin, H., Gottwald, G. and Holm, D. D.} [2001]
     {An integrable shallow water equation
     with linear and nonlinear dispersion.}
     \textit{Phys. Rev. Lett.} \textbf{87} 194501-04.

\bibitem[DuGoHo2003]{DuGoHo2003}
     \textrm{Dullin, H., Gottwald, G. and Holm, D. D.} [2003]
     {Camassa-Holm, Korteweg-de~Vries-5 and other asymptotically
equivalent equations for shallow water waves.}
{\it Fluid Dyn. Res.} {\bf33} (2003) 73Ð95.

\bibitem[DuGoHo2004]{DuGoHo2004}
     \textrm{Dullin, H., Gottwald, G. and Holm, D. D.} [2003]
     {On asymptotically equivalent shallow water wave equations.}
{\it Physica D} {\bf190} (2004) 1-14.

\bibitem[EbMa1970]{EbMa1970}
Ebin, D.G. and  Marsden, J.E. [1970], Groups of diffeomorphisms and
the motion of an incompressible fluid, {\it Ann. of Math.}
\textbf{92}, 102--163.


\bibitem[Er1942]{Er1942} Ertel, H. [1942] Ein Neuer Hydrodynamischer Wirbelsatz, Met. Z., \textbf{59}, 271--281.

\bibitem[FoHo1991]{FoHo1991}
Fordy, A.P.  and Holm, D. D. [1991]
A Tri-Hamiltonian Formulation of the 
Self-Induced Transparency Equations.
{\it Phys. Lett A} {\bf 160} 143--148.

\bibitem[FrHo2001]{FrHo2001}
     \textrm{Fringer, O. B. and Holm, D. D.} [2001]
     {Integrable vs nonintegrable geodesic soliton behavior.}
     \textit{Physica D} \textbf{150} 237-263. 

\bibitem[Fu1996]{Fu1996}
Fuchssteiner, B. [1996]
Some tricks from the symmetry-toolbox for nonlinear equations: generalization of the
CamassaÐHolm equation. {\it Physica D} {\bf95}, 229-243.

\bibitem[GeDo1979]{GeDo1979} 
\textrm{Gelfand, I. M. and Dorfman, I. Ya. R.} [1979]
{\it Funct. Anal. Appl.} {\bf13}, 248.

\bibitem[GiHoKu1982]{GiHoKu1982}
Gibbons, J., Holm, D.~D. and B.~A. Kupershmidt [1982],
Gauge-Invariant Poisson Brackets for Chromohydrodynamics, 
{\it Phys. Lett. A} {\bf 90} 281--283.

\bibitem[GjHo1996]{GjHo1996}
     \textrm{Gjaja, I. and Holm, D. D.}
     {Self-consistent wave-mean flow interaction dynamics and its
     Hamiltonian formulation for a rotating stratified incompressible
     fluid.}  
     {\it Physica D}, {\bf 98} (1996) 343-378.

\bibitem[GuSt1984]{GuSt1984}
Guillemin,  V. and Sternberg, S. [1984], {\em Symplectic Techniques
in Physics}.
\newblock Cambridge University Press.

\bibitem[Ha1984]{Ha1984}
Haine, L. [1984]
The algebraic complete integrability of geodesic flow on SO(4).
{\it Commun, Math, Phys.} {\bf94}, 271-287.

\bibitem[Ho1983]{Ho1983} 
Holm, D.D. [1983]
Magnetic tornadoes: three-dimensional affine motions in
ideal magnetohydrodynamics, {\it Physica D} {\bf 8}, 170-182.

\bibitem[Ho1991]{Ho1991}
     \textrm{Holm, D. D.} [1991]
     {Elliptical vortices and integrable Hamiltonian dynamics of the rotating
shallow-water equations.} 
     {\it J. Fluid Mech.}, {\bf 227} 393-406.

\bibitem[Ho1996a]{Ho1996a}
     \textrm{Holm, D. D.} [1996a]
     {Hamiltonian balance equations.} 
     {\it Physica D}, {\bf 98} 379-414.

\bibitem[Ho1996b]{Ho1996b}
     \textrm{Holm, D. D.} [1996b]
     {The ideal Craik-Leibovich equations.} 
     {\it Physica D}, {\bf 98} (1996) 415-441.

\bibitem[Ho1999]{Ho1999}
     \textrm{Holm, D. D.} [1999] 
     {Fluctuation effects on 3D Lagrangian mean and Eulerian mean fluid
     motion.} \textit{Physica D} \textbf{133}, 215--269.

\bibitem[Ho2002]{Ho2002}
     \textrm{Holm, D. D.} [2002] 
     {Averaged Lagrangians and the mean dynamical effects of
     fluctuations in continuum mechanics.} 
     \textit{Physica D} \textbf{170}, 253--286.

\bibitem[Ho2002a]{Ho2002a}
     \textrm{Holm, D. D.} [2002] 
     {Euler-Poincar\'e dynamics of perfect complex fluids.} 
     In {\it  Geometry, Mechanics, and Dynamics: 
     in honor of the 60th birthday of Jerrold E. Marsden} edited by 
     P. Newton, P. Holmes and A. Weinstein. Springer, pp. 113-167.

\bibitem[Ho2002b]{Ho2002b}
     \textrm{Holm, D. D.} [2002] 
     {Variational principles for Lagrangian averaged fluid dynamics.}
     {\it J. Phys. A: Math. Gen.} {\bf 35} 1-10.

\bibitem[Ho2002c]{Ho2002c}
     \textrm{Holm, D. D.} [2002] 
     {Lagrangian averages, averaged Lagrangians, and the mean effects 
      of fluctuations in fluid dynamics.}
     {\it Chaos} {\bf 12} 518-530.

\bibitem[Ho2005]{Ho2005}
     \textrm{Holm, D. D.} [2005] 
The Euler-Poincar\'e variational framework
for modeling fluid dynamics, 
in {\it Geometric Mechanics and Symmetry: The
Peyresq Lectures}, edited by J. Montaldi and T. Ratiu,Ê London
Mathematical Society Lecture Notes Series 306, Cambridge University Press
(2005).

\bibitem[HoKo1991]{HoKo1991}
Holm, D.~D. and G. Kovacic [1991]
Homoclinic Chaos for Ray Optics in a Fiber. 
{\it Physica D} {\bf 51}, 177--188.

\bibitem[HoKo1992]{HoKo1992}
Holm, D.~D. and G. Kovacic [1992]
Homoclinic Chaos in a Laser-Matter System. 
{\it Physica D} {\bf 56} 270--300.

\bibitem[HoKu1982]{HoKu1982}
Holm, D.~D. and B.~A. Kupershmidt [1982],
Poisson Structures of Superfluids, 
{\it Phys. Lett. A} {\bf 91} 425--430.

\bibitem[HoKu1983]{HoKu1983}
Holm, D.~D. and B.~A. Kupershmidt [1983],
Poisson brackets and Clebsch representations 
for magnetohydrodynamics, multifluid  plasmas, and elasticity, 
{\it Physica D} {\bf 6}, 347--363.

\bibitem[HoKu1988]{HoKu1988}
Holm, D.~D. and B.~A. Kupershmidt [1988],
The Analogy Between Spin Glasses 
and Yang-Mills Fluids,
{\it J. Math Phys.} {\bf 29} 21--30.

\bibitem[HoMa2004]{HoMa2004}
Holm, D.D. and J. E. Marsden [2004]
Momentum maps and measure valued solutions (peakons, filaments,
and sheets) of the Euler-Poincar«e equations
for the diffeomorphism group.
{\it In The Breadth of Symplectic and Poisson Geometry}, (Marsden, J. E.
and T. S. Ratiu, eds) Birkh\"auser Boston, to appear (2004).
http://arxiv.org/abs/nlin.CD/0312048

\bibitem[HoMaRa1986]{HoMaRa1986}
     \textrm{Holm, D. D., Marsden, J. E. and Ratiu, T. S.} [1986] 
     {The Hamiltonian Structure of Continuum Mechanics in Material, Inverse
     Material, Spatial and Convective Representations.}
     In {\it Hamiltonian Structure and
     Lyapunov Stability for Ideal Continuum  Dynamics}, 
     Univ. Montreal Press, pp. 1-124.

\bibitem[HoMaRa1998a]{HoMaRa1998a}
     \textrm{Holm, D. D., Marsden, J. E. and Ratiu, T. S.} [1998a] 
     {The Euler--Poincar\'e equations and semidirect products
     with applications to continuum theories.}
     \textit{Adv. in Math.} \textbf{137}, 1--81.

\bibitem[HoMaRa1998b]{HoMaRa1998b}
     \textrm{Holm, D. D., Marsden, J. E. and Ratiu, T. S.} [1998b] 
     {Euler--Poincar\'e models of ideal fluids with nonlinear dispersion.}
     \textit{Phys. Rev. Lett.} \textbf{80}, 4173--4177.

\bibitem[HoMaRa2002]{HoMaRa2002}
     \textrm{Holm, D. D., Marsden, J. E. and Ratiu, T. S.} [2002]
     {The Euler--Poincar\'{e} equations in geophysical fluid dynamics.} 
     In {\it Large-Scale Atmosphere-Ocean Dynamics 2:
     Geometric Methods and Models}. Edited by J. Norbury and I.
     Roulstone, Cambridge University Press: Cambridge, pp. 251--299. 

\bibitem[HoMaRaWe1985]{HoMaRaWe1985}
     \textrm{Holm, D. D., Marsden, J. E., Ratiu, T. S. and Weinstein, A.}
     [1985]
     {Nonlinear stability of fluid and plasma equilibria.} 
     \textit{Physics Reports} {\bf 123} 1-116.

\bibitem[HoRaTrYo2004]{HoRaTrYo2004}
Holm, D. D., J. T. Rananather, A. Trouv\'e and L. Younes [2004]
Soliton Dynamics in Computational Anatomy.
{\it NeuroImage} {\bf23}, S170-178.
http://arxiv.org/abs/nlin.SI/0411014

\bibitem[HoSt2002]{HoSt2002}
     \textrm{Holm, D. D. and Staley, M. F.} [2002]
     {Wave Structures and Nonlinear Balances in a Family of 1+1
     Evolutionary PDEs.} http://arxiv.org/abs/nlin.CD/0202059.
     Submitted to {\it SIADS.}

\bibitem[HoSt2004]{HoSt2004}
     \textrm{Holm, D. D. and Staley, M. F.} [2004]
     Interaction Dynamics of Singular Wave Fronts.
     Submitted to {\it SIAM J. Appl. Dyn. Syst.}

\bibitem[HoWo1991]{HoWo1991}
Holm, D. D. and K.B. Wolf [1991]
Lie-Poisson Description of Hamiltonian Ray Optics. 
{\it Physica D}  {\bf 51}, 189--199.

\bibitem[HuZh1994]{HuZh1994} 
     \textrm{Hunter, J. K. and Zheng, Y.} [1994]
     {On a completely integrable nonlinear hyperbolic variational
      equation.} 
     \textit{Physica D} {\bf 79}, 361-386.

\bibitem[JoSa98]{JoSa98} Jos\'e, J.V. and Saletan, E.J. [1998],
{\em Classical Dynamics : A Contemporary Approach\/}.
Cambridge University Press.

\bibitem[Jo1998]{Jo1998}
Jost, J. [1998], {\em Riemannian Geometry and Geometric Analysis}.
{\em University Text}.
\newblock Springer-Verlag, second edition.

\bibitem[KaKoSt1978]{KaKoSt1978} 
Kazhdan, D.  Kostant, B., and Sternberg,  S. [1978], Hamiltonian
group actions and dynamical systems of Calogero type, \textit{Comm.
Pure Appl. Math.}, \textbf{31}, 481--508. 

\bibitem[KhMis03]{KhMi2003}
Khesin, B. and Misiolek, G. [2003], Euler equations on homogeneous
spaces and Virasoro orbits, \textit{Adv. in Math.\/}, \textbf{176},
116--144.

\bibitem[Ko1966]{Ko1966} Kostant, B. [1966]
Orbits, symplectic structures and representation theory. {\it Proc.
US-Japan Seminar on Diff. Geom., Kyoto. Nippon Hyronsha, Tokyo} {\bf77}.

\bibitem[KrSc2001]{KrSc2001}
     \textrm{Kruse, H. P., Schreule, J. and Du, W.} [2001] 
     {A two-dimensional version of the CH equation}.
     In \textit{Symmetry and Perturbation Theory: SPT 2001}
     Edited by D. Bambusi, G. Gaeta and M. Cadoni.
     World Scientific: New York, pp 120-127.

\bibitem[La1999]{La1999} 
Lang, S. [1999],
\textit{Fundamentals of Differential Geometry}. Volume 191 of
\textit{Graduate Texts in Mathematics}. Springer-Verlag, New York.

\bibitem[Le2003]{Le2003}
Lee, J.  [2003]
\textit{Introduction to Smooth Manifolds},
Springer-Verlag.

\bibitem[LiMa1987]{LiMa1987} 
Libermann, P. and Marle,  C.-M.  [1987], \textit{Symplectic Geometry and
Analytical Mechanics}. Reidel.

\bibitem[Lie1890]{Lie1890}
Lie, S. {Theorie der Transformationsgruppen. Zweiter Abschnitt}. Teubner.

\bibitem[Man1976]{Man1976} 
Manakov, S.V. [1976]
Note on the integration of Euler's equations of the dynamics of
and $n$-dimensional rigid body.
{\it Funct. Anal. and its Appl.} {\bf 10}, 328--329.

\bibitem[Mar1976]{Mar1976}   Marle, C.-M. [1976],
Symplectic manifolds, dynamical groups, and Hamiltonian
mechanics, in \textit{Differential Geometry and Relativity\/}, 
Cahen, M. and Flato, M. eds., D. Reidel, Boston,  249--269. 

\bibitem[Ma1981]{Ma1981}  Marsden, J.E. [1981],
\textit{Lectures on Geometric Methods in Mathematical Physics\/}.
Volume 37, SIAM, Philadelphia, 1981. 

\bibitem[Ma1992]{Ma1992}
Marsden,  J.E.  [1992], {\em Lectures on Mechanics}. Volume~174 of
{\em London Mathematical Society Lecture Note Series}.
\newblock Cambridge University Press.

\bibitem[MaHu1983]{MaHu1983}
Marsden, J.~E. and T.~J.~R. Hughes [1983], 
{\em Mathematical Foundations of Elasticity}.
\newblock Prentice Hall.
\newblock Reprinted by Dover Publications, NY, 1994.

\bibitem[MaMiOrPeRa2004]{MaMiOrPeRa2004}
Marsden,  J.E., Misiolek, G., Ortega, J.-P., Perlmutter, M., and
Ratiu, T.S.  [2004], {\em Hamiltonian Reduction by Stages}.
Lecture Notes in Mathematics, to appear.
\newblock Springer-Verlag.

\bibitem[MaMoRa90]{MaMoRa1990} 
Marsden, J.E., Montgomery, R., and Ratiu, T.S. [1984],
Reduction, symmetry, and phases in mechanics,  
\textit{Memoirs Amer. Math. Soc.}, \textbf{88}(436), 1--110.

\bibitem[MaRa1994]{MaRa1994}
Marsden,  J.E.  and  Ratiu, T.S. [1994], {\em Introduction to
Mechanics and Symmetry}. Volume~75 of {\em Texts in Applied
Mathematics}, second printing of second edition 2003.
\newblock Springer-Verlag.

\bibitem[MaRa95]{MaRa1995} Marsden, J.E. and  Ratiu,
T.S. [2003], \textit{Geometric Fluid Dynamics}. Unpublished
notes.

\bibitem[MaRa03]{MaRa2003} Marsden, J.E. and  Ratiu,
T.S. [2003], \textit{Mechanics and Symmetry. Reduction Theory}.
In preparation.

\bibitem[MaRaWe84a]{MaRaWe1984a} 
Marsden, J.E., Ratiu, T.S., and Weinstein, A. [1984a], Semidirect products
and reduction in mechanics,
\textit{Trans. Amer. Math. Soc.}, \textbf{281}(1), 147--177.

\bibitem[MaRaWe84b]{MaRaWe1984b} 
Marsden, J.E., Ratiu, T.S., and Weinstein, A. [1984b], Reduction and
Hamiltonian structures on duals of semidirect product Lie algebras,
\textit{Contemporary Math.}, \textbf{28}, 55--100.

\bibitem[MaWe74]{MaWe1974} 
Marsden, J.E. and Weinstein, A. [1983], Reduction of symplectic manifolds
with symmetry, \textit{Rep. Math. Phys.}, \textbf{5}, 121--130.

\bibitem[MaWe1983]{MaWe1983} 
Marsden, J.E. and Weinstein, A. [1983], Coadjoint orbits, vortices and
Clebsch variables for incompressible fluids, \textit{Physica D},
\textbf{7}, 305--323.

\bibitem[MaWil989]{MaWil989}
Mawhin,  J. and Willem, M. [1989], {\em Critical Point Theory and
Hamiltonin Systems}. Volume~74 of {\em Applied
 Mathematical  Sciences}.
\newblock Springer-Verlag, second edition.

\bibitem[McSa1995]{McSa1995}
McDuff,  D. and Salamon,  D. [1995], {\em Introduction to Symplectic
Topology}.
\newblock Clarendon Press.

\bibitem[MeDe1993]{MeDe1993} Melbourne, I. and
Dellnitz, M. [1993], Normal forms for linear Hamiltonian vector
fields commuting with the action of a compact Lie group,
\textit{Proc. Camb. Phil. Soc.}, \textbf{114}, 235--268.

\bibitem[Mi1963]{Mi1963}
Milnor, J. [1963], {\em Morse Theory}.
\newblock Princeton University Press.

\bibitem[MiFo1978]{MiFo1978}
Mishchenko, A. S. and Fomenko, A. T.: Euler equations on finitedimensional
Lie groups, Izv. Acad. Nauk SSSR, Ser. Matem. 42, No.2,
396-415 (1978) (Russian); English translation: Math. USSR-Izv. 12, No.2,
371-389 (1978).

\bibitem[Mi2002]{Mi2002}
Misiolek, G. [2002], Classical solutions of the periodic
Camassa-Holm equation, \textit{Geom. Funct. Anal.\/}, \textbf{12},
1080--1104.

\bibitem[No1918]{No1918}
Noether, E.  [1918] {\it Nachrichten Gesell. Wissenschaft. G\"ottingen}
{\bf2} 235. See also C. H. Kimberling [1972] {\it Am. Math. Monthly}
{\bf79} 136.

\bibitem[Ohk1993]{Ohk1993} Ohkitani, K. [1993] Eigenvalue problems in three-dimensional Euler flows. 
Phys. Fluids A, \textbf{5}, 2570--2572.

\bibitem[Ol2000]{Ol2000}
Olver, P. J.  [2000] {\it Applications of Lie Groups to Differential
Equations}. Springer: New York.

\bibitem[OrRa2004]{OrRa2004}
Ortega, J.-P. and Ratiu, T.S. [2004], {\em Momentum Maps and
Hamiltonian Reduction}. Volume 222 of {\em Progress in Mathematics}.
\newblock Birkh\"auser.

\bibitem[OKh87]{OvKh1987}
Ovsienko, V.Y. and Khesin, B.A. [1987], Korteweg-de Vries superequations as
an Euler equation, \textit{Funct. Anal. Appl.}, \textbf{21},
329--331.

\bibitem[Pa1968]{Pa1968} Palais, R. [1968],
\textit{ Foundations of Global Non-Linear Analysis}, W. A. Benjamin,
Inc., New York-Amsterdam.

\bibitem[Po1901]{Po1901} 
     \textrm{Poincar\'e, H.} [1901]
     {Sur une forme nouvelle des \'{e}quations de la m\'{e}chanique}
     {\it C.R. Acad. Sci.\/} {\bf 132}, 369--371.


\bibitem[Ra1980]{Ra1980} 
Ratiu, T. [1980] The motion of the free
n-dimensional rigid body. {\it Indiana U. Math. J.}
{\bf 29}, 609-627.

\bibitem[RaTuSbSoTe2005]{RaTuSbSoTe2005}
Ratiu, T.S., R. Tudoran, L. Sbano, E. Sousa Dias, G. Terra, A Crash course
in Geometric Mechanics, in {\it Geometric Mechanics and Symmetry: The
Peyresq Lectures}, edited by J. Montaldi and T. Ratiu,Ê London
Mathematical Society Lecture Notes Series 306, Cambridge University Press
(2005).

\bibitem[Ri1860]{Ri1860}
Riemann, B. [1860]
Untersuchungen \"uber die Bewegungen eines fl\"ussigen gleichartigen Ellipsoides.
{\it Abh. d. K\"onigl. Gesell. der Wis. zu G\"ottingen} {\bf9}, 3-36. 

\bibitem[Sc1987]{Sc1987}
Schmid,  R. [1987], {\em Infinite Dimensional Hamiltonian Systems}.
\newblock Bibliopolis.

\bibitem[Se1992]{Se1992}
Serre,  J.-P. [1992], {\em Lie Algebras and Lie Groups}, Volume~1500 of
{\em Lecture Notes in Mathematics}.
\newblock Springer-Verlag.

\bibitem[Sh1998]{Sh1998}
     \textrm{Shkoller, S.} 1998 
     {Geometry and curvature of diffeomorphism groups with $H^1$ metric and
      mean hydrodynamics}
     {\it J. Funct. Anal.,} {\bf160}, 337-365.

\bibitem[Sh2000]{Sh2000}
     \textrm{Shkoller, S.} 2000 
     {Analysis on groups of diffeomorphisms of manifolds with boundary and
      the averaged motion of a fluid}
     {\it J. Differential Geometry,} {\bf55}, 145-191.

\bibitem[Sm1970]{Sm1970} Smale, S.
Topology and Mechanics, {\it Inv. Math.} {\bf10}, 305-331; {\bf11}, 45-64.

\bibitem[So1970]{So1970} Souriau, J. M. [1970]
{\it Structure des Syst\`emes Dynamiques}, Dunod, Paris.

\bibitem[Sp1979]{Sp1979} Spivak, M. [1979],
\textit{Differential Geometry\/}, Volume I. New printing with
corrections. Publish or Perish, Inc. Houston, Texas.

\bibitem[Va1996]{Va1996} Vaisman, I. [1996],
{\it Lectures on the Geometry of Poisson Manifolds}. Volume 118 of
{\it Progress in
Mathematics}, Birkh\"auser.

\bibitem[Wa1983]{Wa1983} Warner, F.W. [1983],
\textit{Foundation of Differentiable Manifolds and Lie Groups}.
Volume 94 of \textit{Graduate Texts in Mathematics},
Springer-Verlag.

\bibitem[We1983]{We1983a} Weinstein, A. [1983a], Sophus Lie and
symplectic geometry,  {\it Exposition Math.\/} {\bf 1}, 95-96.

\bibitem[We1983b]{We1983b} Weinstein, A. [1983b], The local structure of
Poisson manifolds,  {\it Journ. Diff. Geom\/} {\bf 18}, 523--557.

\bibitem[We2002]{We2002}
Weinstein, A. [2002] Geometry of momentum (preprint);
ArXiv:math/SG0208108 v1.

\end{thebibliography}
\end{document}